\documentclass[11pt]{amsart}

\usepackage{evearticle}

\makeatletter
\renewcommand*\env@matrix[1][\arraystretch]{
	\edef\arraystretch{#1}
	\hskip -\arraycolsep
	\let\@ifnextchar\new@ifnextchar
	\array{*\c@MaxMatrixCols c}}
\makeatother

\newcommand{\one}{\mathbf{1}}

\usepackage{etoolbox}
\newcounter{jmstep}
\newcounter{jmstepanchor}
\renewcommand{\thejmstep}{\arabic{jmstep}}

\crefname{jmstep}{Step}{Steps}
\Crefname{jmstep}{Step}{Steps}

\newcommand{\stepn}[2]{
	\refstepcounter{jmstep}
	\stepcounter{jmstepanchor}
	\label{#1}
	\noindent\emph{Step \thejmstep. #2.}
}
\AtBeginEnvironment{proof}{
	\setcounter{jmstep}{0}
}
\newcounter{jmsubstep}[jmstep]
\renewcommand{\thejmsubstep}{\thejmstep\alph{jmsubstep}}
\newcounter{jmsubstepanchor}

\crefname{jmsubstep}{Step}{Steps}
\Crefname{jmsubstep}{Step}{Steps}

\newcommand{\substepn}[2]{
	\refstepcounter{jmsubstep}
	\stepcounter{jmsubstepanchor}
	\label{#1}
	\noindent\emph{Step \thejmsubstep. #2.}
}

\begin{document}

\title[The conformally invariant metric on $\CLE_4$ III: uniqueness]{The conformally invariant metric on $\CLE_4$ III:\\
uniqueness}
\author[E.\ Kammerer, K.\ Kavvadias, J.\ Miller and Y.\ Tian]{Emmanuel Kammerer, Konstantinos Kavvadias, Jason Miller and Yi Tian}

\date{\today}

\maketitle

\begin{abstract}
This paper is the third and final article in a series of papers constructing the canonical conformally invariant metric on the set of loops of the conformal loop ensemble ($\CLE$) with critical parameter $\kappa=4$. The previous two articles construct, as a subsequential limit of the renormalized graph metric on the loops of $\CLE_\kappa$ as $\kappa \downarrow 4$, a conformally invariant, local metric on the loops of a $\CLE_4$ whose metric ball growth from the domain boundary coincides with the uniform exploration of Werner and Wu. In this paper, we establish that this metric is uniquely characterized by its properties, as are its geodesics, and that it is a measurable function of the $\CLE_4$. In particular, we show that the renormalized $\CLE_\kappa$ graph metric converges as $\kappa \downarrow 4$ without passing to a subsequence. A key step in the proof is to show that the metric is determined by the geodesics from each loop to the domain boundary, which are in turn determined by the uniform exploration; this representation will have important applications in future work.
\end{abstract}

\tableofcontents

\setlength{\parindent}{0pt}
\setlength{\parskip}{0.5\baselineskip plus 1pt minus 1pt}

\section{Introduction}

\subsection{Overview}

The \emph{conformal loop ensemble} ($\CLE_\kappa$) with parameter $\kappa \in (8/3, 8)$ is a random collection of non-crossing loops within a simply connected domain in $\BC$ \cite{TreeCLE,CLE}. The $\CLE_\kappa$ is conformally invariant, meaning its law is preserved under conformal automorphisms of the domain.  Its law in an arbitrary simply connected domain is then defined by pushing forward under a conformal map, so that conformal invariance holds for all conformal maps between such domains. CLEs are the loop variant of the Schramm-Loewner evolution (SLE) \cite{s2000sle}.  SLEs and CLEs are conjectured (and in several cases rigorously shown) to describe the scaling limits of interfaces in critical statistical mechanics models on planar lattices \cite{s2001percolation,lsw2004lerw,ss2009dgff,s2010ising} as well as on random planar maps \cite{s2016hc,lsw2017schnyder,gkmw2018active,kmsw2019bipolar,gm2021saw,gm2021percolation}.

\begin{figure}[ht]
	\includegraphics[width=0.49\textwidth]{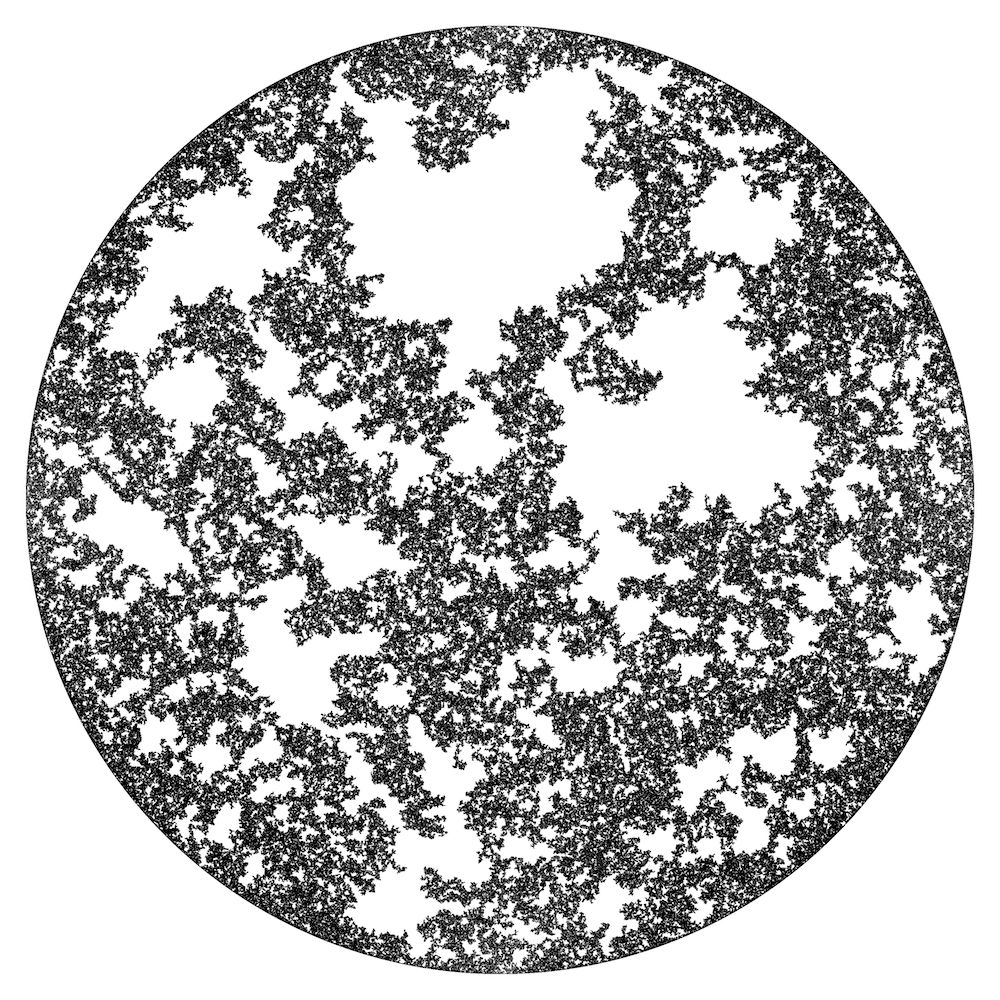}
	\includegraphics[width=0.49\textwidth]{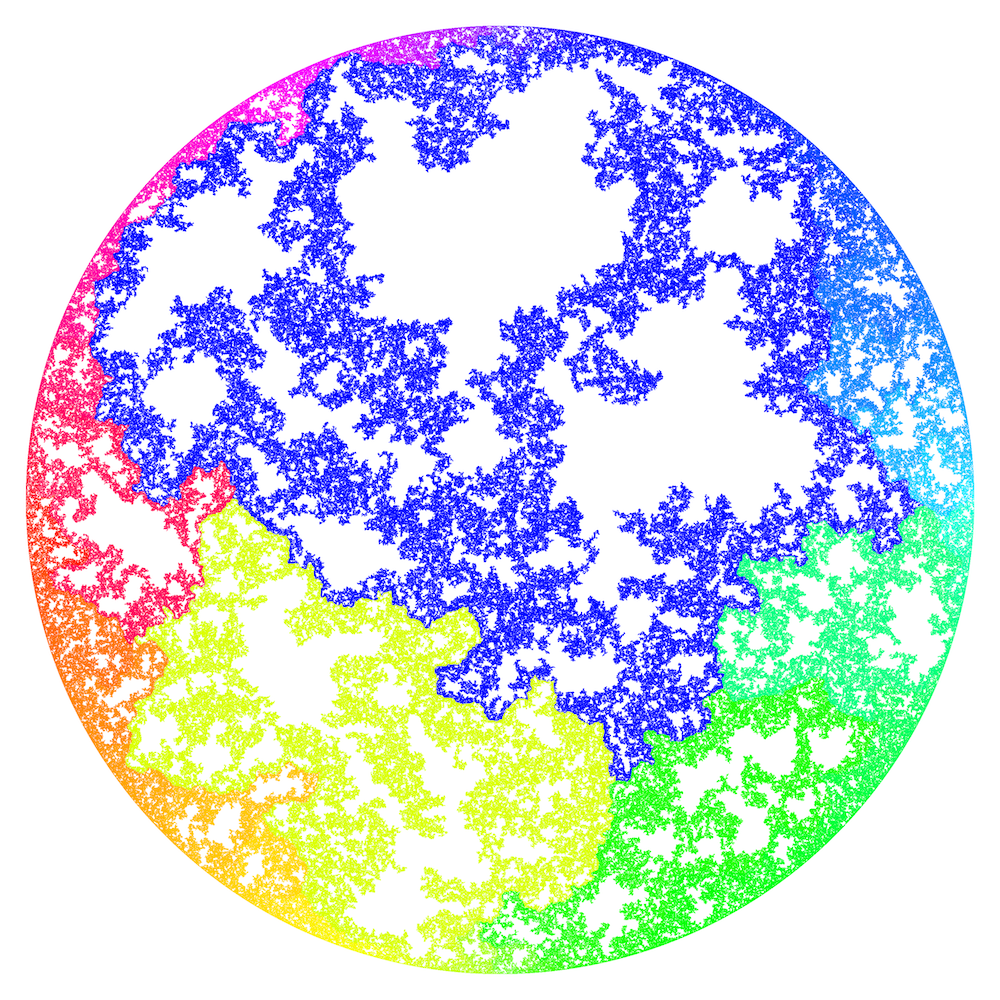}
	\caption{\label{fig:simulations}
		{\bf Left:} A (discrete approximation of a) $\CLE_4$ in the unit disk $\BD$.  {\bf Right:} The loops are colored according to where their unique geodesic to $\partial \BD$ hits $\partial \BD$.  Loops whose geodesics reach $\partial \BD$ at nearby points therefore have similar colors, and loops whose geodesics merge before reaching $\partial \BD$ have the same color.}
	
	\end{figure}

Just like for $\SLE$, the value $\kappa=4$ is critical for $\CLE$.  On the one hand, for $\kappa \in (8/3,4]$, the loops are simple and do not intersect each other or the domain boundary.  On the other hand, for $\kappa \in (4,8)$, the loops are self-intersecting, intersect each other, and also the domain boundary \cite{rs2005basic}.  In this paper, we focus on the critical case $\kappa=4$.  For $\kappa > 4$, since the loops intersect, there is a natural distance that one can define on the set of loops of a $\CLE_\kappa$.  It is simply the graph metric where the vertices of the graph are the loops of the $\CLE_\kappa$ and there is an edge precisely when two loops intersect.  Since a $\CLE_\kappa$ is conformally invariant and this metric is a function of a $\CLE_\kappa$, the metric is trivially conformally invariant.  This construction does not make sense when $\kappa = 4$ since the loops do not intersect each other. Nevertheless, we proved in \cite{kkmt2026cle4_part1, kkmt2026cle4_part2} that such a conformally invariant metric on the loops of a $\CLE_4$ does exist and that it is a geodesic metric. The main aim of this paper is to show that this conformally invariant metric is unique.  See \Cref{fig:simulations} for numerical simulations of the metric.

The program of constructing the $\CLE_4$ metric was initiated by Werner and Wu in \cite{CoInCLEExpl}, in which the distance from any loop in the $\CLE_4$ to the domain boundary was constructed.  Equivalently, they gave a definition of the growth of a metric ball from the domain boundary, where one discovers each loop at a time corresponding to its distance from the domain boundary.  They called this the \emph{uniform exploration} (cf.\ also \cite{LevelLineGFFI,TVSGFF}).  The reason for this terminology is that as one performs the uniform exploration from the domain boundary, the $\CLE_4$ loops are discovered in a Poissonian way, each new loop being rooted at a point chosen according to the harmonic measure on the boundary of the explored region infinitesimally before the time at which it appears.  It was not proved in \cite{CoInCLEExpl} that the uniform exploration corresponds to the growth of the metric balls from the boundary for a metric defined on the entire $\CLE_4$. This is one of the properties of the metric constructed in \cite{kkmt2026cle4_part1,kkmt2026cle4_part2}. It was also conjectured in \cite{CoInCLEExpl} that the uniform exploration is a.s.\ determined by the $\CLE_4$. In this third paper, we prove this conjecture. More precisely, we establish the uniqueness of the metric constructed in \cite{kkmt2026cle4_part1, kkmt2026cle4_part2}, whose metric balls correspond to the uniform exploration, as well as the fact that the metric is a.s.\ determined by the $\CLE_4$.

\begin{theorem}[See \Cref{thm:uniqueness} for a precise statement]
	There exists a unique geodesic metric on the set of loops of the $\CLE_4$ which is local and conformally invariant, and whose metric balls from the domain boundary correspond to the uniform exploration. Moreover, this metric is a measurable function of the $\CLE_4$.
\end{theorem}

The metric constructed in this series of articles was also obtained by Sheffield, Watson, and Wu in unpublished work.  In contrast to the indirect approach taken by Sheffield, Watson, and Wu, we already showed in \cite{kkmt2026cle4_part1} that the metric arises as an appropriately renormalized limit as $\kappa \downarrow 4$ of the adjacency graph metric of $\CLE_\kappa$ loops. Moreover, in \cite{kkmt2026cle4_part2}, we proved the existence of geodesics: this question was not tackled by Sheffield, Watson, and Wu. In this paper, one of the main steps in establishing the uniqueness of the metric will be to show that the $\CLE_4$ distance between arbitrary loops can be represented in terms of the uniform exploration. In particular, we prove that geodesics can be approximated by going ``up and down'' in the geodesics to the boundary. This representation will have important consequences for future work, which we will describe in more detail just below.

We anticipate that the $\CLE_4$ metric will have applications in the study of random planar maps.  Recall that a \emph{planar map} is a graph together with an embedding into the plane such that no two edges cross. Two planar maps are considered to be equivalent if it is possible to deform one into the other using an orientation-preserving homeomorphism of the plane.  In recent years, the asymptotic behavior of random planar maps has received a large amount of attention.  Major results include the Gromov-Hausdorff-Prokhorov convergence of rescaled quadrangulations toward the Brownian sphere \cite{MM06, Mie13, LG13} and of Boltzmann planar maps with large faces toward the $\alpha$-stable carpet or gasket \cite{LGM11, StabCarGas}.  In this general framework, we expect that our random metric space describes the scaling limit of random maps with high degree vertices called $3/2$-stable maps.  These are the dual maps of random planar maps with large faces, where faces of degree $2k$ have a weight of order $c^{k} k^{-2}$ for the appropriate critical constant $c$. They were first studied in the local limit in \cite{BCM18CauchyMaps}. Even though these maps do not admit a scaling limit in the sense of Gromov-Hausdorff or Gromov-Prokhorov \cite{Kam25large32maps}, a non-trivial scaling limit of the distances to the root has been uncovered by one of the authors in \cite{Kam25distancesonCLE4and32maps} which can be identified with a distance from the loops of the $\CLE_4$ to the boundary introduced in \cite{BECriLQG}, and which, in future work, will be shown to agree with the metric constructed here.

As shown in \cite{Kam26GasketsO2maps} (and in \cite{ADSH26volumeOn} in the particular case of quadrangulations, see also \cite{DSHPW25} for the solvable fully-packed case on triangulations), the gasket of a critical $O(2)$ loop-decorated planar map, obtained by forgetting the parts of the map that are inside the outermost loops, is an instance of a $3/2$-stable map. Our expectation is thus consistent with the conjecture that the scaling limit of critical $O(2)$ loop-decorated maps is described by the $\CLE_4$ together with an independent critical Liouville quantum gravity (see e.g.\@ \cite[Conjecture 2.1]{HL24weightedLQG} for a precise conjecture).

The scaling limits of the random maps with large faces mentioned above (which are often called $\alpha$-stable maps for $\alpha \in (1,2)$) are known as the $\alpha$-stable carpets, or gaskets when $\alpha < 3/2$. In this terminology, our metric space can be seen as the ``dual'' of the $3/2$-stable carpet. We also expect that it is possible to construct analogous duals for the $(1/2 + 4/\kappa)$-stable carpets, where $\kappa \in (8/3, 4)$, corresponding to a metric on the loops of a $\CLE_\kappa$, though they would be conformally covariant instead of conformally invariant.

\subsection{Uniqueness of limits of random metric spaces}
Identifying the limit of a sequence of random metric spaces is a central question in the field of random geometry. Sometimes, the limit can be characterized as a function of a stochastic process. The prototypical example is the Brownian tree \cite{Ald91}, which can be obtained by gluing together those pairs of points on the graph of a Brownian excursion that can be connected by a horizontal chord staying below the graph of the excursion. Similarly, the Brownian sphere, defined in \cite{MM06}, can be obtained as a function of the Brownian snake. A similar construction was obtained in \cite{StabCarGas} for the scaling limit of random maps with large faces. In the same spirit, in this paper, we describe the conformally invariant metric on the $\CLE_4$ as a function of the uniform exploration introduced in \cite{CoInCLEExpl} (moreover, we also establish that the uniform exploration is determined by the $\CLE_4$).

Yet, defining a candidate for the limiting metric space is often not sufficient to establish the uniqueness of the limit. A deeper understanding of the geometry is particularly useful for that purpose. In their independent proofs of the uniqueness of the Brownian sphere, Le Gall \cite{LG13} and Miermont \cite{Mie13} both heavily use properties of geodesics. More precisely, \cite{LG13} shows that geodesics can be approximated by going up and down in the tree of geodesics to a distinguished point. An analogous argument was carried out in \cite{StabCarGas} for random maps with large faces.

In the study of random metric spaces having a natural embedding in the complex plane, the approach for proving uniqueness is different. While still relying on properties of geodesics, it makes heavy use of the embedding itself. For example, to prove the uniqueness of the Liouville quantum gravity metric, \cite{ExUniLQG, UniCriSupercriLQGMet} combine properties of geodesics with the independence across scales of the Gaussian free field (GFF) and the bi-Lipschitz equivalence of two conditionally independent samples of the metric given the GFF.  The same pair of ideas (a bi-Lipschitz comparison of two conditionally independent samples of the metric, together with a suitable independence across scales) is also important for the proofs of the uniqueness of the conformally covariant metric on the simple CLE carpet \cite{GeoCLECarp}, of the conformally covariant metric on the non-simple CLE gasket \cite{ExUniCoCoGeoMetNonsimCLEGas}, and of the conformally covariant resistance form, and hence of the canonical Brownian motion, on the non-simple CLE gasket \cite{ExUniBMNonsimCLEGas}.

In the case of the metric on the loops of a $\CLE_4$, one also has a natural embedding into the complex plane.  However, there are significant challenges in making use of this embedding.  First, it is more difficult to compare the metric with the Euclidean metric and its variants.  Second, unlike in the works mentioned above, the metric transforms in a conformally invariant rather than conformally covariant manner.  As a result, the arguments used in the aforementioned works to establish the bi-Lipschitz equivalence of conditionally independent samples do not carry over, leading us to take a different approach.  In order to circumvent these issues, this paper instead borrows ideas from both frameworks: from the latter, we make use of the independence properties of $\CLE_4$ and its coupling with the GFF, and as in the former \cite{LG13, Mie13, StabCarGas}, we prove that geodesics can be approximated by geodesics to the domain boundary, which we shall also refer to as the \emph{root}. This was one of the key ingredients in the identification of the scaling limit of quadrangulations \cite{LG13, Mie13} and of random maps with large faces \cite{StabCarGas}. In the same way, we expect that this property will be crucial in order to obtain the scaling limit of $3/2$-stable maps.

\subsection{Setup}\label{subsec:setup}

Let $U \subsetneq \BC$ be a deterministic simply connected domain. Let $\Gamma_U$ be a non-nested $\CLE_4$ in $U$. We shall write $\Upsilon_U$ for the carpet of $\Gamma_U$, which is the closed subset of $\overline{U}$ obtained by removing the interiors of the loops of $\Gamma_U$ (we choose to include $\partial U$ in $\Upsilon_U$ so that the paths that we will define later can start from $\partial U$). For each open subset $V \subseteq U$, we shall write 
\begin{equation}\label{eq:def-V-star}
    V^\star \defeq V \setminus \overline{\bigcup_{\SCL \in \Gamma_U : \SCL \not\subseteq V} \mathop{\mathrm{int}}(\SCL)}.
\end{equation}
We shall also write $\Gamma_U|_V \defeq \{\SCL \in \Gamma_U : \SCL \subseteq V\}$ for the collection of loops of $\Gamma_U$ that are contained in $V$. For deterministic $z \in U$, we shall set $\SCL(z)$ to be the unique loop of $\Gamma_U$ that surrounds $z$. Note that such a loop exists a.s.

Let $D$ be a metric on $\Gamma_U$. For $\SCL \in \Gamma_U$ and $t \ge 0$, we shall write 
\begin{equation*}
    \SCB_t(\SCL; D) \quad \text{(resp.\ } \SCB_t^-(\SCL; D)\text{)}
\end{equation*}
(or simply $\SCB_t(\SCL)$ (resp.\ $\SCB_t^-(\SCL)$), when there is no danger of confusion) for the closure in $\overline{U}$ of the union of the domains surrounded by $\SCL^\prime$ for $\SCL^\prime \in \Gamma_U$ with $D(\SCL, \SCL^\prime) \le t$ (resp.\ $D(\SCL, \SCL^\prime) < t$). Note that $\SCB_0^-(\SCL; D) = \emptyset$. Note that $\SCB_t(\SCL; D) = \SCB_t^-(\SCL; D)$ unless there exists $\SCL^\prime \in \Gamma_U$ with $D(\SCL, \SCL^\prime) = t$. 

We will also talk about the distance between a loop and a subset of $\overline{U}$, for example the boundary of the domain. For that reason, we define the distance from a loop to a subset of $\overline{U}$ and the metric ball from a subset as follows. Let $A,B \subseteq \overline U$. Then we shall write
\begin{itemize}
    \item $D(\SCL, A) \defeq \inf\{t \ge 0 : \SCB_t(\SCL; D) \cap A \neq \emptyset\}$ for all $\SCL \in \Gamma_U$;
    \item $\SCB_t(A; D)$ (resp.\ $\SCB_t^-(A; D)$) for the closure of the union of $A$ and the domains surrounded by $\SCL$ for $\SCL \in \Gamma_U$ with $D(\SCL, A) \le t$ (resp.\ $D(\SCL, A) < t$) for all $t \ge 0$, except for $\SCB_0^-(A; D)\defeq \emptyset$;
    \item $D(A, B) \defeq \inf\{t \ge 0 : \SCB_t(A; D) \cap B \neq \emptyset\}$
\end{itemize}
(with the convention that $\inf\emptyset = \infty$).

Next, let us introduce a notion of continuous path for the metric space $(\Gamma_U, D)$. Since the loops of $\Gamma_U$ are drawn in $U$, it is convenient to define such a path as a continuous path taking its values in $\Upsilon_U$. Intuitively, this path should take its values in the set of loops. However, since the loops of the $\CLE_4$ do not intersect, this is not possible. This is why we write the following definition. Let $P \colon [0, 1] \to \Upsilon_U$ be a continuous path. Whenever we regard $P$ as a subset of $\overline U$, we shall identify it with its image $P([0,1])$. Then we say that $P$ is \emph{admissible} if $\bigl(\bigcup_{\SCL \in \Gamma_U} \SCL\bigr) \cap P$ is dense in $P$. We shall refer to
\begin{equation*}
    \len(P; D) \defeq \sup_{\SCL_0, \SCL_1, \ldots, \SCL_n} \sum_{j = 1}^n D(\SCL_{j - 1}, \SCL_j),
\end{equation*}
as the \emph{$D$-length} of $P$, where the supremum is taken over all sequences $\SCL_0, \ldots, \SCL_n$ of (not necessarily distinct) loops of $\Gamma_U$ that $P$ hits in chronological order. Let $A, B \subseteq X \subseteq \overline U$. Then, we shall write 
\begin{equation*}
    D(A, B; X) \defeq \inf_{\substack{P \subseteq X\\P \colon A \to B}} \len(P; D),
\end{equation*}
where $P$ ranges over all admissible paths in $X$ from $A$ to $B$ (with the convention that $\inf\emptyset = \infty$). The mapping $D(\bullet, \bullet;X)$ is called the internal metric in $X$; it may take the value $+\infty$ if no admissible path in $X$ connects $A$ to $B$. Let us emphasize that $D(A,B)$ and $D(A,B;\overline U)$ are defined differently: the former is given in terms of metric balls and the latter as an infimum of $D$-lengths of admissible paths.

We shall write $U_\BQ \defeq U \cap \BQ^2$. By abuse of notation, we shall also denote by $D$ the mapping
\begin{equation*}
    D \colon U_\BQ \times U_\BQ \to \BR, \quad (x, y) \mapsto D(\SCL(x), \SCL(y)). 
\end{equation*}
We shall equip $\BR^{U_\BQ \times U_\BQ}$ with the product topology. Note that $\BR^{U_\BQ \times U_\BQ}$ is a Polish space. 

We denote by $\mathscr{Loop}(U)$ the metric space consisting of loops in $U$, and equipped with the metric given by
\begin{equation*}
    d_{\mathscr{Loop}(U)}(\SCL_1, \SCL_2) \defeq \inf_{P_1, P_2} \sup_{t \in \BR/\BZ} \left\lvert P_1(t) - P_2(t)\right\rvert, 
\end{equation*}
where $P_1 \colon \BR/\BZ \to \SCL_1$ (resp.\ $P_2 \colon \BR/\BZ \to \SCL_2$) ranges over all parameterizations of $\SCL_1$ (resp.\ $\SCL_2$). We shall write $\mathscr{LoopEns}(U)$ for the metric space consisting of all closed discrete subsets of $\mathscr{Loop}(U)$, and equipped with the Hausdorff metric induced by $d_{\mathscr{Loop}(U)}$.

\subsection{Main results}
In \cite{kkmt2026cle4_part1}, we proved that subsequential limits of the $\CLE_\kappa$ together with its rescaled graph distance as $\kappa \downarrow 4$ give rise to a weak geodesic $\CLE_4$ metric coupling. Then, in \cite{kkmt2026cle4_part2}, we proved that every weak geodesic $\CLE_4$ metric coupling in fact admits geodesics; more precisely, any such coupling is a geodesic metric coupling according to the following definition.
\begin{definition}\label{def:weak_axioms}
	We define a \emph{geodesic CLE$_4$ metric coupling} to be a family of couplings
	\begin{equation*}
		D = \left\{\left( \Gamma_U, D_{\Gamma_U}^U \right)\right\}_U
	\end{equation*}
	where $U \subsetneq \BC$ ranges over all simply connected domains and $\Gamma_U$ is a non-nested CLE$_4$ in $U$, such that the following conditions are satisfied for every such $U$ and all deterministic choices of the auxiliary data below:
	\begin{enumerate}[label=(\Roman*), ref=\Roman*]
		\item\label{it:weak_axiom_geodesic} {\bfseries (Geodesic metric)} $D_{\Gamma_U}^U$ is a.s.\ a metric on $\Gamma_U$ such that the following hold.
		\begin{enumerate}[label=(\roman*), ref=\roman*]
			\item\label{it:axiom_geodesic_path} For each $\SCL_1, \SCL_2 \in \Gamma_U$, there exists an admissible path $P$ (which we refer to as a \emph{$D_{\Gamma_U}^U$-geodesic}) connecting $\SCL_1$ and $\SCL_2$ such that $D_{\Gamma_U}^U(\SCL_1, \SCL_2) = \len(P; D_{\Gamma_U}^U)$. The same is true with a deterministic connected arc (resp.\ two deterministic and disjoint connected arcs) of $\partial U$ in place of $\SCL_1$ (resp.\ $\SCL_1$ and $\SCL_2$).
			\item\label{it:axiom_geodesic_cadlag} Let $\phi\colon U \to \BD$ be a conformal mapping. For each $\SCL \in \Gamma_U$, the process $(\phi(\SCB_t(\SCL; D^U_{\Gamma_U}){)})_{t\ge 0}$ is c\`adl\`ag for the Hausdorff distance (with the convention of \Cref{subsec:setup}), and the same is true with a deterministic connected arc of $\partial U$ in place of $\SCL$.
			\item\label{it:axiom_geodesic_twoball} For all $t, s\ge 0$, if $\SCB_t(\SCL_1; D^U_{\Gamma_U}) \cap \SCB_s(\SCL_2; D^U_{\Gamma_U}) \neq \emptyset$, then $D^U_{\Gamma_U} (\SCL_1, \SCL_2) \le t+s$. The same is true with a deterministic connected arc (resp.\ two deterministic and disjoint connected arcs) of $\partial U$ in place of $\SCL_1$ (resp.\ $\SCL_1$ and $\SCL_2$).
		\end{enumerate}
		Note that $\len(P; D_{\Gamma_U}^U) \ge D_{\Gamma_U}^U(\SCL_1,\SCL_2)$ for every admissible path $P$ connecting $\SCL_1$ and $\SCL_2$, by the triangle inequality along the chains in the definition of $\len$; a path attaining equality is thus length-minimizing, which is what the term geodesic refers to. 
		\item\label{it:weak_axiom_locality} {\bfseries (Locality)} Let $V \subset U$ be a deterministic simply connected subdomain. Let $I \subset \partial U$ be a deterministic connected arc. Let $\{V_j\}_j$ be the connected components of $V^\star$. Then, conditionally on the $\sigma$-algebra generated by
		\begin{multline}\label{eq:weak_axiom_locality}
			\bigl\{\SCL \in \Gamma_U : \SCL \not\subset V\bigr\}, \quad \bigl\{\SCB_t(\SCL; D_{\Gamma_U}^U) : \SCL \in \Gamma_U, \ \SCL \not\subset V, \ t \in [0, D_{\Gamma_U}^U(\SCL, \partial V^\star)]\bigr\}, \\
			\text{and} \quad \bigl\{\SCB_t(I; D_{\Gamma_U}^U) : t \in [0, D_{\Gamma_U}^U(I, \partial V^\star)]\bigr\},
		\end{multline}
		the $(\Gamma_U|_{V_j}, D_{\Gamma_U}^U(\bullet, \bullet; V_j))$'s are independent and their conditional laws are those of $(\Gamma_{V_j}, D_{\Gamma_{V_j}}^{V_j})$, respectively. 
		\item\label{it:weak_axiom_conformal_invariance} {\bfseries (Conformal invariance)} Let $\phi \colon U \to \phi(U)$ be a deterministic conformal mapping. Then,
		\begin{equation*}
			\left(\phi(\Gamma_U), \left( D_{\Gamma_U}^U\left(\phi^{-1}(\SCL_1), \phi^{-1}(\SCL_2) \right) \right)_{\SCL_1, \SCL_2 \in \phi(\Gamma_U)} \right)
			\overset{(\mathrm{d})}{=} \left( \Gamma_{\phi(U)}, \left(D_{\Gamma_{\phi(U)}}^{\phi(U)}(\SCL_1, \SCL_2) \right)_{\SCL_1, \SCL_2 \in \Gamma_{\phi(U)}} \right).
		\end{equation*}
		The equality in law is understood as in \Cref{def:weak_axioms}.\item\label{it:axiom_uniform_exploration_coupling} {\bfseries (Uniform exploration)} The collection $\{(\SCL, D_{\Gamma_U}^U(\SCL, \partial U))\}_{\SCL \in \Gamma_U}$ has the law of a uniform exploration of $\Gamma_U$ (cf.~\cite{CoInCLEExpl}), in the sense of \Cref{def:weak_axioms}.
	\end{enumerate}
\end{definition}
Note that, in the above definition, the metric $D^U_{\Gamma_U}$ is not required to be a function of $\Gamma_U$. We will actually show that any geodesic $\CLE_4$ metric coupling can be described as a uniquely defined measurable function of the $\CLE_4$. We thus introduce the notion of a strong $\CLE_4$ metric, which is essentially a geodesic $\CLE_4$ metric coupling such that the metric $D^U_{\Gamma_U}$ is a measurable function of $\Gamma_U$. The coupling $(\Gamma_U, D^U_{\Gamma_U})$ can then be recovered from $\Gamma_U$ together with the measurable map $D^U$. More precisely, let us define a strong $\CLE_4$ metric in the following way.
\begin{definition}\label{def:axioms}
    By a (strong) $\CLE_4$ metric, we mean a family of \textbf{deterministic} measurable mappings 
    \begin{equation*}
        D = \left\{D^U \colon \mathscr{LoopEns}(U) \to \BR^{U_\BQ \times U_\BQ},\quad  \Gamma_U \mapsto D_{\Gamma_U}^U\right\}_U,
    \end{equation*}
    where $U \subsetneq \BC$ ranges over all simply connected domains, such that the following conditions are satisfied whenever $U \subsetneq \BC$ is a deterministic simply connected domain and $\Gamma_U$ is a non-nested $\CLE_4$ in $U$:
    \begin{enumerate}[label=(\Roman*), ref=\Roman*]
        \item\label{it:axiom_geodesic} {\bfseries (Geodesic metric)} The mapping $D_{\Gamma_U}^U$ a.s.\ determines a metric on $\Gamma_U$ such that for each $\SCL_1, \SCL_2 \in \Gamma_U$, there exists an admissible path $P$ (which we shall refer to as a \emph{$D_{\Gamma_U}^U$-geodesic}) connecting $\SCL_1$ and $\SCL_2$ such that $D_{\Gamma_U}^U(\SCL_1, \SCL_2) = \len(P; D_{\Gamma_U}^U)$. Moreover, for all conformal mapping $\phi \colon U \to \BD$, the process $(\phi(\SCB_t(\SCL; D^U_{\Gamma_U})))_{t\ge 0}$ is c\`adl\`ag with respect to the Hausdorff distance, and, for all $t, s\ge 0$, if $\SCB_t(\SCL_1; D^U_{\Gamma_U}) \cap \SCB_s(\SCL_2; D^U_{\Gamma_U}) \neq \emptyset$, then $D^U_{\Gamma_U} (\SCL_1, \SCL_2) \le t+s$. The same is true with a deterministic connected arc (resp.\ two deterministic disjoint connected arcs) of $\partial U$ in place of $\SCL_1$ (resp.\ $\SCL_1$ and $\SCL_2$).
        \item\label{it:axiom_locality} {\bfseries (Locality)} Let $V \subseteq U$ be a deterministic simply connected subdomain. Then, a.s., for each connected component $V_j$ of $V^\star$, 
        \begin{equation*}
            D_{\Gamma_U}^U(\SCL_1, \SCL_2; V_j) = D_{\Gamma_U|_{V_j}}^{V_j}(\SCL_1, \SCL_2), \quad \forall \SCL_1, \SCL_2 \in \Gamma_U|_{V_j}.
        \end{equation*}
        \item\label{it:axiom_conformal_invariance} {\bfseries (Conformal invariance)} Let $\phi \colon U \to \phi(U)$ be a deterministic conformal mapping. Then, a.s., 
        \begin{equation*}
            D_{\Gamma_U}^U(\SCL_1, \SCL_2) = D_{\phi(\Gamma_U)}^{\phi(U)}(\phi(\SCL_1), \phi(\SCL_2)), \quad \forall \SCL_1, \SCL_2 \in \Gamma_U. 
        \end{equation*}

        \item\label{it:axiom_uniform_exploration} {\bfseries (Uniform exploration)} The collection $\{(\SCL, D_{\Gamma_U}^U(\SCL, \partial U))\}_{\SCL \in \Gamma_U}$ has the law of a uniform exploration of $\Gamma_U$ (cf.\ \cite{CoInCLEExpl}).
    \end{enumerate}
\end{definition}

\begin{remark}
	We record the following two properties.
    \begin{itemize}
    	\item In Axiom~\eqref{it:axiom_locality}, the internal metric $D_{\Gamma_U}^U(\bullet, \bullet; V_j)$ is a.s.\ determined by $V_j$ and $\Gamma_U|_{V_j}$. 
        
        \item By Axiom~\eqref{it:axiom_conformal_invariance}, $D$ is uniquely determined by $D^U$ for any single fixed simply connected domain $U \subsetneq \BC$. 
    \end{itemize}
\end{remark}

Next, we state the main results of the paper. We prove that any geodesic $\CLE_4$ metric coupling in the sense of \Cref{def:weak_axioms} coincides a.s.\ with a strong $\CLE_4$ metric in the sense of \Cref{def:axioms} and that such a strong $\CLE_4$ metric is unique.

\begin{theorem}[Uniqueness of the metric]\label{thm:uniqueness}
	\begin{enumerate}
		\item\label{thm:uniqueness-existence} Let $D$ be a geodesic $\CLE_4$ metric coupling. Then there is a strong $\CLE_4$ metric $\widetilde{D}$ such that for every simply connected open subset $U\subsetneq \BC$, $D_{\Gamma_U}^U = \widetilde D_{\Gamma_U}^U$ a.s.
		\item\label{thm:uniqueness-uniqueness} If $D$ and $\widetilde D$ are two strong $\CLE_4$ metrics, then for each deterministic simply connected domain $U \subsetneq \BC$ we have $D_{\Gamma_U}^U = \widetilde D_{\Gamma_U}^U$ a.s.
	\end{enumerate}
\end{theorem}

One of the main consequences of the existence of the strong $\CLE_4$ metric is the following corollary.  It states that the uniform $\CLE_4$ exploration introduced in  \cite{CoInCLEExpl} is a.s.\ determined by the associated $\CLE_4$, which answers \cite[Question~3]{CoInCLEExpl}.  Moreover, it states that if the $\CLE_4$ is coupled with a family of $\BCLE_4(\rho)$ processes for $\rho \in (-2,0)$ as in Section~\ref{subsec:gff_labelled_cle_4} (see also \cite[Section~2.7]{kkmt2026cle4_part1} for more details), then the $\BCLE_4(\rho)$'s are a.s.\ determined by the associated $\CLE_4$.

\begin{corollary}[Measurability]\label{thm:uniform-exploration-measurability}
    Let $U \subsetneq \BC$ be a deterministic simply connected domain. Then:
    \begin{enumerate}
        \item\label{thm:uniform-exploration-measurability-0} The uniform exploration of $\Gamma_U$ is a.s.\ determined by $\Gamma_U$.
        \item\label{thm:uniform-exploration-measurability-1} Let $\Gamma_U$ and a family of $\BCLE_4(\rho)$'s for $\rho \in (-2, 0)$ be coupled as in Section~\ref{subsec:gff_labelled_cle_4}. Then the $\BCLE_4(\rho)$'s are a.s.\ determined by $\Gamma_U$. 
    \end{enumerate}
\end{corollary}

Another main consequence is the uniqueness of the scaling limit of the $\CLE_\kappa$ equipped with the graph distance. The existence of subsequential scaling limits was obtained in  \cite{kkmt2026cle4_part1}. Let $\kappa \in (4, 8)$. Let $\Gamma_U^\kappa$ be a non-nested $\CLE_\kappa$ in $U$. We shall write $D_{\Gamma_U^\kappa}^U \colon \Gamma_U^\kappa \times \Gamma_U^\kappa \to \BR$ for the graph distance associated with $\Gamma_U^\kappa$, where $\SCL_1, \SCL_2 \in \Gamma_U^\kappa$ are adjacent if and only if $\SCL_1 \cap \SCL_2 \neq \emptyset$. For deterministic $z \in U$, let $\SCL^\kappa(z)$ be defined in the same manner as $\SCL(z)$ but with $\Gamma_U^\kappa$ in place of $\Gamma_U$. We shall write
	\begin{equation*}
		1/\ka_\kappa \defeq \BP\lbrack\SCL^\kappa(z) \cap \partial U \neq \emptyset\rbrack = 1 - \frac{\sin(\pi(\kappa/4 + 8/\kappa))}{\sin(\pi(\kappa/4 - 1))}
	\end{equation*}
	(the above probability was computed in \cite{BdryTouchingNonsimCLE}; its exact form will not be important for this work).

\begin{theorem}[Convergence of non-simple CLE graph distance]\label{thm:convergence_nonsimple_cle}
	The following is true for every simply connected domain $U \subsetneq \BC$.  The random metrics $\ka_\kappa^{-1} D_{\Gamma_U^\kappa}^U$ converge in law, with respect to the product topology of $\BR^{U_\BQ \times U_\BQ}$, to $D_{\Gamma_U}^U$ as $\kappa \downarrow 4$, where $D$ is the strong $\CLE_4$ metric of \Cref{thm:uniqueness}.
\end{theorem}

Finally, we state the last main results of the paper concerning the behavior of geodesics of the strong $\CLE_4$ metric.  The first of these states that, a.s., any two $\CLE_4$ loops are connected by a unique geodesic. The same holds for geodesics between a deterministic boundary arc and a loop, and between two deterministic boundary arcs.

\begin{theorem}[Uniqueness of geodesics]\label{thm:uniqueness-of-geodesics}
    Let $U \subsetneq \BC$ be a deterministic simply connected domain. Then, a.s., for each $\SCL_1, \SCL_2 \in \Gamma_U$, the $D_{\Gamma_U}^U$-geodesic connecting $\SCL_1$ and $\SCL_2$ is unique in the following sense: if $P_1$ and $P_2$ are $D_{\Gamma_U}^U$-geodesics connecting $\SCL_1$ and $\SCL_2$, then they visit the same set of loops, that is, $\{\SCL \in \Gamma_U : \SCL \cap P_1 \neq \emptyset\} = \{\SCL \in \Gamma_U : \SCL \cap P_2 \neq \emptyset\}$. The same is true with a deterministic connected arc (resp.\ two deterministic disjoint connected arcs) of $\partial U$ in place of $\SCL_1$ (resp.\ $\SCL_1$ and $\SCL_2$).
\end{theorem}

Another result, which we obtain in the proof of \Cref{thm:uniqueness}, states that geodesics between two loops can be approximated using geodesics to the boundary. More precisely, the distance between two loops $\SCL_1$ and $\SCL_2$ can be approximated using a sequence of loops $\SCL_1=\SCL'_0, \ldots, \SCL'_n=\SCL_2$ and some paths $P_i$ between $\SCL_{i-1}'$ and $\SCL_i'$ for each $1 \le i \le n$ such that for any loop $\SCL$ visited by $P_i$ at some time $t \in [0,1]$, there is a geodesic from $\SCL$ to $\partial U$ containing either the time-reversal of $P_i\vert_{[0,t]}$ or the path $P_i\vert_{[t,1]}$. See \Cref{thm:main_result} for more details.

\subsection{Outline}

Next, we will give an outline of the content of the rest of the paper.

\smallskip
\noindent\emph{\Cref{section:preliminaries}: Preliminaries.} We give some background on Gaussian free fields (GFF), Schramm-Loewner evolutions (SLE) with force points, conformal loop ensembles (CLE), level lines and two-valued sets of the GFF, the labeled $\CLE_4$, couplings between the GFF and $\CLE_4$, the Carathéodory topology, and some useful estimates; in particular, four-arm estimates, estimates for extremal distance, and estimates for the $\CLE_4$ metric.

\smallskip
\noindent\emph{\Cref{sec:uniqueness_and_measurability_of_the_metric}: Uniqueness of the metric: genealogy, main result, and proof strategy.} In this section, we define a genealogy on the loops of the $\CLE_4$ using only the uniform exploration. Then, we define a distance $\widetilde{D}$ on the $\CLE_4$ that only relies on pieces of geodesics to the boundary. More precisely, the metric $\widetilde{D}$ is obtained by going up and down in the tree of geodesics to the boundary. We then state two key results, namely \Cref{thm:main_result,prop D tilde does not depend on the metric}. The first one states that the metric of any geodesic $\CLE_4$ metric coupling coincides a.s.\ with $\widetilde{D}$, and the second one shows that the geodesics to the boundary (which we also call the root) are determined by the above-mentioned genealogy: a loop $\SCL(y)$ is a descendant of a loop $\SCL(x)$ in the genealogy if and only if $\SCL(x)$ lies on a geodesic from $\partial U$ to $\SCL(y)$. Combined together, \Cref{thm:main_result,prop D tilde does not depend on the metric} entail that any geodesic $\CLE_4$ metric coupling is uniquely determined by the associated uniform exploration (see \Cref{thm metric determined by the exploration}). We then explain the proof strategy in order to obtain \Cref{thm:main_result,prop D tilde does not depend on the metric}.

\smallskip
\noindent\emph{\Cref{sec:crossings_of_annuli}: Crossings of annuli.} We construct an exploration of the crossings of an annulus by $\CLE_4$ loops, and we deduce from it in \Cref{lem:crossings_determined_by_the_field} that these crossings are a.s.\ determined by the restriction of the GFF to the annulus. We consider the remaining connected components of the annulus, i.e.,~the components which are left once the crossings of the annulus and the excursions into it from its inner and outer boundaries have been removed. Next, we prove in \Cref{prop:conformal_rectangles} that, conditionally on these components and with conditional probability at least $1-2^{-bR}$, where $b>0$ is a universal constant and $2^{-R}$ (resp.\ $1/2$) is the inner (resp.\ outer) radius of the annulus, there exist conformal rectangles included in these connected components whose left and right sides are contained in crossings of the annulus made by loops and whose $\CLE_4$-metric distance from top to bottom is larger than their $\CLE_4$-metric distance from left to right.

\smallskip
\noindent\emph{\Cref{sec:good_annuli}: Good annuli.} We establish that with very high probability, there are many annuli which are \emph{good} in the sense that they are crossed by exactly two loops of the $\CLE_4$ and that any geodesic from the inner boundary of the annulus to $\partial \BD$ has to intersect one of these two loops. The conformal rectangles constructed in the previous section are the main ingredient in the proof of the second property.

\smallskip
\noindent\emph{\Cref{sec:covering_geodesics-by_good_annuli}: Covering geodesics by good annuli.} This section is devoted to the proof of  \Cref{thm:main_result,prop D tilde does not depend on the metric}. The main idea is to cover geodesics using the good annuli obtained in the previous section and argue, using these annuli, that geodesics can be approximated by concatenations of pieces of geodesics to the boundary.

\smallskip
\noindent\emph{\Cref{sec:measurability}: Measurability.} We prove in \Cref{lem:measurability} that any geodesic $\CLE_4$ metric coupling, and hence also the associated uniform exploration, is a measurable function of the $\CLE_4$. This completes the proof of \Cref{thm:uniqueness}, and it also completes the proofs of \Cref{thm:uniform-exploration-measurability} and, together with the main results of \cite{kkmt2026cle4_part1, kkmt2026cle4_part2}, of \Cref{thm:convergence_nonsimple_cle}.

\smallskip
\noindent\emph{\Cref{section:uniqueness-of-geodesics}: Uniqueness of geodesics.} We will prove \Cref{thm:uniqueness-of-geodesics}, which is restated as \Cref{thm:uniqueness-geodesics}.

\smallskip
\noindent\emph{Appendix~\ref{sec:cle_4_four_arm_exponents}: Nested $\CLE_4$ four-arm exponents.} We will prove estimates on four-arm events for nested $\CLE_4$, extending the results in \cite{U2C4ASimCLE}. These estimates are used in \Cref{sec:good_annuli} to prove \Cref{prop:bounded_cle_crossings}, which bounds the number of crossings of an annulus by loops of the nested $\CLE_4$.

\subsection*{Acknowledgements}

E.K.\ acknowledges the support of a Research Fellowship from Emmanuel College, Cambridge. J.M.\ received support from the ERC Consolidator Grant ARPF (Horizon Europe UKRI G120614). Y.T.\ was supported by a Cambridge International Scholarship from the Cambridge Trust.  K.K. was supported by the Simons Collaboration Grant
\emph{Probabilistic Paths to Quantum Field Theory}. We thank Wendelin Werner for stimulating discussions at an early stage of this work. E.K.\ also thanks Juhan Aru and Ellen Powell for insightful discussions on this topic before starting this work.

\section{Preliminaries}\label{section:preliminaries}

\subsection{Gaussian free fields}
\label{subsec:gff}

We briefly recall some elementary properties of the Gaussian free field (GFF); see \cite{ss2007gff} for a review. Let $G \subseteq \BC$ be a simply connected domain with harmonically non-trivial boundary and let $H_0^1(G)$ denote the Hilbert space closure of $C_0^{\infty}(G)$ with respect to the Dirichlet inner product defined by setting, for all $f, g \in C_0^{\infty}(G)$,
\begin{align*}
(f,g)_{\nabla} \defeq\frac{1}{2\pi} \int_G \nabla f(z) \cdot \nabla g(z) dz.
\end{align*}
The zero-boundary GFF $h$ is the random distribution defined by
\begin{align}\label{eqn:gff_series}
h = \sum_{n \geq 1} \alpha_n \phi_n,
\end{align}
where $(\phi_n)_{n \geq 1}$ is a $(\cdot , \cdot)_{\nabla}$-orthonormal basis of $H_0^1(G)$ and $(\alpha_n)_{n \geq 1}$ is a sequence of independent $\mathcal{N}(0,1)$-distributed random variables. The law of $h$ does not depend on the choice of orthonormal basis $(\phi_n)_{n \geq 1}$.

A GFF on $G$ is not a random function but rather a random variable taking values in the space of distributions (also called generalized functions) on $G$; a.s., the series in~\eqref{eqn:gff_series} does not converge in $H_0^1(G)$.  Nevertheless, the above series converges a.s.\ in the space of distributions on $G$. We say that a random distribution $h$ on $G$ is a GFF with boundary conditions given by a function $f$ defined on $\partial G$ if $h$ can be expressed as the sum of a zero-boundary GFF on $G$ and the harmonic extension of $f$ from $\partial G$ to $G$.

The following are two important properties of the GFF that we will frequently use.
	\begin{itemize}
		\item \emph{The Markov property.} If $h$ is a zero-boundary GFF on $G$ and $V \subseteq G$ is open, then we can write $h = h_{V}^0 + h_{V}$, where $h_{V}^0$ is a zero-boundary GFF on $V$ (that is, independent zero-boundary GFFs, one in each connected component of $V$) and $h_{V}$ is a distribution on $G$ which is harmonic on $V$. Moreover, $h_{V}^0$ and $h_{V}$ are independent.
		\item \emph{Conformal invariance.} If $\phi \colon G \to G'$ is a conformal transformation and $h$ is a GFF on $G$, then $h \circ \phi^{-1}$ is a GFF on $G'$. 
	\end{itemize}
	There is also a variant of the Markov property where one thinks of conditioning on the values of the GFF on a random set, in the same vein as for the strong Markov property and stopping times for Brownian motion. Specifically, a random closed set $A \subseteq G$ coupled with $h$ is said to be \emph{local} for $h$ \cite{ContourLineGFF} if we can write $h = h_A^0 + h_A$ where $h_A$ is a distribution on $G$ which is harmonic on $G \setminus A$ and, given $(A,h_A)$, the field $h_A^0$ has the law of a GFF on $G \setminus A$ with zero boundary conditions, independently in each connected component of $G \setminus A$.

\subsection{Chordal $\SLE_{\kappa}(\underline{\rho})$ processes}
\label{subsec:chordal_sle}

The chordal $\SLE_\kappa(\underline{\rho})$ processes are variants of $\SLE_\kappa$ where one keeps track of extra marked points \cite{CoRestr}. To define them, we fix $\kappa>0$ and integers $\ell, r \geq 0$, and we write $[a,b]_\BZ \defeq [a,b] \cap \BZ$ for the set of integers in an interval $[a,b]$. We let $\underline{x}_L = (x_{\ell,L},\ldots,x_{1,L})$ and $\underline{x}_R = (x_{1,R},\ldots,x_{r,R})$, where $x_{\ell,L} < \cdots < x_{1,L} \leq 0^-$ and $0^+ \leq$ $x_{1,R} < \cdots < x_{r,R}$ (we allow the innermost force points to be located at $0^-$ and $0^+$, i.e.\ immediately to the left and to the right of the seed of the curve), and let $\underline{\rho}_L = (\rho_{1,L},\ldots,\rho_{\ell,L})$ and $\underline{\rho}_R = (\rho_{1,R},\ldots,\rho_{r,R})$, where $\rho_{j,q} \in \BR$ for $q \in \{L,R\}$ and $j \in [1, N_q]_\BZ$ with $N_L = \ell$, $N_R = r$. Let $(g_t)$ denote the solution to the ODE
\begin{align}\label{eqn:loewner_ode}
	\partial_t g_t(z) = \frac{2}{g_t(z) - W_t}, \quad g_0(z) = z \in \BH,
\end{align}
where $W$ is a solution to
\begin{align}\label{eqn:multiforce_point_sde}
	&\mathrm{d} W_t = \sum_{j=1}^{\ell} \frac{\rho_{j,L}}{W_t - V_t^{j,L}} \mathrm{d}t + \sum_{j=1}^r \frac{\rho_{j,R}}{W_t - V_t^{j,R}} \mathrm{d}t + \sqrt{\kappa} \mathrm{d}B_t, \quad W_0 = 0,\notag \\
	&\mathrm{d}V_t^{j,q} = \frac{2}{V_t^{j,q} - W_t} \mathrm{d}t, \quad V_0^{j,q} = x_{j,q}, \  j \in [1, N_q]_\BZ, \ q \in \{L,R\},
\end{align}
where $B$ is a standard Brownian motion. For each $z \in \BH$, the solution $g_t(z)$ to~\eqref{eqn:loewner_ode} is defined up to the time $\tau(z)$ at which $g_t(z) - W_t$ first hits $0$, and we write $\BH_t \defeq \{z \in \BH : \tau(z) > t\}$ for the domain of $g_t$.

For any value of $\underline{\rho}_L$ and $\underline{\rho}_R$, it is clear that~\eqref{eqn:multiforce_point_sde} has a unique strong solution until the first time $t$ such that $W_t=V^{j,q}_t$ for some $q \in \{L,R\}$ and $j \in [1, N_q]_\BZ$. Actually, it was shown in \cite{IG1,LevelLineGFFI} that, when $\sum_{j=1}^k \rho_{j,L} > -2$ for all $k \in [1,\ell]_\BZ$ and $\sum_{j=1}^k \rho_{j,R} > -2$ for all $k \in [1,r]_\BZ$, there exists a solution to~\eqref{eqn:multiforce_point_sde} defined for all times $t\geq 0$, and for this solution the set of $t \geq 0$ for which $W_t=V^{j,q}_t$ for some $q \in \{L,R\}$ and $j \in [1, N_q]_\BZ$ a.s.\ has zero Lebesgue measure. The uniqueness in law of such a solution is also proved in \cite{IG1}. Moreover, it was shown in \cite{IG1,LevelLineGFFI} that there a.s.\ exists a continuous curve $\eta$ such that the domain $\BH_t$ of $g_t$ is given by the unbounded connected component of $\BH \setminus \eta([0,t])$, for all $t\ge 0$. The curve $\eta$ is called a chordal $\SLE_{\kappa}(\underline{\rho})$ in $\BH$ from~$0$ to~$\infty$. (Processes with $\sum_{j=1}^k \rho_{j,q} \le  -2$ are not produced by~\eqref{eqn:multiforce_point_sde} and have to be constructed by other means; their continuity was later proved in \cite{CLEPerc,ms2019lightcone}.)

A chordal $\SLE_{\kappa}(\underline{\rho})$ process in a proper simply connected domain $G \subseteq \BC$ is defined as the image of a chordal $\SLE_{\kappa}(\underline{\rho})$ on $\BH$ under a conformal transformation mapping $\BH$ onto $G$ which sends the starting point $0$ to the desired starting point on $\partial G$ and the target point $\infty$ to the desired target point on $\partial G$. The force points on $\partial G$ are then defined as the images of the corresponding force points on $\partial \BH$. Since the conformal automorphisms of $\BH$ fixing $0$ and $\infty$ are the scalings $z \mapsto az$ with $a > 0$, and the law of a chordal $\SLE_\kappa(\underline{\rho})$ is invariant under such scalings modulo time parameterization, this does not depend on the choice of conformal transformation.

\subsection{Conformal loop ensembles}
\label{subsec:cle}

Now, we will briefly review CLE and refer the reader to \cite{CLE,CLEPerc,TreeCLE} for more details.

CLE in $\BD$ is a countable collection $\Gamma$ of non-nested loops $(\gamma_j)_{j \in J}$ in $\overline{\BD}$ that satisfies the following properties.
\begin{enumerate}[label=(\Roman*), ref=\Roman*]
	\item\label{it:cle_conf_inv} \textbf{Conformal invariance}: For any M\"obius transformation $\phi \colon \BD \to \BD$, the laws of $\Gamma$ and $\phi(\Gamma)$ are the same. This makes it possible to define CLE on any proper simply connected domain $G \subseteq \BC$ as the image of a CLE on $\BD$ under a conformal transformation mapping $\BD$ onto $G$.
	\item\label{it:cle_dmp} \textbf{Domain Markov property}: For any non-trivial simply connected domain $W \subseteq \BD$, we let $W^\star$ be defined as in~\eqref{eq:def-V-star}, with $\BD$, $\Gamma$ and $W$ in place of $U$, $\Gamma_U$ and $V$; that is, $W^\star$ is the set obtained by removing from $W$ the closure of the union of the interiors of the loops of $\Gamma$ which are not contained in $W$. Then, conditionally on the loops of $\Gamma$ which intersect $\BD \setminus W$ (which, we note, determines $W^\star$), the law of $\Gamma$ restricted to each component $V$ of $W^\star$ is that of a CLE in $V$. Moreover, the restrictions of $\Gamma$ to distinct components of $W^\star$ are independent.
\end{enumerate}

It was shown in \cite{CLE,TreeCLE} that for each $\kappa \in (8/3,8)$ there exists a $\CLE_\kappa$ where the loops can be constructed using a branching variant of $\SLE_\kappa$ (specifically, a tree built out of a collection of $\SLE_\kappa(\kappa-6)$ processes starting from the same point). The loops in the $\CLE_\kappa$ are then $\SLE_\kappa$-type curves. From the perspective of the branching tree construction, it is not obvious that the loops of the $\CLE_\kappa$ are independent of the choice of root, are continuous curves, and are locally finite. In the case that $\kappa \in (8/3,4]$, this was proved in \cite{CLE} as a consequence of the representation of the $\CLE_\kappa$ loops as the outer boundaries of clusters of the Brownian loop-soup; see also \cite{CLEPerc} for another argument. In the case where $\kappa \in (4,8)$, this was proved as a consequence of the continuity of $\SLE_\kappa(\kappa-6)$ \cite{IG1}, the reversibility of $\SLE_\kappa$ \cite{IG3}, and the continuity of space-filling $\SLE_\kappa$ \cite{IG4}.

If $\kappa \in (8/3,4]$, we have that $\CLE_{\kappa}$ consists of disjoint simple loops which a.s.\ do not intersect the domain boundary. However, if $\kappa \in (4,8)$, then the loops in the $\CLE_{\kappa}$ are non-simple and can intersect each other and the domain boundary. In the present work, we will focus on the $\kappa = 4$ case, where the loops do not intersect each other or the domain boundary but \emph{nearly} do. It was further shown in \cite{CLE} that if one has an ensemble of loops which satisfies~\eqref{it:cle_conf_inv} and~\eqref{it:cle_dmp}, together with the additional assumption that the loops are simple and do not intersect the boundary or each other, then it must be a $\CLE_\kappa$ for some $\kappa \in (8/3,4]$.

\subsection{Level lines of the GFF}
\label{subsec:level_lines_of_gff}

Let $\lambda = \pi / 2$ and let $h$ be an instance of the GFF on $\BH$ with boundary conditions given by $-\lambda$ (resp.\ $\lambda$) on $\BR_-$ (resp.\ $\BR_+$). Then, it was shown in \cite{ContourLineGFF} that a zero level line $\eta$ of $h$ from $0$ to $\infty$ can be made sense of in a rigorous way and has the law of a chordal $\SLE_4$ from $0$ to $\infty$. Moreover, $\eta$ is characterized by the property that for all $t \geq 0$, the conditional law of $h$ given $\eta|_{[0,t]}$ is that of a GFF on $\BH \setminus \eta([0,t])$ with boundary conditions given by $-\lambda$ (resp.\ $\lambda$) on the left (resp.\ right) side of $\eta([0,t])$ and on $\BR_-$ (resp.\ $\BR_+$).

The results in \cite{ContourLineGFF} were generalized in \cite{LevelLineGFFI} as follows. Fix a collection of weights $(\underline{\rho}_L ; \underline{\rho}_R)$ and force points 
$(\underline{x}_L ; \underline{x}_R)$, and let $h$ be a GFF on $\BH$ with boundary conditions given by $-\lambda( 1+ \sum_{i=0}^j \rho_{i,L})$ on $[x_{j+1,L}, x_{j,L})$ for $j \in [0,\ell]_\BZ$ and $\lambda (1+\sum_{i=0}^j \rho_{i,R})$ on $[x_{j,R}, x_{j+1,R})$ for $j \in [0,r]_\BZ$, where $\rho_{0,L} = \rho_{0,R} = 0$, $x_{0,L} = 0^-$, $x_{\ell + 1,L} = -\infty$, $x_{0,R} = 0^+$, and $x_{r+1,R} = \infty$.  Assume moreover that $\sum_{j=1}^{k} \rho_{j,L} > -2$ for all $k \in [1,\ell]_\BZ$ and $\sum_{j=1}^{k} \rho_{j,R} > -2$ for all $k \in [1,r]_\BZ$.  Then, the zero level line $\eta$ of $h$ from $0$ to $\infty$ can be rigorously defined. Moreover, it has the law of a chordal $\SLE_4(\underline{\rho}_L ; \underline{\rho}_R)$ and it is a measurable function of $h$. Furthermore, the level line is characterized by the property that the conditional law of $h$ given $\eta|_{[0,t]}$ is that of a GFF on $\BH \setminus \eta([0,t])$ with boundary conditions given by $-\lambda$ (resp.\ $\lambda$) on the left (resp.\ right) side of $\eta([0,t])$ and the same boundary values as $h$ on $\partial \BH$.

For all $u \in \BR$, the level line of $h$ with height $u$ from $0$ to $\infty$ is defined as the level line of $h-u$ from $0$ to $\infty$. Then, we have the following interaction rules for level lines with different heights proved in \cite{LevelLineGFFI}.

\begin{theorem}[{\cite[Theorems~1.1.4 and 1.1.5]{LevelLineGFFI}}]
	Let $h$ be a GFF on $\BH$ with piecewise constant boundary conditions that change only finitely many times. For all $u \in \BR$ and $x \in \partial \BH$, we denote by $\gamma_u^x$ the level line of $h$ with height $u$ starting from $x$ and ending at $\infty$. Fix $x_2 \leq x_1$ and $u_1, u_2 \in \BR$.
	\begin{enumerate}
		\item If $u_2 < u_1$, then $\gamma_{u_2}^{x_2}$ a.s.\ stays to the left of $\gamma_{u_1}^{x_1}$.
		\item If $u_2 = u_1$, then $\gamma_{u_2}^{x_2}$ may intersect $\gamma_{u_1}^{x_1}$ and, upon intersecting, the two curves merge and never separate.
		\item If $u_1 - u_2 \geq 2 \lambda$, then $\gamma_{u_1}^{x_1}$ and $\gamma_{u_2}^{x_2}$ do not intersect each other a.s.
	\end{enumerate}
\end{theorem}

\subsection{The coupling between $\CLE_4$ and the GFF}
\label{subsec:nested_cle_gff}

Now we describe the coupling between a zero-boundary GFF and a nested $\CLE_4$ that we are going to consider.
Let $G \subseteq \BC$ be a proper simply connected domain.  Recall that the nested $\CLE_4$ $\overline{\Gamma}$ on $G$ is obtained from the non-nested $\CLE_4$ $\Gamma$ on $G$ by the following inductive procedure.  Let $\Gamma^1\defeq\Gamma$ and inductively,  we assume that $n \in \BN$ and a collection of loops $\Gamma^n$ in $G$ (which we call the $n$-th generation loops) has been defined.  Then conditional on $\{\Gamma^k\}_{1 \leq k \leq n}$,  for each loop $\SCL$ in $\Gamma^n$,  we let $\Gamma_{\SCL}$ be a non-nested $\CLE_4$ in the open region $\mathop{\mathrm{int}}(\SCL)$ which is surrounded by $\SCL$,  and we take these CLEs to be conditionally independent given $\{\Gamma^k\}_{1 \leq k \leq n}$.  We then let $\Gamma^{n+1} \defeq \bigcup_{\SCL \in \Gamma^n} \Gamma_{\SCL}$.  Finally,  we let
\begin{align*}
	\overline{\Gamma}\defeq\bigcup_{n=1}^{\infty} \Gamma^n.
\end{align*}

Let us now review the coupling between $\overline{\Gamma}$ and a zero-boundary GFF $\Psi$ on $G$ that we are going to use (see e.g.\ \cite[Section~4]{aru2019bounded}).  Let $\{X_{\SCL}\}_{\SCL \in \overline{\Gamma}}$ be random variables indexed by the $\CLE_4$ loops which are conditionally independent given $\overline{\Gamma}$,  with 
\begin{align*}
	\BP[X_{\SCL} = 1 \mid \overline{\Gamma}] = \BP[X_{\SCL} = -1 \mid \overline{\Gamma}] = \frac{1}{2}.
\end{align*}
For $n \in \BN$,  we let $\Psi_n$ be the piecewise constant function on $G$ which is defined Lebesgue almost everywhere on $G$ and for each $n$-th generation loop $\SCL \in \Gamma^n$,  satisfies 
\begin{align*}
	\Psi_n|_{\mathop{\mathrm{int}}(\SCL)} = 2\lambda \sum_{k=1}^n X_{\SCL^{(k)}}
\end{align*}
where $\SCL^{(n)} = \SCL$ and, for $k = 1,\ldots,n-1$, $\SCL^{(k)}$ is the unique loop in $\Gamma^k$ which disconnects $\SCL$ from $\partial G$.  Recall also that $\lambda = \frac{\pi}{2}$.  Then we have that $\Psi_n$ converges a.s.\ in the distributional sense, as $n \to \infty$, to a zero-boundary GFF $\Psi$ on $G$.

In the coupling described above, we have that $\Psi$ and $\left(\overline{\Gamma}, \{X_{\SCL}\}_{\SCL \in \overline{\Gamma}}\right)$ are a.s.\ given by measurable functions of each other.  Furthermore, the following Markov property holds.  For each $n \in \BN$,  the conditional law of $(\Psi,\overline{\Gamma})$ given
\begin{align*}
	\bigcup_{k=1}^n \Gamma^k \quad \text{and} \quad \left\{X_{\SCL} :  \SCL \in \bigcup_{k=1}^n \Gamma^k \right\}
\end{align*}
is described as follows.  Let $\{\Psi_{\SCL}\}_{\SCL \in \Gamma^n}$ be conditionally independent zero-boundary GFFs on the domains $\mathop{\mathrm{int}}(\SCL)$ surrounded by the loops in $\Gamma^n$.  Then we have that
\begin{align*}
	\Psi|_{\mathop{\mathrm{int}}(\SCL)} = \Psi_{\SCL} + 2\lambda \sum_{k=1}^n X_{\SCL^{(k)}} \quad \text{for all} \quad \SCL \in \Gamma^n,
\end{align*}
where the $\SCL^{(k)}$ are as above.  Moreover,  the set $\overline{\Gamma}|_{\mathop{\mathrm{int}}(\SCL)}$ of loops of $\overline{\Gamma}$ which are contained in $\mathop{\mathrm{int}}(\SCL)$ is the same as the set of \emph{level loops} of the field $\Psi_{\SCL}$, that is, the loops of the nested $\CLE_4$ which is coupled with $\Psi_{\SCL}$ as described above.

\subsection{Labeled $\CLE_4$}
\label{subsec:labelled_cle_4}

Labeled $\CLE_4$ in $\BD$, also called the uniform exploration of the $\CLE_4$ on $\BD$, was constructed in \cite{CoInCLEExpl} and corresponds to a Markovian exploration of $\CLE_4$ loops, with labels on each loop that keep track of the time when each loop is discovered. We will describe briefly the construction and refer to \cite{CoInCLEExpl} and \cite{LevelLineGFFI} for more details.

To construct the exploration, we define the measure $M$ as the image measure of $\mu \otimes \omega$ by the mapping $(\ell, x) \mapsto x \cdot \ell$, where $\omega$ denotes the harmonic measure on $\partial \BD$ as seen from $0$, and $\mu$ denotes the $\SLE_4$ bubble measure on $\BD$ pinned at $1$ introduced in \cite{CLE}. We also let $(\ell_t)_{t \geq 0}$ denote a Poisson point process (PPP) with intensity measure given by $M$ times the Lebesgue measure on $\BR_+$. This process defines an ordering on the (pinned) loops, also often called bubbles, where each loop $\ell$ is equipped with time label $t_{\ell}$.

The exploration targeted at $0$ is defined as follows. Let $\tau_0$ denote the first time that the PPP discovers a loop that surrounds $0$ and for all $\varepsilon \in (0,1)$, we let $\ell_{t_1^{\varepsilon}},\ldots, \ell_{t_{n_{\varepsilon}}^{\varepsilon}}$ denote the loops discovered until time $\tau_0$ with diameter at least $\varepsilon$, such that $t_1^{\varepsilon}<\cdots<t_{n_{\varepsilon}}^{\varepsilon}$. 
We define $D_0 \defeq \BD$ and inductively let $D_i$ denote the connected component of $D_{i-1} \setminus \phi_i(\ell_{t_i^{\varepsilon}})$ containing $0$, where $\phi_i$ is the conformal mapping from $\BD$ onto $D_{i-1}$ such that $\phi_i(0) = 0$ and $\phi_i'(0) > 0$. Then, it was shown in \cite{CoInCLEExpl} that the boundary of $D_{n_{\varepsilon}}$ converges as $\varepsilon \to 0$ to a loop which has the same law as the loop in a $\CLE_4$ encircling zero. The loop surrounding zero is then equipped with the time label $\tau_0$ and we say that it is discovered at time $\tau_0$. More generally, if we set $D^\varepsilon_t = D_i$ when $t_{i-1}^\varepsilon \le t<t_i^\varepsilon$ for all $t \le \tau_0$, where $t_0^\varepsilon \defeq 0$, then $D^\varepsilon_t$ converges in the Carath\'eodory sense to a simply connected open domain $C_t(0)$, which is called the connected component of the unexplored region containing the origin.

It was shown in \cite[Lemma~6]{CoInCLEExpl} that the measure $M$ is invariant under M\"obius transformations of $\BD$. Therefore, by conformal invariance, the exploration targeted at any fixed point $z \in \BD$ is defined by mapping the exploration targeted at $0$ via a M\"obius transformation of $\BD$ which maps $0$ to $z$.

The labeled $\CLE_4$ exploration is then defined by coupling the exploration processes targeted at any two fixed points $z,w \in \BD$ together so that they agree up until the first time that they are separated,  and then evolve independently. Each loop $\SCL$ of the $\CLE_4$ can then be labeled by the time $t_\SCL$ at which it is discovered. For all $z \in \BD$, we denote by $C_t(z)$ the connected component of the unexplored region containing $z$. The explored region is the compact set 
\[\SCB_t(\partial \BD) \defeq \overline{\BD} \setminus \bigcup_{z \in \BD_\BQ} C_t(z).\]
We use here the same symbol $\SCB_t$ as for the metric balls of \Cref{def:axioms}.  Let us emphasize that this is not an assumption: that the explored region at time $t$ coincides with the metric ball of radius $t$ from $\partial \BD$ is precisely what Axiom~\eqref{it:axiom_uniform_exploration} of \Cref{def:axioms} asserts.
Since every rational point is encircled by a loop and since $C_t(z)$ is open, we can also write $\overline{C_t(z)} = \overline{\bigcup_{x \in C_t(z) \cap \BQ^2} \mathop{\mathrm{int}}(\SCL(x)) }$.

\subsection{Two-valued sets of the GFF}
\label{subsec:two_valued_level_sets}

Next, we introduce the two-valued (local) sets of a GFF $h$ on $\BD$ that we are going to consider. Recall from \Cref{subsec:gff} the definition of a local set of $h$, introduced in \cite{ContourLineGFF}. We refer to \cite{aru2019bounded} and \cite{TVSGFF} for more detailed expositions of two-valued sets of the GFF.

\begin{definition}
	Let $a,b>0$ and let $h$ be a GFF in a proper simply connected domain $G \subseteq \BC$. We say that a local set $A$ of $h$ is a \emph{two-valued local set} of $h$ of levels $-a$ and $b$ if it satisfies the following properties.
	\begin{enumerate}
		\item \label{it:boundary_conditions}
		If $h_{A}$ denotes the harmonic extension to $G\setminus A$ of the boundary values of $h$ on $A$, then we have that $h_{A} \in \{-a,b\}$ a.s.
		\item \label{it:thin_local_set}
		For any smooth test function $f \in C_0^{\infty}(G)$, a.s.,
		\begin{align*}
			(h,f) = (h-h_{A} , f) + \int_{G \setminus A} h_{A}(x) f(x) \mathrm{d}x.
		\end{align*}
		\item \label{it:finitely_many_components} The set $A \cup \partial G$ has a finite number of connected components.
	\end{enumerate}
\end{definition}

Suppose that $h$ has zero boundary conditions. It was shown in \cite{aru2019bounded} that if $a + b \geq 2\lambda$, then such a set exists and it is a.s.\ determined by $h$; we denote it by $\mathbb{A}_{-a,b}$, and by $\mathbb{A}_{-a,b}(O)$ when we wish to record the domain $O$ in which it is constructed. Moreover, $\mathbb{A}_{-a,b}$ is unique in the sense that if $A'$ is another local set of $h$ satisfying~\eqref{it:boundary_conditions}, \eqref{it:thin_local_set} and~\eqref{it:finitely_many_components}, then $A' = \mathbb{A}_{-a,b}$ a.s.

Let us now briefly review the construction of $\mathbb{A}_{-a,-a+2\lambda}$ for $a \in (0,2\lambda)$, which will be needed for our purposes. We follow \cite[Section~3.2]{TVSGFF}. 

To begin with, we let $\eta$ denote the level line of $h+a-\lambda$ from $-\ri$ to $\ri$ and set $A^1 \defeq \eta([0,\infty])$. In each connected component $O$ of $\BD \setminus A^1$ lying to the left of $\eta$, we let $x$ and $y$ be the two endpoints of the arc $\partial O \cap \partial \BD$, and we assume that $(-\ri ,x,y)$ are ordered in the counterclockwise way. We then let $\eta^O$ denote the level line of $h|_{O} + a -\lambda$ from $x$ to $y$.

As for the connected components $O$ of $\BD \setminus A^1$ lying to the right of $\eta$, we perform an analogous procedure except for exploring the level line of $h|_O + a - \lambda$ from $y$ to $x$. We also denote those level lines by $\eta^O$.

In the next step of the exploration, we let $A^2$ be the closure of the union of $A^1$ with the $\eta^O$'s for each connected component $O$ of $\BD \setminus A^1$. In the connected components of $\BD \setminus A^2$ whose boundaries are contained in $A^2$, we stop the exploration. Note that the boundary conditions of $h$ on each of the aforementioned components are a.s.\ constant and given by either $-a$ or $2\lambda - a$. In the rest of the components, we iterate as before to obtain $A^n$. We note that if $O$ is a connected component of $\BD \setminus A^2$ such that $\partial O \cap \partial \BD \neq \emptyset$, then we have that the boundary conditions of $h|_{O}$ are given by $0$ on $\partial O \cap \partial \BD$, and $-a$ (resp.\ $2\lambda-a$) on $\partial O \setminus \partial \BD$ if $O$ lies to the left (resp.\ right) of $A^1$. Then, we set $\mathbb{A}_{-a,-a+2\lambda}\defeq\overline{\bigcup_{n \in \BN} A^n}$.

\begin{remark}\label{rem:target_invariance}
	The uniqueness of the two-valued set $\mathbb{A}_{-a,b}$ implies that the starting and ending points of the level lines and the order in which they were sampled to produce $\mathbb{A}_{-a,b}$ do not matter. This boils down to the reversibility of the $\SLE_4(\underline{\rho})$ processes.
\end{remark}

\begin{remark}\label{remark coupling TVS BCLE}
	The set $\mathbb{A}_{-a,-a+2\lambda}$ constructed above is equal to the range of a clockwise boundary conformal loop ensemble with $\kappa=4$ and parameter $\rho$ (i.e., a $\cwBCLE_4(\rho)$ process) with $\rho = -a/\lambda$ introduced in \cite[Section~7]{CLEPerc}. Recall that the loops in the complement of the $\cwBCLE_4(\rho)$ with a clockwise (resp.\ counterclockwise) orientation are its ``true'' (resp.\ false) loops and the true loops intersect the domain boundary as an $\SLE_4(\rho)$ process.  By varying $a \in (0,2\lambda)$, we obtain $\cwBCLE_4(\rho)$ for all $\rho \in (-2,0)$. The loops with boundary condition $2\lambda -a$ in $\mathbb{A}_{-a,-a+2\lambda}$ correspond to the true loops of the $\cwBCLE_4(\rho)$ and the loops with boundary condition $-a$ in $\mathbb{A}_{-a,-a+2\lambda}$ correspond to the false loops of the $\cwBCLE_4(\rho)$.
\end{remark}

\subsection{Coupling between the GFF and the labeled $\CLE_4$}
\label{subsec:gff_labelled_cle_4}

Now, we are ready to construct a coupling that we will use between a zero-boundary GFF $h$ on $\BD$ and a labeled $\CLE_4$ $\{(\ell , t_{\ell})\}_{\ell \in \Gamma}$. We will follow the construction given in \cite[Section~6]{TVSGFF}.

Fix $r \in (0,2\lambda)$. For all $j \in \BN$, we define sets $\BB_r^j$ iteratively as follows. First, we set $\BB_r^1\defeq\mathbb{A}_{-r,-r+2\lambda}$. In each connected component $O$ of $\BD \setminus \BB_r^1$ such that the boundary conditions of $h|_O$ are given by $-r$, we explore the set $\mathbb{A}_{-r,-r+2\lambda}(O)$ which is defined in the same way as $\mathbb{A}_{-r,-r+2\lambda}$, but with the zero-boundary GFF $h|_O - h_{\BB_r^1}$ in place of $h$, where we recall that $h_{\BB_r^1}$ is the harmonic extension to $\BD \setminus \BB_r^1$ of the boundary values of $h$ on $\BB_r^1$. Then, we let $\BB_r^2$ denote the closure of the union of the sets explored. Note that the boundary conditions of $h$ on the connected components of $\BD \setminus \BB_r^2$ are constant and lie in $\{2\lambda-r,2\lambda-2r,-2r\}$.

Suppose that we have constructed $\BB_r^j$ for $j \in \BN$. Then, in each connected component $O$ of $\BD \setminus \BB_r^j$ such that the boundary conditions of $h|_O$ are given by $-jr$, we explore the set $\mathbb{A}_{-r,-r+2\lambda}(O)$ defined in the same way as $\mathbb{A}_{-r,-r+2\lambda}$ but with the zero-boundary GFF $h|_O - h_{\BB_r^j}$ in place of $h$. We define $\BB_r^{j+1}$ as the closure of the union of $\BB_r^j$ and the sets explored.

We set $\BB_r\defeq\overline{\bigcup_{j \in \BN} \BB_r^j}$ and $\BB_0\defeq\overline{\bigcup_{n \in \BN} \BB_{2^{-n}}}$. To every loop $\ell$ of $\BB_0$ (i.e., whose interior $\mathop{\mathrm{int}}(\ell)$ is a connected component of $\BD \setminus \BB_0$), we associate a label $t_{\ell}$ such that the boundary conditions of $h|_{\mathop{\mathrm{int}}(\ell)}$ are given by $2\lambda - t_{\ell}$. Then, we have the following.

\begin{proposition}[{\cite[Proposition~6.6]{TVSGFF}}]
	\label{prop coupling TVS uniform explo}
	The collection of loops in $\BB_0$ has the law of a non-nested $\CLE_4$ $\Gamma$ on $\BD$ and $\{(\ell,t_{\ell})\}_{\ell \in \Gamma}$ has the law of a labeled $\CLE_4$.
\end{proposition}

Besides, let us recall from \cite{kkmt2026cle4_part1} that for all $r \in (0, 2\lambda)$, one can couple $(\SCB_s(\partial \BD))_{s\ge 0}$ with a $\cwBCLE_4(-r/\lambda)$ so that all the domains encircled by false (counterclockwise) loops of the $\cwBCLE_4(-r/\lambda)$ are exactly the connected components of the unexplored region $\BD\setminus \SCB_r(\partial \BD)$ whose boundary intersects $\partial \BD$. More precisely, we proved the following result.
\begin{lemma}[{\cite[Lemma~2.9]{kkmt2026cle4_part1}}]\label{lemma BCLE loops}
	Consider the uniform exploration $(\SCB_s(\partial \BD))_{s\ge 0}$, let $r \in (0, 2\lambda)$ and consider the $\cwBCLE_4(-r/\lambda)$ coupled with $(\SCB_s(\partial \BD))_{s\ge 0}$ as above. Then, the interiors of the false $\cwBCLE_4(-r/\lambda)$ loops are exactly the connected components of $\BD \setminus \SCB_r(\partial \BD)$ whose boundary intersects $\partial \BD$.
\end{lemma}

\subsection{Four-arm estimates}

Next,  we mention one useful four-arm estimate for bi-chordal $\SLE_4$ which we are going to use in the last section. It estimates the probability that both of the two branches of the bi-chordal $\SLE_4$ intersect a small Euclidean ball centered at an interior point of the domain where the bi-chordal $\SLE_4$ is defined. Lemma~\ref{lem:bichordal-bulk-4A} follows immediately from \cite[Theorem~1.1]{zhan2020two}.

\begin{lemma}\label{lem:bichordal-bulk-4A}
    Let $z_1, z_2, z_3, z_4$ be points on $\partial\BD$ ordered in the  counterclockwise way.  Let also $\{\eta_{12}, \eta_{34}\}$ be a bi-chordal $\SLE_4$ in $\BD$ with (interior) link pattern $\{\{z_1, z_2\}, \{z_3, z_4\}\}$.  Then,  there exists a universal constant $C \in (0,\infty)$ such that the following holds for all $\varepsilon \in (0,1)$.  The probability that both $\eta_{12}$ and $\eta_{34}$ intersect $B_\varepsilon(0)$ is at most $C \varepsilon^2$.
\end{lemma}

\begin{proof}
The claim in the statement of the lemma follows from \cite[Theorem~1.1]{zhan2020two}.
\end{proof}

\subsection{Estimates for extremal distance}
We record here some estimates related to extremal distance. The first one controls the probability that a planar Brownian motion hits the left or right side of a very thin horizontal rectangle.
\begin{lemma}[{\cite[Proposition~{3.69}]{lawler2008conformally}}]
		\label{lem:extremal_length}
		For $L>0$ we set $R_L \defeq(0,L) \times (0,\pi)$ and let $\partial_1 = [0,\ri\pi]$ and $\partial_2 = \partial_{2,L} = [L,L+\ri\pi]$.  Let $\tau = \tau_{R_L}$ denote the first time that a complex Brownian motion exits $R_L$ and set
		\begin{align*}
			&f_1(z) = f_{1,L}(z) = 2 \min\{\BP_z[B_{\tau} \in \partial_1] ,  \BP_z[B_{\tau} \in \partial_2]\}
		\end{align*}
		and set $\Theta(R_L ; \partial_1 , \partial_2) = \sup\{ f_1(z) : z \in R_L\}$, which we abbreviate to $\Theta(L)$.  Then $\Theta(L)$ is a continuous, strictly decreasing function of $L$ such that 
		\begin{align*}
			\Theta(L) = \left(\frac{8}{\pi}\right) e^{-L/2} + O(e^{-L}) \quad \text{as} \quad L \to \infty.
		\end{align*}
		Moreover the supremum in the definition of $\Theta$ is attained when $z$ is the center $(\frac{L}{2}) + \ri (\frac{\pi}{2})$.
\end{lemma}
The next lemma gives an upper bound for the extremal distance between two compact subsets of the unit disk.
\begin{lemma}[{\cite[Lemma~4.2]{kkmt2026cle4_part2}}]\label{lem:Beurling}
		There is a universal constant $C > 0$ such that the following is true: Let $K_1, K_2 \subset \overline\BD$ be two disjoint, connected, simply connected, and compact subsets,  such that $\partial \BD \not\subseteq K_1 \cup K_2$. Then the extremal distance in $\BD$ between $K_1$ and $K_2$ is at most 
		\begin{equation}\label{eq:Beurling}
			\left.C\middle/\log\left(\frac{\diam(K_1) \wedge \diam(K_2)}{\dist(K_1, K_2)} \vee 1\right)\right., 
		\end{equation}
		where $\diam$ and $\dist$ denote Euclidean diameter and distance, respectively. 
\end{lemma}

\subsection{Locality property for metric balls}

In the present subsection, we record the fact that metric balls satisfy a locality property which is a consequence of the results from \cite{kkmt2026cle4_part2}. Roughly speaking, this means that conditionally on a metric ball, the loop configuration in each complementary component is a non-nested $\CLE_4$, and the internal metric on that component is precisely the metric associated with this $\CLE_4$.

\begin{lemma}[{\cite[Lemma~2.10]{kkmt2026cle4_part2}}]\label{lem:metric-ball-locality}
	Let $U \subsetneq \BC$ be a simply connected domain. Let $z \in U$ be deterministic. Let $\tau$ be a stopping time for the filtration generated by 
	\begin{equation*}
		\left(\SCB_t(\SCL(z); D_{\Gamma_U}^U), \Gamma_U|_{\SCB_t(\SCL(z); D_{\Gamma_U}^U)}\right)_{t \ge 0}. 
	\end{equation*}
	Then the following is true: Let $V \subset U$ be a deterministic simply connected domain. Write $\{V_j\}_j$ for the connected components of $V^\star$ (recall~\eqref{eq:def-V-star} for the definition of $V^\star$).  Let $(\Gamma_U, D_{\Gamma_U}^U, \{D_{\Gamma_U|_{V_j}}^{V_j}\}_j)$ be coupled as in \Cref{def:weak_axioms}, Axiom~\eqref{it:weak_axiom_locality} (locality). Then, on the event that $\SCB_\tau(\SCL(z); D_{\Gamma_U}^U) \cap V^\star = \emptyset$, conditionally on $\SCB_\tau(\SCL(z); D_{\Gamma_U}^U)$, on $V^\star$ and on the $\sigma$-algebra generated by~\eqref{eq:weak_axiom_locality}, the $(\Gamma_U|_{V_j}, D_{\Gamma_U|_{V_j}}^{V_j})$'s are independent and their conditional laws are those of $(\Gamma_{V_j}, D_{\Gamma_{V_j}}^{V_j})$, respectively.
	The same is true with a deterministic connected arc of $\partial U$ in place of $\SCL(z)$. 
\end{lemma}

\subsection{Estimates for the $\CLE_4$ metric}\label{section:estimates}

In this last subsection, we record a few useful results that we proved in \cite{kkmt2026cle4_part2}. We established them for a weak geodesic $\CLE_4$ metric coupling. But since we also showed in \cite[Theorem~1.4]{kkmt2026cle4_part2} that any geodesic $\CLE_4$ metric coupling is a weak geodesic $\CLE_4$ metric coupling, they hold for any geodesic $\CLE_4$ metric coupling. In this subsection, $D$ is a geodesic $\CLE_4$ metric coupling as defined in \Cref{def:weak_axioms}.

Let us first recall a useful bound for the distance from top to bottom across a conformal rectangle.

\begin{proposition}[{\cite[Proposition~3.2]{kkmt2026cle4_part2}}]\label{prop:distance_across_rectangle}
	For each $r > 0$, write
	\begin{equation*}
		V_r \defeq \left\{z \in \BC : 0 < \Re(z) < r, \ 0 < \Im(z) < 1\right\};
	\end{equation*}
	write $T_r$ (resp.\ $B_r$) for the top (resp.\ bottom) side of $V_r$; let $\Gamma_{V_r}$ be a non-nested $\CLE_4$ in $V_r$. (By convention, $V_r$ is of conformal modulus $1/r$.) Then there is a deterministic constant $c \in (0, 1)$ such that 
	\begin{equation*}
		\BP\!\left\lbrack D_{\Gamma_{V_r}}^{V_r}(T_r, B_r) \le 1/r\right\rbrack \ge c
	\end{equation*}
	for all sufficiently large $r > 0$. 
\end{proposition}
Moreover, a key estimate in \cite{kkmt2026cle4_part2} was the following estimate for the probability that two metric balls stopped when they hit a small Euclidean ball do not intersect each other.
\begin{proposition}[{\cite[Proposition~4.10]{kkmt2026cle4_part2}}]\label{prop:hitting-two-metric-balls}
	Let $\Gamma_\BD$ be a non-nested $\CLE_4$ in $\BD$. Let $\delta > 0$. Then for each $x, y, z \in \BD$ with $\lvert x - z\rvert \wedge \lvert y - z\rvert \ge \delta$, 
		\begin{equation}\label{eq:hitting-two-metric-balls}
			\BP\!\left\lbrack\SCB^{B_\varepsilon(z)}(\SCL(x); D_{\Gamma_\BD}^\BD) \cap \SCB^{B_\varepsilon(z)}(\SCL(y); D_{\Gamma_\BD}^\BD) = \emptyset\right\rbrack \le \varepsilon^{2 + o(1)} \quad \text{as } \varepsilon \to 0, 
		\end{equation}
		at a rate depending only on $\delta$, where $\SCB^{B_\varepsilon(z)}(\SCL; D_{\Gamma_\BD}^\BD) \defeq \SCB_{D_{\Gamma_\BD}^\BD(\SCL, B_\varepsilon(z))}(\SCL; D_{\Gamma_\BD}^\BD)$. The same is true with two deterministic disjoint connected arcs of $\partial\BD$ in place of $\SCL(x)$ and $\SCL(y)$. (In this case, the rate of~\eqref{eq:hitting-two-metric-balls} depends only on the deterministic arcs and is uniform over all $z \in \BD$.)
\end{proposition}

Next, let us state a key lemma. For all $a>0$, let $V_a \defeq (0,a)\times (0,1)$. The lemma below states that with positive probability, the $D^{V_a}_{\Gamma_{V_a}}$-distance from top to bottom is at least twice the $D^{V_a}_{\Gamma_{V_a}}$-distance from left to right.
\begin{lemma}[{\cite[Lemma~5.1]{kkmt2026cle4_part2}}]\label{lem:genealogy}
		Fix $a>0$ and let $\Gamma_{V_a}$ be a $\CLE_4$ in $V_a$.  Then, there exists $r_0>0$ such that with positive probability we have that 
		\begin{align*}
			D_{\Gamma_{V_a}}^{V_a}(\{0\} \times [0,1] ,  \{a\} \times [0,1]) < r_0 < 2r_0 \leq D_{\Gamma_{V_a}}^{V_a}([0,a] \times \{0\} ,  [0,a] \times \{1\}).
		\end{align*}
		More precisely, the following holds with positive probability. Let us divide $V_a$ into three rectangles $G_1= (0,a) \times (4/5, 1) $, $G_2= (0,a) \times (1/5, 4/5)$ and $G_3=  (0,a) \times (0, 1/5)$. Let $L_2=\{0\} \times [1/5,4/5]$ be the left boundary of $G_2$ and let $L=\{0\}\times [0, 1]$ be the left boundary of $V_a$. Let $\partial_\mathrm{R} \SCB_{r_0}(L; D^{V_a}_{\Gamma_{V_a}})$ be the right boundary of the metric ball $\SCB_{r_0}(L; D^{V_a}_{\Gamma_{V_a}})$, which we define as the intersection of $\SCB_{r_0}(L; D^{V_a}_{\Gamma_{V_a}})$ with the closure of the union of the connected components $C$ of $V_a \setminus \SCB_{r_0}(L; D^{V_a}_{\Gamma_{V_a}})$ such that $\overline{C} \cap (\{a\}\times (0,1)) \neq \emptyset$. Then,
			\begin{enumerate}[label=(\alph*)]
				\item\label{it:lemma crossing a} We have $D_{\Gamma_{V_a}}^{V_a}(L_2 ,  \{a\} \times [1/5,4/5]) < r_0$;
				\item\label{it:lemma crossing b} The right boundary $\partial_\mathrm{R} \SCB_{r_0}(L; D^{V_a}_{\Gamma_{V_a}})$ intersects $L_2$;
				\item\label{it:lemma crossing c} We have $D_{\Gamma_{V_a}}^{V_a}([0,a] \times \{0\}, [0,a] \times \{1/5\})\ge r_0$ and $D_{\Gamma_{V_a}}^{V_a}([0,a] \times \{4/5\}, [0,a] \times \{1\})\ge r_0$.
			\end{enumerate} 
		
\end{lemma}

Let us end this subsection by recalling a property of geodesics, which states that they do not touch the boundary of the domain.

\begin{corollary}[{\cite[Corollary~6.5]{kkmt2026cle4_part2}}]\label{cor:geodesic_avoid_boundary}
	Almost surely, for each $\SCL_1, \SCL_2 \in \Gamma_{\BD}$, any $D_{\Gamma_{\BD}}^{\BD}$-geodesic connecting $\SCL_1$ and $\SCL_2$ does not touch $\partial{\BD}$. 
	A similar statement holds with two deterministic disjoint connected arcs of $\partial\BD$ in place of $\SCL_1$ and $\SCL_2$ (in which case the $D_{\Gamma_\BD}^\BD$-geodesic touches $\partial\BD$ only at its endpoints).
\end{corollary}

\section{Uniqueness of the metric: genealogy, main result, and proof strategy}
\label{sec:uniqueness_and_measurability_of_the_metric}

In this section and in Sections~\ref{sec:crossings_of_annuli}, \ref{sec:good_annuli} and~\ref{sec:covering_geodesics-by_good_annuli}, we will show that any geodesic $\CLE_4$ metric can be recovered from the uniform exploration. This result is the key ingredient toward proving \Cref{thm:uniqueness}. Let us state this result more precisely.
\begin{theorem}\label{thm metric determined by the exploration}
	Let $U$ be a simply connected domain and $(\SCB_t(\partial U))_{t\ge 0}$ be a uniform exploration of $\Gamma_U$ from the boundary. Let $D$ be a geodesic $\CLE_4$ metric coupling and assume that for all $t\ge 0$,
	\[
	\SCB_t(\partial U) = \SCB_t(\partial U; D^U_{\Gamma_U}).
	\]
	Then, the distance $D^U_{\Gamma_U}$ is a uniquely defined measurable function of $(\SCB_t(\partial U))_{t\ge 0}$.
\end{theorem}
In particular, the above result already shows that the geodesic $\CLE_4$ metric coupling is unique in law. To prove Theorem~\ref{thm metric determined by the exploration}, we first define a genealogy of loops using the uniform exploration. Then, we prove that the geodesic $\CLE_4$ metric coupling is uniquely determined by the distances along geodesics to the boundary and we conclude by showing that the geodesics to the boundary are indeed described by the genealogy of loops. By conformal invariance, we focus on $U=\BD$.

\subsection{Definition of a genealogy from the uniform exploration}\label{section:genealogy-of-loops}

Let $D$ be a geodesic $\CLE_4$ metric coupling. Consider a loop $\SCL(x)$ of the $\CLE_4$ discovered at some time $t_x\ge 0$ (note that a.s.\ $t_x>0$, since the loops of a $\CLE_4$ do not touch $\partial \BD$) by the uniform exploration $(\SCB_t(\partial \BD; D_{\Gamma_\BD}^\BD))_{t\ge 0}$. Let $V$ be the connected component of $\BD \setminus \SCB_{t_x}(\partial \BD; D_{\Gamma_\BD}^\BD)$ such that $\SCL(x)\subseteq \partial V$. Note that such a component exists and is unique.

Let $t_x<t_1< \cdots< t_k$ for $k \ge 1$
. 
The set $\SCB_{t_1}(\partial \BD; D_{\Gamma_\BD}^\BD)\cap V$ has countably many connected components, which are denoted by $(C_{1,j})_{j\ge 1}$ (with the convention $C_{1,j}=\emptyset$ if there are less than $j$ connected components). Similarly, for all $2 \le i\le k$, let $(C_{i,j})_{j\ge 1}$ be the family of connected components of $\SCB_{t_i}(\partial \BD; D_{\Gamma_\BD}^\BD) \cap V \setminus \SCB_{t_{i-1}}(\partial \BD; D_{\Gamma_\BD}^\BD)$. Throughout, the components $C_{i,j}$ are understood to be taken relative to the sequence $t_x<t_1<\cdots<t_k$.

We say that a component $C_{1,j}$ for $j\ge 1$ is a \emph{descendant} of $\SCL(x)$ if $\overline{C_{1,j}} \cap \overline{ \partial V \setminus \SCL(x)}= \emptyset$. Then, by induction, for all $2 \le i \le k$, we say that $C_{i,j}$ is a descendant of $\SCL(x)$ if $\overline{C_{i,j}}$ does not intersect the union of the $C_{i-1,j'}$'s which are not descendants of $\SCL(x)$. We then write $\SCL(x) \preceq C_{i,j}$ whenever $C_{i,j}$ is a descendant of $\SCL(x)$.

\begin{figure}[h]
	\centering
	\includegraphics[scale=0.8]{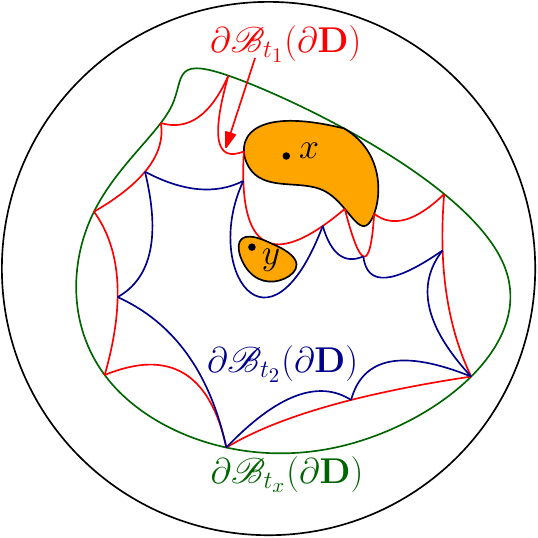}
	\caption{Definition of the genealogy. Here, the loop $\SCL(y)$ is a descendant of $\SCL(x)$.}
	\label{fig:def_genealogy}
\end{figure}

Let $\SCL(y)$ be another loop. Assume that $\SCL(y)\subseteq V$. Let $t_y \defeq D_{\Gamma_\BD}^\BD(\partial \BD, \SCL(y))$. Then, we say that $\SCL(y)$ is a descendant of $\SCL(x)$, and we write $\SCL(x)\preceq \SCL(y)$ if there exist $t_x<t_1< \cdots <t_k$ for some $k\ge 1$ (with $t_{k-1}<t_y \le t_{k}$) such that
the connected component $C_{k,j}$ which contains $\SCL(y)$ is a descendant of $\SCL(x)$. If $\SCL(y) \not\subseteq V$, we declare that $\SCL(x) \npreceq \SCL(y)$.  By convention, we also set $\SCL(x) \preceq \SCL(x)$. See \Cref{fig:def_genealogy} for an example.

Equivalently, a point $y\in V$ is a descendant of $\SCL(x)$ if there exist $k\ge 1$ and $t_x< t_1< \cdots< t_{k}$ such that $y$ is in one of the associated connected components $C_{k,j}$ which is a descendant of $\SCL(x)$, for some $j\ge 1$. A loop $\SCL(y)$ surrounding $y$ is a descendant of $\SCL(x)$ if and only if $y$ is a descendant of $\SCL(x)$.

One can see that the set of descendants of $\SCL(x)$ is
\begin{equation*}
\bigcup_{k\ge 1} \bigcup_{t_x<t_1< \cdots<t_{k} } \bigcup_{\substack{j\ge1 \\ C_{k,j} \succeq \SCL(x)}} C_{k,j}.
\end{equation*}
Moreover, a loop $\SCL(y)$ is a descendant of $\SCL(x)$ if and only if it is included in the above union. Note that we can choose $t_{k}$ so that $t_y=t_k$.

\begin{lemma}
	The relation $\preceq$ defines a partial order on $\Gamma$.
\end{lemma}
\begin{proof}
	Reflexivity and antisymmetry are obvious. For transitivity, let $\SCL(x) \preceq \SCL(y) \preceq \SCL(z)$. Let $t_w= D_{\Gamma_\BD}^\BD(\partial\BD, \SCL(w))$ for all $w \in \{x,y,z\}$. Let $t_x<t_1<\cdots< t_{k}$ be a sequence associated with $\SCL(x) \preceq \SCL(y)$ with $t_{k}= t_y$ and let $t_y< s_1< \cdots< s_{\ell}$ be a sequence associated with $\SCL(y) \preceq \SCL(z)$ with $s_{\ell}= t_z$. Let $V$ (resp.\ $V'$) denote the connected component associated with $\SCL(x)$ (resp.\ $\SCL(y)$) as above.  Since $\SCL(y)$ is a descendant of $\SCL(x)$, we have $V' \subseteq V$, so that every connected component taken relative to $V'$ is contained in a unique connected component taken relative to $V$ at the same time.  Concatenating the two sequences yields $t_x<t_1<\cdots<t_k=t_y<s_1<\cdots<s_{\ell}=t_z$, and it follows by induction on $\ell$ that the component of the associated family relative to $V$ which contains $\SCL(z)$ is a descendant of $\SCL(x)$.  Hence $\SCL(x) \preceq \SCL(z)$.
\end{proof}

\begin{lemma}\label{lemma the genealogy of the uniform exploration is contained in the geodesics to the root}
	Let $x, y \in \BD_\BQ$. Assume that $\SCL(y)$ is a descendant of $\SCL(x)$. Then, 
	a.s., every geodesic $P$ from $\partial \BD$ to $\SCL(y)$ intersects $\SCL(x)$.
\end{lemma}

\begin{proof}
	\stepn{step:genealogy-geo-setup}{Setup and notation}

	\substepn{step:ggs-times-V}{The times $t_x$, $t_y$ and the component $V$} Let $t_x\defeq D^\BD_{\Gamma_\BD}(\partial \BD, \SCL(x))$ and $t_y \defeq D_{\Gamma_\BD}^\BD(\partial \BD ,\SCL(y))$. Let $V$ be the connected component of $\BD \setminus \SCB_{t_x}(\partial \BD; D_{\Gamma_\BD}^\BD)$ such that $\SCL(x) \subseteq \partial V$, which as noted above exists and is unique.
	
	\substepn{step:ggs-families}{The families of components and the geodesic $P$} Since $\SCL(y)$ is a descendant of $\SCL(x)$, there exists an integer $k\ge 1$, $t_x<t_1 <\cdots <t_{k}$ and two families of connected components $(C_{i,j})_{1 \le i \le k, 1\le j \le m_i}$ and $(C'_{i,j})_{1 \le i \le k, 1\le j \le m'_i}$ with $m_i, m'_i\in \BN \cup \{\infty\}$ such that:
	\begin{itemize}
		\item For all $1\le i \le k$, the components $C_{i,j}$ for $1\le j \le m_i$ and $C'_{i,j}$ for $1 \le j \le m'_i$ are the connected components of $\SCB_{t_{i}}(\partial \BD; D_{\Gamma_\BD}^\BD) \cap V \setminus \SCB_{t_{i-1}}(\partial \BD; D_{\Gamma_\BD}^\BD)$ such that for all $j \le m_i$,
		\begin{equation}\label{eq descendant}
			\overline{C_{i,j}}   \cap  \bigcup_{j'=1}^{m'_{i-1}} C'_{i-1, j'}= \emptyset,
		\end{equation}
		where by convention $m_0= m'_0= 1$, $C'_{0,1}$ is $\overline{\partial V \setminus \SCL(x)}$ and $C_{0,1}$ is the closure of the domain encircled by the loop $\SCL(x)$.
		\item We have $m_{k}=1$ and  $\SCL(y)\subseteq C_{k,1}$.
	\end{itemize}
		\begin{figure}[h]
		\centering
		\includegraphics[scale=0.8]{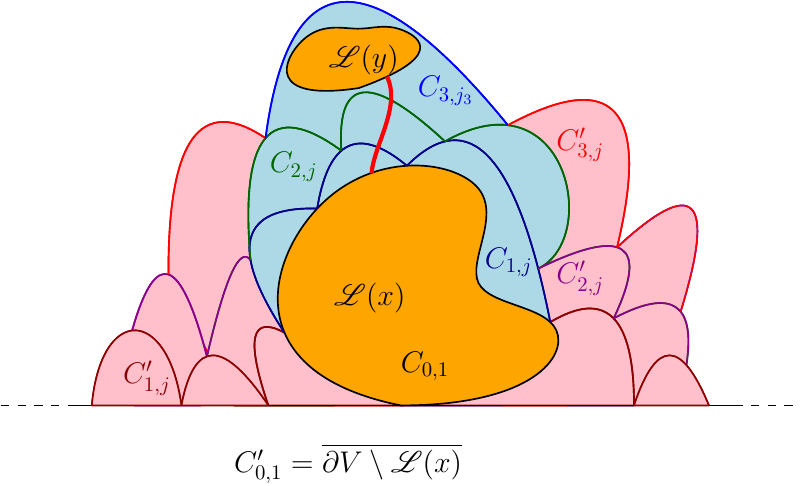}
		\caption{Illustration of the proof of \Cref{lemma the genealogy of the uniform exploration is contained in the geodesics to the root}. The components $C'_{i,j}$ are colored in pink (except $C'_{0,1}$) and the components $C_{i,j}$ are colored in blue. The part of the geodesic $P$ to $\partial \BD$ between $\SCL(y)$ and $\SCL(x)$ is drawn in red. The boundary $\partial \BD$ is not represented for the sake of simplicity.}
		\label{fig:genealogy_implies_genealogy}
	\end{figure}
	See \Cref{fig:genealogy_implies_genealogy}. By convention, we write $t_0=t_x$. Let $P\colon [0, D^\BD_{\Gamma_\BD}( \partial \BD, \SCL(y))] \to \overline{\BD}$ be a $D^\BD_{\Gamma_\BD}$-geodesic from $\partial \BD$ to $\SCL(y)$.
	
	\stepn{step:genealogy-geo-induction}{The downward induction and the case $i=k$} Let us prove by downward induction that for all $1 \le i \le k$ there exists a unique $j_i \le m_i$ such that $P ((t_{i-1}, t_{i}]) \subseteq C_{i, j_i}$. 
	
	This is clearly true for $i=k$ because $\SCL(y)\subseteq C_{k,1}$ and $P ((t_{k-1}, t_k])$ is connected, contains $\SCL(y)$ and is included in $\SCB_{t_k}(\partial \BD; D_{\Gamma_\BD}^\BD) \setminus \SCB_{t_{k-1}}(\partial \BD; D_{\Gamma_\BD}^\BD)$.
	
	\stepn{step:genealogy-geo-induction-step}{The induction step} Next, let us assume that, for some $2 \le i \le k$, there exists a unique $j_i \le m_i$ such that $P((t_{i-1}, t_{i}]) \subseteq C_{i, j_i}$. Note that $P ((t_{i-2}, t_{i-1}])$ is a connected set. 
	
	Then, note that $P((t_{i-1}, t_{i-1}+ \varepsilon]) \subseteq C_{i,j_i}$ for all $\varepsilon \in (0, t_i-t_{i-1})$, so that $P(t_{i-1}) \in \overline{C_{i,j_i}}$ by the continuity of $P$ at $t_{i-1}$. Since $t_{i-1} \in (t_{i-2}, t_{i-1}]$, it follows that
	\begin{equation}\label{eq crossing geodesic}
		P((t_{i-2}, t_{i-1}]) \cap \overline{C_{i,j_i}} \neq \emptyset.
	\end{equation}
	
	Moreover, $P((t_{i-2}, t_{i-1}])$ is a connected subset of $\SCB_{t_{i-1}}(\partial \BD; D_{\Gamma_\BD}^\BD) \setminus \SCB_{t_{i-2}}(\partial \BD; D_{\Gamma_\BD}^\BD)$, which can be written as the disjoint union of the connected components $C_{i-1,j}$ for $1 \le j \le m_{i-1}$ and $C'_{i-1, j}$ for $1 \le j \le m'_{i-1}$. So either 
	\begin{equation*}
	P ((t_{i-2}, t_{i-1}]) \subseteq C_{i-1, j}
	\end{equation*}
	for some $j \le m_{i-1}$ or 
	\begin{equation*}
	P ((t_{i-2}, t_{i-1}]) \subseteq C'_{i-1, j}
	\end{equation*}
	for some $j\leq m'_{i-1}$. The second case is excluded by~\eqref{eq descendant} and~\eqref{eq crossing geodesic}. The induction step is complete.
	
	\stepn{step:genealogy-geo-conclusion}{Conclusion of the proof} In order to end the proof, we use the same argument: $P((t_x, t_x+ \varepsilon]) \subseteq C_{1,j_1}$ for all $\varepsilon \in (0, t_1-t_x)$, so that $P(t_x) \in \overline{C_{1,j_1}}$ by the continuity of $P$ at $t_x$, and hence
	\begin{equation}\label{eq crossing 2}
		P([0, t_x])   \cap \overline{C_{1, j_1}} \neq \emptyset.
	\end{equation}
	Furthermore, by~\eqref{eq descendant} we have $\overline{C_{1,j_1}}\cap \overline{\partial V \setminus \SCL(x)}=\emptyset$, so that
	\begin{equation*}\overline{C_{1,j_1}} \cap \SCB_{t_x}^- (\partial \BD; D_{\Gamma_\BD}^\BD)= \emptyset.\end{equation*} 
	Thus $\overline{C_{1,j_1}}$ and $\SCB_{t_x}^-(\partial \BD; D_{\Gamma_\BD}^\BD)$ are disjoint. Moreover, $P([0,t_x])$ is connected, is contained in $\SCB_{t_x}(\partial \BD; D_{\Gamma_\BD}^\BD)$, intersects $\SCB_{t_x}^-(\partial \BD; D_{\Gamma_\BD}^\BD)$ (since $P(0) \in \partial \BD \subseteq \SCB_{t_x}^-(\partial \BD; D_{\Gamma_\BD}^\BD)$), and intersects $\overline{C_{1,j_1}}$ by~\eqref{eq crossing 2}.  Consequently, $P([0,t_x])$ cannot be contained in $\overline{C_{1,j_1}} \cup \SCB_{t_x}^-(\partial \BD; D_{\Gamma_\BD}^\BD)$.  Since $\SCB_{t_x}(\partial \BD; D_{\Gamma_\BD}^\BD) \setminus (\overline{C_{1,j_1}} \cup \SCB_{t_x}^-(\partial \BD; D_{\Gamma_\BD}^\BD)) \subseteq \SCL(x)$, we conclude that $P([0, t_x]) \cap \SCL(x) \neq \emptyset$.
\end{proof}

\subsection{Main ingredients and proof strategy}
\label{subsec:outline_of_the_proof}

Now, let us describe the strategy of the proof of \Cref{thm metric determined by the exploration} that will be the object of Sections~\ref{sec:crossings_of_annuli}, \ref{sec:good_annuli} and \ref{sec:covering_geodesics-by_good_annuli}. Let $D$ be a geodesic $\CLE_4$ metric coupling. Let $\Gamma$ be a non-nested $\CLE_4$ in $\BD$. To lighten notation, in the present section, we shall write $D_\Gamma = D_\Gamma^\BD$.

\subsubsection{Statement of the two main intermediate results.}
For any loops $\SCL_1, \SCL_2 \in \Gamma \cup \{\partial \BD\}$,  we set
\begin{align*}
	D_{\Gamma}^0(\SCL_1 ,  \SCL_2) = \inf_{P : \SCL_1 \to \SCL_2} \len(P ; D_\Gamma),
\end{align*}
where the infimum ranges over all admissible paths from $\SCL_1$ to $\SCL_2$ such that for all $t \in [0,1]$ with $P(t) \in \SCL$ for some $\SCL \in \Gamma \cup \{\partial \BD\}$,  there exists a $D_\Gamma$-geodesic from $\SCL$ to $\partial \BD$ containing either the time-reversal of $P|_{[0,t]}$ or $P|_{[t,1]}$.  Moreover,  we set
\begin{align*}
	\widetilde D_\Gamma(\SCL_1,\SCL_2) = \inf_{\SCL_0',\SCL_1',\ldots,\SCL_n'} \left\{ \sum_{j=1}^n D_{\Gamma}^0(\SCL_{j-1}' ,  \SCL_j') \right\},
\end{align*}
where $\SCL_1 = \SCL_0',\SCL_1',\ldots,\SCL_n' = \SCL_2$ ranges over all sequences of loops of $\Gamma \cup \{\partial \BD\}$.

Then, we have the following main ingredients for the proof of \Cref{thm metric determined by the exploration}. The first one states that the distance between two loops can be obtained by going up and down in the tree of geodesics to the boundary of the domain.
\begin{proposition}\label{thm:main_result}
	It is a.s.\ the case that
	\begin{align*}
		D_\Gamma(\SCL_1 ,  \SCL_2) = \widetilde D_\Gamma(\SCL_1,\SCL_2) \quad \text{for all} \quad \SCL_1 ,  \SCL_2 \in \Gamma.
	\end{align*}
\end{proposition}

In passing, we will also check the converse inclusion of \Cref{lemma the genealogy of the uniform exploration is contained in the geodesics to the root}, which will be the second main ingredient to get \Cref{thm metric determined by the exploration}.
\begin{proposition}\label{prop D tilde does not depend on the metric}
	Almost surely, for all $x,y \in \BD_\BQ$, the loop $\SCL(y)$ is a descendant of $\SCL(x)$ for $\preceq$ if and only if $\SCL(x)$ intersects a geodesic from $\partial \BD$ to $\SCL(y)$. In particular, since $\preceq$ was constructed from the uniform exploration, $\widetilde{D}_\Gamma$ is a uniquely defined measurable function of the uniform exploration from the boundary $(\SCB_t(\partial \BD; D_\Gamma))_{t\ge 0}$.
\end{proposition}

\begin{proof}[Proof of \Cref{thm metric determined by the exploration} assuming \Cref{thm:main_result,prop D tilde does not depend on the metric}]
    By \Cref{thm:main_result}, the distance $D_\Gamma$ is a uniquely defined measurable function of the geodesics to the root and of the distances to the boundary. But, by \Cref{prop D tilde does not depend on the metric}, we also know that the geodesics are a uniquely defined measurable function of the distances to the boundary. Thus, $D_\Gamma$ is a uniquely defined measurable function of the distances to the boundary.
\end{proof}

\subsubsection{Outline of the proof.}
Now we give an outline of the main ideas behind the proof of \Cref{thm:main_result}. The proof of \Cref{prop D tilde does not depend on the metric} will also follow as a by-product of these ideas and is given at the end of  \Cref{sec:covering_geodesics-by_good_annuli}.

Fix $x ,  y \in \BD_\BQ$ distinct points and let $\SCL$ (resp.\ $\SCL'$) be the loop in $\Gamma$ surrounding $x$ (resp.\ $y$).  Suppose that we are working on the event that $\SCL \neq \SCL'$ and let $P$ be a $D_\Gamma$-geodesic connecting $\SCL$ with $\SCL'$.  The main idea of the proof is to find Euclidean balls $B_{2^{-k_1}}(z_1),\ldots,{B_{2^{-k_{N-1}}}(z_{N-1})}$,  and a collection of pairwise distinct loops $\SCL_1,\ldots,\SCL_N$ such that the points $z_1,\ldots,z_{N-1}$ lie on a sufficiently dense Euclidean grid of points and the following properties are satisfied.
\begin{enumerate}
	\item \label{it:good_loops_in_an_annulus}
	$\SCL = \SCL_1,\SCL' = \SCL_N$,  and $P$ intersects $\SCL_j$ before intersecting $\SCL_{j+1}$ for all $1 \leq j \leq N-1$.  Moreover,  we have that for some $R>1$ large but fixed,  the loops $\SCL_j,\SCL_{j+1}$ both cross the annulus
	\begin{align*}
		A_j{\defeq}{B_{2^{-k_j}}(z_j)} \setminus \overline{B_{2^{-k_j - R}}(z_j)} \quad \text{for all} \quad 1 \leq j \leq N-1.
	\end{align*}
	\item \label{it:geodesic_between_loops}
	For all $1 \leq j \leq N$,  we let $s_j$ (resp.\ $t_j$) be the first (resp.\ last) time that $P$ intersects $\SCL_j$.  Then,  we have for all $1 \leq j \leq N-1$,  that $P([t_j,s_{j+1}]) \subseteq B_{2^{-k_j}}(z_j)$,  and the following is true for every loop $\ell \in \Gamma$ that $P|_{(t_j,s_{j+1})}$ intersects.  Let $t \in (t_j,s_{j+1})$ be such that $P(t) \in \ell$.  Then,  we have that either the time-reversal of $P|_{[t_j,t]}$ or $P|_{[t,s_{j+1}]}$ is contained in some $D_\Gamma$-geodesic connecting $\ell$ with $\partial \BD$.
\end{enumerate}

Let us now explain why conditions~\eqref{it:good_loops_in_an_annulus} and~\eqref{it:geodesic_between_loops} imply that $D_\Gamma(\SCL,\SCL') = \widetilde D_\Gamma(\SCL,\SCL')$.  Since $x,y \in \BD_\BQ$ were arbitrary, and since a.s.\ every loop of $\Gamma$ is of the form $\SCL(w)$ for some $w \in \BD_\BQ$,  this will complete the proof of \Cref{thm:main_result}.

First,  we note that it is easy to see using the fact that $P$ is a geodesic and an admissible path that $P \cap \SCL_j = P([s_j,t_j])$ for all $1\leq j \leq N$ (see \Cref{lem:intersection_is_a_segment}),  and hence we have that
\begin{align*}
	D_\Gamma(\SCL,\SCL') = \len(P ; D_\Gamma) = \sum_{j=1}^{N-1} \len(P|_{[t_j,s_{j+1}]} ; D_\Gamma).
\end{align*}
Moreover,  the definition of $D_\Gamma^{0}$ implies that 
\begin{align*}
	D_{\Gamma}^0(\SCL_j ,  \SCL_{j+1}) \leq \len(P|_{[t_j,s_{j+1}]} ; D_\Gamma) \quad \text{for all} \quad 1 \leq j \leq N-1
\end{align*}
and hence
\begin{align*}
	\widetilde D_\Gamma(\SCL,\SCL') \leq \sum_{j=1}^{N-1} D_{\Gamma}^0(\SCL_j,\SCL_{j+1}) \leq \sum_{j=1}^{N-1} \len(P|_{[t_j,s_{j+1}]} ; D_\Gamma) = D_\Gamma(\SCL,\SCL').
\end{align*}
Clearly,  we have that $D_\Gamma(\SCL,\SCL') \leq \widetilde D_\Gamma(\SCL,\SCL')$ and so $D_\Gamma(\SCL,\SCL') = \widetilde D_\Gamma(\SCL,\SCL')$.

Therefore,  the main difficulty in proving \Cref{thm:main_result} is to prove that condition~\eqref{it:geodesic_between_loops} is satisfied.  In particular,  we would like that for all $1 \leq j \leq N-1$,  the following hold.
\begin{enumerate}
    \setcounter{enumi}{2}
	\item \label{it:exactly_two_loops_crossing}
	$\SCL_j,\SCL_{j+1}$ are the only loops in $\Gamma$ that cross $A_j$.
	\item \label{it:geodesic_to_the_boundary_passes_through_the_loops}
	Any $D_\Gamma$-geodesic connecting the inner boundary of $A_j$ with $\partial \BD$,  has to intersect either $\SCL_j$ or $\SCL_{j+1}$.
\end{enumerate}

We would like to prove that with very high probability,  we can find many annuli of the form $B_{2^{-k}}(z) \setminus \overline{B_{2^{-R-k}}(z)}$ for $z \in \BD,  k \in \BN$,  and $R>1$ large but fixed,  for which conditions~\eqref{it:exactly_two_loops_crossing} and ~\eqref{it:geodesic_to_the_boundary_passes_through_the_loops} hold.  We will call these annuli that satisfy conditions~\eqref{it:exactly_two_loops_crossing} and ~\eqref{it:geodesic_to_the_boundary_passes_through_the_loops} \emph{good}.  The first step will be to prove that conditions~\eqref{it:exactly_two_loops_crossing} and ~\eqref{it:geodesic_to_the_boundary_passes_through_the_loops} hold with high probability for the annulus $B_{1/2}(z) \setminus \overline{B_{2^{-R}}(z)}$ provided that $R>1$ is large.  In particular,  we will show in \Cref{prop:good_annulus_macroscopic_scale} that the above probability is at least $1 - 2^{-\alpha R}$ for some universal constant $\alpha > 2$.

The second step will be to use the independence across scales property of the GFF (see \cite[Lemma~3.1]{LocMetGFF}) to argue that for fixed $z \in \BD$,  we have off an event whose probability decays to zero as $K \to \infty$ exponentially in $K$,  that there exists $1 \leq k \leq K$ such that the annulus $B_{2^{-k}}(z) \setminus \overline{B_{2^{-k-R}}(z)}$ satisfies conditions~\eqref{it:exactly_two_loops_crossing} and ~\eqref{it:geodesic_to_the_boundary_passes_through_the_loops} (see \Cref{prop:good_annulus_arbitrary_scale}).   Then, we can take a union bound over all $z$ lying in a suitably dense grid of points.  However, in order to be able to apply the independence across scales property argument, we need to generate the crossings of annuli by loops in $\Gamma$ in a way which is local with respect to some suitable GFF on $\BD$.  For that purpose,  we are going to use a coupling $(\Psi,\Gamma)$ between a zero-boundary GFF $\Psi$ on $\BD$ and $\Gamma$ described e.g.\ in \cite{aru2019bounded}.  In particular,  we will show that for a fixed annulus $A \subseteq \BD$,  the collection of crossings of $A$ made by loops in $\Gamma$ is a.s.\ determined by the restriction of $\Psi$ to $A$ (see \Cref{lem:crossings_determined_by_the_field}).

That condition~\eqref{it:exactly_two_loops_crossing} holds with high probability will follow from two main estimates.  The first estimate states that there exists a universal constant $b>2$ such that for $R>1$ sufficiently large,  we have off an event with probability at most $2^{-bR}$ that there are at most four crossings of $B_{1/2}(0) \setminus \overline{B_{2^{-R}}(0)}$ made by loops in $\Gamma$ (see \Cref{prop:bounded_cle_crossings}).  The second estimate roughly states that if $A \subseteq \BD$ is an annulus of size $2^{-R}$ and we grow two metric balls with respect to $D_\Gamma$ starting from two distinct loops in $\Gamma$ and stopped at the first time that they intersect $A$,  then it is very likely that their $D_\Gamma$-distance is very small (see \Cref{prop:crossing_distances_of_small_annuli_are_small} and \Cref{lem:crossing_distances_of_small_annuli_are_uniformly_small}).  This will imply that there have to be at least two loops in $\Gamma$ crossing $A$ since otherwise,  the $D_\Gamma$-distance between the inner and outer boundaries of $A$ would not be sufficiently small.  Combining the two estimates then yields condition~\eqref{it:exactly_two_loops_crossing}: a loop which crosses $A$ does so at least twice, so the bound of at most four crossings together with the presence of at least two crossing loops forces there to be exactly two such loops, each crossing $A$ exactly twice.

That condition~\eqref{it:geodesic_to_the_boundary_passes_through_the_loops} holds with sufficiently high probability will follow from a more complicated argument (see \Cref{prop:conformal_rectangles,lem:geodesic_intersects_loops}).  The argument goes roughly as follows. Fix $R>1$ large and consider the annulus $A{\defeq}{B_{1/2}(0)} \setminus \overline{B_{2^{-R}}(0)}$.  Let $\Gamma_A^{\mathrm{out}}$ (resp.\ $\Gamma_A^{\mathrm{in}}$) denote the collection of segments of loops in $\Gamma$ which make an excursion from the outer (resp.\ inner) boundary of $A$.  Let also $\Gamma_A$ denote the collection of loops in $\Gamma$ which cross $A$.  Suppose that we are working on the event that the region bounded by the loop in $\Gamma$ containing $0$ does not contain $B_{1/2}(0)$,  and let $U_1,\ldots,U_N$ denote the connected components of the set
\begin{align*}
	A \setminus \overline{\bigcup_{\ell \in \Gamma_A^{\mathrm{in}} \cup \Gamma_A^{\mathrm{out}}} \mathop{\mathrm{int}}(\ell) \cup \bigcup_{\ell \in \Gamma_A} \ell}.
\end{align*}

We will show in \Cref{prop:conformal_rectangles} that there exists a universal constant $b>0$ such that for all sufficiently large $R>1$, the following holds. Conditional on $U_1,\ldots,U_N$ and on the event that the field $\Psi|_{U_j}$ has zero boundary conditions for a fixed $j \in \{1,\ldots,N\}$, there exist, with conditional probability at least $1 - 2^{-bR}$, conformal rectangles $R_{j,m}$ for $1 \leq m \leq 8$ satisfying the following properties:
\begin{enumerate}
    \setcounter{enumi}{4}
	\item \label{it:conformal_rectangles_inside_annuli}
	Each conformal rectangle $R_{j,m}$ is contained in $U_j$, and its left (resp.\ right) boundary is contained in a crossing of $A$ made by a loop in $\Gamma$.
	\item \label{it:geodesic_to_the_boundary_hits_the_two_loops}
	No loop in $\Gamma$ crosses $R_{j,m}$ from its top to its bottom boundary, and the $D_\Gamma$-distance between the left and right sides of $R_{j,m}$ is strictly less than the $D_\Gamma$-distance between its top and bottom sides.
\end{enumerate}

Suppose that conditions~\eqref{it:conformal_rectangles_inside_annuli} and \eqref{it:geodesic_to_the_boundary_hits_the_two_loops} hold. If $P$ is a $D_\Gamma$-geodesic crossing $A$ which does not intersect any loop of $\Gamma$ crossing $A$, then there must exist $j \in \{1,\ldots,N\}$ such that $P$ crosses the conformal rectangle $R_{j,m}$ from its top to its bottom side for all $m = 1,\ldots,8$. However, condition~\eqref{it:geodesic_to_the_boundary_hits_the_two_loops} implies that we can replace this segment of $P$ with another admissible path having the same endpoints but strictly smaller $D_\Gamma$-length. Indeed, we can replace each top-to-bottom crossing $I$ of $R_{j,m}$ made by $P$ with the concatenation of a subsegment of $I$ and a $D_{\Gamma}$-geodesic connecting the left and right sides of $R_{j,m}$. Doing this for all $1 \leq m \leq 8$ allows us to replace the portion of $P$ crossing all eight rectangles with an admissible path having the same endpoints but strictly smaller $D_\Gamma$-length. This contradicts the geodesic property of $P$. It follows that $P$ must intersect a loop in $\Gamma$ which crosses~$A$.

\section{Crossings of annuli}
\label{sec:crossings_of_annuli}

Throughout this section,  we let $\Psi$ be a zero-boundary GFF on $\BD$ and let $\overline{\Gamma}$ be a nested $\CLE_4$ on $\BD$ such that the pair $(\Psi,\overline{\Gamma})$ is coupled as in \Cref{subsec:nested_cle_gff}.   We also let $\Gamma$ denote the collection of outermost loops in $\overline{\Gamma}$,  and we fix $R>1$ and let $A \defeq B_{1/2}(0) \setminus \overline{B_{2^{-R}}(0)}$ denote the corresponding annulus. Our goal in the present section is to study the collection of segments of loops in $\overline{\Gamma}$ which are contained in $A$,  except for their two endpoints,  which lie on $\partial A$.   First we will prove the following measurability statement. We say that a loop $\SCL$ in $\overline{\Gamma}$ has height $(2k+1)\lambda$,  where $k \in \BZ$,  if the boundary values of $\Psi$ on one side of $\SCL$ are equal to $2k\lambda$ and on the other side are equal to $(2k+2)\lambda$.  The height of a segment of a loop is by definition the height of the loop of which it is a part.

\begin{lemma}\label{lem:crossings_determined_by_the_field}
	The collection of segments of loops in $\overline{\Gamma}$ with height in $\{-\lambda,\lambda\}$ which are contained in $A$,  except for their two endpoints,  which lie on $\partial A$,  is a.s.\ determined by the restriction of the field $\Psi$ to $A$.
\end{lemma}

Before we prove \Cref{lem:crossings_determined_by_the_field},  we will give a brief overview of its proof and introduce the setup that we are going to use.

The main idea of the proof of \Cref{lem:crossings_determined_by_the_field} is the following.  First we will construct a deterministic,  countable collection of conformal rectangles $(\SCR_n)_{n \in \BN}$ satisfying the following property.  Let $\SCL \colon [0,1] \to \overline{A}$ be a segment of a loop of $\overline{\Gamma}$ with height in $\{-\lambda,\lambda\}$ such that $\{\SCL(0),\SCL(1)\} \subseteq \partial A$ and $\SCL((0,1)) \subseteq A$.  Then for all $0 < s < t < 1$,  there exists $n \in \BN$ such that $\SCL([s,t])$ is part of a crossing of $\SCR_n \cap A$ of height in $\{-\lambda,\lambda\}$ from its top to its bottom boundary.   Next we will use the exploration procedure introduced in \cite[Section~5]{kavvadias2022conformal}, which recovers the crossings of $\SCR_n \cap A$ in a way which is measurable with respect to $\Psi|_{\SCR_n \cap A}$.  Since $s$ and $t$ were arbitrary,  letting $s \downarrow 0$ and $t \uparrow 1$ we obtain that $\SCL$ is a.s.\ determined by $\Psi|_{A}$.

\subsection{Constructing tubes containing the crossings}
\label{subsubsec:crossings_inside_tubes}

Now we proceed with the construction of the conformal rectangles mentioned above.  Fix distinct points $x,y \in \partial \BD$ and $0 < \delta < 2^{-R-2}|x-y|$.  Let $x_{\delta}^\rL$ be the point on the clockwise arc of $\partial \BD$ from $x$ to $y$ such that $|x_{\delta}^\rL - x| = \delta$ and let $x_{\delta}^\rR$ be the point on the counterclockwise arc of $\partial \BD$ from $x$ to $y$ such that $|x_{\delta}^\rR - x| = \delta$.  Similarly,  we let $y_{\delta}^\rL$ be the point on the clockwise arc of $\partial \BD$ from $y$ to $x$ such that $|y_{\delta}^\rL - y| = \delta$ and let $y_{\delta}^\rR$ be the point on the counterclockwise arc of $\partial \BD$ from $y$ to $x$ such that $|y_{\delta}^\rR - y| = \delta$;  each of these four points is uniquely determined since $\delta < |x-y|/2$.

\begin{figure}[h]
	\centering
	$\vcenter{\hbox{\includegraphics[scale=0.8]{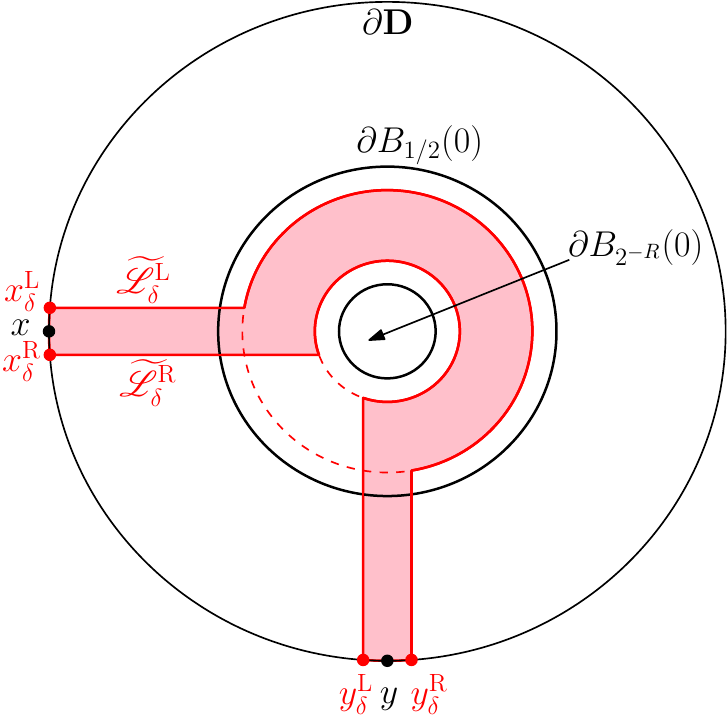}}}$
	\ \
	$\vcenter{\hbox{\includegraphics[scale=0.8]{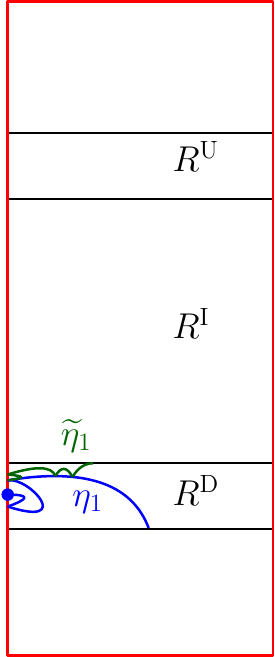}}}$
	\caption{Exploration of the crossings. On the left-hand side, the tube $G$ constructed using the curves $\widetilde{\SCL}_{\delta}^\rL$ and $\widetilde{\SCL}_{\delta}^\rR$. On the right-hand side, an example of the first stage of the exploration of $R^\rD$.}
	\label{fig:constructing_tubes}
\end{figure}

Moreover, we define simple curves $\widetilde{\SCL}_{\delta}^\rL,\widetilde{\SCL}_{\delta}^\rR$ as follows.  First we let $I_{\delta}^\rL(x)$ (resp.\ $I_{\delta}^\rR(x)$) be the segment starting from $x_{\delta}^\rL$ (resp.\ $x_{\delta}^\rR$) which is parallel to $[x,0]$ and points toward the interior of $\BD$,  stopped at the first time that it hits $\partial B_{1/2-\delta}(0)$ (resp.\ $\partial B_{2^{-R}+\delta}(0)$).  Also we let $I_{\delta}^\rL(y)$ (resp.\ $I_{\delta}^\rR(y)$) be the segment starting from $y_{\delta}^\rL$ (resp.\ $y_{\delta}^\rR$) which is parallel to $[y,0]$ and points toward the interior of $\BD$,  stopped at the first time that it hits $\partial B_{2^{-R}+\delta}(0)$ (resp.\ $\partial B_{1/2-\delta}(0)$).  Note that the Euclidean distance from $I_{\delta}^\rL(x)$ and $I_{\delta}^\rR(x)$ to $[x,0]$ (resp.\ from $I_{\delta}^\rL(y)$ and $I_{\delta}^\rR(y)$ to $[y,0]$) is equal to $\delta \sqrt{1-\delta^2/4} \leq \delta$.  Then the curve $\widetilde{\SCL}_{\delta}^\rL$ follows $I_{\delta}^\rL(x)$ up until it hits $\partial B_{1/2-\delta}(0)$,  then follows $\partial B_{1/2-\delta}(0)$ in the clockwise direction up until it hits $I_{\delta}^\rR(y)$,  and finally it follows $I_{\delta}^\rR(y)$ up until it hits $y_{\delta}^\rR$.  Similarly the curve $\widetilde{\SCL}_{\delta}^\rR$ follows $I_{\delta}^\rR(x)$ up until it hits $\partial B_{2^{-R}+\delta}(0)$,  then follows $\partial B_{2^{-R}+\delta}(0)$ in the clockwise direction up until it hits $I_{\delta}^\rL(y)$,  and finally it follows $I_{\delta}^\rL(y)$ up until it hits $y_{\delta}^\rL$.  The assumption that $\delta < 2^{-R-2}|x-y|$ ensures that $\widetilde{\SCL}_{\delta}^\rL$ and $\widetilde{\SCL}_{\delta}^\rR$ are indeed simple and that they are disjoint.  Let $G$ be the conformal rectangle bounded by $\widetilde{\SCL}_{\delta}^\rL ,  \widetilde{\SCL}_{\delta}^\rR$,  the clockwise arc of $\partial \BD$ from $x_{\delta}^\rR$ to $x_{\delta}^\rL$,  and the counterclockwise arc of $\partial \BD$ from $y_{\delta}^\rL$ to $y_{\delta}^\rR$,  with the four marked boundary points $x_{\delta}^\rR ,  x_{\delta}^\rL ,  y_{\delta}^\rR ,  y_{\delta}^\rL$ listed in cyclic order. See the left-hand side of \Cref{fig:constructing_tubes}.  The countable collection of conformal rectangles $(\SCR_n)_{n \in \BN}$ referred to above consists of the tubes $G$ obtained in this way as $x,y$ range over a fixed countable dense subset of $\partial \BD$ and $\delta$ ranges over the rational numbers in $(0,2^{-R-2}|x-y|)$.

\subsection{Exploration of the crossings}
\label{subsubsec:crossings_exploration}

Let $\phi$ be the unique conformal map from $G$ onto the rectangle $R{\defeq}(0,1) \times (0,H)$,  where $H \in (0,\infty)$, such that $\phi(x_{\delta}^\rR) = 0,  \phi(x_{\delta}^\rL) = \ri H,  \phi(y_{\delta}^\rR) = 1+\ri H$,  and $\phi(y_{\delta}^\rL) = 1$.  Fix $\xi \in (0,H/10)$ and set
\begin{align*}
	R^{\rD} = (0,1) \times (2\xi,3\xi), \quad R^{\rU} = (0,1) \times (H-3\xi,H-2\xi),  \quad \text{and} \quad R^{\text{I}} = (0,1) \times (3\xi ,  H-3\xi).
\end{align*}

We will perform a level line exploration of the field $\Psi|_G \circ \phi^{-1} - \lambda$  in the rectangles $R^{\rD},R^{\rU}$ from their left to their right boundaries.  Here we recall that,  conditionally on $\Psi|_{\BD \setminus G}$,  the field $\Psi|_G$ is a GFF on $G$ whose boundary values are equal to $0$ on the two arcs of $\partial G$ which are contained in $\partial \BD$,  that is,  on the two sides of $G$ which are mapped by $\phi$ to the left and to the right side of $R$.  The purpose of the exploration is to generate the crossings of $A$  in a way which is measurable with respect to $\Psi$.  The exploration is very similar to the level line exploration introduced in \cite[Section~5]{kavvadias2022conformal}. 

Now we proceed with the description of the exploration.  We fix $u \in (-2\lambda,0)$ with $2\lambda/u \in \BN$ (e.g., $u \defeq -\lambda$) and let $\eta_1$ be the level line of $\Psi|_G \circ \phi^{-1}  - \lambda$ of height $u$, i.e.\@ the level line of the field $\Psi|_G \circ \phi^{-1}  - \lambda-u$, started from the midpoint of the left boundary of $R^{\rD}$ (resp.\ $R^{\rU}$).  The  choice of $u \in (-2\lambda, 0)$ guarantees that the level line $\eta_1$ is well-defined (in the sense that the weights of the associated force points are larger than $-2$) and that $\eta_1$ intersects the parts of $\partial R^{\rD}$ (resp.\ $\partial R^{\rU}$) immediately to the left and to the right of its starting point (see \cite[Lemma~2.1]{miller2017intersections}).  If $\eta_1$ exits $R^{\rD}$ (resp.\ $R^{\rU}$) on its right boundary,  then the whole procedure terminates,  as explained below.  Otherwise there are two possible outcomes:

\begin{enumerate}
	\item $\eta_{1}$ exits $R^{\rD}$ (resp.\ $R^{\rU}$) on its bottom boundary.  If this occurs,  then we start to explore the level line $\widetilde{\eta}_{1}$ of $-\Psi|_G \circ \phi^{-1} + \lambda+u$ from the point of the left boundary of $R^{\rD}$ (resp.\ $R^{\rU}$) with the largest imaginary part among the points which have been visited by the exploration.  If $\widetilde{\eta}_{1}$ hits the top boundary of $R^{\rD}$ (resp.\ $R^{\rU}$) before hitting either its right or bottom boundary,  then the first stage of the exploration is concluded (as depicted on the right-hand side of \Cref{fig:constructing_tubes}) and we then iterate the exploration as described below.  If $\widetilde{\eta}_{1}$ hits the bottom boundary of $R^{\rD}$ (resp.\ $R^{\rU}$) before hitting either its right or its top boundary,  then we explore the level line of $\Psi|_G \circ \phi^{-1} - \lambda$ of height $u$ from the point of the left boundary of $R^{\rD}$ (resp.\ $R^{\rU}$) with the largest imaginary part among the points which have been visited by the exploration at this point in time.  We repeat this step until the newly generated level line of height $u$ hits the top boundary of $R^{\rD}$ (resp.\ $R^{\rU}$) before hitting either its bottom or its right boundary, which concludes the first stage of the exploration.
	
	\item $\eta_{1}$ exits  $R^{\rD}$ (resp.\ $R^{\rU}$) on its top boundary.  If this occurs,  then we start to explore the level line $\widetilde{\eta}_{1}$ of $-\Psi \vert_G\circ \phi^{-1} + \lambda+u$ from the  point on the left boundary of $R^{\rD}$ (resp.\ $R^{\rU}$) which has the smallest imaginary part among the points which have been visited by the exploration.  If $\widetilde{\eta}_{1}$ hits the bottom boundary of $R^{\rD}$ (resp.\ $R^{\rU}$) before hitting either its right or its top boundary,  then the first stage of the exploration is concluded and next we iterate the exploration as described below.  If $\widetilde{\eta}_{1}$ hits the top boundary of $R^{\rD}$ (resp.\ $R^{\rU}$) before hitting either the bottom or right boundary,  then we explore the level line of $\Psi|_G \circ \phi^{-1} - \lambda$ of height $u$ from the  point of the left boundary of $R^{\rD}$ (resp.\ $R^{\rU}$) which has the smallest imaginary part among the points which have been visited by the exploration at this point in time.  We repeat this step until the newly generated level line of height $u$ hits the bottom boundary of $R^{\rD}$ (resp.\ $R^{\rU}$) before hitting either its top or its right boundary,  concluding the first stage of the exploration.  
\end{enumerate}

We note that if we discover a level line which hits the right boundary of $R^{\rD}$ (resp.\ $R^{\rU}$) in the procedure described above,  then the entire procedure terminates.

Let $K_1$ be the union of the left boundary of $R^{\rD}$ (resp.\ $R^{\rU}$) and the set discovered by the first stage of the exploration.   If the procedure has not been terminated,  we let $U_{1}$ be the connected component of $R^{\rD} \setminus K_1$ (resp.\ $R^{\rU} \setminus K_1$) which has the right boundary of $R^{\rD}$ (resp.\ $R^{\rU}$) as part of its boundary and we let $L_{1} = \overline{\partial U_{1} \setminus \partial R^{\rD}}$ (resp.\ $\overline{\partial U_{1} \setminus \partial R^{\rU}}$).  

We note that $L_{1}$ consists of two level lines,  one of $\Psi|_G \circ \phi^{-1} - \lambda-u$ and one of $-\Psi|_G \circ \phi^{-1} + \lambda+u$,  which both start from the left boundary of $R^{\rD}$ (resp.\ $R^{\rU}$) and one of which exits the rectangle on its top boundary while the other exits on its bottom boundary.  This implies that the boundary conditions of $\Psi|_G \circ \phi^{-1} - \lambda$ along $L_{1}$ are constant and equal to either $-\lambda + u$ or $\lambda + u$.  We say that $L_{1}$ is a crossing of $R^{\rD}$ (resp.\ $R^{\rU}$) from its top to its bottom boundary of height $u$,  where heights are measured with respect to the field $\Psi|_G \circ \phi^{-1} - \lambda$.  Equivalently,  the boundary values of $\Psi|_G \circ \phi^{-1}$ on the two sides of $L_1$ are equal to $u$ and to $2\lambda+u$
.

We define the successive stages of the exploration inductively as follows.  Suppose that $j$ stages have been completed and that the exploration has not yet discovered a level line which hits the right boundary of $R^{\rD}$ (resp.\ $R^{\rU}$).  Let $K_j$ be the union of the left boundary of $R^{\rD}$ (resp.\ $R^{\rU}$) and the set discovered by the exploration after its first $j$ stages,  let $U_j$ be the connected component of $R^{\rD} \setminus K_j$ (resp.\ $R^{\rU} \setminus K_j$) which has the right boundary of $R^{\rD}$ (resp.\ $R^{\rU}$) as part of its boundary,  and let $L_{j} \defeq \overline{\partial U_j \setminus \partial R^{\rD}}$ (resp.\ $\overline{\partial U_j \setminus \partial R^{\rU}}$) be the $j$th crossing and $u_j$ its height. For the iteration step,  we make the following definition.

\begin{enumerate}
	\item If the boundary conditions of $\Psi|_G \circ \phi^{-1} - \lambda$ on $L_{j}$ are equal to $\lambda + u_j$,  then we let $\eta_{j+1}$ be the level line of $\Psi|_G \circ \phi^{-1} - \lambda$ of height $u_{j+1} = u_j + u$ starting from the rightmost point of $L_{j}$ with imaginary part equal to $5\xi/2$ (resp.\ $H - 5 \xi / 2$).
	\item If the boundary conditions of $\Psi|_G \circ \phi^{-1} - \lambda$ on $L_{j}$ are equal to $-\lambda + u_j$,  then we let $\eta_{j+1}$ be the level line of $\Psi|_G \circ \phi^{-1} - \lambda$ of height $u_{j+1} = u_j - u$ starting from the rightmost point on $L_{j}$ with imaginary part equal to $5\xi/2$ (resp.\ $H - 5\xi / 2$).
\end{enumerate}

Then we repeat the steps described in the first stage of the exploration by replacing $R^{\rD}$ (resp.\ $R^{\rU}$) by $U_j$,  the left boundary of $R^{\rD}$ (resp.\ $R^{\rU}$) by $L_{j}$, and $\eta_{1}$ by $\eta_{j+1}$ as follows.  In the choice of where to start the next level line,  we replace the point of intersection with the largest (resp.\ smallest) imaginary part by the last point with respect to the upward (resp.\ downward) exploration of $L_{j}$ which was intersected by the previous level line.  If any of these level lines exits the rectangle $R^{\rD}$ (resp.\ $R^{\rU}$) on its right boundary,  then we stop the exploration at that point.  From this exploration,  we obtain a sequence of crossings $(L_{j})_j$ of $R^{\rD}$ (resp.\ $R^{\rU}$) from its top to its bottom boundary with the property that $L_{j} \cap L_{j+1} \neq \emptyset$ for each $j$.  Note that so that the crossing $L_j$ is part of a loop of $\overline{\Gamma}$ precisely when $u_j \in 2\lambda \BZ$. Moreover, we have the following.

\begin{lemma}\label{lem:finite_number_of_crossings}
	The explorations in $R^{\rD}$ and in $R^{\rU}$ both terminate after finitely many stages a.s.
\end{lemma}

\begin{proof}
	The statement of the lemma follows from the same argument as the one used to prove \cite[Lemma~5.6]{kavvadias2022conformal}.
\end{proof}

For $q \in \{\rD,\rU\}$,  we let $\widetilde{\eta}_i^q$ for $i=1,\ldots,n_q$ be the crossings of $R^q$ with height in $\{-\lambda,\lambda\}$ which are discovered by the exploration,  ordered from left to right, and write $\widetilde{\eta}_0^q$ (resp.\ $\widetilde{\eta}_{n_q+1}^q$) for the left (resp.\ right) side of $R^q$.  Note that each of these crossings will continue into $R^{\text{I}}$ and either form a crossing of $R^{\text{I}}$ or make an excursion into $R^{\text{I}}$ which returns to $R^q$.  We let $\widetilde{\eta}_1,\ldots,\widetilde{\eta}_{\widetilde{n}}$ be the aforementioned crossings of $R^{\text{I}}$,  ordered from left to right,  $\widetilde{\eta}_{\widetilde{n}+2},\ldots,\widetilde{\eta}_{\widetilde{n}+\widetilde{m}+1}$ be the level line excursions into $R^{\text{I}}$ from either $R^{\rD}$ or $R^{\rU}$ which return to the rectangle from which they started,  with height in $\{-\lambda,\lambda\}$,  and let $\widetilde{\eta}_0$ (resp.\ $\widetilde{\eta}_{\widetilde{n}+1}$) be the left (resp.\ right) side of $R^{\text{I}}$.

Now we are ready to prove \Cref{lem:crossings_determined_by_the_field}.

\begin{proof}[Proof of \Cref{lem:crossings_determined_by_the_field}]
	Suppose that we have the setup described in the previous paragraphs.  First we note that the level line interaction rules (see \cite[Theorem~1.1.4]{LevelLineGFFI}) imply that any segment of a loop in $\overline{\Gamma}$ with height in $\{-\lambda,\lambda\}$ which either crosses $\phi^{-1}(R^{\text{I}})$ from its top to its bottom boundary or makes an excursion into $\phi^{-1}(R^{\text{I}})$ from one of these two boundaries back to that same boundary has to be one of the level line segments $\phi^{-1}(\widetilde{\eta}_1),\ldots,\phi^{-1}(\widetilde{\eta}_{\widetilde{n}}),\phi^{-1}(\widetilde{\eta}_{\widetilde{n}+2}),\ldots,\phi^{-1}(\widetilde{\eta}_{\widetilde{n}+\widetilde{m}+1})$.  Moreover, since level lines are locally determined by the field (see \cite[Lemma~2.3.15]{LevelLineGFFI}),  we obtain that the sub-collection of those segments of loops in $\overline{\Gamma}$ which are contained in $A$ is a.s.\ determined by $\Psi|_A$.  Note that by letting $\xi \to 0$ and $\delta \to 0$ along countable sequences and then varying $x,y$ over a fixed countable and dense subset of $\partial \BD$,  we can recover all of the segments of loops in $\overline{\Gamma}$ which stay in $A$.  Combining the above completes the proof of the lemma.
\end{proof}

Next we let $\Gamma_A^{\mathrm{out}}$ (resp.\ $\Gamma_A^{\mathrm{in}}$) be the collection of segments of loops in $\Gamma$ (i.e., outermost loops in $\overline{\Gamma}$) with height in $\{-\lambda,\lambda\}$ which make an excursion into $A$ from the outer (resp.\ inner) boundary of $A$.  Let also $\Gamma_A$ denote the collection of segments of loops in $\overline{\Gamma}$ with height in $\{-\lambda,\lambda\}$ which cross $A$.  Given a loop $\SCL$ of $\overline{\Gamma}$,  we write $\mathop{\mathrm{int}}(\SCL)$ for the region enclosed by $\SCL$.  Then the local finiteness of $\overline{\Gamma}$ (see \cite{CLE}) implies that the complement in $A$ of the set
\begin{align*}
	\overline{\bigcup_{\SCL \in \Gamma_A^{\mathrm{out}} \cup \Gamma_A^{\mathrm{in}}} \mathop{\mathrm{int}}(\SCL)  \cup \bigcup_{\SCL \in \Gamma_A} \SCL }
\end{align*}
has finitely many connected components which we denote by $U_1,\ldots,U_N$.

We will be interested in those $U_j$ which are not surrounded by any loop in $\overline{\Gamma}$ and which are such that $\partial U_j$ contains a crossing of $A$ made by a loop of $\overline{\Gamma}$ with height in $\{\pm \lambda\}$.  First we will show the following.

\begin{lemma}\label{lem:boundary_part_of_loops}
	Fix $j \in \{1,\ldots,N\}$.  If 
	\begin{align*}
		\partial U_j \not \subseteq \overline{\bigcup_{\SCL \in \Gamma_A^{\mathrm{out}} \cup \Gamma_A^{\mathrm{in}} \cup \Gamma_A} \SCL},
	\end{align*}
	then we have that $U_j$ is surrounded by an outermost loop in $\overline{\Gamma}$.
\end{lemma}

\begin{proof}
	First we note that the construction of the coupling between $\Psi$ and $\overline{\Gamma}$ explained in \Cref{subsec:nested_cle_gff} implies that the outermost loops in $\overline{\Gamma}$ have height in $\{-\lambda,\lambda\}$.  Moreover for every fixed point $z \in \BD$,  we have that it is a.s.\ the case that there exists an outermost loop $\ell$ in $\overline{\Gamma}$ such that $\ell$ surrounds $z$.  Applying this to each point of a fixed countable dense subset of $\partial A$,  we obtain that it is a.s.\ the case that every such point is surrounded by an outermost loop in $\overline{\Gamma}$.  Suppose now that $z$ is a point of $\partial U_j$ which does not lie in $\overline{\bigcup_{\SCL \in \Gamma_A^{\mathrm{out}} \cup \Gamma_A^{\mathrm{in}} \cup \Gamma_A} \SCL}$.  Then $z \in \partial A$,  and since $\overline{\bigcup_{\SCL \in \Gamma_A^{\mathrm{out}} \cup \Gamma_A^{\mathrm{in}} \cup \Gamma_A} \SCL}$ is closed,  the same is true for all of the points of $\partial A$ in a neighborhood of $z$;  in particular,  we can take $z$ to belong to the countable dense set fixed above.  Therefore there exists a loop $\ell$ in $\overline{\Gamma}$ of height in $\{-\lambda,\lambda\}$ such that $z$ is surrounded by $\ell$.  But then the region which is enclosed by $\ell$ has to intersect $U_j$.  Since $\ell \cap U_j = \emptyset$ by the definition of $U_j$,  we obtain that $U_j$ is contained in the region which is enclosed by $\ell$.  In particular,  we have that $U_j$ is surrounded by an outermost loop in $\overline{\Gamma}$ and so this completes the proof of the lemma.
\end{proof}

\Cref{lem:boundary_part_of_loops} implies that we will only be interested in connected components whose boundaries are contained in the closure of the union of the segments of loops in $\overline{\Gamma}$ with height in $\{-\lambda,\lambda\}$.  Next we will show that if $U_j$ is a connected component as above such that
\begin{align*}
	\partial U_j  \subseteq \overline{\bigcup_{\SCL \in \Gamma_A^{\mathrm{out}} \cup \Gamma_A^{\mathrm{in}} \cup \Gamma_A} \SCL}
\end{align*}
and the boundary conditions of $\Psi|_{U_j}$ on $\partial U_j$ are not given by the zero function,  then $U_j$ is surrounded by an outermost loop in $\overline{\Gamma}$.

\begin{lemma}\label{lem:zero_boundary_conditions}
	Fix $j \in \{1,\ldots,N\}$ such that
	\begin{align*}
		\partial U_j  \subseteq \overline{\bigcup_{\SCL \in \Gamma_A^{\mathrm{out}} \cup \Gamma_A^{\mathrm{in}} \cup \Gamma_A} \SCL}
	\end{align*}
	and such that $\partial U_j$ contains a crossing of $A$ made by a loop of $\overline{\Gamma}$ with height in $\{\pm \lambda\}$. If the boundary conditions of the field $\Psi|_{U_j}$ on $\partial U_j$ are not given by the zero function, then we have that $U_j$ is contained in the region enclosed by an outermost loop in $\overline{\Gamma}$.
\end{lemma}

\begin{proof}
	Suppose that the assumptions in the statement of the lemma hold.  Then the set $\partial U_j$ contains exactly two arcs $I_j^1$ and $I_j^2$ which both belong to $\Gamma_A$.  We have the following two cases.

	\stepn{step:zbc-case-where-boundary}{The case where the boundary conditions of $\Psi|_{U_j}$ on $\partial U_j \setminus (I_j^1 \cup I_j^2)$ are not equal to zero}
	In that case,  there exists a segment of a loop $\SCL$ in $\Gamma_A^{\mathrm{out}} \cup \Gamma_A^{\mathrm{in}}$ such that $\SCL \subseteq \partial U_j$ and the boundary conditions of $\Psi|_{U_j}$ on $\SCL$ are non-zero.  Since $\SCL$ is part of some loop $\ell$ in $\overline{\Gamma}$ with height in $\{-\lambda,\lambda\}$,  we obtain that the boundary conditions of $\Psi|_{U_j}$ on $\SCL$ are constant and given by either $2\lambda$ or $-2\lambda$ (recall also the coupling between $\Psi$ and $\overline{\Gamma}$ described in \Cref{subsec:nested_cle_gff}).  
	
	If $U_j \cap \mathop{\mathrm{int}}(\ell) \neq \emptyset$,  then the construction of $U_j$ implies that $U_j \subseteq \mathop{\mathrm{int}}(\ell)$,  where we recall that $\mathop{\mathrm{int}}(\ell)$ denotes the region enclosed by $\ell$.  Thus the component $U_j$ is surrounded by an outermost loop in $\overline{\Gamma}$.  If $U_j \cap \mathop{\mathrm{int}}(\ell) = \emptyset$,  then the boundary condition of $\Psi\vert_{U_j}$ on $\ell$ is non-zero. This means that $\ell$ cannot be an outermost loop. So there exists a loop $\overline{\ell}$ in $\overline{\Gamma}$ which is an outermost loop and surrounds $\ell$.  Then it follows that $U_j \cap \mathop{\mathrm{int}}(\overline{\ell}) \neq \emptyset$ and so $\overline{\ell}$ surrounds $U_j$.

	\stepn{step:zbc-case-where-boundary-2}{The case where the boundary conditions of $\Psi\vert_{U_j}$ on $I_j^1 \cup I_j^2$ are not equal to zero}
	Let $\ell_j^1$ (resp.\ $\ell_j^2$) be the loop in $\overline{\Gamma}$ containing $I_j^1$ (resp.\ $I_j^2$).  Then the boundary conditions of $\Psi|_{U_j}$ are equal to $2\lambda$ or to $-2\lambda$ on at least one of $I_j^1$ and $I_j^2$.  Without loss of generality,  we can assume that the former case holds since a similar argument will work for the latter case.  Moreover we order the crossings of $A$ made by loops in $\overline{\Gamma}$ in counterclockwise order starting from $I_j^1$ and assume without loss of generality that $I_j^2$ lies to the right of $I_j^1$.
	
	Suppose that the boundary conditions of $\Psi|_{\mathop{\mathrm{int}}(\ell_j^1)}$ on $\ell_j^1$ are given by zero.  Then the loop $\ell_j^1$ is not an outermost loop and so we let $\overline{\ell}_j^1$ be the outermost loop in $\overline{\Gamma}$ surrounding $\ell_j^1$.  If $\overline{\ell}_j^1$ does not surround $\ell_j^2$,  we must have that $\overline{\ell}_j^1$ separates $I_j^1$ from $I_j^2$ and so $\overline{\ell}_j^1$ makes a crossing of $A$ lying between $I_j^1$ and $I_j^2$.  But the crossings made by $\overline{\ell}_j^1$ have height in $\{-\lambda,\lambda\}$ and so we obtain a contradiction since there are no crossings in $\Gamma_A$ intersecting $U_j$.  Hence we must have that $\overline{\ell}_j^1$ surrounds $\ell_j^2$ and so $U_j \subseteq \mathop{\mathrm{int}}(\overline{\ell}_j^1)$.

	Suppose that the boundary conditions of $\Psi|_{\mathop{\mathrm{int}}(\ell_j^1)}$ on $\ell_j^1$ are given by either $2\lambda$ or $-2\lambda$.  Then the right side of $I_j^1$ (when seen as a collection of prime ends in $\partial U_j$) is contained in $\partial \mathop{\mathrm{int}}(\ell_j^1)$, which implies that $U_j \cap \mathop{\mathrm{int}}(\ell_j^1) \neq \emptyset$.  Hence we obtain that $U_j$ is surrounded by $\ell_j^1$.

	Combining everything,  we obtain that $U_j$ is surrounded by an outermost loop in $\overline{\Gamma}$ in both cases and so this completes the proof of the lemma.
\end{proof}

\begin{remark}\label{rem:no_crossings}
	We note that if the boundary conditions of $\Psi|_{U_j}$ are given by the zero function,  then the conditional law given $\partial U_j$ of the collection of loops in $\overline{\Gamma}$ which are contained entirely in $U_j$ is that of a nested $\CLE_4$ on $U_j$. Moreover, if there are no crossings of the annulus by loops of $\overline{\Gamma}$ with height in $\{-\lambda,\lambda\}$ and $\widehat{\Gamma}$ denotes the collection of loops in $\overline{\Gamma}$ with height in $\{-\lambda,\lambda\}$ and which intersect $\partial A \cup (2^{-R} ,  \frac{1}{2})$ (where the latter denotes the radial slit $\{ t : 2^{-R} < t < \frac{1}{2}\} \subseteq \BR$),  we have that the set
	\begin{align*}
		U_1 \defeq  A \setminus \overline{\bigcup_{\SCL \in \widehat{\Gamma}} \mathop{\mathrm{int}}(\SCL)}
	\end{align*}
	is simply connected and the conditional law given $\partial U_1$ of the loops in $\Gamma$ which are contained entirely in $U_1$ is that of a non-nested $\CLE_4$ on $U_1$,  where $\Gamma$ denotes the collection of outermost loops in $\overline{\Gamma}$. More generally, let $\SP_A$ be the (countable) collection of polygonal paths connecting the outer boundary to the inner boundary of $A$ and made of finitely many segments whose endpoints are chosen in a fixed arbitrary countable dense subset of $A$. On the event that there are no crossings of $A$ by loops of $\overline{\Gamma}$ with height in $\{-\lambda, \lambda\}$, for all $\gamma \in \SP_A$, we denote by $\widehat{\Gamma}_{\gamma}$ the collection of loops in $\overline{\Gamma}$ wih height in $\{-\lambda, \lambda\}$ which intersect $\partial A \cup \gamma$. Then the set
	\[
	U_\gamma \defeq A \setminus \overline{ \bigcup_{\SCL \in \widehat{\Gamma}_\gamma} \mathrm{int}(\SCL)}
	\]
	is simply connected and the conditional law given $\partial U_\gamma$ of the loops in $\Gamma$ which are contained in $U_\gamma$ is that of a non-nested $\CLE_4$ on $U_\gamma$. The rectangles $U_\gamma$ will not be used in the present section but will be useful in \Cref{sec:good_annuli}.
\end{remark}

\begin{figure}[h]
	\centering
	$\vcenter{\hbox{\includegraphics[scale=0.8]{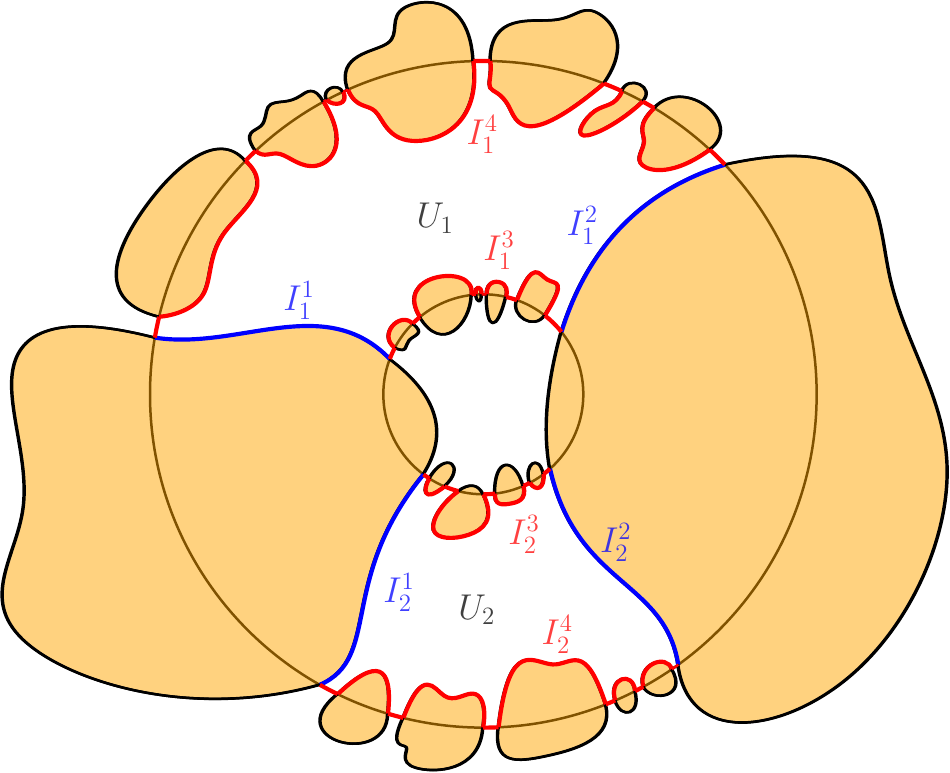}}}$

	\caption{A typical configuration for the connected components $U_j$ which will be used in \Cref{prop:conformal_rectangles}. In this example, exactly two loops cross the annulus, and each of them crosses it exactly twice, hence creating two simply connected regions $U_1$ and $U_2$ in the complement.}
	\label{fig:constructing_conformal rectangles1}
\end{figure}

\begin{remark}\label{rem:crossings_general_case}
	Next we note that we can perform the above exploration of crossings of $A$ in exactly the same way except that,  for fixed $n \in \BZ$, we explore the crossings of $A$ made by loops of $\overline{\Gamma}$ with height in $\{(2n-1) \lambda ,  (2n+1) \lambda\}$.  Again we let $\Gamma_{n,A}^{\mathrm{out}}$ (resp.\ $\Gamma_{n,A}^{\mathrm{in}}$) denote the collection of segments of loops in $\overline{\Gamma}$ with height in $\{(2n-1) \lambda ,  (2n+1) \lambda\}$  which make an excursion into $A$ from the outer (resp.\ inner) boundary of $A$.  Let also $\Gamma_{n,A}$ denote the collection of segments of loops in $\overline{\Gamma}$ with height in $\{(2n-1) \lambda ,  (2n+1) \lambda\}$  which cross $A$.  Then the local finiteness of $\overline{\Gamma}$ implies that the complement in $A$ of the set
	\begin{align*}
		\overline{\bigcup_{\SCL \in \Gamma_{n,A}^{\mathrm{out}} \cup \Gamma_{n,A}^{\mathrm{in}}} \mathop{\mathrm{int}}(\SCL)  \cup \bigcup_{\SCL \in \Gamma_{n,A}} \SCL }
	\end{align*}
	has finitely many connected components which we denote again by $U_1,\ldots,U_N$. 
	
	Moreover, arguing as in \Cref{lem:boundary_part_of_loops,lem:zero_boundary_conditions},  we obtain that the following is true.  Fix $j \in \{1,\ldots,N\}$ and suppose that either $\partial U_j$ is not contained in the set $\overline{\bigcup_{\SCL \in \Gamma_{n,A}^{\mathrm{out}} \cup \Gamma_{n,A}^{\mathrm{in}} \cup \Gamma_{n,A}} \SCL}$,  or it is a subset of $\overline{\bigcup_{\SCL \in \Gamma_{n,A}^{\mathrm{out}} \cup \Gamma_{n,A}^{\mathrm{in}} \cup \Gamma_{n,A}} \SCL}$ and the boundary conditions of $\Psi|_{U_j}$ on $\partial U_j$ are not equal to $2 \lambda n$,  and that in both cases $\partial U_j$ contains a crossing of $A$ made by a loop of $\overline{\Gamma}$ with height in $\{(2n-1)\lambda ,  (2n+1)\lambda\}$.  Then the component $U_j$ is surrounded by an outermost loop in $\overline{\Gamma}$.  Furthermore, for a component $U_j$ such that the boundary conditions of $\Psi|_{U_j}$ on $\partial U_j$ are given by $2\lambda n$,  we have that the conditional law given $\partial U_j$ of the collection of loops in $\overline{\Gamma}$ which are contained entirely in $U_j$ is that of a nested $\CLE_4$ on $U_j$.  Also the set $\partial U_j$ consists of four marked boundary arcs $I_j^1,I_j^2,I_j^3,I_j^4$  with the following properties.  The arcs $I_j^1,I_j^2$ intersect both the inner and the outer boundary of $A$ and belong to $\Gamma_{n,A}$.  The arcs $I_j^3,I_j^4$ are disjoint and are contained in $\overline{\bigcup_{\SCL \in \Gamma^\mathrm{in}_{n,A} }\SCL}$ and $\overline{\bigcup_{\SCL \in \Gamma^\mathrm{out}_{n,A} }\SCL}$, respectively. See \Cref{fig:constructing_conformal rectangles1} for an example.	
	
	In the case that there are no crossings of the annulus by loops with height in $\{(2n-1)\lambda ,  (2n+1)\lambda\}$,  we have that the set
	\begin{align*}
		U_{1}{\defeq}A \setminus \overline{\bigcup_{\SCL \in \widehat{\Gamma}_n} \mathop{\mathrm{int}}(\SCL)}
	\end{align*}
	is simply connected,  where $\widehat{\Gamma}_n$ denotes the collection of loops in $\overline{\Gamma}$ with height either $(2n+1)\lambda$ or $(2n-1) \lambda$, and which intersect $\partial A \cup (2^{-R},\frac{1}{2})$.  Then we have that $\partial U$ consists of four disjoint arcs $I_1,I_2,I_3,I_4$ with the property that both $I_1$ and $I_2$ are contained in the closure of the union of the loops in $\widehat{\Gamma}_n$ which intersect $(2^{-R},\frac{1}{2})$ and $I_3$ (resp.\ $I_4$) is contained in the closure of the union of the loops in $\widehat{\Gamma}_n$ which intersect the inner (resp.\ outer) boundary of $A$. (Indeed, we may choose $I_1$ and $I_2$ to be the maximal arcs that intersect $(2^{-R},\frac{1}{2})$ and are contained, except for their endpoints, in the interior of the annulus, and let $I_3$ and $I_4$ be their complementary components.) Similarly to the case where there are crossings of the annulus by loops,  we have that the conditional law given $\partial U$ of the loops in $\overline{\Gamma}$ which are contained entirely in $U$ is that of a nested $\CLE_4$ on $U$. We define in exactly the same way $U_\gamma$ for a path $\gamma \in \SP_A$ by replacing $(2^{-R}, \frac{1}{2})$ with $\gamma$ in the above definition of the conformal rectangle $U$. The rectangles $U_\gamma$ will not be used in the present section but will be useful in \Cref{sec:good_annuli}.
\end{remark}

\subsection{Constructing conformal rectangles}
\label{subsec:conformal_rectangles}

In the present subsection,  we are going to construct pairwise disjoint conformal rectangles which are contained in the components $U_j$ for $j \in \{1,\ldots,N\}$ introduced in \Cref{subsubsec:crossings_exploration} for which the field $\Psi|_{U_j}$ has zero-boundary conditions,  and such that the conformal rectangles satisfy properties~\eqref{it:conformal_rectangles_inside_annuli} and~\eqref{it:geodesic_to_the_boundary_hits_the_two_loops} introduced in \Cref{subsec:outline_of_the_proof}.  In particular,  we will prove the following result.

\begin{lemma}\label{prop:conformal_rectangles}
	
    Suppose that we have the same setup as in \Cref{rem:no_crossings,rem:crossings_general_case},  where we explore crossings of the nested $\CLE_4$ $\overline{\Gamma}$ on $\BD$ with height in $\{-\lambda ,  \lambda\}$.  Then,  there exist universal constants $b>0$, $R_0>0$,  and $B>0$,  such that the following is true for all $R \geq R_0$.  Conditional on $U_1,\ldots,U_N$,  the following holds with conditional probability at least $1-2^{-b R}$.  Fix $j \in \{1,\ldots,N\}$ such that the boundary conditions of $\Psi|_{U_j}$ on $\partial U_j$ are given by the zero function and such that $\partial U_j$ contains a crossing of $A$ made by a loop of $\overline{\Gamma}$ with height in $\{\pm \lambda\}$.  Suppose also that we are working on the event that $I_j^3 \cap \partial B_{2^{-2R/3}}(0) = \emptyset$ and $I_j^4 \cap \partial B_{2^{-R/3}}(0) = \emptyset$.  Then,  there exists a collection of conformal rectangles $R_{j,m}$ contained in $U_j$ for $1\leq m \leq 8$,  satisfying the following properties (see \Cref{fig:constructing_conformal rectangles2} for an illustration).
	
	\begin{enumerate}
		\item \label{it:disconnecting_conformal_rectangles}
		$R_{j,m+1}$ disconnects $R_{j,m}$ from $I_j^3$ for all $m = 1,2,\ldots,7$ in the following sense.  Any path in $U_j$ which connects $I_j^3$ with $I_j^4$ and intersects $R_{j,m}$,  has to cross $R_{j,m+1}$ from its top boundary to its bottom boundary.  Moreover,  the left (resp.\ right) boundary of $R_{j,m}$ is contained in $I_j^1$ (resp.\ $I_j^2$).
		\item \label{it:left_to_right_crossings_shorter_than_top_to_bottom}
		For $1 \leq m \leq 8$,  we let $\Gamma_{j,m}$ denote the collection of outermost loops in $\overline{\Gamma}|_{U_j}$ with height in $\{-\lambda,\lambda\}$ which intersect $\partial R_{j,m}$.  Then,  there are no loops in $\Gamma_{j,m}$ intersecting both the top and bottom marked boundary arcs of $\partial R_{j,m}$ and the $D^{U_{j,m}}_{\Gamma|_{U_{j,m}}}$-distance between the left / right marked boundary arcs of $\partial U_{j,m}$ is strictly less than $B$ and the  $D^{U_{j,m}}_{\Gamma|_{U_{j,m}}}$-distance between the top / bottom marked boundary arcs of $\partial U_{j,m}$ is at least $B$,  where $U_{j,m}$ denotes the connected component of $R_{j,m} \setminus \overline{\bigcup_{\SCL \in \Gamma_{j,m}} \SCL}$ whose boundary intersects both the left and right marked boundary arcs of $R_{j,m}$. For later use, we also require the following more precise properties. Here and below,  given a conformal rectangle $W$ with four marked boundary arcs,  we write $L(W)$ (resp.\ $R(W)$) for its left (resp.\ right) marked boundary arc.  One can divide the rectangle $U_{j,m}$ from top to bottom into three conformal rectangles $U^1_{j,m},U^2_{j,m}, U^3_{j,m}$ as in \Cref{fig:constructing_conformal rectangles2} such that {
			\begin{enumerate}
				\item \label{condition a genealogy} The $D^{U_{j,m}}_{\Gamma|_{U_{j,m}}}$-distance from the left boundary $L(U^2_{j,m})$ of $U^2_{j,m}$ to the right boundary of $U^2_{j,m}$ is less than $B$;
				\item \label{condition b genealogy} The right boundary (as defined in \Cref{lem:genealogy}) of the metric ball $\SCB_B(L(U_{j,m}), D^{U_{j,m}}_{\Gamma|_{U_{j,m}}})$ intersects the left boundary $L(U^2_{j,m})$;
				\item \label{condition c genealogy} The $D^{U_{j,m}}_{\Gamma|_{U_{j,m}}}$-distance from the bottom side to the top side of $U^1_{j,m}$ (resp.\ $U^3_{j,m}$) is at least $B$.
			\end{enumerate}}
	\end{enumerate}

	\begin{figure}[h]
		\centering
		$\vcenter{\hbox{\includegraphics[scale=0.8]{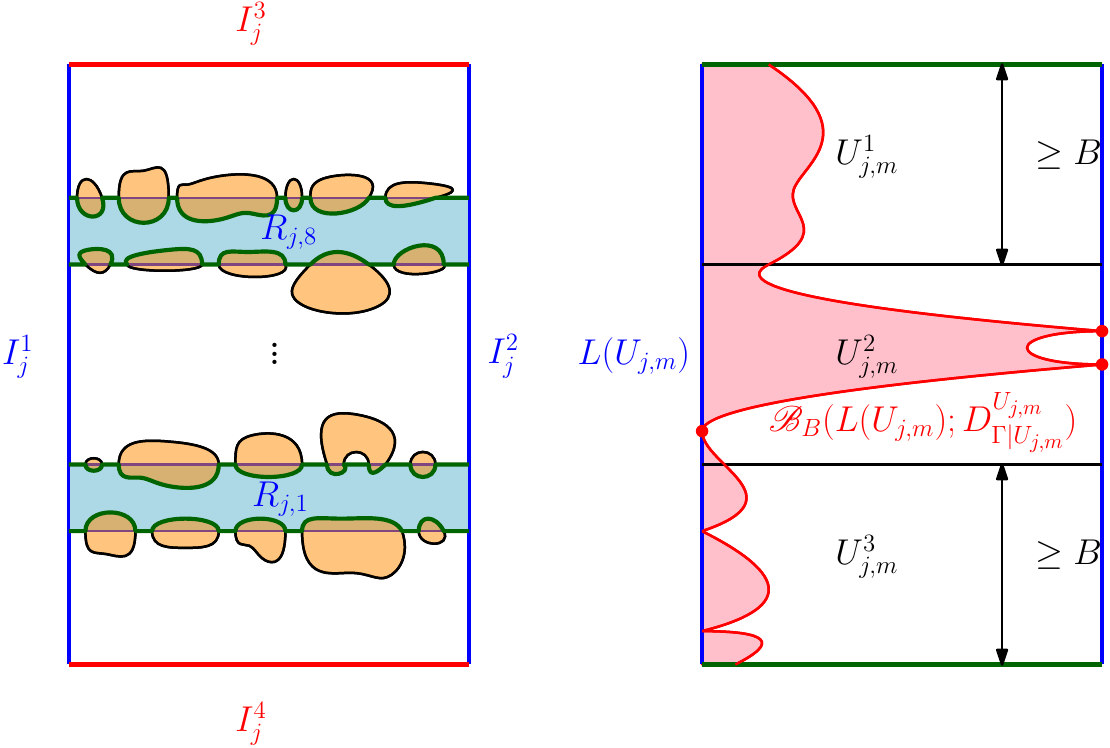}}}$

		\caption{Illustration of \Cref{prop:conformal_rectangles}. On the left-hand side, we draw the connected component $U_j$ with its associated boundaries. The top boundary $I^3_j$ is associated with the inner boundary of the annulus while the bottom boundary $I^4_j$ is associated with the outer boundary of the annulus. The left and right boundaries $I^1_j$ and $I^2_j$ correspond to crossings by loops. The conformal rectangles $R_{j,1}, \ldots, R_{j,8}$ are drawn in blue. The top and bottom sides of the conformal rectangles $U_{j, 1}, \ldots, U_{j,8}$ are drawn in green. On the right-hand side, one can see a conformal rectangle $U_{j,m}$ and its division into three sub-rectangles $U^1_{j,m}, U^2_{j,m}$ and $U^3_{j,m}$ satisfying Conditions~\eqref{condition a genealogy}, \eqref{condition b genealogy} and~\eqref{condition c genealogy}.}
		\label{fig:constructing_conformal rectangles2}
	\end{figure}
	
	Moreover,  in the case that there are no crossings of the annulus by loops,  there exists a collection of conformal rectangles $R_m$ for $1\leq m \leq 8$,  with four distinguished marked boundary arcs satisfying conditions~\eqref{it:disconnecting_conformal_rectangles} and~\eqref{it:left_to_right_crossings_shorter_than_top_to_bottom} with the component $U$ in place of $U_j$ and the arcs $\{I_1,I_2,I_3,I_4\}$ in place of $\{I_j^1,I_j^2,I_j^3,I_j^4\}$.
\end{lemma}

Now let us give a brief outline of the proof of \Cref{prop:conformal_rectangles}.  Suppose that we have the same setup as in \Cref{prop:conformal_rectangles} and fix $j \in \{1,\ldots,N\}$ so that the boundary conditions of $\Psi|_{U_j}$ are given by the zero function.  We let $\phi_j$ denote the unique conformal mapping from the simply connected domain $U_j$ onto the Euclidean rectangle $V_j{\defeq}(0,L_j) \times (0,1)$ for some $L_j \in (0,\infty)$ such that $\phi_j(I_j^1 \cap \partial^{\mathrm{out}} A) = 0$, $\phi_j(I_j^2 \cap \partial^{\mathrm{out}} A) = L_j$, $\phi_j(I_j^2 \cap \partial^{\mathrm{in}} A) = L_j + \ri$,  and $\phi_j(I_j^1 \cap \partial^{\mathrm{in}} A) = \ri$.  We assume that the arc $I_j^2$ lies to the right of $I_j^1$ in the counterclockwise direction.  We also assume that $I_j^3$ (resp.\ $I_j^4$) lies on the inner (resp.\ outer) boundary of $A$ and that $I_j^3 \cap \partial B_{2^{-2 R / 3}}(0) = \emptyset$ and $I_j^4 \cap \partial B_{2^{-R/3}}(0) = \emptyset$.

The first step will be to divide $U_j$ into $\lfloor R \rfloor +1$ conformal rectangles $G_{j,1},\ldots,G_{j,\lfloor R \rfloor +1}$ which satisfy property~\eqref{it:disconnecting_conformal_rectangles} in the statement of \Cref{prop:conformal_rectangles}.  Next,  we will prove in \Cref{lem:rectangles_are_not_wide} that the above conformal rectangles are not too wide in the sense that all of them have conformal modulus at most $1/R$ provided $R>1$ is sufficiently large.

The next step will be to argue that with high probability,  a positive fraction of the $G_{j,n}$'s have the property that there is no loop in $\Gamma$ intersecting both their top and bottom boundaries.  This will be achieved in \Cref{lem:level_line_hits_right_side_with_positive_prob} using \Cref{lem:rectangles_are_not_wide} by arguing that there exists a universal constant $p \in (0,1)$ such that with conditional probability at least $p$,  a level line of $\Psi$ starting from some point on the left boundary of $G_{j,n}$ intersects the right boundary before intersecting either its top or bottom boundaries.  This will prevent the loops in $\Gamma$ from making a crossing of $G_{j,n}$ between its top and bottom boundaries due to the level line interaction rules (see \cite[Theorem~1.1.4]{LevelLineGFFI}).

Finally,  we will show in \Cref{lem:distance_from_left_to_right_less_than_distance_from_bottom_to_top} using \Cref{lem:genealogy} that the following holds with high probability.  Among the conformal rectangles $G_{j,n}$ for which there is no loop in $\Gamma$ making a crossing between their top and bottom boundaries,  a positive fraction of them satisfies the following property.  The $D_\Gamma$-distances between their left and right boundaries are strictly smaller than the $D_\Gamma$-distances between their top and bottom boundaries.  We will deduce \Cref{prop:conformal_rectangles} by combining \Cref{lem:rectangles_are_not_wide,lem:level_line_hits_right_side_with_positive_prob,lem:distance_from_left_to_right_less_than_distance_from_bottom_to_top}.

We start by showing in the next lemma that the conformal rectangle $V_j$ is not too wide.

\begin{lemma}\label{lem:rectangles_are_not_wide}
	Suppose that we have the setup described above.  Then, there exist  universal constants $a \in (0,1)$,   $R_0>0$ such that for all $R \geq R_0$,  we have that $L_j \leq 1 / (aR)$.
\end{lemma}

\begin{proof}[Proof of \Cref{lem:rectangles_are_not_wide}]
	We will argue by contradiction.  We fix $a \in (0,1),  R>1$ (to be chosen later in the proof in a universal way) and assume that $L_j > 1 / (aR)$.  Let $z_j$ denote the center of $V_j$.  Then,  by applying rotation and scaling,  and using the conformal invariance of the complex Brownian motion,  we obtain that the probability that a complex Brownian motion starting from $z_j$ exits $V_j$ on $[\ri,L_j + \ri]$ (resp.\ $[0,L_j]$) is equal to the probability that a complex Brownian motion starting from the center of $(0,L_j^{-1}) \times (0,1)$ first exits $(0,L_j^{-1}) \times (0,1)$ on $[0,\ri]$ (resp.\ $[L_j^{-1} ,  L_j^{-1} +\ri]$).  But \Cref{lem:extremal_length} implies that there exists a universal constant $\widetilde{c}>0$ such that the latter probability is at least $\widetilde{c} e^{-\pi / (2L_j)}$.  Since we have already assumed that $L_j > 1 / (aR)$,  we obtain that the latter probability is at least
	\begin{align}\label{eqn:1}
		\widetilde{c} e^{-\pi{aR} / 2 }.
	\end{align}

	Suppose first that $|\phi_j^{-1}(z_j)| \geq 2^{-R/2}$.  Since we have assumed that $I_j^3 \cap \partial B_{2^{-2R/3}}(0) = \emptyset$,  the Beurling estimate implies that there exists some universal constant $C \in (0,\infty)$ such that the probability that a complex Brownian motion starting from $\phi_j^{-1}(z_j)$ exits $U_j$ on $I_j^3$ is at most $C 2^{-R / 12}$.  Thus,  combining~\eqref{eqn:1} with the conformal invariance of the complex Brownian motion,  we obtain that     \begin{align*}
	\widetilde{c} e^{-\pi a R / 2} \leq C 2^{-R / 12}.
	\end{align*}
	But the latter inequality leads to a contradiction if $a \in (0,1)$ is sufficiently small (in a universal way),  and $R$ is sufficiently large.
	
	Similarly,  if $|\phi_j^{-1}(z_j)| < 2^{-R/2}$,  we have that the probability that a complex Brownian motion starting from $\phi_j^{-1}(z_j)$ exits $U_j$ on $I_j^4$ is at most $C 2^{-R/12}$ for some universal constant $C \in (0,\infty)$,  since $I_j^4 \cap \partial B_{2^{-R/3}}(0) = \emptyset$.  Therefore,  by arguing as in the previous case and using~\eqref{eqn:1},  we obtain a contradiction if $a \in (0,1)$ is sufficiently small (in a universal way) and $R>1$ is sufficiently large.  This completes the proof of the lemma.	
\end{proof}

Next,  we set 
\begin{align*}
	V_{j,k} = (0 ,  L_j) \times \left(\frac{k}{R} ,  \frac{k+1}{R}\right)
\end{align*}
for all $k = 0,1,\ldots,\lfloor R \rfloor - 1$,  and we partition $V_{j,k}$ into $24$ rectangles $V_{j,k}^m$ of the form
\begin{align*}
	V_{j,k}^m = (0 ,  L_j) \times \left(\frac{k}{R} + \frac{m}{24 R} ,  \frac{k}{R} + \frac{m+1}{24 R} \right)
\end{align*}
for $m=0,1,\ldots,23$.  We also fix $\mu \in (0,\lambda)$ and for $m=0,1,\ldots,23$,  we let $a_{j,k}^m$ be the midpoint of the left boundary of $V_{j,k}^m$ and let $\eta_{j,k}^m$ be the level line starting from $a_{j,k}^m$ of height $\mu$ of the field $\Psi|_{U_j} \circ \phi_j^{-1}$.  

We will show in the next lemma that with positive conditional probability,  the curve $\eta_{j,k}^m$ intersects the right boundary of $V_{j,k}^m$ before intersecting either its top or bottom boundaries.

\begin{lemma}\label{lem:level_line_hits_right_side_with_positive_prob}
	Fix $j \in \{1,\ldots,N\}$,  $k \in \{0,1,\ldots,\lfloor R \rfloor - 1\}$ and $m \in \{0,1,\ldots,23\}$.  There exists a constant $p \in (0,1)$ depending only on $\mu$ and the constants $a,R_0$ of \Cref{lem:rectangles_are_not_wide} such that the following holds.  Conditional on $U_1,\ldots,U_N$ and $\eta_{j,\widetilde{k}}^{\widetilde{m}}$ for $(\widetilde{k},\widetilde{m}) \neq (k,m)$,  and on the event that the boundary conditions of $\Psi|_{U_j}$ are given by the zero function,  and that $I_j^3 \cap \partial B_{2^{-2 R /3}}(0) = \emptyset,  I_j^4 \cap \partial B_{2^{-R /3}}(0) = \emptyset$,  we have that the curve $\eta_{j,k}^m$ first intersects the right boundary of the rectangle 
	\begin{align*}
		(0 ,  L_j) \times \left(\vert a_{j,k}^m \vert - \frac{1}{48 R} ,  \vert a_{j,k}^m\vert + \frac{1}{48 R} \right)
	\end{align*}
	before hitting either its top or bottom boundaries,  with conditional probability at least $p$.
\end{lemma}

\begin{proof}

	First, we recall that conditional on $U_1,\ldots,U_N$ and $\eta_{j,\widetilde{k}}^{\widetilde{m}}$ for $(\widetilde{k} ,  \widetilde{m})$ such that either $\widetilde{k} \neq k$ or $\widetilde{k} = k$ and $\widetilde{m} \neq m$,  and on the event that $\Psi|_{U_j}$ has zero-boundary conditions,  we have that the following holds.  The conditional law of the curve $\eta_{j,k}^m$ stopped at the first time that it intersects either the top,  bottom,  or right boundary of the rectangle
	\begin{align*}
		(0 ,  L_j) \times \left( \vert a_{j,k}^m \vert- \frac{1}{48 R} ,  \vert a_{j,k}^m\vert + \frac{1}{48 R} \right)
	\end{align*}
	is that of an SLE$_4(\rho^\rL ; \rho^\rR)$ process in $V_j$ starting from $a_{j,k}^m$ with the force points located at $(a_{j,k}^m)^-$ and $(a_{j,k}^m)^+$ respectively,  and such that $\rho^\rL = -1 + \frac{\mu}{\lambda},  \rho^\rR = -1 - \frac{\mu}{\lambda}$ and the latter process is stopped at the corresponding time.  	
	By rotating, scaling and translating so that $a_{j,k}^m$ is mapped to $0$, we get an SLE$_4(\rho^\rL ; \rho^\rR)$ process in $(\vert a_{j,k}^m \vert R -R, \vert a_{j,k}^m \vert R ) \times (0, L_j R)$ starting from $0$ with the force points located at $0^-$ and $0^+$ respectively, and the process is stopped at the first time it intersects the top, left or right boundary of $(\pm 1/ 48) \times (0, L_j R)$. 	
	Then, notice that since $L_j R \le 1/a$, the (conditional) probability that $\eta^m_{j,k}$ first intersects the right boundary of the rectangle
	\begin{align*}
		(0 ,  L_j) \times \left( \vert a_{j,k}^m \vert - \frac{1}{48 R} ,  \vert a_{j,k}^m \vert + \frac{1}{48 R} \right)
	\end{align*}
	before hitting either its top or bottom boundaries is at least the probability that an SLE$_4(\rho^\rL ; \rho^\rR)$ process in $(\vert a_{j,k}^m \vert R -R, \vert a_{j,k}^m \vert R ) \times (0, 1/a)$ starting from $0$ with the force points located at $0^-$ and $0^+$ respectively, which is stopped at the first time it intersects the top, left or right boundary of $(\pm 1/ 48) \times (0, 1/a)$, hits the top boundary of $(\pm 1/ 48) \times (0, 1/a)$ before hitting the left or right boundaries. 	
	Let $\mathcal{Z}_{j,k}^m$ denote the Radon--Nikodym derivative of the law of this stopped SLE$_4(\rho^\rL ; \rho^\rR)$ process in $(\vert a_{j,k}^m \vert R -R, \vert a_{j,k}^m \vert R) \times (0, L_j R)$ with respect to that of an SLE$_4(\rho^\rL ; \rho^\rR)$ process on the half-plane starting from $0$ with the force points located at $0^-$ and $0^+$ stopped at the first time it intersects the top, left or right boundary of $(\pm 1/ 48) \times (0, 1/a)$. This Radon-Nikodym derivative exists thanks to \cite[Lemma~2.7]{miller2017intersections}. By \cite[Lemma~2.8]{miller2017intersections}, we know that $\mathcal{Z}_{j,k}^m$ is a.s.\ bounded from below and above by a universal constant that depend only on $a,R_0$. Therefore, the claim of the lemma follows from combining with \cite[Lemma~2.5]{miller2017intersections}.\end{proof}

Now,  we fix $k \in \{0,1,\ldots,\lfloor R \rfloor - 1\}$ and suppose that we are working on the event that $\eta_{j,k}^m$ intersects the right boundary of 
\begin{align*}
	(0 ,  L_j) \times \left( \vert a_{j,k}^m \vert - \frac{1}{48 R} ,  \vert a_{j,k}^m \vert+ \frac{1}{48 R} \right)
\end{align*}
before intersecting either its top or bottom boundaries,  for all 
\begin{align*}
	m \in [0, 23]_\BZ \setminus \{3q+1: q \in [0,7]_\BZ \}
\end{align*}

Recall that combining the level line interaction rules (see \cite[Theorem~1.1.4]{LevelLineGFFI}) with the coupling between the nested $\CLE_4$ and the GFF described in \Cref{subsec:nested_cle_gff},  we obtain that there are no loops in $\phi_j(\overline{\Gamma}|_{U_j})$ intersecting both the bottom and top boundaries of $V_{j,k}^{m-1}$,  or both the bottom and top boundaries of $V_{j,k}^{m+1}$ for all
\begin{align*}
	m \in \{3q+1: q \in [0,7]_\BZ \}.
\end{align*}
Let $G_{j,k}^m$ denote the connected subdomain of $V_j$ lying between the outermost loops in $\phi_j(\overline{\Gamma}|_{U_j})$ with height in $\{-\lambda,\lambda\}$ which intersect the bottom boundary of $V_{j,k}^{m-1}$ and the top boundary of $V_{j,k}^{m+1}$ for $m$ as above.  Note that $V_{j,k}^{m} \subseteq G_{j,k}^m$ and let $f_{j,k}^m$ be the conformal map from $G_{j,k}^m$ onto a rectangle $(0,b_{j,k}^m) \times (0,1)$ such that the bottom (resp.\ top) left corner of $\partial G_{j,k}^m$ is mapped to $0$ (resp.\ $\ri$) and the bottom (resp.\ top) right corner of $\partial G_{j,k}^m$ is mapped to $b_{j,k}^m$ (resp.\ $b_{j,k}^m + \ri$).

\begin{lemma}\label{lem:distance_from_left_to_right_less_than_distance_from_bottom_to_top}
	Suppose that we have the setup described above and let $a,R_0$ be the constants in \Cref{lem:rectangles_are_not_wide}.  Fix $j \in \{1,\ldots,N\}$,  $k \in \{0,1,\ldots,\lfloor R \rfloor -1\}$ and $m \in \{3q+1 : q \in [0,7]_\BZ\}$.  Then,  there exist constants $p \in (0,1) ,  B \in (0,\infty)$ depending only on $a,R_0$ such that the following is true.  We condition on $U_1,\ldots,U_N$ and the outermost loops in $\phi_j(\overline{\Gamma}|_{U_j})$ with height in $\{-\lambda ,  \lambda\}$ which intersect either the bottom boundary of $V_{j,k}^{m-1}$ or the top boundary of $V_{j,k}^{m+1}$,  and we assume that we are working on the event that the field $\Psi|_{U_j}$ has zero-boundary conditions,  and that $I_j^3 \cap \partial B_{2^{-2R/3}}(0) = \emptyset,  I_j^4 \cap \partial B_{2^{-R/3}}(0) = \emptyset$.  Then,  with conditional probability at least $p$,  the following hold.  First,  we have that
	\begin{align*}
			D^{f^m_{j,k}(G^m_{j,k})}_{f^m_{j,k}(\phi_j(\Gamma|_{U_j})\vert_{G^m_{j,k}})}
		([0,\ri] ,  [b_{j,k}^m ,  b_{j,k}^m + \ri]) < B \leq 
				D^{f^m_{j,k}(G^m_{j,k})}_{f^m_{j,k}(\phi_j(\Gamma|_{U_j})\vert_{G^m_{j,k}})}
			([0,b_{j,k}^m] ,  [\ri ,  \ri + b_{j,k}^m]).
	\end{align*}
	 More precisely, the rectangle $f^m_{j,k}(G^m_{j,k}) = (0,b^m_{j,k}) \times (0,1)$ can be divided from top to bottom into three conformal rectangles $G^{m,1}_{j,k}$, $G^{m,2}_{j,k}$ and $G^{m,3}_{j,k}$ such that $L(G_{j,k}^{m,i}) \subseteq L(f_{j,k}^m(G_{j,k}^m))$ and $R(G_{j,k}^{m,i}) \subseteq R(f_{j,k}^m(G_{j,k}^m))$ for all $i \in \{1,2,3\}$,  and such that the following hold. 
	\begin{enumerate}[label=({\alph*})]
		\item \label{condition a genealogy1} The $D^{f^m_{j,k}(G^m_{j,k})}_{f^m_{j,k}(\phi_j(\Gamma|_{U_j})\vert_{G^m_{j,k}})}$-distance from the left boundary $L(G^{m,2}_{j,k})$ of $G^{m,2}_{j,k}$ to its right boundary $R(G^{m,2}_{j,k})$ is smaller than $B$;
		\item \label{condition b genealogy1} The right boundary (as defined in \Cref{lem:genealogy}) of the metric ball $\SCB_B(L(f^m_{j,k}(G^{m}_{j,k})),\allowbreak D^{f^m_{j,k}(G^m_{j,k})}_{f^m_{j,k}(\phi_j(\Gamma|_{U_j})\vert_{G^m_{j,k}})})$ intersects $L(G^{m,2}_{j,k})$;
		\item \label{condition c genealogy1} The $D^{f^m_{j,k}(G^m_{j,k})}_{f^m_{j,k}(\phi_j(\Gamma|_{U_j})\vert_{G^m_{j,k}})}$-distance from the bottom boundary to the top boundary of $G^{m,1}_{j,k}$ (resp.\ $G^{m,3}_{j,k}$) is at least $B$.
		\end{enumerate}
	The same statement holds if we condition instead on all of the outermost loops in $\phi_j(\overline{\Gamma}|_{U_j})$ with height in $\{-\lambda,\lambda\}$ which intersect the bottom boundary of $V_{j,k'}^{m'-1}$ or the top boundary of $V_{j,k'}^{m'+1}$ for some $k'$ and some $m' \in \{3q+1 : q \in [0,7]_\BZ\}$,  since the proof below only uses that,  conditionally on the collection in question,  the loops of $\overline{\Gamma}$ contained in $G^m_{j,k}$ form a nested $\CLE_4$ in $G^m_{j,k}$.  It is in this form that we shall use the lemma.
\end{lemma}

\begin{proof}
	\stepn{step:dist-lr-bt-setup}{Setup and the conditional law of $\Gamma_{j,k}^m$} Fix $j \in \{1,\ldots,N\}$ such that the field $\Psi|_{U_j}$ has zero-boundary conditions and set $\Gamma_{j,k}^m \defeq	
	f^m_{j,k}(\phi_j(\Gamma|_{U_j})\vert_{G^m_{j,k}})	$.  Then, the conditional law of $\Gamma_{j,k}^m$ can be described as follows.  Let $\widetilde{\Gamma}$ be a non-nested $\CLE_4$ on $\BD$ and let $g_{j,k}^m$ be the unique conformal map from $\BD$ onto $(0,b_{j,k}^m) \times (0,1)$ which maps $-1$ to $0$ and $0$ to the center of $(0,b_{j,k}^m) \times (0,1)$.  Similarly,  we let $g$ be the unique conformal map from $\BD$ onto $(0,24{/a}) \times (0,1)$ that maps $-1$ to $0$ and $0$ to the center of $(0,24{/a}) \times (0,1)$.  Then,  we can set 
	\begin{align*}
		\Gamma^m_{j,k}		= g_{j,k}^m(\widetilde{\Gamma}).
	\end{align*}
	We emphasize that we use the same sample $\widetilde{\Gamma}$ in $g_{j,k}^m(\widetilde{\Gamma})$ and in $g(\widetilde{\Gamma})$ below.  In this way,  the metrics $D^{f^m_{j,k}(G^m_{j,k})}_{\Gamma^m_{j,k}}$ and $D^{g(\BD)}_{g(\widetilde{\Gamma})}$ are both push-forwards of $D^{\BD}_{\widetilde{\Gamma}}$ under a conformal map,  which is what will allow us to compare them in~\eqref{eqn:metrics_inclusion}.

	\stepn{step:dist-lr-bt-arc-inclusions}{Inclusions of the boundary arcs} We claim that
	\begin{align}\label{eqn:arc_inclusion_1}
		g^{-1}([0,\ri]) \subseteq (g_{j,k}^m)^{-1}([0,\ri]),  \quad g^{-1}([24{/a} ,  24{/a} + \ri]) \subseteq  (g_{j,k}^m)^{-1}([b_{j,k}^m ,  b_{j,k}^m + \ri])
	\end{align}
	and
	\begin{align}\label{eqn:arc_inclusion_2}
		(g_{j,k}^m)^{-1}([\ri,\ri+b_{j,k}^m]) \subseteq g^{-1}([\ri,\ri+24{/a} ]),  \quad (g_{j,k}^m)^{-1}([0,b_{j,k}^m]) \subseteq g^{-1}([0,24{/a} ]).
	\end{align}
	Indeed,  first we note that since $V_{j,k}^m \subseteq G_{j,k}^m$ and these two conformal rectangles have the same left and right marked boundary arcs,  the monotonicity of extremal distance under inclusion together with \Cref{lem:rectangles_are_not_wide} gives that
	\begin{align}\label{eqn:bounded_width}
		b_{j,k}^m \leq 24 R  L_j \leq 24{/a}.
	\end{align}
	Hence, combining~\eqref{eqn:bounded_width} with \Cref{lem:extremal_length},  we obtain that the probability that a complex Brownian motion starting from the center of $(0,b_{j,k}^m) \times (0,1)$ exits $(0,b_{j,k}^m) \times (0,1)$ on $[0,\ri]$ is equal to the probability that it exits $(0,b_{j,k}^m) \times (0,1)$ on $[b_{j,k}^m ,  b_{j,k}^m + \ri]$,  and it is at least the probability that a complex Brownian motion starting from the center of $(0,24{/a}) \times (0,1)$ exits $(0,24{/a}) \times (0,1)$ on $[0,\ri]$,  and the latter is equal to the probability that it exits $(0,24{/a}) \times (0,1)$ on $[24{/a} ,  24{/a} + \ri]$.  Similarly,  the probability that a complex Brownian motion starting from the center of $(0,b_{j,k}^m) \times (0,1)$ exits the rectangle on $[\ri,b_{j,k}^m + \ri]$ (resp.\ $[0,b_{j,k}^m]$) is at most the probability that a complex Brownian motion starting from the center of $(0,24{/a}) \times (0,1)$ exits the rectangle on $[\ri,24{/a} + \ri]$ (resp.\ $[0,24{/a}]$).  Therefore,  we obtain~\eqref{eqn:arc_inclusion_1} and~\eqref{eqn:arc_inclusion_2} by combining with the conformal invariance of complex Brownian motion.
	
	\stepn{step:dist-lr-bt-comparison}{Comparison of the distances} Therefore,  combining~\eqref{eqn:arc_inclusion_1} and~\eqref{eqn:arc_inclusion_2} with the conformal invariance of the metric $D_\Gamma$ (see Axiom~\eqref{it:axiom_conformal_invariance}),  we obtain that a.s.,
	\begin{align}\label{eqn:metrics_inclusion}
		&
			D_{\Gamma^m_{j,k}}^{f^m_{j,k}(G^m_{j,k})}
				([0,\ri] ,  [b_{j,k}^m ,  b_{j,k}^m + \ri]) \leq D^{g(\BD)}_{g(\widetilde{\Gamma})}([0,\ri] ,  [24{/a} ,  24{/a} + \ri]) \quad \text{and} \\ \notag
		&D^{g(\BD)}_{g(\widetilde{\Gamma})}([\ri,\ri+24{/a}] ,  [0,24{/a}]) \leq 
			D_{\Gamma^m_{j,k}}^{f^m_{j,k}(G^m_{j,k})}
				([\ri,\ri+b_{j,k}^m] ,  [0,b_{j,k}^m]).
	\end{align}
	\Cref{lem:genealogy} implies that we may choose the constants $p \in (0,1)$ and $B \in (0,\infty)$ in the statement of the lemma,  which depend only on $a$,  in such a way that with probability at least $p$,  we have that
	\begin{align}\label{eqn:unconditional_comparison_of_distances}
		D^{g(\BD)}_{g(\widetilde{\Gamma})}([0,\ri] ,  [24{/a} ,  24{/a} + \ri]) < B \leq D^{g(\BD)}_{g(\widetilde{\Gamma})}([\ri,\ri+24{/a}] ,  [0,24{/a} ]).
	\end{align}
	Here we take $B \defeq r_0$,  where $r_0$ is the constant from \Cref{lem:genealogy}.  

	\stepn{step:dist-lr-bt-conclusion}{Conclusion of the proof} Thus,  the proof of the lemma is complete by combining~\eqref{eqn:metrics_inclusion} with~\eqref{eqn:unconditional_comparison_of_distances} (the fact that Conditions \ref{condition a genealogy1}, \ref{condition b genealogy1}, and \ref{condition c genealogy1} hold also comes from \Cref{lem:genealogy}: the three sub-rectangles of $(0,24/a) \times (0,1)$ are sent to $(0, b^m_{j,k})\times (0,1)$ using $g^m_{j,k} \circ g^{-1}$,  the rectangles $G_1,G_2,G_3$ of \Cref{lem:genealogy} corresponding to $G^{m,1}_{j,k},G^{m,2}_{j,k},G^{m,3}_{j,k}$ respectively). More precisely, for all $i \in \{1, 2, 3\}$, we set $G^{m,i}_{j,k} \defeq (g^m_{j,k} \circ g^{-1})(G_i)$. However, the sides of $G^{m,i}_{j,k} $ do not correspond to the image of the sides of $G_i$, but we let the left side $L(G^{m,i}_{j,k} )$ (resp.\@ right side $R(G^{m,i}_{j,k} )$) of $G^{m,i}_{j,k} $ be the intersection of $\overline{G^{m,i}_{j,k} }$ with $[0, \ri]$ (resp.\@ $[b^m_{j,k}, b^m_{j,k} + \ri]$). The inclusions \eqref{eqn:arc_inclusion_1} and \eqref{eqn:arc_inclusion_2} show that 
	\[
	g^{-1}(L(G_i)) \subseteq (g^m_{j,k})^{-1}(L(G^{m,i}_{j,k} )) \qquad \text{and} \qquad
	g^{-1}(R(G_i)) \subseteq (g^m_{j,k})^{-1}(  R(G^{m,i}_{j,k} )),
	\]
	where $L(G_i)$ and $R(G_i)$ are the left and right sides of $G_i$. This enables us to check Conditions \ref{condition a genealogy1}, \ref{condition b genealogy1}, and \ref{condition c genealogy1} from \Cref{lem:genealogy}.
\end{proof}

Now we are ready to prove \Cref{prop:conformal_rectangles}.

\begin{proof}[Proof of \Cref{prop:conformal_rectangles}]
    \stepn{step:conf-rect-setup}{Setup and the two properties} \Cref{lem:level_line_hits_right_side_with_positive_prob,lem:distance_from_left_to_right_less_than_distance_from_bottom_to_top} imply that there exist two constant{s} $p \in (0,1)$ and $B>0$ depending only on the constant{s} $R_0{,a}$ in \Cref{lem:rectangles_are_not_wide} such that the following is true for all $k \in \{0,1,\ldots,\lfloor R \rfloor - 1\}$.  Conditionally on $U_1,\ldots,U_N$,  the outermost loops in $\overline{\Gamma}$ with height in $\{-\lambda,\lambda\}$ that intersect the bottom and top boundaries of $\phi_j^{-1}(V_{j,k})$,  and on the event that the field $\Psi|_{U_j}$ has zero-boundary conditions for a fixed $j \in \{1,\ldots,N\}$,  and that $I_j^3 \cap \partial B_{2^{-2R/3}}(0) = \emptyset,  I_j^4 \cap \partial B_{2^{-R/3}}(0) = \emptyset$,  we have that with probability at least $p$,  the following hold for all $m \in \{3q+1: q \in [0,7]_\BZ \}$,  where we have decreased $p$ if necessary so as to take a union over the finitely many values of $m$ (each of the two lemmas above bounds the conditional probability of the corresponding event given all of the others,  so that the simultaneous bound follows upon replacing $p$ by a suitable power of $p$).

	\begin{enumerate}[label=({\alph*})]
		\item \label{it:disjoint_set_of_loops}
		There is no loop in $\phi_j(\overline{\Gamma}|_{U_j})$ with height in $\{-\lambda,\lambda\}$ which intersects both the bottom boundary of $V_{j,k}^{m-1}$ and the top boundary of $V_{j,k}^{m+1}$.
		
		\item \label{it:comparison_of_distances}
		\begin{align*}
			D^{f^m_{j,k}(G^m_{j,k})}_{f_{j,k}^m(\phi_j(\Gamma|_{U_j})|_{G_{j,k}^m})}([0,i] ,  [b_{j,k}^m ,  b_{j,k}^m + i]) < B \leq D^{f^m_{j,k}(G^m_{j,k})}_{f_{j,k}^m(\phi_j(\Gamma|_{U_j})|_{G_{j,k}^m})}([0,b_{j,k}^m] ,  [i ,  i + b_{j,k}^m]),
		\end{align*}
		and more precisely Conditions \ref{condition a genealogy1},\ref{condition b genealogy1} and \ref{condition c genealogy1} of \Cref{lem:distance_from_left_to_right_less_than_distance_from_bottom_to_top} hold.	\end{enumerate}

	\stepn{step:conf-rect-rectangles}{Choice of $k$ and definition of the rectangles $R_{j,m}$} Recall that $G_{j,k}^m$ denotes the connected subdomain of $V_j$ lying between the outermost loops in $\phi_j(\overline{\Gamma}|_{U_j})$ with height in $\{-\lambda,\lambda\}$ which intersect the bottom boundary of $V_{j,k}^{m-1}$ and the top boundary of $V_{j,k}^{m+1}$.  Then, using the independence stemming from Axiom \ref{it:axiom_locality} (locality) for property~\ref{it:comparison_of_distances},  and,  for property~\ref{it:disjoint_set_of_loops},  the fact that \Cref{lem:level_line_hits_right_side_with_positive_prob} bounds the conditional probability given all of the other level lines (note that the $\eta_{j,k}^m$ all have the same height $\mu$,  so that they can merge with each other and are not independent) and hence bounds it given the level lines associated with $V_{j,k'}$ for $k' < k$,  we obtain that the number of $k \in \{0,1,\ldots,\lfloor R \rfloor - 1\}$ for which properties~\ref{it:disjoint_set_of_loops} and~\ref{it:comparison_of_distances} both hold stochastically dominates a binomial random variable with parameters $\lfloor R\rfloor$ and $p$.  Therefore,  there exists a constant $b>0$ depending only on $p$ such that off an event with probability at most $2^{-b R}$,  we have that there exists $k \in \{0,1,\ldots,\lfloor R \rfloor - 1\}$ such that properties~\ref{it:disjoint_set_of_loops} and~\ref{it:comparison_of_distances} hold.

For $k$ as in the previous paragraph,  we set $R_{j,m}$ for all $1 \leq m \leq 8$,  to be the image under $\phi_j^{-1}$ of the rectangle whose left (resp.\@ right) boundary is contained in the left (resp.\@ right) boundary of $V_j$,  and its bottom (resp.\@ top) boundary is given by the bottom (resp.\@ top) boundary of $V_{j,k}^{3(m-1)}$ (resp.\@ $V_{j,k}^{3(m-1)+2}$).  Also, we set $U_{j,m} = \phi_j^{-1}(G_{j,k}^{3(m-1)+1})$ for all $1 \leq m \leq 8$. Then,  it is clear that property~\ref{it:disjoint_set_of_loops} holds by the choice of the $R_{j,m}$'s and $U_{j,m}$'s. 

\stepn{step:conf-rect-conclusion}{Conclusion of the proof} Moreover,  combining with the conformal invariance property of the metric $D_\Gamma$ (Axiom~\eqref{it:axiom_conformal_invariance}),  we obtain that property~\ref{it:comparison_of_distances} holds as well.  It remains to deduce the two properties in the statement of the lemma.  Property~\eqref{it:disconnecting_conformal_rectangles} is immediate from the choice of the $R_{j,m}$'s,  since the $G_{j,k}^{3(m-1)+1}$ are stacked one above the other inside $V_{j,k}$ and each of them has the left (resp.\ right) side of $V_j$ as part of its boundary,  and since the top of $V_j$ corresponds to $I_j^3$.  Property~\eqref{it:left_to_right_crossings_shorter_than_top_to_bottom} then follows from properties~\ref{it:disjoint_set_of_loops} and~\ref{it:comparison_of_distances} together with the conformal invariance of $D_\Gamma$,  the sub-rectangles $U_{j,m}^1,U_{j,m}^2,U_{j,m}^3$ being the images of $G^{3(m-1)+1,1}_{j,k},G^{3(m-1)+1,2}_{j,k},G^{3(m-1)+1,3}_{j,k}$ under $\phi_j^{-1} \circ (f^{3(m-1)+1}_{j,k})^{-1}$. 

This completes the proof of the lemma.
\end{proof}

\begin{remark}\label{rem:conformal_rectangles}
Recall from \Cref{rem:crossings_general_case} that we can perform the exploration of the crossings of $A$ made by loops in $\overline{\Gamma}$ with height in $\{(2n-1) \lambda ,  (2n+1) \lambda\}$ for fixed $n \in \BZ$ in exactly the same way as in the case $n=0$.  Let $U_1,\ldots,U_N$ be the corresponding connected components obtained by these crossings (or excursions into $A$).  Note that for a component $U_j$ such that the boundary conditions of $\Psi$ on $\partial U_j$ are equal to $2\lambda n$,  we have that the conditional law, given $\Psi|_{\BD \setminus U_j}$, of the collection of loops in $\overline{\Gamma}$ which are contained entirely in $U_j$ is the same as that of a nested $\CLE_4$ on $U_j$,  and that (for $n \neq 0$) $U_j$ is necessarily surrounded by an outermost loop in $\overline{\Gamma}$.  Therefore,  arguing exactly as in the proof of \Cref{prop:conformal_rectangles} (applied to the field $\Psi|_{U_j} - 2\lambda n$, which has zero boundary conditions on $\partial U_j$ and whose level loops of height in $\{-\lambda,\lambda\}$ are precisely the loops of $\overline{\Gamma}$ that have height in $\{(2n-1)\lambda, (2n+1)\lambda\}$),  we obtain the same statement for these components $U_j$, with the same constants $b,R_0$,  and $B$;  in particular,  these constants do not depend on $n$.
\end{remark}

\section{Good annuli}
\label{sec:good_annuli}

\subsection{Outline and setup}
\label{subsubsec:outline_and_setup1}

In the present section,  we are going to prove that with very high probability,  there are many annuli which are \emph{good},  i.e.,  which satisfy conditions~\eqref{it:exactly_two_loops_crossing} and~\eqref{it:geodesic_to_the_boundary_passes_through_the_loops} introduced in \Cref{subsec:outline_of_the_proof}.

Throughout the present section,  we let $\Psi$ be a zero-boundary GFF on $\BD$ and let $\overline{\Gamma}$ be a nested $\CLE_4$ on $\BD$ such that the pair $(\Psi,\overline{\Gamma})$ is coupled as in \Cref{subsec:nested_cle_gff},  we let $\Gamma$ denote the collection of outermost loops in $\overline{\Gamma}$,  and,  as in \Cref{sec:crossings_of_annuli},  we fix $R>1$ and let $A \defeq B_{1/2}(0) \setminus \overline{B_{2^{-R}}(0)}$.

First,  we will show that an event $E(\Psi)$,  which encodes the information about $A$ that we shall need in order to deduce that $A$ is \emph{good},  occurs with probability tending to $1$ as $R \to \infty$.  This is the content of the following lemma. Recall from \Cref{rem:no_crossings} the definition of $U_\gamma$ for $\gamma \in \SP_A$ in the case where no loop in $\overline{\Gamma}$ with height in $\{-\lambda,\lambda\}$ crosses $A$. We let $\phi_\gamma\colon U_\gamma \to V_\gamma \defeq (0, L_\gamma) \times(0,1)$ be the unique conformal mapping from $U_\gamma$ that sends the left, right, top and bottom sides of $U_\gamma$ to the respective sides of $V_\gamma \defeq (0, L_\gamma) \times(0,1)$.

\begin{lemma}\label{prop:good_annulus_macroscopic_scale}
	
    There exist universal constants $\alpha>2$ and $R_0 > 1$ such that for all $R \geq R_0$,  there exists a constant $c_R>0$ depending only on $R$ such that the following is true.  Suppose that we have the same setup as in \Cref{prop:conformal_rectangles} and let $U_1,\ldots,U_N$ denote the connected components of the complement in $A$ of the union of the crossings of $A$ made by loops in $\overline{\Gamma}$ with height in $\{-\lambda,\lambda\}$ together with the interiors of the excursions into $A$ from the inner and outer boundary of $A$ made by such loops,  as in \Cref{rem:crossings_general_case} with $n=0$.  In the case that no loop in $\overline{\Gamma}$ with height in $\{-\lambda,\lambda\}$ crosses $A$,  we instead set $N=1$ and let $U_1$ be as in \Cref{rem:no_crossings}. Let $E(\Psi)$ denote the event that the following hold.
	\begin{enumerate}
		\item \label{it:at_most_four_crossings}
		There are at most four crossings of $A$ made by loops in $\overline{\Gamma}$ with height in $\{-\lambda,\lambda\}$.
		\item \label{it:distances_across}
		For every $j \in \{1,\ldots,N\}$ such that the boundary conditions of $\Psi|_{U_j}$ are given by the zero function,  we have that (recall the definitions of $\phi_j$,  $V_j$ and $L_j$ from \Cref{subsec:conformal_rectangles})
		\begin{align*}
			D^{V_j}_{\phi_j(\Gamma|_{U_j})}((0,L_j) \times \{0\} ,  (0,L_j) \times \{1\}) \geq c_R,
		\end{align*}
		and in the case where there is no crossing of $A$ by loops with height in $\{-\lambda, \lambda\}$, we have for all $\gamma \in \SP_A$,
		\begin{align*}
				D^{V_\gamma}_{\phi_\gamma(\Gamma|_{U_\gamma})}((0,L_\gamma) \times \{0\} ,  (0,L_\gamma) \times \{1\}) \geq c_R.
		\end{align*}
		
		Moreover,  in the case that there are exactly four crossings of $A$ made by loops in $\overline{\Gamma}$ with height in $\{-\lambda,\lambda\}$,  there exists for each such $j$ a collection of conformal rectangles $R_{j,m} \subseteq U_j$ for $1 \leq m \leq 8$ for which properties~\eqref{it:disconnecting_conformal_rectangles} and~\eqref{it:left_to_right_crossings_shorter_than_top_to_bottom} in the statement of \Cref{prop:conformal_rectangles} hold,  with $B$ the universal constant appearing there.
	\end{enumerate}
	Then,  we have that $\BP[E(\Psi)] \geq 1 - 2^{-\alpha R}$.
\end{lemma}

\begin{remark}\label{rem:good_event_arbitrary_scale}
	Note that \Cref{lem:crossings_determined_by_the_field} combined with the coupling between $\Psi$ and $\overline{\Gamma}$ described in \Cref{subsec:nested_cle_gff} together with Axiom~\ref{it:axiom_locality} (locality)
	 implies that the event $E(\Psi)$ introduced in \Cref{prop:good_annulus_macroscopic_scale} is a.s.\ determined by $\Psi|_A$ together with the metrics $D^{V_j}_{\phi_j(\Gamma\vert_{U_j})}$ for $1 \leq j \leq N$,  which are conditionally independent given $\Psi|_A$.  Hence, if $\widetilde{\Psi}$ is a GFF on $\BD$ whose boundary conditions are piecewise constant and change only a finite number of times,  then since the laws of $\Psi|_A$ and $\widetilde{\Psi}|_A$ are mutually absolutely continuous (see \cite[Proposition~3.4]{IG1}),  we can define the event $E(\widetilde{\Psi})$ in the exact same way that we defined $E(\Psi)$ but with the field $\Psi$ replaced by $\widetilde{\Psi}$.  
	
	Moreover,  for all $z \in \BD$ and $k \in \BN$ such that $\overline{B_{2^{-k}}(z)} \subseteq \BD$,  we can define the event $E_{z,k}(\Psi)$ in the same way that we defined $E(\Psi)$ except that we explore the crossings of the annulus $B_{2^{-k-1}}(z) \setminus \overline{B_{2^{-k-R}}(z)}$ 
made by the loops in $\overline{\Gamma}$ with height in $\{-\lambda,\lambda\}$ instead of the corresponding crossings of~$A$.  Note that the case $z=0$,  $k=0$ corresponds to the annulus $A$ used in the definition of $E(\Psi)$.  We emphasize that the notation $E_{z,k}(\Psi)$ suppresses the dependence on $R$,  which will always be clear from the context.  Finally,  when the event is considered for a field whose boundary conditions are not identically equal to $0$,  the heights $\{-\lambda,\lambda\}$ above are to be replaced by $\{(2n-1)\lambda ,  (2n+1)\lambda\}$ for the relevant $n \in \BZ$;  see \Cref{rem:crossings_general_case,rem:good_annulus_macroscopic_scale}.
\end{remark}

Next, using \Cref{prop:good_annulus_macroscopic_scale,rem:good_event_arbitrary_scale} combined with the independence across scales property of the GFF (see \Cref{subsec:outline_of_the_proof}),  we will show that there are many \emph{good} annuli with sufficiently high probability.  This is the main result of the present subsection and we state it in the following lemma.

\begin{lemma}\label{prop:good_annulus_arbitrary_scale}
	Let $\alpha>2$ and $R_0>1$ be the constants in \Cref{prop:good_annulus_macroscopic_scale}.  Suppose that we have the same setup as in \Cref{prop:good_annulus_macroscopic_scale} and \Cref{rem:good_event_arbitrary_scale}.  Fix a compact set $F \subseteq \BD$ and let $\nu \in (0,1)$ be such that $\alpha > \max\{2(1+\nu) ,  2 / (1-\nu)^2\}$;  such a $\nu$ exists since $\alpha>2$.
	Fix also $R_1 \geq R_0$. 
	Then,  it is a.s.\ the case that there exist $R>R_1$ and $K_0 \in \BN$ such that for all $K \geq K_0$ and all $z \in \left(\frac{1}{100}2^{-(1+\nu) R K}\BZ\right)^2 \cap F$, there exists $k \in [\nu  K R,  (1+\nu)  K R -R ]_\BZ$ such that $E_{z, k}(\Psi)$ occurs.
\end{lemma}

Let us now give a brief outline of the proof of \Cref{prop:good_annulus_macroscopic_scale}.  Condition~\eqref{it:at_most_four_crossings} in the statement of \Cref{prop:good_annulus_macroscopic_scale} will follow from the following lemma, which states that the probability of having more than four crossings of $A$ made by loops in $\overline{\Gamma}$ is sufficiently small.  Its proof will be given in Section~\ref{sec:cle_4_four_arm_exponents}.

\begin{lemma}\label{prop:bounded_cle_crossings}
	Let $\overline{\Gamma}$ be a nested $\CLE_4$ on $\BD$, let $F \subseteq \BD$ be a compact set,  and set $d \defeq \dist(F ,  \partial \BD)>0$.  Then,  there exists a universal constant $b_0 \in (0,1)$ such that the following is true.  For each $\epsilon \in (0,1)$,  there exists $\delta \in (0,1)$ depending only on $\epsilon$ and $d$ such that for all $z \in F$ and all $0 < r_1 < r_2 < d/2$ with $r_1 / r_2 \leq \delta$,  the probability that there are at least four crossings of $B_{r_2}(z) \setminus \overline{B_{r_1}(z)}$ made by loops in $\overline{\Gamma}$ is at most $(r_1 / r_2)^{2 - \epsilon}$.  Also,  the probability that there are more than four crossings of $B_{r_2}(z) \setminus \overline{B_{r_1}(z)}$ made by loops in $\overline{\Gamma}$ is at most $(r_1 / r_2)^{2 + b_0}$.
\end{lemma}

As for condition~\eqref{it:distances_across} in the statement of \Cref{prop:good_annulus_macroscopic_scale},  we would like to apply \Cref{prop:conformal_rectangles} to deduce the existence of conformal rectangles satisfying the desired properties.  In order to do so,  we have to verify the hypotheses of \Cref{prop:conformal_rectangles} that $I_j^3 \cap \partial B_{2^{-2R/3}}(0) = \emptyset$ and $I_j^4 \cap \partial B_{2^{-R/3}}(0) = \emptyset$.  Since $I_j^3$ (resp.\ $I_j^4$) is contained in the union of the excursions into $A$ from the inner (resp.\ outer) boundary of $A$,  this amounts to knowing that,  with sufficiently high probability,  the following holds.  On the event that there are exactly four crossings of $A$ made by loops in $\overline{\Gamma}$ with height in $\{-\lambda,\lambda\}$,  there is no crossing of either $B_{1/2}(0) \setminus \overline{B_{2^{-R / 3}}(0)}$ or $B_{2^{-2R/3}}(0) \setminus \overline{B_{2^{-R}}(0)}$ made by a loop in $\overline{\Gamma}$ with height in $\{-\lambda,\lambda\}$ other than the ones arising from the four crossings of $A$.  This will be achieved in \Cref{prop:exactly_four_crossings}.

In the proof of \Cref{prop:exactly_four_crossings},  we are going to use the independence across scales property of the GFF (see \Cref{lem:M_good_event_general,lem:R_N_derivative_bounds_general}) to argue that there are many values of $k \in [2,R]_\BZ$ such that 
the following holds.  Almost surely,  the conditional law,  given the values of $\Psi$ on $\BD \setminus B_{2^{-k}}(0)$,  of $\Psi - \lfloor \Psi_{2^{-k}}(0) / (2\lambda) \rfloor 2 \lambda$ restricted to $B_{2^{-k-1}}(0)$ is mutually absolutely continuous with respect to the law of the restriction to $B_{2^{-k-1}}(0)$ of a zero-boundary GFF on $B_{2^{-k}}(0)$,  and the corresponding Radon--Nikodym derivatives in both directions have conditional moments of all orders which are bounded by deterministic constants independent of $k$ and $R$.  Therefore,  by combining these Radon--Nikodym derivative bounds with \Cref{prop:bounded_cle_crossings} applied to a nested $\CLE_4$ on $B_{2^{-k}}(0)$ coupled with a zero-boundary GFF on $B_{2^{-k}}(0)$,  and arguing in the same way for the crossings of $B_{2^{-2 R / 3}}(0) \setminus \overline{B_{2^{-R}}(0)}$,  we will obtain \Cref{prop:exactly_four_crossings}.  Combining the latter with \Cref{prop:bounded_cle_crossings,prop:conformal_rectangles} will then complete the proof of \Cref{prop:good_annulus_macroscopic_scale}.

\subsection{Radon--Nikodym derivatives for restrictions of GFFs} Let us start by stating Radon--Nikodym derivative estimates for restrictions of zero-boundary GFFs. Let $F$ be a compact subset of $\BD$ and set $d\defeq \dist(F, \partial \BD)>0$. Let $k_0\defeq \left\lfloor \frac{\log (d^{-1})}{\log (2)}\right\rfloor +1$ (so that $2^{-k_0}<d$,  and hence $\overline{B_{2^{-k}}(z)} \subseteq \BD$ for all $z \in F$ and $k \geq k_0$) and let $z \in F$. Note that the Markov property of $\Psi$ implies that for all $k \ge k_0$, we can write $\Psi = \Psi_{z,k}^0 + \kh_{z,k}$,  where $\Psi_{z,k}^0$ is a zero-boundary GFF on $B_{2^{-k}}(z)$ and $\kh_{z,k}$ is a random distribution on $\BD$ which,  when restricted to $B_{2^{-k}}(z)$,  is given by a harmonic function and is independent of $\Psi_{z,k}^0$. Let $\mathcal{F}_{z,k}$ denote the $\sigma$-algebra generated by the values of $\Psi$ on $\BD \setminus B_{2^{-k}}(z)$,  and note that $\kh_{z,k}$ is $\mathcal{F}_{z,k}$-measurable.

Fix a constant $M \in (0,\infty)$. We say that the Euclidean ball $B_{2^{-k}}(z)$ for $k \ge k_0$ is $M$-good for $\Psi$ if
\begin{align*}
	\sup\left\{|\kh_{z,k}(w) - \kh_{z,k}(z)| : w \in B_{31 \cdot 2^{-k} / 32}(z)\right\} \leq M
\end{align*}
and we denote this event by $E^M_{z,k}$;  note that $E^M_{z,k} \in \mathcal{F}_{z,k}$.

\begin{lemma}\label{lem:M_good_event_general}
	Fix $a \in (0,\infty),  b \in (0,1)$.  Then there exist constants $c_0,M>0$ depending only on $F,a$ and $b$ such that the following holds for all $K \geq k_0$ and $z \in F$.  Off an event with probability at most $c_0 e^{-a K}$,  there are at least $b K$ values of $k \in [k_0,K]_{\BZ}$ for which $E_{z,k}^M$ occurs.
\end{lemma}
\begin{proof}
	It follows from the exact same argument used to prove  \cite[Proposition~4.3]{GeoLQGnSLE}.  
\end{proof}

Moreover, we have the following Radon--Nikodym estimates. For all $z \in \BD$ and $r>0$ such that $\overline{B_r(z)} \subset \BD$, let us denote by $\Psi_r(z)$ the circle average of $\Psi$ on $\partial B_r(z)$.

\begin{lemma}\label{lem:R_N_derivative_bounds_general}
	For all $z \in F$ and $k\ge k_0$, a.s.\ the conditional law given $\mathcal{F}_{z,k}$ of $\Psi - \lfloor \Psi_{2^{-k}}(z)/(2\lambda) \rfloor 2 \lambda$ restricted to $B_{15 \cdot 2^{-k} / 16}(z)$ is mutually absolutely continuous with respect to the law of $\Psi_{z,k}^0$ restricted to $B_{15 \cdot 2^{-k} / 16}(z)$.  Moreover let $\mathcal{Z}_{z,k}$ (resp.\ $\mathcal{W}_{z,k}$) be the Radon--Nikodym derivative of the former with respect to the latter (resp.\ latter with respect to former).  Then,  for all $p \in \BR$,  there exists a constant $C_{p,M} \in (0,\infty)$ depending only on $p$ and $M$ such that,  a.s.\ on the event $E_{z,k}^M$,  we have that
	\begin{align*}
		&\BE\!\left[ \mathcal{Z}_{z,k}(\Psi_{z,k}^0|_{B_{15 \cdot 2^{-k} / 16}(z)})^p \mid \mathcal{F}_{z,k} \right] \leq C_{p,M}\\
		&\BE\!\left[  \mathcal{W}_{z,k}( (\Psi  - \lfloor \Psi_{2^{-k}}(z)/(2\lambda) \rfloor 2 \lambda) |_{B_{15 \cdot 2^{-k} / 16}(z)})^p \mid \mathcal{F}_{z,k} \right] \leq C_{p,M}.
	\end{align*}
\end{lemma}

\begin{proof}
	The above result follows from the same argument used to prove \cite[Lemma~4.1]{GeoLQGnSLE}.
\end{proof}

\subsection{Proving \Cref{prop:good_annulus_macroscopic_scale}}
\label{subsubsec:good_annuli_macroscopic_scale}

Now we proceed with the proof of \Cref{prop:good_annulus_macroscopic_scale}. We will apply \Cref{lem:M_good_event_general,lem:R_N_derivative_bounds_general} in the particular case $F= \overline{B_{1/2}(0)}$, $z=0$. Note that here $k_0=2$.

We state and prove the main input in the proof of \Cref{prop:good_annulus_macroscopic_scale} (together with \Cref{prop:conformal_rectangles}).

\begin{lemma}\label{prop:exactly_four_crossings}
	Suppose that we have the same setup as in \Cref{prop:good_annulus_macroscopic_scale}.  Then,  there exist universal constants $b>2,R_0'>1$ such that the following holds for all $R \geq R_0'$ (we use a different symbol here since $R_0$ already denotes the constant appearing in \Cref{prop:good_annulus_macroscopic_scale}).  Off an event with probability at most $2^{-b R}$,  we have that,  on the event that there are exactly four crossings of $A$ made by loops in $\overline{\Gamma}$ with height in $\{-\lambda,\lambda\}$,  every crossing of $B_{1/2}(0) \setminus \overline{B_{2^{-R/3}}(0)}$ and every crossing of $B_{2^{-2R/3}}(0) \setminus \overline{B_{2^{-R}}(0)}$ made by a loop in $\overline{\Gamma}$ with height in $\{-\lambda,\lambda\}$ is part of one of these four crossings of $A$.
\end{lemma}

\begin{proof}
	\stepn{step:four-cross-constants}{Choice of constants} Let $b_0 \in (0,1)$ be the constant in \Cref{prop:bounded_cle_crossings}.  Fix $c \in (0,1/3)$ and let $\epsilon \in (0,1/3-c)$ be small enough and $q\in (1,\infty)$ close enough to $1$ so that
	\begin{align}\label{eqn:choice_of_constants}
		\frac{1}{q} \left((2-\epsilon) (1 - c - \epsilon) + (2+b_0) c\right) > 2
	\end{align}
	Let also $a \in (0,1)$ be such that $1-a < \epsilon$.  Then,  \Cref{lem:M_good_event_general} applied with $b$ there equal to the current $a$ and with the constant $a$ of \Cref{lem:M_good_event_general} there taken sufficiently large implies that we can choose a large constant $M>0$ such that for all $R>1$ sufficiently large,  we have off an event with probability at most $2^{-3R}$ that there are at least $aR$ values of $k \in [k_0,R]_\BZ$ (recall that $k_0=2$ here) for which $E_{0,k}^M$ occurs.  For the remainder of the proof,  we assume that this event occurs.

	\stepn{step:four-cross-good-scales}{Good scales and crossing events} Let $k_f$ (resp.\ $k_l$) denote the smallest integer $k$ in $[cR ,  (c+\epsilon) R]$ (resp.\ $[(1-c-\epsilon)R ,  (1-c) R - 1]$) such that $E_{0,k}^M$ occurs;  note that these two windows are disjoint and contained in $[k_0,R]_\BZ$ for $R$ large,  since $\epsilon<1/3-c$ and $c<1/3$ give $c+\epsilon<1/3<1-c-\epsilon$.  Each window contains at least $\epsilon R - 2$ integers,  whereas the number of $k \in [k_0,R]_\BZ$ for which $E^M_{0,k}$ fails is at most $(1-a)R<\epsilon R$.  Since $1-a<\epsilon$ strictly,  it follows that,  provided $R$ is sufficiently large,  $k_f$ and $k_l$ are well defined,  and by construction $k_f \leq (c+\epsilon) R$ and $k_l \geq (1-c-\epsilon) R$.  We emphasize that both $k_f$ and $k_l$ are stopping times for the filtration $(\mathcal{F}_{0,k})_{k \geq k_0}$,  which is increasing in $k$;  this is the reason for taking the smallest,  rather than the largest,  good scale in the second window.

	Let $E_1$ (resp.\ $E_2$) denote the event that there are at least six crossings of $B_{1/4}(0) \setminus \overline{B_{2^{-cR}}(0)}$ (resp.\ $B_{2^{-(1-c) R}}(0) \setminus \overline{B_{2^{-R}}(0)}$) made by loops in $\overline{\Gamma}$.  Similarly,  we let $E_3$ (resp.\ $E_4$) denote the event that there are at least four crossings of $B_{2^{-k_f-2}}(0) \setminus \overline{B_{2^{-R}}(0)}$ (resp.\ $B_{1/4}(0) \setminus \overline{B_{2^{-k_l}}(0)}$) made by loops in $\overline{\Gamma}$.  The reason for the factors $1/4$ and $2^{-k_f-2}$ is that \Cref{prop:bounded_cle_crossings} requires the outer radius to be strictly smaller than half of the distance from the center to the boundary of the ambient domain,  and that \Cref{lem:R_N_derivative_bounds_general} controls the field only on $B_{15 \cdot 2^{-k}/16}(0)$;  since these modifications change the moduli of the annuli involved by an additive constant,  they affect the estimates below only through the implicit constants.  Note that $E_1$ (resp.\ $E_4$) is $\mathcal{F}_{0,k_f}$ (resp.\ $\mathcal{F}_{0,k_l}$)-measurable,  by \Cref{lem:crossings_determined_by_the_field} and since $2^{-k_f} \leq 2^{-cR}$ (resp.\ $2^{-k_l} \geq 2^{-(1-c)R}$),  while $E_3$ (resp.\ $E_2$) is determined by the restriction of $\Psi$ to $B_{2^{-k_f-2}}(0)$ (resp.\ $B_{2^{-(1-c)R}}(0)$). Note also that the level lines of $\Psi_{0,k}^0 + 2\lambda \lfloor \Psi_{2^{-k}}(0) / (2\lambda) \rfloor$ with height in $\{-\lambda,\lambda\}$ correspond to the level lines of $\Psi_{0,k}^0$ with height in $\{-(2n - 1)\lambda ,  -(2n+1) \lambda\}$,  where $n = \lfloor \Psi_{2^{-k}}(0) / (2\lambda) \rfloor$.  In particular,  the nested $\CLE_4$ coupled with $\Psi_{0,k}^0 + 2\lambda n$ consists of the same loops as the one coupled with $\Psi_{0,k}^0$,  only with their heights shifted by $2\lambda n$,  so that the events $E_2$ and $E_3$,  which do not refer to the heights of the loops,  are unaffected by this shift.

	\stepn{step:four-cross-prob-bound}{Bounding $\BP[E_1 \cap E_3]$ and $\BP[E_2 \cap E_4]$} Let $p \in (1,\infty)$ be such that $1/p+1/q=1$,  where $q$ is as above. Then, using H\"older's inequality,  and writing $E_3(\varphi)$ for the event $E_3$ defined in terms of the nested $\CLE_4$ coupled with the field $\varphi$,
		\begin{align*}
			\BP\left[E_1 \cap E_3 \right]&=
			\BE\left[  \one_{E_1} \BP \left[    E_3 \mid  \mathcal{F}_{0,k_f}\right] \right]\\
			&=\BE \left[ \one_{E_1}  \BE\left[ \left. \one_{E_3}\left(\Psi^0_{0,k_f} + \lfloor \Psi_{2^{-k_f}}(0)/(2\lambda) \rfloor 2\lambda\right) \mathcal{Z}_{0,k_f} (\Psi_{0,k_f}^0|_{B_{15 \cdot 2^{-k_f} / 16}(0)})\right\vert \mathcal{F}_{0,k_f} \right]\right]\\
			&\le \BE \left[ \one_{E_1} \BP\left[ \left. E_3\left(\Psi^0_{0,k_f} + \lfloor \Psi_{2^{-k_f}}(0)/(2\lambda) \rfloor 2\lambda\right) \right\vert \mathcal{F}_{0,k_f} \right]^{1/q} C_{p,M}^{1/p} \right]\\
			&\le C_{p,M}^{1/p} \BP[E_1] O\big(2^{-(2-\epsilon)(1-c-\epsilon)R/q}\big)\\
			&\le O(2^{-(1/q)((2+b_0) c + (1-c-\epsilon)(2-\epsilon))R}),
		\end{align*}
		where in the third line we used \Cref{lem:R_N_derivative_bounds_general},  which applies since $E^M_{0,k_f}$ occurs by the definition of $k_f$,  in the fourth line we used \Cref{prop:bounded_cle_crossings} (applied in $B_{2^{-k_f}}(0)$,  which is licit by scale invariance) together with the fact that $k_f \leq (c+\epsilon)R$,  and in the last line we used \Cref{prop:bounded_cle_crossings} again to bound $\BP[E_1] = O(2^{-(2+b_0)cR})$ and then $1/q<1$. Similarly, using the fact that $k_l \geq (1-c-\epsilon)R$, we deduce that
		\begin{align}\label{eqn:bounded_prob}
			\max\{\BP[E_1 \cap E_3] ,  \BP[E_2 \cap E_4]\} \leq O(2^{-(1/q)((2+b_0) c + (1-c-\epsilon)(2-\epsilon))R}),
		\end{align}
		if $R>1$ is sufficiently large,  with the implicit constant depending only on $p,q$ and $M$.

	\stepn{step:four-cross-conclusion}{Conclusion of the proof} Note that if there are exactly four crossings of $A$ made by loops in $\overline{\Gamma}$ with height in $\{-\lambda,\lambda\}$ and there is a further crossing of $B_{1/2}(0) \setminus \overline{B_{2^{-R/3}}(0)}$ or $B_{2^{-2R/3}}(0) \setminus \overline{B_{2^{-R}}(0)}$ which does not cross $A$,  then,  since $c<1/3$, either $E_1 \cap E_3$ or $E_2 \cap E_4$ occurs.  Indeed,  in the first case the four crossings of $A$ each restrict to a crossing of $B_{1/4}(0) \setminus \overline{B_{2^{-cR}}(0)}$ and of $B_{2^{-k_f-2}}(0) \setminus \overline{B_{2^{-R}}(0)}$,  while the further crossing restricts to a crossing of $B_{1/4}(0) \setminus \overline{B_{2^{-cR}}(0)}$ because $cR \leq R/3$;  since a loop which crosses an annulus does so an even number of times,  the further crossing is accompanied by at least one more,  so that $E_1$ occurs and $E_3$ occurs as well.  The second case is symmetric,  using that $(1-c)R \geq 2R/3$.  Hence, the proof of the lemma is complete by combining~\eqref{eqn:choice_of_constants} with~\eqref{eqn:bounded_prob},  and by adding the probability $2^{-3R}$ of the event excluded at the beginning of the proof,  which shows that the total probability in question is at most $2^{-bR}$ for some universal $b>2$ provided $R$ is sufficiently large.
\end{proof}

Now we are ready to prove \Cref{prop:good_annulus_macroscopic_scale}.

\begin{proof}[Proof of \Cref{prop:good_annulus_macroscopic_scale}]
	First,  we note that combining \Cref{prop:bounded_cle_crossings} with \Cref{prop:exactly_four_crossings},  we obtain that there exist universal constants $\alpha>2, R_0 >1$ such that the following holds for all $R \geq R_0$.  Off an event with probability at most $2^{-\alpha R}$,  we have that there are at most four crossings of $A$ made by loops in $\overline{\Gamma}$ with height in $\{-\lambda,\lambda\}$.  Moreover,  on the event that there are exactly four crossings of $A$ made by loops in $\overline{\Gamma}$ with height in $\{-\lambda,\lambda\}$,  we have that there are no other such crossings either of $B_{1/2}(0) \setminus \overline{B_{2^{-R/3}}(0)}$ or $B_{2^{-2R/3}}(0) \setminus \overline{B_{2^{-R}}(0)}$ which do not cross $A$.  The first of these two assertions is precisely condition~\eqref{it:at_most_four_crossings}.

	Note that if the event in the previous paragraph occurs and there are exactly four crossings of $A$ made by loops in $\overline{\Gamma}$ with height in $\{-\lambda,\lambda\}$ (so that the arcs $I_j^3,I_j^4$ are defined),  then we have that
	\begin{align*}
		I_j^3 \cap \partial B_{2^{-2R/3}}(0) = I_j^4 \cap \partial B_{2^{-R/3}}(0) = \emptyset \quad \text{for all} \quad j \in \{1,\ldots,N\}.
	\end{align*}
	Since there are at most four crossings of $A$,  we have that $N \leq 4$,  so that we may apply \Cref{prop:conformal_rectangles} for each $j \in \{1,\ldots,N\}$ separately and take a union bound over $j$.  This gives the second assertion of condition~\eqref{it:distances_across},  after possibly taking $R_0$ to be larger in a universal way.  As for the first assertion,  we have a.s.\ that
	\begin{align*}
		D^{V_j}_{\phi_j(\Gamma|_{U_j})}((0,L_j) \times \{0\} ,  (0,L_j) \times \{1\}) > 0 \quad \text{for all} \quad j \in \{1,\ldots,N\},
	\end{align*}
	and therefore,  for each fixed $R \geq R_0$,  we may choose $c_R>0$ small enough that the probability that the above distance is smaller than $c_R$ for some $j \in \{1,\ldots,N\}$ is at most $2^{-\alpha R}/2$. In the case where there is no crossing of loops of $\overline{\Gamma}$ with height in $\{-\lambda, \lambda\}$, by applying Axiom \eqref{it:weak_axiom_locality} (locality) to $B_{1/2}(0)^\star$, we see that the $D^\BD_{\Gamma}$-distance $d$ between the inner and the outer boundary of $A$ is a.s.\ strictly positive. We then notice that $D^{V_\gamma}_{\phi_\gamma(\Gamma \vert_{U_\gamma})}((0, L_\gamma ) \times \{0\} , (0, L_\gamma) \times \{1\})$ is at least $d$ (since any admissible path from the top to bottom boundary of $U_\gamma$ is an admissible path from the inner to the outer boundary of $A^\star$).
		
	Combining this with the estimates of the previous paragraphs and replacing $\alpha$ by a slightly smaller universal constant which is still larger than $2$ completes the proof of the lemma.
\end{proof}

\begin{remark}\label{rem:good_annulus_macroscopic_scale}
Arguing exactly as in the proof of \Cref{prop:good_annulus_macroscopic_scale} and making use of the observation of \Cref{rem:conformal_rectangles},  we obtain that the statement of \Cref{prop:good_annulus_macroscopic_scale} still holds (with the same constants $\alpha,R_0$,  and $c_R$) if,  for any fixed $n \in \BZ$,  everywhere in its statement we replace the set of heights $\{-\lambda,\lambda\}$ by $\{(2n-1)\lambda, (2n+1)\lambda\}$,  the components $U_1,\dots,U_N$ by the corresponding ones defined in \Cref{rem:conformal_rectangles},  and the condition that the boundary conditions of $\Psi|_{U_j}$ are given by the zero function by the condition that they are equal to $2\lambda n$.
\end{remark}

\subsection{Proving \Cref{prop:good_annulus_arbitrary_scale}}
\label{subsubsec:good_annuli_arbitrary_scales}

Now we will prove \Cref{prop:good_annulus_arbitrary_scale}.  This will follow from combining \Cref{prop:good_annulus_macroscopic_scale} with the independence across scales argument used in the proof of \Cref{prop:exactly_four_crossings}.  More precisely,  we will use  \Cref{lem:M_good_event_general,lem:R_N_derivative_bounds_general} for Euclidean balls centered at points of a fixed compact subset $F$ of $\BD$. Recall that $\Psi_r(z)$ denotes the circle average of the field $\Psi$ on $\partial B_r(z)$.

\begin{proof}[Proof of \Cref{prop:good_annulus_arbitrary_scale}]
	\stepn{step:good-ann-arb-setup}{Setup and choice of constants} Suppose that we have the same setup as in \Cref{lem:M_good_event_general,lem:R_N_derivative_bounds_general} and let $\alpha>2$ and $R_0>1$ be the constants in the statement of \Cref{prop:good_annulus_macroscopic_scale}.  We fix $b \in (2,\alpha),  \delta \in (0,1)$ such that $b(1-\delta - \nu) > 2(1+\nu)$.  Such constants exist:  since $(1+\nu)(1-\nu)^2 \leq 1$,  the assumption $\alpha > 2/(1-\nu)^2$ gives $\alpha > 2(1+\nu)/(1-\nu)$,  so that we may first pick $b \in (2(1+\nu)/(1-\nu) ,  \alpha)$,  which automatically satisfies $b>2$,  and then pick $\delta>0$ sufficiently small.  Throughout the proof we will take $R>1$ to be a sufficiently large constant;  in particular we take $R > \max\{R_0,R_1\}$.

	\stepn{step:good-ann-arb-good-scales}{The good scales $j_m$} Let $M>0$ and $c_0>0$ be the constants obtained by applying \Cref{lem:M_good_event_general} with $a$ there equal to $4$ and $b$ there equal to $1-\delta/2$,  so that $M$ and $c_0$ depend only on $\delta$ and $F$.  Fix $z \in F$ and let $j_1$ be the first $k \in \BN,  k \geq k_0$ such that the event $E_{z,k}^M$ occurs.  Inductively,  given $j_1,\ldots,j_m$,  we let $j_{m+1}$ be the first $k \geq j_m + R$ such that the event $E_{z,k}^M$ occurs,  with the convention that $j_{m+1}=\infty$ if there is no such $k$.  Note that each $j_m$ is a stopping time for the filtration $(\mathcal{F}_{z,k})_{k \geq k_0}$,  which is increasing in $k$.  We fix $K \geq k_0$ and set
	\begin{align*}
		m_K = \max\{m \in \BN : j_m \leq R(K-1)\},
	\end{align*}
	with the convention that $m_K=0$ if there is no such $m$.
	 Note that,  since $j_{m+1}$ is the first good scale which is at least $j_m+R$,  every $k$ for which $E^M_{z,k}$ occurs lies in $[j_m ,  j_m+R)$ for some $m$.  Hence,  if $m_K < (1-\delta) K$,  then the number of $k \in [k_0 , R K]_\BZ$ for which the event $E_{z,k}^M$ occurs is at most $(m_K+1)(R+1) \leq ((1-\delta) K + 1)(R+1) < (1-\delta/2) RK$ for $K \in \BN$ sufficiently large.  Therefore, combining with \Cref{lem:M_good_event_general},  we obtain that off an event with probability at most $c_0 e^{-4RK}$,  we have that $m_K \geq (1-\delta)K$.  For the remainder of the proof,  we assume that we are working on the above event.
	
	\stepn{step:good-ann-arb-estimate}{The conditional estimate at a good scale}

	\substepn{step:gae-comparison}{Comparison with the zero-boundary field} Now,  we choose $q \in (1,\infty)$ such that $2 < b < \alpha / q$ (which is possible since $b<\alpha$) and let $p \in (1,\infty)$ be such that $1 / p + 1/q = 1$.  Let $\widetilde{\Psi}_{z,j_m} \defeq \Psi^0_{z,j_m}$ denote the zero-boundary GFF on $B_{2^{-j_m}}(z)$ appearing in the Markov decomposition of $\Psi$ at scale $j_m$,  which is independent of $\mathcal{F}_{z,j_m}$.  Then,  \Cref{lem:R_N_derivative_bounds_general} combined with H\"older's inequality implies that there exists a constant $C_{p,M}' \in (0,\infty)$ depending only on $p,M$,  and $F$ such that
	\begin{align*}
		\BP[(E_{z,j_m}(\Psi))^c \mid \mathcal{F}_{z,j_m}] \leq C_{p,M}' (\BP[(E_{z,j_m}(\widetilde{\Psi}_{z,j_m} + 2\lambda n_m))^c \mid \mathcal{F}_{z,j_m}])^{1/q}\\
		\text{for all} \quad m \in \BN \quad \text{a.s.\ on $\{j_m<\infty\}$},
	\end{align*}
	where $n_m = \lfloor \Psi_{2^{-j_m}}(z) / (2\lambda) \rfloor$.  Here we used that $E^M_{z,j_m}$ occurs by the definition of $j_m$,  that the annulus $B_{2^{-j_m-1}}(z) \setminus \overline{B_{2^{-j_m-R}}(z)}$ associated with $E_{z,j_m}$ is contained in $B_{15 \cdot 2^{-j_m}/16}(z)$,  so that \Cref{lem:R_N_derivative_bounds_general} applies,  and that,  by \Cref{rem:good_event_arbitrary_scale},  the event $E_{z,j_m}$ is determined by the restriction of the field to that annulus together with metrics whose conditional law given the field is the same under both of the measures being compared.  Note that \Cref{prop:good_annulus_macroscopic_scale} and \Cref{rem:good_annulus_macroscopic_scale} applied on $B_{2^{-j_m}}(z)$,  together with the scale and conformal invariance of the construction, imply that possibly by taking $R$ to be larger,  we can assume that
	\begin{align*}
		\BP[(E_{z,j_m}(\widetilde{\Psi}_{z,j_m}{+ 2\lambda n_m}))^c{ \mid  \mathcal{F}_{z,j_m} }] \leq 2^{-\alpha R}
	\end{align*}
	and hence
	\begin{align*}
		\BP[(E_{z,j_m}(\Psi))^c \mid \mathcal{F}_{z,j_m} ] \leq C_{p,M}' 2^{-\alpha R / q} \quad \text{a.s.}
	\end{align*}

	\substepn{step:gae-domination}{The bound $2^{-bR}$ and stochastic domination} Therefore,  possibly by taking again $R$ to be larger (in a way that depends only on $p,q,\alpha,\delta,b,M$ and $F$),  we can assume that
	\begin{align*}
		\BP[(E_{z,j_m}(\Psi))^c \mid \mathcal{F}_{z,j_m}] \leq 2^{-b R } 
	\end{align*}
	a.s.\ on the event that $j_m < \infty$. In particular, the number $X_K$ of $m \in [\nu K +1, (1-\delta)K ]_\BZ$, such that $j_m=\infty$ or $E_{z,j_m}(\Psi)$ occurs stochastically dominates a binomial random variable with success probability $1-2^{-bR}$ and $\lfloor (1-\delta-\nu)K\rfloor -1$ trials.  Indeed,  since $j_{m+1} \geq j_m + R$,  the annulus $B_{2^{-j_m-1}}(z) \setminus \overline{B_{2^{-j_m-R}}(z)}$ associated with $E_{z,j_m}(\Psi)$ is disjoint from $B_{2^{-j_{m+1}}}(z)$,  meaning that $E_{z,j_m}(\Psi)$ is measurable with respect to the $\sigma$-algebra generated by $\mathcal{F}_{z, j_{m+1}}$ together with the metrics associated with the components of the annuli explored at the scales $j_1,\ldots,j_m$ (recall from \Cref{rem:good_event_arbitrary_scale} that $E_{z,k}(\Psi)$ is not determined by the field alone).  The conditional bound above may therefore be applied iteratively to each index $m$,  given the full history up to time $j_m$,  which includes the outcomes of the events $E_{z,j_i}(\Psi)$ for all $i < m$. 

	\substepn{step:gae-counting}{Counting the good scales} Moreover, a.s.\ on the event that
	\begin{align*}
		|\{k_0 \leq k \leq RK : E_{z,k}^M \,\,\text{occurs}\}| \geq (1-\delta / 2) RK,
	\end{align*}
	we have that the random variable
	\begin{align*}
		|\{ m \in [\nu K + 1 ,  m_K]_\BZ : E_{z,j_m}(\Psi) \,\,  \text{occurs} \}|
	\end{align*}
	is at least $X_K$
		. In particular,  we have that the probability that there is no $m \in [\nu K + 1 ,  m_K]_\BZ$ such that $E_{z,j_m}(\Psi)$ occurs is at most $(2^{-bR})^{T_K} \leq 2^{bR} 2^{-b(1-\delta-\nu)RK}$, plus the previous $c_0 e^{-4RK}$.  Moreover,  since $j_m \geq j_1 + (m-1)R \geq (m-1)R$ and $j_{m_K} \leq R(K-1)$,  every such $m$ satisfies $j_m \in [\nu KR ,  RK - R]_\BZ \subseteq [\nu KR ,  (1+\nu)KR - R]_\BZ$,  as required in the statement of the lemma.  Therefore,  since $b(1-\delta-\nu) >2(1+\nu)$ and since $4/\log 2 >2(1+\nu)$,  the proof of the lemma is complete by taking a union bound over all $z \in F \cap \left(\tfrac{1}{100}2^{-(1+\nu) RK} \BZ\right)^2$,  of which there are $O(2^{2(1+\nu)RK})$ many,  the factor $2^{bR}$ above being harmless since it does not depend on $K$,  summing over $K$ and applying the Borel--Cantelli lemma.\end{proof}

\section{Covering geodesics by good annuli}
\label{sec:covering_geodesics-by_good_annuli}

In the present section,  we will prove \Cref{thm:main_result} and at the end of this section we will also obtain \Cref{prop D tilde does not depend on the metric}.  As explained in \Cref{subsec:outline_of_the_proof},  the main idea behind the proof of \Cref{thm:main_result} is to cover any $D_\Gamma$-geodesic $P$ connecting any two loops in $\Gamma$ by \emph{good} annuli.

As a first step,  we will use \Cref{prop:good_annulus_arbitrary_scale} to find annuli $A_1,\ldots,A_{N-1}$ satisfying the properties of the events $E_{z,k}(\Psi)$ introduced in \Cref{rem:good_event_arbitrary_scale},  in such a way that $P$ is covered by the Euclidean balls bounded by the outer boundaries of the $A_j$'s together with the loops in $\Gamma$ which $P$ intersects.  However,  \Cref{prop:good_annulus_arbitrary_scale} requires that the annuli be centered at points of a grid contained in a deterministic compact set $F \subseteq \BD$,  and so we will need that $P \subseteq F$ for some such $F$.  Since the image of $P$ is a compact subset of $\overline{\BD}$,  the latter follows from the fact that $P$ does not intersect $\partial \BD$,  which was shown in \Cref{cor:geodesic_avoid_boundary}.

Next,  we would like to show that we can choose the annuli $A_1,\ldots,A_{N-1}$ so that there are exactly two loops in $\Gamma$ crossing $A_j$,  for all $1 \leq j \leq N-1$ (see condition~\eqref{it:exactly_two_loops_crossing} in \Cref{subsec:outline_of_the_proof}).  Note first that condition~\eqref{it:at_most_four_crossings} in the definition of the events $E_{z,k}(\Psi)$ implies that there are at most two such loops,  since a loop which crosses an annulus has to cross it at least twice.  On the other hand, the definition of the events $E_{z,k}(\Psi)$ does not exclude the case where we have only one loop in $\Gamma$ or no loops in $\Gamma$ crossing the annulus.  Nevertheless,  condition~\eqref{it:distances_across} in the definition of the $E_{z,k}(\Psi)$'s implies that in both of these cases,  the geodesic $P$ has to travel $D_\Gamma$-distance at least $c_R$ in order to cross the annulus,  where $c_R>0$ is the constant from \Cref{prop:good_annulus_macroscopic_scale}.  We will show in \Cref{lem:crossing_distances_of_small_annuli_are_uniformly_small} (see also \Cref{prop:crossing_distances_of_small_annuli_are_small}) that with very high probability,  the latter does not happen.  Thus,  we will deduce that there are exactly two loops in $\Gamma$ crossing $A_j$ and that each of them crosses $A_j$ exactly twice,  for all $1 \leq j \leq N-1$.

Now,  we note that the existence of the conformal rectangles introduced in \Cref{prop:conformal_rectangles} (which are indexed by the connected components of the complement in $A_j$ of the loops crossing $A_j$,  and not by the annuli themselves) implies that the following hold for all $1 \leq j \leq N-1$ (see \Cref{lem:geodesic_intersects_loops}).  Whenever the geodesic $P$ makes a crossing of $A_j$,  it has to intersect at least one of the two loops in $\Gamma$ crossing $A_j$.  The same argument applies to any $D_\Gamma$-geodesic connecting $\partial \BD$ with a loop in $\Gamma$ which is separated from $\partial \BD$ by the two loops crossing $A_j$ together with the fifth of the conformal rectangles contained in each of the two components of $A_j$ cut out by these loops (the rectangles being numbered starting from the one closest to the outer boundary of $A_j$;  the reason for taking the fifth one is explained in the proof of \Cref{thm:main_result}).  This is the case in particular for every loop in $\Gamma$ which $P$ intersects while making a crossing of $A_j$,  which is the only situation in which we shall use it.  Hence, combining with the fact that replacing an initial segment of a $D_\Gamma$-geodesic by another $D_\Gamma$-geodesic with the same endpoints again produces a $D_\Gamma$-geodesic,  we will deduce that (after possibly discarding some of the loops and of the annuli from the chain,  so that $P$ intersects both of the loops crossing $A_j$ and so that the second of these two loops is also the first of the two loops crossing $A_{j+1}$) the part $P_j$ of $P$ between the last time it intersects the first such loop and the first time it intersects the second satisfies the following property.  Let $t \in [0,1]$ be such that $P_j(t) \in \ell$ for some $\ell \in \Gamma$ (no other case has to be considered in the definition of $D_\Gamma^0$,  since $P$ does not intersect $\partial \BD$).  Then,  there exists a $D_\Gamma$-geodesic from $\ell$ to $\partial \BD$ containing (up to monotone reparameterization) either the time-reversal of $P_j|_{[0,t]}$ or $P_j|_{[t,1]}$,  where we have parameterized $P_j$ by $[0,1]$;  for $t \in \{0,1\}$ the assertion is immediate,  since any path contained in a single loop has zero $D_\Gamma$-length.  This will complete the proof of \Cref{thm:main_result}.

\subsection{Some useful results}
Recall that as explained at the beginning of the present section,  we need that with very high probability,  the following is true.  Whenever the geodesic $P$ and its time-reversal both enter a small Euclidean ball before intersecting each other,  the $D_\Gamma$-length of the part of $P$ which is left after removing the two corresponding initial segments has to be very small.  In particular,  whenever $P$ makes a crossing of a small Euclidean annulus,  the $D_\Gamma$-length of that crossing has to be very small.  This is the content of the following lemma.

\begin{lemma}\label{lem:crossing_distances_of_small_annuli_are_uniformly_small}
	Let $F \subseteq \BD$ be a compact set and fix $x,y \in F$. Then,  there exist a constant $c_2 >0$ (depending only on $x,y$,  and $F$) and,  for each $c>0$,  constants $R_0 = R_0(x,y,F,c)>1$ and $c_1 = c_1(x,y,F,c)>1$,  such that the following is true for all $c>0$ and all $R \geq R_0$. Let $\SCL$ (resp.\ $\SCL'$) denote the loop in $\Gamma$ surrounding $x$ (resp.\ $y$) and suppose that we are working on the event that $\SCL \neq \SCL'$. Let $P$ be a $D_\Gamma$-geodesic connecting $\SCL$ to $\SCL'$ viewed as a path from $\SCL$ to $\SCL'$,  and let $\widetilde{P}$ denote the time-reversal of $P$. Off an event with probability at most $c_1 e^{-c_2 c R}$,  we have that for all $z \in F$,  if the paths $P$ and $\widetilde{P}$ intersect $\overline{B_{2^{-R}}(z)}$ before intersecting each other, {and if we let $t$ (resp.\ $\widetilde{t}$) be the first time that $P$ (resp.\ $\widetilde{P}$) intersects $\overline{B_{2^{-R}}(z)}$, } then we have that $t \leq 1 - \widetilde{t}$ and that $\len(P|_{[t , 1-\widetilde{t}]} ; D_\Gamma)$,  i.e.\ the $D_\Gamma$-length of the part of $P$ which is left after removing $P|_{[0,t]}$ and the time-reversal of $\widetilde{P}|_{[0,\widetilde{t}]}$,  is strictly less than $c$.
\end{lemma}

The main ingredient of the proof of \Cref{lem:crossing_distances_of_small_annuli_are_uniformly_small} is the following lemma.

\begin{lemma}\label{prop:crossing_distances_of_small_annuli_are_small}
	Suppose that we have the same setup as in \Cref{lem:crossing_distances_of_small_annuli_are_uniformly_small}. Let $\epsilon>0$.	 Then,  there exist constants $c_1,R_0>1$ (depending only on $x,y$,  $\epsilon$,  and $F$) and $c_2>0$ (depending only on $x,y$,  and $F$) such that the following is true for all $R \geq R_0$ and all $z \in F$. The probability that both $P$ and $\widetilde{P}$ intersect $\overline{B_{2^{-R}}(z)}$ before intersecting each other is at most $c_1 2^{-(2-\epsilon)R}$.  Moreover,  for all $c>0$ and conditioned on the event that both $P$ and $\widetilde{P}$ intersect $\overline{B_{2^{-R}}(z)}$ before intersecting each other,  the following holds off an event with conditional probability at most $c_1 e^{-c_2 c}$.  Let $t$ (resp.\ $\widetilde{t}$) be the first time that $P$ (resp.\ $\widetilde{P}$) intersects $\overline{B_{2^{-R}}(z)}$.  Then,  the distance between $P|_{[0,t]}$ and $\widetilde{P}|_{[0,\widetilde{t}]}$ with respect to $D_\Gamma$ (equivalently,  since $P$ is a $D_\Gamma$-geodesic,  the quantity $\len(P|_{[t , 1-\widetilde{t}]} ; D_\Gamma)$) is strictly less than $c / R$.
\end{lemma}

\begin{proof}
    This follows immediately from Lemma~\ref{lem:Beurling} and \Cref{prop:distance_across_rectangle,prop:hitting-two-metric-balls} (here, we use that a geometric random variable has an exponential tail). 
\end{proof}

\begin{proof}[Proof of \Cref{lem:crossing_distances_of_small_annuli_are_uniformly_small}]
	Let $F'$ denote the closed $\dist(F,\partial \BD)/2$-neighborhood of $F$,  which is a compact subset of $\BD$.  The claim in the statement of the lemma follows from taking a union bound over all $w \in \left(\frac{2^{-R}}{100} \BZ\right)^2 \cap F'$ and applying \Cref{prop:crossing_distances_of_small_annuli_are_small} with $F'$ in place of $F$,  with $R-1$ in place of $R$,  with $c(R-1)$ in place of $c$,  and with $\epsilon>0$ to be chosen just below,  and noting that for all $z \in F$,  there exists $w \in \left(\frac{2^{-R}}{100} \BZ\right)^2 \cap F'$ such that $\overline{B_{2^{-R}}(z)}$ $\subseteq B_{2^{-(R-1)}}(w)$.  (We work with $F'$ rather than with $F$ since the grid points which are closest to a given $z \in F$ need not lie in $F$.)  Let us also note that the first time at which $P$ (resp.\ $\widetilde{P}$) intersects $\overline{B_{2^{-(R-1)}}(w)}$ is at most the first time at which it intersects $\overline{B_{2^{-R}}(z)}$,  so that the piece of $P$ associated with $z$ is contained in the piece of $P$ associated with $w$ and it therefore suffices to bound the $D_\Gamma$-length of the latter.  Finally,  the number of grid points is $O(2^{2R})$ while the probability that both $P$ and $\widetilde{P}$ intersect $\overline{B_{2^{-(R-1)}}(w)}$ before intersecting each other is at most $c_1 2^{-(2-\epsilon)(R-1)}$,  so that the union bound leaves us with a factor $O(2^{\epsilon R})$;  taking $\epsilon \defeq c_2 c / (2 \log 2)$ (this is the reason why the constant $c_1$ in the statement of the lemma is allowed to depend on $c$) and renaming the constants,  the result follows.
\end{proof}

For the remainder of the present section,  we fix distinct points $x,y \in \BD_\BQ$ and let $\SCL$ (resp.\ $\SCL'$) denote the loop in $\Gamma$ surrounding $x$ (resp.\ $y$).  We assume that we are working on the event that $\SCL \neq \SCL'$ and let $P$ be a $D_\Gamma$-geodesic connecting $\SCL$ to $\SCL'$.  We parameterize $P$ by $[0,1]$ and assume that $P(0) \in \SCL ,  P(1) \in \SCL'$.

Combining \Cref{prop:good_annulus_arbitrary_scale,cor:geodesic_avoid_boundary} with \Cref{lem:crossing_distances_of_small_annuli_are_uniformly_small},  we obtain that the following conditions hold.

\begin{corollary}\label{cor:conditions_a_b_c}
	Suppose that we have the same setup as in \Cref{lem:crossing_distances_of_small_annuli_are_uniformly_small} and \Cref{prop:good_annulus_arbitrary_scale}. Then, it is a.s.\ the case that there exist $R>1,  K \in \BN$ and a compact set $F \subseteq \BD$ such that 
\begin{enumerate}[label=(\alph*), ref=\alph*]
	\item \label{it:existence_of_good_annuli}
	$P \subseteq F$ and for all $z \in \left(\frac{1}{100}2^{-(1+\nu) RK}\BZ\right)^2 \cap F$,  there exists $\nu R K \leq k \leq (1+\nu) R K -R$ such that $E_{z,k}(\Psi)$ occurs,  where $\nu$ is as in \Cref{prop:good_annulus_arbitrary_scale}.
	\item \label{it:distances_across_annuli_are_small}
	Let $\widetilde{P}$ denote the time-reversal of $P$.  Then,  for all $z \in F$,  if both $P$ and $\widetilde{P}$ intersect $\overline{B_{2^{-\nu K R}}(z)}$ before they intersect each other,  we have that the following holds.  Let $t$ (resp.\ $\widetilde{t}$) denote the first time that $P$ (resp.\ $\widetilde{P}$) intersects $\overline{B_{2^{-\nu K R}}(z)}$.  Then,  $\len(P|_{[t , 1-\widetilde{t}]} ; D_\Gamma)$ is strictly less than $c_R$,  where recall that $c_R$ is the constant in condition~\eqref{it:distances_across} in \Cref{prop:good_annulus_macroscopic_scale}.
	\item \label{it:loops_not_too_small}
	We have that $\diam(\SCL) \geq 100 \cdot 2^{-R}$,  $\diam(\SCL') \geq 100 \cdot 2^{-R}, \dist(\SCL,\SCL') \geq 100 \cdot 2^{-R}$,   $\dist(P(0) ,  P \setminus \SCL) \geq 100 \cdot 2^{-R}$ and $\dist(P(1) ,  P \setminus \SCL^{\prime}) \geq 100 \cdot 2^{-R}$.

\end{enumerate}
\end{corollary}
\begin{proof}
	
	We start by modifying $P$ near its starting and ending points.  Let $t_1$ (resp.\@ $s_1$) denote the last (resp.\@ first) time that $P$ intersects $\SCL$ (resp.\@ $\SCL^{\prime}$),  fix a point $w \in \SCL \setminus \{P(t_1)\}$ (resp.\@ $w' \in \SCL^{\prime} \setminus \{P(s_1)\}$),  and prepend to $P$ an arc of $\SCL$ (resp.\@ $\SCL^{\prime}$) from $w$ (resp.\@ $w'$) to $P(0)$ (resp.\@ $P(1)$).  Since any path contained in a single loop has zero $D_{\Gamma}$-length,  the resulting path is still a $D_{\Gamma}$-geodesic from $\SCL$ to $\SCL^{\prime}$,  and we denote it again by $P$.  As $P \setminus \SCL \subseteq P([t_1,1])$ (resp.\@ $P \setminus \SCL^{\prime} \subseteq P([0,s_1])$),  as $P([t_1,1])$ (resp.\@ $P([0,s_1])$) is compact and $w \notin P([t_1,1])$ (resp.\@ $w' \notin P([0,s_1])$),  we obtain that $\mathrm{dist}(P(0) ,  P \setminus \SCL) > 0$ (resp.\@ $\mathrm{dist}(P(1) ,  P \setminus \SCL^{\prime}) > 0$) a.s.
	
	Let $(F_n)_{n\ge 1}$ be a deterministic sequence of compact subsets of $\BD$ such that for all $n\ge 1$, we have $\BP[P \subset F_n \mid \SCL \neq \SCL'] \ge 1-1/n^2$. Such a sequence does exist thanks to \Cref{cor:geodesic_avoid_boundary} together with the fact that $\SCL$ is a.s.\ a compact subset of $\BD$. Moreover,  we may and do take each $F_n$ to be a closed Euclidean disk centered at the origin and contained in $\BD$,  and we fix deterministic constants $\epsilon_n>0$ such that the event $\{P \subseteq F_n\}$ above can be replaced by the event that $P$ is contained in the set of points of $F_n$ lying at distance at least $\epsilon_n$ from $\partial F_n$.  On the latter event,  $F_n$ contains the Euclidean $\epsilon_n$-neighborhood of $P$,  which will be used in the construction of the chains of loops below. Let also $(R_{0,n})_{n\ge 1}$ be a deterministic sequence of positive numbers such that
\begin{align*}
	&\BP\left[ \min \left( \diam(\SCL),  \diam(\SCL') , \dist(\SCL,\SCL'),  \dist(P(0) ,  P \setminus \SCL),  \dist(P(1),P \setminus \SCL^{\prime}) \right)\right.\\
	&\qquad\qquad\left.\geq 100 \cdot 2^{-R_{0,n}} \mid \SCL \neq \SCL' \right] \\
	&\ge 1-\frac{1}{n^2}.
\end{align*}
	Such a sequence exists since the four quantities appearing above are a.s.\ positive on the event that $\SCL \neq \SCL'$.  Note that,  since we shall have $R>R_{0,n}$ below,  the above event implies condition~\eqref{it:loops_not_too_small}.
	Next, by \Cref{prop:good_annulus_arbitrary_scale} applied with $R_1 = R_{0,n}$ and with $F_n$ in place of $F$, there exist a deterministic $R(n)>R_{0,n}$ and an a.s.\ finite random variable $K_0(n) \in \BN$ such that for all $K \ge K_0(n)$, for all $z \in \left(\frac{1}{100}2^{-(1+\nu) R(n)K}\BZ\right)^2 \cap F_n$,  there exists $\nu R(n) K \leq k \leq (1+\nu) R(n) K -R(n)$ such that $E_{z,k}(\Psi)$ occurs.  Enlarging $K_0(n)$ if necessary,  we may and do also assume that $\nu K_0(n) \geq 1$ and that $2^{-\nu K_0(n) R(n)} \leq \epsilon_n$.  (That $R(n)$ can be taken to be deterministic is what the proof of \Cref{prop:good_annulus_arbitrary_scale} gives,  and this is what allows us to apply the Borel--Cantelli lemma below for a fixed value of $R$.)
	But, using \Cref{lem:crossing_distances_of_small_annuli_are_uniformly_small}, the Borel--Cantelli lemma and using that for all deterministic $R>R_{0,n}$, we have
	\[\sum_{K\ge 1} c_1 e^{-c_2 c_{R} \nu K R }< \infty,\]
	we obtain in particular that there exists a random $K(n) \ge K_0(n)$ such that~\eqref{it:distances_across_annuli_are_small} holds with $R,K$ replaced by $R(n),K(n)$.
	Finally, since $\sum_{n \geq 1} 1/n^2 < \infty$,  applying the Borel--Cantelli lemma once more, we obtain the desired result.\end{proof}

Conditions~\eqref{it:existence_of_good_annuli}--\eqref{it:loops_not_too_small} will play a crucial role in constructing chains of loops in $\Gamma$ and \emph{good} annuli which cover the geodesic $P$.  Before proceeding to the explicit construction of the aforementioned chains of loops and annuli,   we state and prove two useful lemmas.  The first lemma (\Cref{lem:intersection_is_a_segment}) states that the intersection of $P$ with any loop in $\Gamma$ is a segment of $P$.  The second lemma (\Cref{lem:geodesic_intersects_loops}) states that whenever the geodesic $P$ crosses an annulus for which the event $E_{z,k}(\Psi)$ occurs and there are two loops in $\Gamma$ crossing the annulus,  we have that $P$ intersects at least one of the loops in $\Gamma$ crossing the annulus.

\begin{lemma}\label{lem:intersection_is_a_segment}
	Let $\ell$ be a loop in $\Gamma$ such that $P \cap \ell \neq \emptyset$.  Then,  there exist times $0 \leq s \le t \leq 1$ such that $P \cap \ell  = P([s,t])$.
\end{lemma}

\begin{proof}
	Let $\ell \in \Gamma$ be such that $P \cap \ell \neq \emptyset$ and let $s$ (resp.\ $t$) be the first (resp.\ last) time that $P$ intersects $\ell$.  By the definitions of $s$ and $t$,  we have that $P \cap \ell \subseteq P([s,t])$,  and so it suffices to show that $P([s,t]) \subseteq \ell$.

	We first claim that $\len(P|_{[s,t]} ; D_\Gamma) = 0$.  Indeed,  since $P(s),P(t) \in \ell$,  we can replace $P|_{[s,t]}$ by an arc of $\ell$ from $P(s)$ to $P(t)$.  The resulting path is admissible,  it connects $\SCL$ with $\SCL'$,  and its $D_\Gamma$-length is equal to $\len(P ; D_\Gamma) - \len(P|_{[s,t]} ; D_\Gamma)$ since any path which is contained in a single loop has zero $D_\Gamma$-length.  As $P$ is a $D_\Gamma$-geodesic,  the claim follows.

	Suppose now that $\ell' \in \Gamma$ is a loop which $P|_{[s,t]}$ intersects.  Then the definition of $\len(\cdot ; D_\Gamma)$ implies that
	\begin{align*}
		D_\Gamma(\ell,\ell') \leq \len(P|_{[s,t]} ; D_\Gamma) = 0,
	\end{align*}
	and hence $\ell' = \ell$ since $D_\Gamma$ is a metric on $\Gamma$.  Therefore every loop in $\Gamma$ which $P([s,t])$ intersects is equal to $\ell$. Now, assume by contradiction that there is $u\in (s, t)$ such that $P(u)\not\in \ell$. By admissibility of $P$, we know that $P\cap \bigcup_{L \in \Gamma} L$ is dense in $P$. Since $\ell$ is closed and $P(u) \not\in \ell$, we deduce that there exists a sequence $t_n \in [0, 1]$ and loops $\ell_n \in \Gamma\setminus \{\ell\}$ such that $P(t_n) \in \ell_n$ and $P(t_n) \to P(u)$. By taking a subsequence, we may assume that $t_n \to u'\in [0, 1]$. In particular, $P(u')= P(u)$. One can see that $s<u'<t$ because otherwise, we can remove from $P$ the piece of path $P([u', u])$ or $P([u, u'])$ and obtain an admissible path with shorter $D_\Gamma$-length. But then, for $n$ large enough, we have $t_n \in (s,t)$, which is absurd since we have shown that every loop in $\Gamma$ which intersects $P([s,t])$ is equal to $\ell$. Therefore, $P([s,t])\subset \ell$. 
	This completes the proof of the lemma.
\end{proof}

\begin{lemma}\label{lem:geodesic_intersects_loops}
	Suppose that we have the same setup as in \Cref{cor:conditions_a_b_c},  so that in particular $R$ and $P$ are as fixed above.  Let $z \in \BD ,  k \in \BN$ be such that $\overline{B_{2^{-k}}(z)} \subseteq \BD$ and $E_{z,k}(\Psi)$ occurs.  Suppose also that the collection of loops in $\overline{\Gamma}$ with height in $\{-\lambda,\lambda\}$ which cross $B_{2^{-k-1}}(z) \setminus \overline{B_{2^{-k-R}}(z)}$

is non-empty and that there are at least two loops in $\overline{\Gamma}$ with height in $\{-\lambda,\lambda\}$ which cross $B_{2^{-k-1}}(z) \setminus \overline{B_{2^{-k-R}}(z)}$.  Let $\ell_1,\ell_2$ denote the elements of the above collection (recall that,  since $E_{z,k}(\Psi)$ occurs,  it consists of at most two loops; we use the letters $\ell_1,\ell_2$ rather than $\SCL_1,\SCL_2$ in order to avoid a clash with the chain of loops constructed in \Cref{subsec:constructing_chains}).  Then,  wherever the geodesic $P$ makes a crossing of $B_{2^{-k-1}}(z) \setminus \overline{B_{2^{-k-R}}(z)}$ (in particular,  whenever it makes a crossing of the larger annulus $B_{2^{-k}}(z) \setminus \overline{B_{2^{-k-R}}(z)}$),  the portion of $P$ which makes that crossing must intersect either $\ell_1$ or $\ell_2$.  More generally,  the same conclusion holds whenever $P$ crosses one of the connected components $U_j$ of the complement in $B_{2^{-k-1}}(z) \setminus \overline{B_{2^{-k-R}}(z)}$ of the crossings of that annulus made by loops in $\overline{\Gamma}$ with height in $\{-\lambda,\lambda\}$ together with the interiors of the excursions into it made by such loops (we use the notation of \Cref{prop:good_annulus_macroscopic_scale}),  between the two marked boundary arcs $I_j^3$ and $I_j^4$ of $\partial U_j$.
\end{lemma}

\begin{proof}
	\stepn{step:geo-loops-setup}{Setup} Suppose that $P$ makes a crossing of  $B_{2^{-k}}(z) \setminus \overline{B_{2^{-k-R}}(z)}$.  We will assume that it crosses the annulus from its outer to its inner boundary since the exact same argument works whenever it makes a crossing from its inner to its outer boundary.

	\stepn{step:geo-loops-crossing-loops}{The crossing loops $\ell_1$ and $\ell_2$} Let $0 \leq s < t \leq 1$ be such that $P(s) \in \partial B_{2^{-k}}(z),  P(t) \in \partial B_{2^{-k-R}}(z),  P((s,t)) \subseteq  B_{2^{-k}}(z) \setminus \overline{B_{2^{-k-R}}(z)}$,  and $P([s,t])$ is part of the crossing that $P$ makes.  Note that $P|_{[s,t]}$ contains in particular a crossing of the smaller annulus $B_{2^{-k-1}}(z) \setminus \overline{B_{2^{-k-R}}(z)}$ from its outer to its inner boundary.  Note also that the definition of $E_{z,k}(\Psi)$ implies that there are at most four crossings of $B_{2^{-k-1}}(z) \setminus \overline{B_{2^{-k-R}}(z)}$ made by loops in $\overline{\Gamma}$ with height in $\{-\lambda,\lambda\}$,  and hence at most two such loops,  since a loop which crosses an annulus has to cross it at least twice.  Consequently,  $\ell_1$ and $\ell_2$ are the only loops in $\overline{\Gamma}$ with height in $\{-\lambda,\lambda\}$ which cross $B_{2^{-k-1}}(z) \setminus \overline{B_{2^{-k-R}}(z)}$.

	\stepn{step:geo-loops-avoiding-case}{The case where $P|_{[s,t]}$ avoids $\ell_1$ and $\ell_2$} Suppose that $P|_{[s,t]}$ does not intersect either $\ell_1$ or $\ell_2$.  Then,  using the notation of \Cref{prop:good_annulus_macroscopic_scale},  we obtain that there exists $j \in \{1,\ldots,N\}$ such that the field $\Psi|_{U_j}$ has zero-boundary conditions and $P|_{[s,t]}$ makes a crossing of $U_j$ from $I_j^4$ to $I_j^3$.  Indeed,  $P|_{[s,t]}$ has to cross some component $U_j$ from $I^4_j$ to $I^3_j$ since it does not intersect $\ell_1$ or $\ell_2$;  moreover,  if the boundary conditions of $\Psi|_{U_j}$ were not equal to zero,  then \Cref{lem:zero_boundary_conditions} would imply that $U_j \subseteq \mathop{\mathrm{int}}(\ell_0)$ for some outermost loop $\ell_0$ in $\overline{\Gamma}$,  i.e.\ for some $\ell_0 \in \Gamma$;  this is impossible,  since $\mathop{\mathrm{int}}(\ell_0)$ is an open set which does not intersect any loop in $\Gamma$ (the loops in $\Gamma$ being pairwise non-nested),  so that the admissibility of $P$ forces $P \cap \mathop{\mathrm{int}}(\ell_0) = \emptyset$.  Then,  condition~\eqref{it:disconnecting_conformal_rectangles} in \Cref{prop:conformal_rectangles} (and,  more precisely,  the fact that the left and right boundaries of $R_{j,m}$ are contained in $I_j^1$ and $I_j^2$,  so that $R_{j,m}$ separates $I_j^3$ from $I_j^4$ in $U_j$) implies that $P|_{[s,t]}$ has to cross all of the topological rectangles $R_{j,m}$ for $1 \leq m \leq 4$ from the bottom to their top boundaries.  In particular,  the path $P|_{[s,t]}$ has to intersect all of the paths $P_{j,m}$ for $1 \leq m \leq 4$,  where $P_{j,m}$ denotes an admissible path in $U_{j,m}$ connecting the left and right marked boundary arcs of $R_{j,m}$ with $\len(P_{j,m} ; D_\Gamma)<B$,  the existence of which is guaranteed by condition~\eqref{it:left_to_right_crossings_shorter_than_top_to_bottom} in \Cref{prop:conformal_rectangles} and by the definition of the internal metric $D^{U_{j,m}}_{\Gamma|_{U_{j,m}}}$.  Let $\tau_{j,m}$ denote the first time after $s$ that $P$ intersects $P_{j,m}$.  Then,  the path $P|_{[\tau_{j,1},\tau_{j,4}]}$ must cross from bottom to top both $R_{j,2}$ and $R_{j,3}$,  and hence condition~\eqref{it:left_to_right_crossings_shorter_than_top_to_bottom} in \Cref{prop:conformal_rectangles} implies that the $D_\Gamma$-length of $P|_{[\tau_{j,1},\tau_{j,4}]}$ is at least $2 B$ (here we use that $R_{j,2}$ and $R_{j,3}$ are disjoint and that the $D_\Gamma$-length of a path contained in $U_{j,m}$ is at least the $D^{U_{j,m}}_{\Gamma|_{U_{j,m}}}$-distance between its endpoints).  However,  concatenating the paths $P_{j,1}$ and $P_{j,4}$ with the loops $\ell_1,\ell_2$ (more precisely,  we follow $P_{j,1}$ from $P(\tau_{j,1})$ until we reach $I_j^1$,  we then travel along the loop containing $I_j^1$ until we reach the endpoint of $P_{j,4}$ on $I_j^1$,  which costs no $D_\Gamma$-length,  and we finally follow $P_{j,4}$ until $P(\tau_{j,4})$) and noting that the $D_\Gamma$-lengths of $P_{j,1}$ and $P_{j,4}$ are strictly less than $B$ (see condition~\eqref{it:left_to_right_crossings_shorter_than_top_to_bottom} in \Cref{prop:conformal_rectangles}),  we obtain an admissible path from $P(\tau_{j,1})$ to $P(\tau_{j,4})$ whose $D_\Gamma$-length is strictly less than $2B$.  This implies that we can construct an admissible path connecting $\SCL$ to $\SCL'$ whose $D_\Gamma$-length is strictly less than that of $P$.  But this contradicts the fact that $P$ is a $D_\Gamma$-geodesic.  It follows that $P|_{[s,t]}$ must intersect either $\ell_1$ or $\ell_2$.  Finally,  we note that the argument of the previous paragraph used nothing about $P|_{[s,t]}$ beyond the fact that it crosses the component $U_j$ between $I_j^3$ and $I_j^4$,  and that the component crossed by $P$ automatically has zero boundary conditions,  by the argument given there.  This proves the second assertion of the lemma as well.  This completes the proof of the lemma.
\end{proof}

\subsection{Constructing chains of loops}\label{subsec:constructing_chains}
Now we are going to use conditions~\eqref{it:existence_of_good_annuli}--\eqref{it:loops_not_too_small} from \Cref{cor:conditions_a_b_c} and \Cref{lem:intersection_is_a_segment,lem:geodesic_intersects_loops} in order to construct sequences of loops $\SCL_1,\ldots,\SCL_N$ and Euclidean balls 
\begin{align*}
	B_{2^{-k_1}}(z_1),\ldots, {B_{2^{-k_{N-1}}}(z_{N-1})}
\end{align*}
satisfying the following properties.  First,  we have that $\SCL = \SCL_1,\SCL' = \SCL_N$,  and the geodesic $P$ intersects $\SCL_j$ before intersecting $\SCL_{j+1}$ for all $1 \leq j \leq N-1$.  Secondly,  we have that $\SCL_{j+1}$ makes a crossing of $B_{2^{-k_j-1}}(z_j) \setminus \overline{B_{2^{-k_j - R}}(z_j)}$ and the event $E_{z_j,k_j}(\Psi)$ occurs for all $1 \leq j \leq N-1$.  Thirdly,  we have that $P \cap B_{2^{-k_j-R}}(z_j) \neq \emptyset$ for all $1 \leq j \leq N-1$,  and
\begin{align*}
	P \subseteq \left(\bigcup_{j=1}^{N-1} B_{2^{-k_j}}(z_j)\right) \cup \left(\bigcup_{j=1}^N \SCL_j\right).
\end{align*}
Throughout the construction,  we work on the almost sure event given by \Cref{cor:conditions_a_b_c},  and we assume without loss of generality that $\nu K \geq 1$ (so that $2^{-\nu K R} \leq 2^{-R}$) and that $F$ contains the Euclidean $2^{-\nu K R}$-neighborhood of $P$.  Both of these were arranged in the proof of \Cref{cor:conditions_a_b_c},  and the latter guarantees that the grid points used below can indeed be taken in $F$.

\setcounter{jmstep}{0}\stepn{step:chains-init}{Initialization of the construction} 

\substepn{step:ci-first-annulus}{The first loop and the annulus $A_1$} First,  we set $\SCL_1 = \SCL$.  Note that \Cref{lem:intersection_is_a_segment} implies that there exists $t_1 \in [0,1]$ such that $P \cap \SCL_1 = P([0,t_1])$.  Let $z_1 \in \left(\frac{1}{100} 2^{-(1+\nu) R K }\BZ\right)^2 \cap F$ be such that $|z_1 - P(t_1)| \leq 2^{-(1+\nu) R K } / 50$,  where we break ties in some arbitrary but fixed measurable way.  Let also $k_1$ denote the smallest integer $\nu K R \leq k \leq (1+\nu) R K - R$ such that the event $E_{z_1,k}(\Psi)$ occurs;  such an integer exists by condition~\eqref{it:existence_of_good_annuli}.  Note that $2^{-k_1-R} \geq 2^{-(1+\nu) R K}$ and hence $P(t_1) \in B_{2^{-k_1-R}}(z_1)$,  and that $2^{-k_1} \leq 2^{-\nu K R} \leq 2^{-R}$.  Note that condition~\eqref{it:loops_not_too_small} implies that $\SCL$ has to make a crossing of the annulus $B_{2^{-R}}(z_1) \setminus \overline{B_{2^{-(1+\nu) R K}}(z_1)}$ since $P(t_1) \in B_{2^{-(1+\nu) R K}}(z_1)$ and since $\diam(\SCL) \geq 100 \cdot 2^{-R}$,  so that $\SCL$ cannot be contained in $B_{2^{-R}}(z_1)$.  In particular,  it has to make a crossing of $A_1{\defeq}{B_{2^{-k_1-1}}(z_1)} \setminus \overline{B_{2^{-k_1 -R}}(z_1)}$,  since the latter annulus is contained in the former one.  Let us emphasize that we take the outer radius of $A_1$ to be $2^{-k_1-1}$ and not $2^{-k_1}$,  so that $A_1$ is exactly the annulus appearing in the definition of the event $E_{z_1,k_1}(\Psi)$ (see \Cref{rem:good_event_arbitrary_scale});  in this way the connected components and the conformal rectangles provided by that event are the ones associated with $A_1$ itself,  and the loops crossing $A_1$ are the ones which the event controls.  Moreover,  $P$ has to exit $B_{2^{-R}}(z_1)$ after time $t_1$,  since $P(1) \in \SCL'$ while $\dist(\SCL,\SCL') \geq 100 \cdot 2^{-R}$ by condition~\eqref{it:loops_not_too_small} and $P(t_1) \in \SCL \cap B_{2^{-(1+\nu)RK}}(z_1)$.

\substepn{step:ci-second-loop}{Exactly one further loop crosses $A_1$} We now show that there exists exactly one more loop in $\Gamma$ which crosses $A_1$.  First,  we note that the definition of $E_{z_1,k_1}(\Psi)$ implies that there are at most four crossings of $B_{2^{-k_1-1}}(z_1) \setminus \overline{B_{2^{-k_1-R}}(z_1)}$ made by loops in $\Gamma$,  i.e.\ of $A_1$,  and hence that at most two loops in $\Gamma$ cross $A_1$,  since a loop which crosses $A_1$ does so at least twice.  Since $\SCL$ crosses $A_1$,  it therefore suffices to rule out the case that $\SCL$ is the unique loop in $\Gamma$ crossing $A_1$.  So suppose that this is the case.  Then,  there exists a connected component $U_1$ of the set
\begin{align*}
	A_1 \setminus \overline{\bigcup_{\ell \in \Gamma_{A_1}^{\mathrm{out}} \cup \Gamma_{A_1}^{\mathrm{in}}} \mathop{\mathrm{int}}(\ell) \cup \bigcup_{\ell \in \Gamma_{A_1}} \ell }
\end{align*}
such that the field $\Psi|_{U_1}$ has zero-boundary conditions and $P|_{[t_1,1]}$ makes a crossing of $U_1$ from its top to its bottom boundary (that the boundary conditions can be taken to be equal to zero follows as in the proof of \Cref{lem:geodesic_intersects_loops},  since an admissible path cannot enter the interior of a loop in $\Gamma$).  Then,  the definition of $E_{z_1,k_1}(\Psi)$ (see condition~\eqref{it:distances_across} in \Cref{prop:good_annulus_macroscopic_scale}) implies that
\begin{align*}
	D^{V_1}_{\phi_1(\Gamma|_{U_1})}((0,L_1) \times \{0\} ,  (0,L_1) \times \{1\}) \geq c_R,
\end{align*}
where $\phi_1 \colon U_1 \to V_1 = (0,L_1) \times (0,1)$ is the conformal map from \Cref{subsec:conformal_rectangles}.
In particular,  we have that the $D_\Gamma$-length of the part of $P$ between the time $t_1$ and the first time after $t_1$ that $P$ intersects $\partial B_{2^{-k_1-1}}(z_1)$ is at least $c_R$ (here we use the conformal invariance and the locality of $D_\Gamma$).  Let now $t$ (resp.\ $\widetilde{t}$) denote the first time that $P$ (resp.\ $\widetilde{P}$) intersects $\overline{B_{2^{-\nu K R}}(z_1)}$.  Since $P(t_1) \in B_{2^{-k_1-R}}(z_1) \subseteq B_{2^{-\nu K R}}(z_1)$,  we have that $t \leq t_1 \leq 1-\widetilde{t}$,  and since $B_{2^{-k_1}}(z_1) \subseteq B_{2^{-\nu K R}}(z_1)$, the part of $P$ between the time $t_1$ and the first time after $t_1$ that $P$ intersects $\partial B_{2^{-k_1-1}}(z_1)$ is contained in $P|_{[t , 1-\widetilde{t}]}$.  Hence $\len(P|_{[t , 1-\widetilde{t}]} ; D_\Gamma) \geq c_R$,  which contradicts condition~\eqref{it:distances_across_annuli_are_small} from \Cref{cor:conditions_a_b_c}.  Hence, there exists another loop $\SCL_2$ in $\Gamma$,  distinct from $\SCL$,  which crosses $A_1$ and,  by the above, there are no other loops in $\Gamma$ crossing $A_1$.  Since each of $\SCL$ and $\SCL_2$ crosses $A_1$ at least twice while there are at most four crossings in total,  each of them crosses $A_1$ exactly twice and the total number of crossings of $A_1$ is equal to four.

\substepn{step:ci-geodesic-hits}{The geodesic intersects $\SCL_2$ after time $t_1$} Note that $P|_{(t_1,1]}$ does make a crossing of $A_1$,  since $P(t_1) \in B_{2^{-k_1-R}}(z_1)$ while $P$ exits $B_{2^{-R}}(z_1) \supseteq B_{2^{-k_1-1}}(z_1)$ after time $t_1$.  \Cref{lem:geodesic_intersects_loops} implies that the path $P|_{(t_1,1]}$ must intersect either $\SCL_1$ or $\SCL_2$ while making its first crossing of $A_1$ after time $t_1$.  But since $P((t_1,1]) \cap \SCL_1 = \emptyset$ by the definition of $t_1$,  we obtain that $P$ has to intersect $\SCL_2$ while making that crossing,  and in particular after time $t_1$.

\stepn{step:chains-induction}{The inductive step of the construction}

\substepn{step:cind-hypothesis}{The inductive hypothesis and the annulus $A_m$} Next,  we assume that we have defined loops $\SCL_1,\ldots,\SCL_m$ in $\Gamma$ (where $m \geq 2$) and Euclidean balls $B_{2^{-k_1}}(z_1),\ldots,B_{2^{-k_{m-1}}}(z_{m-1})$ satisfying the following properties.  We have that $P$ intersects $\SCL_j$ before intersecting $\SCL_{j+1}$ for all $1 \leq j \leq m-1$.  Also,  the loop $\SCL_{j+1}$ makes a crossing of $A_j{\defeq}{B_{2^{-k_j-1}}(z_j)} \setminus \overline{B_{2^{-k_j-R}}(z_j)}$ and the event $E_{z_j,k_j}(\Psi)$ occurs for all $1 \leq j \leq m-1$.  Moreover,  we have that $P \cap B_{2^{-k_j - R}}(z_j) \neq \emptyset$ for all $1 \leq j \leq m-1$.  Note that \Cref{lem:intersection_is_a_segment} implies that the loops $\SCL_1,\ldots,\SCL_m$ are pairwise distinct and that $t_1 < t_2 < \cdots < t_m$,  where $t_j$ denotes the last time that $P$ intersects $\SCL_j$;  indeed,  $P$ intersects each of these loops along a single segment and it intersects them in this order,  and these segments are disjoint since the loops in $\Gamma$ are.

If $\SCL_m = \SCL'$,  then we stop the procedure.  Suppose that $\SCL_m \neq \SCL'$.  Let $t_m$ be the last time that $P$ intersects $\SCL_m$ and let $z_m \in \left(\frac{1}{100}2^{-(1+\nu) R K}\BZ\right)^2 \cap F$ be such that $|z_m - P(t_m)| \leq 2^{-(1+\nu) R  K } /  50$,  where we break ties in some arbitrary but fixed measurable way.  We let $k_m$ denote the smallest integer $\nu R K \leq k \leq (1+\nu) R K - R$ for which $E_{z_m,k}(\Psi)$ occurs,  so that,  as in \Cref{step:chains-init},  $P(t_m) \in B_{2^{-k_m-R}}(z_m)$ and $2^{-k_m} \leq 2^{-R}$,  and set
\begin{align*}
	A_m \defeq B_{2^{-k_m-1}}(z_m) \setminus \overline{B_{2^{-k_m - R}}(z_m)}.
\end{align*}

\substepn{step:cind-at-least-one}{At least one loop crosses $A_m$} As in \Cref{step:chains-init},  we will show that there are exactly two loops in $\Gamma$ which cross $A_m$ and that the total number of crossings is equal to four.  Indeed,  first we note that the definition of $E_{z_m,k_m}(\Psi)$ implies that there are at most two loops in $\Gamma$ crossing $A_m$ (exactly as in \Cref{step:chains-init}).  Suppose that there are no crossings of $A_m$ made by loops in $\Gamma$. Let $\SP_{A_m}$ be a countable family of polygonal paths defined as in \Cref{rem:no_crossings} with $A_m$ instead of $A$. Then, by \Cref{rem:no_crossings}, for all $\gamma \in \SP_{A_m}$, the set
\begin{align*}
	U_{\gamma} \defeq A_m \setminus \overline{\bigcup_{\ell \in \widehat{\Gamma}_{\gamma}} \mathop{\mathrm{int}}(\ell)}
\end{align*}
is simply connected,  the field $\Psi|_{U_{\gamma}}$ has zero-boundary conditions,  and the definition of $E_{z_m,k_m}(\Psi)$ implies that
\begin{align*}
	D^{V_{\gamma}}_{\phi_{\gamma}(\Gamma|_{U_{\gamma}})}((0,L_{\gamma}) \times \{0\} ,  (0,L_{\gamma}) \times \{1\}) \geq c_R,
\end{align*}
where recall that $\widehat{\Gamma}_{\gamma}$ denotes the collection of loops in $\overline{\Gamma}$ with height in $\{-\lambda,\lambda\}$ which intersect $\partial A_m \cup {\gamma}$
.  Note that $t_m>t_1$ and hence $P(t_m) \notin \SCL$.  Hence, since $P(t_m) \in B_{2^{-(1+\nu) R K}}(z_m)$,  we obtain by combining with condition~\eqref{it:loops_not_too_small} that $P(0) \notin \overline{B_{2^{-k_m-1}}(z_m)}$ (indeed,  $P(t_m) \in P \setminus \SCL$,  so that $|P(0)-P(t_m)| \geq 100 \cdot 2^{-R}$,  while $|P(t_m)-z_m| \leq 2^{-(1+\nu)RK}/50$ and $2^{-k_m-1} \leq 2^{-R}$).   It follows that $P|_{[0,t_m]}$ makes a crossing of $A_m$ from its outer to its inner boundary. 
Moreover,  we have that $P(1) \notin B_{2^{-\nu R K}}(z_m)$.  Indeed,  $P(t_m) \in P \setminus \SCL^{\prime}$ since $P(t_m) \in \SCL_m$ and we have assumed that $\SCL_m \neq \SCL^{\prime}$,  and so condition~\eqref{it:loops_not_too_small} in Corollary~\ref{cor:conditions_a_b_c} implies that $|P(1) - P(t_m)| \ge 100 \cdot 2^{-R}$,  while $|P(t_m)-z_m| \leq 2^{-(1+\nu) RK} / 50$ and $2^{-\nu RK} \leq 2^{-R}$.  Let now $t$ (resp.\@ $\widetilde{t}$) denote the first time that $P$ (resp.\@ $\widetilde{P}$) intersects $\overline{B_{2^{-\nu RK}}(z_m)}$. Exactly as in the first step, the two crossings of $A_m$ that $P|_{[0,t_m]}$ makes are contained in $P|_{[t,1-\widetilde{t}]}$. Moreover, since the path $P|_{[0,t_m]}$ crosses $A_m$ twice, there exists a polygonal path $\gamma \in \SP_{A_m}$ such that the first crossing is a crossing of $U_\gamma$, whence $\mathrm{len}(P|_{[t,1-\widetilde{t}]} ; D_{\Gamma}) \ge c_R$,  which contradicts condition~\eqref{it:distances_across_annuli_are_small}.  Therefore,  there exists at least one loop $\widetilde{\SCL}$ in $\Gamma$ which crosses $A_m$.

\substepn{step:cind-tilde-loop}{Properties of the loop $\widetilde{\SCL}$} Suppose that $\widetilde{\SCL}$ is the unique loop in $\Gamma$ that crosses $A_m$.  Then,  $P|_{[0,t_m]}$ must intersect $\widetilde{\SCL}$ while making a crossing of $A_m$.  Indeed, if it was not the case, as explained in the proof of Lemma~\ref{lem:geodesic_intersects_loops} (using the notation of \Cref{prop:good_annulus_macroscopic_scale}),  there exists $j \in \{1,\cdots,N\}$ such that the field $\Psi|_{U_j}$ has zero boundary conditions and $P|_{[0,t_m]}$ makes a crossing of $U_j$ from $I_j^4$ to $I_j^3$,  where we also have that $I_j^1 \cup I_j^2 \subseteq \widetilde{\SCL}$.  Thus,  by condition~\eqref{it:distances_across} in Lemma~\ref{prop:good_annulus_macroscopic_scale},  we obtain that the $D_{\Gamma}$-length of the above crossing is at least $c_R$,  and so arguing as in the previous paragraph gives that $\mathrm{len}(P|_{[t,1-\widetilde{t}]} ; D_{\Gamma}) \ge c_R$,  which contradicts condition~\eqref{it:distances_across_annuli_are_small}.  It follows that $P|_{[0,t_m]}$ intersects $\widetilde{\SCL}$ while making a crossing of $A_m$.

Since $t_m$ is the last time that $P$ intersects $\SCL_m$,  we obtain that either $\widetilde{\SCL} = \SCL_m$ or $P$ intersects $\widetilde{\SCL}$ before intersecting $\SCL_m$.  In both cases,  Lemma~\ref{lem:intersection_is_a_segment} implies that $P$ does not intersect $\widetilde{\SCL}$ after time $t_m$,  since $P \cap \widetilde{\SCL}$ is a single segment of $P$.  Let us emphasize that this fact does not depend on the number of loops crossing $A_m$; it will be used again below,  after we have ruled out the case that $\widetilde{\SCL}$ is the only such loop.  Let us also note,  again independently of that number,  that $\widetilde{\SCL} \neq \SCL^{\prime}$: otherwise,  $P$ would intersect $\SCL^{\prime}$ at some time at most $t_m$,  and since $P \cap \SCL^{\prime}$ is a segment of $P$ which contains the time $1$,  we would get $P(t_m) \in \SCL^{\prime}$ and hence $\SCL_m = \SCL^{\prime}$,  which we have excluded.

\substepn{step:cind-exactly-two}{Exactly two loops cross $A_m$ and the loop $\SCL_{m+1}$} Recall that $P(1) \notin B_{2^{-\nu R K}}(z_m)$ and $P(t_m) \in B_{2^{-(1+\nu) R K}}(z_m)$.  Hence,  we have that $P|_{[t_m,1]}$ makes a crossing between the inner and outer boundaries of $A_m$ and we have that $t < 1-\widetilde{t}$,  where recall that $t$ (resp.\@ $\widetilde{t}$) is the first time that $P$ (resp.\@ $\widetilde{P}$) intersects $\partial B_{2^{-\nu R K}}(z_m)$.  In particular,  we have that $P|_{(t_m,1)}$ must cross one of the connected components of the set 
\begin{align*}
	A_m \setminus \overline{\bigcup_{\ell \in \Gamma_{A_m}^{\mathrm{out}} \cup \Gamma_{A_m}^{\mathrm{in}}} \mathop{\mathrm{int}}(\ell) \cup \bigcup_{\ell \in \Gamma_{A_m}} \ell },
\end{align*}
and so the $D_{\Gamma}$-length of the above crossing is at least $c_R$ by condition~\eqref{it:distances_across} in Lemma~\ref{prop:good_annulus_macroscopic_scale}.  Moreover,  the crossing of $A_m$ that $P|_{[t_m,1]}$ makes is contained in $P|_{[t,1-\widetilde{t}]}$,  and so $\mathrm{len}(P|_{[t,1-\widetilde{t}]} ; D_{\Gamma}) \ge c_R$,  which contradicts condition~\eqref{it:distances_across_annuli_are_small}.  It follows that there are exactly two loops in $\Gamma$ which cross $A_m$,  and hence exactly four crossings of $A_m$,  since each of these two loops crosses $A_m$ at least twice while there are at most four crossings.  We let $\SCL_{m+1}$ denote the one which is distinct from $\widetilde{\SCL}$,  and we claim that $P$ intersects $\SCL_{m+1}$ after time $t_m$ and before exiting $\overline{B_{2^{-k_m-1}}(z_m)}$.  Indeed,  since $P(1) \notin B_{2^{-\nu R K}}(z_m)$,  we have that $P|_{(t_m,1]}$ must exit $\overline{B_{2^{-k_m-1}}(z_m)}$ and so it makes a crossing of $A_m$ from its inner to its outer boundary,  and Lemma~\ref{lem:geodesic_intersects_loops} implies that it has to intersect $\widetilde{\SCL}$ or $\SCL_{m+1}$ while doing so; as the former is impossible,  the claim follows in this case.  Note also that Lemma~\ref{lem:intersection_is_a_segment} implies that $\SCL_{m+1} \notin \{\SCL_1,\cdots,\SCL_m\}$,  since $P$ intersects each $\SCL_j$ along a single segment which ends at time $t_j \leq t_m$.  If $\SCL_{m+1} = \SCL^{\prime}$,  then we stop the procedure; otherwise,  we continue inductively as described above.

\stepn{step:chains-conclusion}{Termination of the construction and conclusion} Note that all of the loops $\SCL_{m+1}$ constructed in \Cref{step:chains-init,step:chains-induction} cross the annuli $A_m$ and hence have Euclidean diameters at least $2^{-(1+\nu) R K}$ (a loop crossing $A_m$ has diameter at least $2^{-k_m-1}-2^{-k_m-R} \geq 2^{-k_m-2} \geq 2^{-(1+\nu)RK}$,  where we used that $R \geq 2$ and that $k_m \leq (1+\nu)RK-R$),  and the same is true of $\SCL_1 = \SCL$,  which crosses $A_1$.  Moreover,  the loops $\SCL_1,\SCL_2,\ldots$ are pairwise distinct and they all intersect $P \subseteq F$.  Therefore,  the a.s.\ local finiteness of $\Gamma$ implies that the procedure described in the previous paragraphs stops after finitely many steps $N$.  Moreover,  if $\SCL_N \neq \SCL'$,  then the inductive argument described above implies that we can find a loop $\SCL_{N+1} \in \Gamma \setminus \{\SCL_1,\ldots,\SCL_N\}$ intersected by $P$ such that $\SCL_{N+1}$ crosses the annulus $B_{2^{-k_N-1}}(z_N) \setminus \overline{B_{2^{-k_N - R}}(z_N)}$.  But this contradicts the choice of $N$ and so we must have that $\SCL_N = \SCL'$.

It remains to check the three properties announced at the beginning of the present subsection.  The first two are immediate from the construction,  and we have already noted that $P(t_j) \in B_{2^{-k_j-R}}(z_j)$ for all $1 \leq j \leq N-1$.  As for the covering property,  let $s_j$ (resp.\ $t_j$) denote the first (resp.\ last) time that $P$ intersects $\SCL_j$,  so that $P \cap \SCL_j = P([s_j,t_j])$ by \Cref{lem:intersection_is_a_segment},  and note that $s_1 = 0$,  that $t_N = 1$ and that $s_1 \leq t_1 \leq s_2 \leq \cdots \leq s_N \leq t_N$.  By construction,  $P$ intersects $\SCL_{j+1}$ after time $t_j$ and before exiting $\overline{B_{2^{-k_j-1}}(z_j)}$,  and hence $P([t_j,s_{j+1})) \subseteq B_{2^{-k_j-1}}(z_j) \subseteq B_{2^{-k_j}}(z_j)$ while $P(s_{j+1}) \in \SCL_{j+1}$.  Since
\begin{align*}
	P = \left(\bigcup_{j=1}^{N-1} P([t_j , s_{j+1}])\right) \cup \left(\bigcup_{j=1}^{N} P([s_j,t_j])\right),
\end{align*}
the desired inclusion follows.

\subsection{Going up and down in the geodesic tree} Now we are ready to prove \Cref{thm:main_result}.  

\begin{proof}[Proof of \Cref{thm:main_result}]
	\stepn{step:main-result-setup}{Setup and notation} Suppose that we have the same setup as in \Cref{subsec:constructing_chains} and let $\SCL_1,\ldots,\SCL_N$ be the loops and let $A_1,\ldots,A_{N{-1}}$ be the annuli constructed there.  For all $1 \leq i \leq N$,  we let $s_i$ (resp.\ $t_i$) be the first (resp.\ last) time that $P$ intersects $\SCL_i$.  Note that the construction of the loops $\SCL_1,\ldots,\SCL_N$ described above implies that $\SCL_i$ is one of the loops in $\Gamma$ crossing the annulus $A_{i-1}$ and that there is exactly one more loop $\widetilde{\SCL}_i$ in $\Gamma$ which crosses $A_{i-1}$ as well,  for all $2 \leq i \leq N$.  Moreover,  the path $P$ intersects $\widetilde{\SCL}_i$ before intersecting $\SCL_i$ and $|P(t_{i-1}) - z_{i-1}| \leq 2^{-(1+\nu) R K } / 50$.  
	
	\stepn{step:main-result-region}{The region $G$ between $\widetilde{\SCL}_i$ and $\SCL_i$} Fix $2 \leq i \leq N$.  Note that since the event $E_{z_{i-1},k_{i-1}}(\Psi)$ occurs,  we have that both $\widetilde{\SCL}_i$ and $\SCL_i$ cross $A_{i-1}$ exactly twice,  which implies that there are exactly two connected components $U_1,U_2$ of the set
	\begin{align*}
		A_{i-1} \setminus \overline{\bigcup_{\ell \in \Gamma_{A_{i-1}}^{\mathrm{out}} \cup \Gamma_{A_{i-1}}^{\mathrm{in}}} \mathop{\mathrm{int}}(\ell) \cup \bigcup_{\ell \in \Gamma_{A_{i-1}}} \ell}
	\end{align*}
	such that the field $\Psi|_{U_j}$ has zero-boundary conditions for $j=1,2$.  Furthermore,  \Cref{lem:geodesic_intersects_loops} implies that $P|_{(t_{i-1},1]}$ cannot cross $A_{i{-}1}$ without intersecting either $\widetilde{\SCL}_i$ or $\SCL_i$.  Thus,  combining with the definition of the topological rectangles $R_{j,m}$ for $j=1{,}2,  m = 1,2,\ldots,8$ obtained in condition~\eqref{it:disconnecting_conformal_rectangles} in \Cref{prop:conformal_rectangles},  we have that $P([t_{i-1},s_i])$ is contained in the bounded connected component $G$ of the set
	\begin{align*}
		\BC \setminus \left(\mathop{\mathrm{int}}(\widetilde{\SCL}_i) \cup \mathop{\mathrm{int}}(\SCL_i) \cup R_{1,5} \cup R_{2,5}\right).
	\end{align*}
	Let us emphasize that we use here the fifth of the eight rectangles in each of the two components,  and not the first one.  The reason is that,  by condition~\eqref{it:disconnecting_conformal_rectangles} in \Cref{prop:conformal_rectangles},  any path which starts in $G$,  does not intersect $\widetilde{\SCL}_i \cup \SCL_i$,  and leaves $A_{i-1}$ through its outer boundary,  has to cross each of the four rectangles $R_{j,1},\ldots,R_{j,4}$ for some $j \in \{1,2\}$ between their top and bottom boundaries;  this is exactly what the argument of \Cref{lem:geodesic_intersects_loops} requires,  and it will be used below.  Note also that $\widetilde{\SCL}_i$ and $\SCL_i$ together with $R_{1,5}$ and $R_{2,5}$ form a topological annulus,  so that the complement above has a unique bounded connected component.
	Indeed, $P(t_{i-1})$ is in $G$,  since $P(t_{i-1}) \in B_{2^{-k_{i-1}-R}}(z_{i-1})$ and since the part of $B_{2^{-k_{i-1}-R}}(z_{i-1})$ which is not surrounded by $\widetilde{\SCL}_i$ or $\SCL_i$ is contained in $G$. Moreover, $P$ intersects $\widetilde{\SCL}_i$ and then $\SCL_i$,  it cannot enter $\mathop{\mathrm{int}}(\widetilde{\SCL}_i)$ or $\mathop{\mathrm{int}}(\SCL_i)$ since it is admissible, and it cannot cross $R_{1,5}$ (or $R_{2, 5}$) between its top and bottom boundaries, otherwise one can build a shorter path going from $\widetilde{\SCL}_i$ to $\SCL_i$ (recall from condition~\eqref{it:left_to_right_crossings_shorter_than_top_to_bottom} in \Cref{prop:conformal_rectangles} that the $D_\Gamma$-distance between the left and right sides of that rectangle,  which lie on $\widetilde{\SCL}_i$ and on $\SCL_i$ respectively,  is strictly less than $B$,  while the distance between its top and bottom sides is at least $B$).
		
	\stepn{step:main-result-reduction}{Reduction to the case $\widetilde{\SCL}_i = \SCL_{i-1}$} {Let us first show that we can choose $\SCL_1, \ldots, \SCL_N$ and $A_1, \ldots, A_{N-1}$ so that $\widetilde{\SCL}_i = \SCL_{i-1}$. In the case $i=2$, using the statement~\eqref{it:loops_not_too_small} of \Cref{cor:conditions_a_b_c}, the loop $\SCL_1$ is too large to be contained in $G$ (indeed $G \subseteq B_{2^{-k_1-1}}(z_1)$,  whose diameter is at most $2^{-R}$,  while $\diam(\SCL_1) \geq 100 \cdot 2^{-R}$), so that we have $\SCL_1 \in \{ \widetilde{\SCL}_2, \SCL_2\}$ (a loop which meets $\overline{G}$ and is neither $\widetilde{\SCL}_2$ nor $\SCL_2$ would be contained in $G$,  since it cannot cross the rectangles) and thus $\SCL_1= \widetilde{\SCL}_2$. In the case $i>2$, we know that $P$ comes from outside of $A_{i-1}$ (since $\SCL_1$ is too large to be contained in $G$) and that $P$ hits $\widetilde{\SCL}_i$ before hitting $\SCL_i$. Assume by contradiction that $\widetilde{\SCL}_i \neq \SCL_{i-1}$. Since $\vert P(t_{i-1}) -z_{i-1} \vert \le 2^{-(1+\nu)RK}/50$, we know that $\SCL_{i-1}$ intersects $\overline{G}$. Since $\SCL_{i-1} \not\in \{\widetilde{\SCL}_i, \SCL_i\}$, we have actually $\SCL_{i-1} \subset G$. But since $P$ comes from outside of $A_{i-1}$ and does not cross the rectangles $R_{1, 5}, R_{2, 5}$, we see that $P$ hits $\widetilde{\SCL}_i$ before hitting $\SCL_{i-1}$ (indeed,  in order to reach $\SCL_{i-1} \subset G$ it has to enter $G$,  and the boundary of $G$ consists of arcs of $\widetilde{\SCL}_i$ and $\SCL_i$ together with parts of $R_{1,5}$ and $R_{2,5}$;  recall also that $P$ hits $\widetilde{\SCL}_i$ before $\SCL_i$). Then, there are two cases:
	\begin{itemize}
		\item If $\widetilde{\SCL}_{i-1}= \widetilde{\SCL}_i$ or $\widetilde{\SCL}_{i-1}$ is hit after $\widetilde{\SCL}_{i}$, then we can modify the sequences $\SCL_1, \ldots, \SCL_N$ and $A_1, \ldots, A_{N-1}$ by removing $\SCL_{i-1}$ and $A_{i-2}$. The new sequences will still satisfy the above properties. See the left-hand side of \Cref{fig:chain_of_annuli_pathological}.
		\item If $\widetilde{\SCL}_{i-1}$ is hit before $\widetilde{\SCL}_{i}$, then $\widetilde{\SCL}_{i-1}$ is outside $G$. Using the fact that the loops are disjoint, we see that the rectangles $R_{1,5}, R_{2,5}$ associated with $A_{i-1}$ and the rectangles $R'_{1,5}, R'_{2,5}$ associated in the same way with $A_{i-2}$ must cross, which contradicts condition~\eqref{it:left_to_right_crossings_shorter_than_top_to_bottom} in \Cref{prop:conformal_rectangles} (a crossing of two such rectangles would force one of the loops $\widetilde{\SCL}_{i-1},\SCL_{i-1},\widetilde{\SCL}_i,\SCL_i$,  which contain their left and right sides,  to intersect both the top and the bottom boundary of one of them). See the right-hand side of \Cref{fig:chain_of_annuli_pathological}.
	\end{itemize}
	\begin{figure}[h]
		\centering
		$\vcenter{\hbox{\includegraphics[scale=0.8]{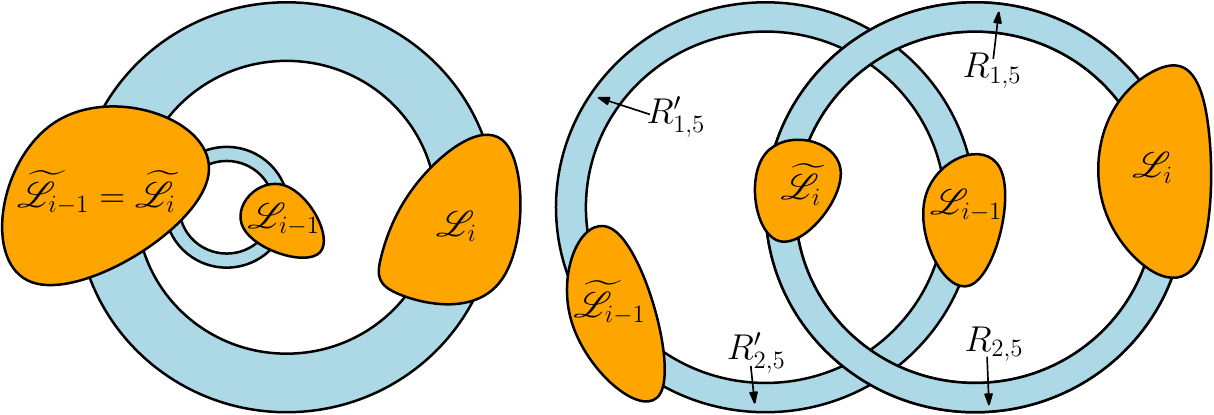}}}$
		\caption{The two cases we get rid of in the construction of the chain of loops. On the left, the case $\widetilde{\SCL}_{i-1}= \widetilde{\SCL}_i$. On the right, the case where $\widetilde{\SCL}_{i-1}$ is hit before $\widetilde{\SCL}_{i}$}
		\label{fig:chain_of_annuli_pathological}
	\end{figure}
	In that way we ensure that $\widetilde{\SCL}_i = \SCL_{i-1}$.
	}
	
	\stepn{step:main-result-comparison}{Comparison of $D_\Gamma^0$ and $D_\Gamma$} Fix $t \in (t_{i-1} ,  s_i)$ such that $P(t) \in \ell$ for some $\ell \in \Gamma$.  We will show that either the time-reversal of $P|_{[t_{i-1},t]}$ or $P|_{[t,s_i]}$ is part of a $D_\Gamma$-geodesic connecting $\ell$ with $\partial \BD$.  Note first that necessarily $\ell \notin \{\widetilde{\SCL}_i ,  \SCL_i\}$,  because $t_{i-1}$ (resp.\ $s_i$) is the last (resp.\ the first) time that $P$ intersects $\widetilde{\SCL}_i = \SCL_{i-1}$ (resp.\ $\SCL_i$).  As for the two endpoints $t \in \{t_{i-1} , s_i\}$,  which also have to be considered in the definition of $D_\Gamma^0$,  the assertion is immediate,  since any path contained in a single loop has zero $D_\Gamma$-length and hence a $D_\Gamma$-geodesic from that loop to $\partial \BD$ may be taken to start at $P(t)$.  Let now $Q$ be a $D_\Gamma$-geodesic connecting $\ell$ with $\partial \BD$.  First,  we note that $\ell \subseteq \overline{G}$
	 since no loop in $\Gamma$ can cross any of the topological rectangles $R_{j,m}$.  Since $Q$ has to leave $G$ in order to reach $\partial \BD$,  and since it would then have to cross the four rectangles $R_{j,1},\ldots,R_{j,4}$ for some $j \in \{1,2\}$ between their top and bottom boundaries if it were to avoid $\widetilde{\SCL}_i$ and $\SCL_i$ (this is the reason for our choice of $R_{1,5}$ and $R_{2,5}$ in the definition of $G$),  the argument in \Cref{lem:geodesic_intersects_loops} implies that $Q$ has to intersect either $\widetilde{\SCL}_i$ or $\SCL_i$. Therefore, we obtain that we can modify $Q$ so that it agrees with either the time-reversal of $P|_{[t_{i-1},t]}$ or $P|_{[t,s_i]}$ stopped at the first time that it intersects one of the loops in $\{\widetilde{\SCL}_i,\SCL_i\}$ (these two paths are $D_\Gamma$-geodesics from $\ell$ to $\SCL_{i-1}$ and to $\SCL_i$ respectively,  being sub-paths of the $D_\Gamma$-geodesic $P$ whose endpoints lie on loops,  which have zero $D_\Gamma$-diameter). The modified path $\widetilde{Q}$ will still be a geodesic from $\ell$ to $\partial \BD$ (indeed, in a geodesic metric space, if $\gamma$ is a geodesic from $x$ to $z$ and goes through a point $y$ and if $\eta$ is a geodesic from $x$ to $y$, if we replace the part of $\gamma$ from $x$ to $y$ by $\eta$, then the obtained path is a geodesic). Using the fact that $\widetilde{\SCL}_i = \SCL_{i-1}$, it follows that
	\begin{align}\label{eqn:main_loop_inequality}
		D_{\Gamma}^0(\SCL_{i-1},\SCL_i) \leq D_\Gamma(\SCL_{i-1},\SCL_i)
	\end{align}
	by the definition of $D_\Gamma^0$:  indeed,  the above shows that $P|_{[t_{i-1},s_i]}$ is one of the admissible paths over which the infimum defining $D_\Gamma^0(\SCL_{i-1},\SCL_i)$ is taken,  and $\len(P|_{[t_{i-1},s_i]} ; D_\Gamma) = D_\Gamma(\SCL_{i-1},\SCL_i)$ since $P$ is a $D_\Gamma$-geodesic and the loops have zero $D_\Gamma$-diameter.
	
	\stepn{step:main-result-conclusion}{Conclusion of the proof} Combining~\eqref{eqn:main_loop_inequality} with the definition of $\widetilde D_\Gamma$ and the fact that $P$ is a $D_\Gamma$-geodesic connecting $\SCL$ with $\SCL'$,  we obtain that
	\begin{align*}
		\widetilde D_\Gamma(\SCL,\SCL') \leq \sum_{j=1}^{N-1} D_{\Gamma}^0(\SCL_j,\SCL_{j+1}) \leq \sum_{j=1}^{N-1} D_\Gamma(\SCL_j,\SCL_{j+1}) = \sum_{j=1}^{N-1} \len(P|_{[t_j ,  s_{j+1}]} ; D_\Gamma) = D_\Gamma(\SCL,\SCL').
	\end{align*}
	Here,  the last equality follows from \Cref{lem:intersection_is_a_segment} together with the fact that any path contained in a single loop has zero $D_\Gamma$-length,  so that $\len(P ; D_\Gamma) = \sum_{j=1}^{N-1} \len(P|_{[t_j , s_{j+1}]} ; D_\Gamma)$.
	This completes the proof of the theorem since clearly we have that $D_\Gamma(\SCL,\SCL') \leq \widetilde D_\Gamma(\SCL,\SCL')$, the points $x,y \in \BD_{\BQ}$ were arbitrary,  and a.s.\ every loop in $\Gamma$ is of the form $\SCL(w)$ for some $w \in \BD_\BQ$.
\end{proof}

\subsection{The genealogy is given by the uniform exploration}
\begin{figure}[h]
	\centering
	$\vcenter{\hbox{\includegraphics[scale=0.8]{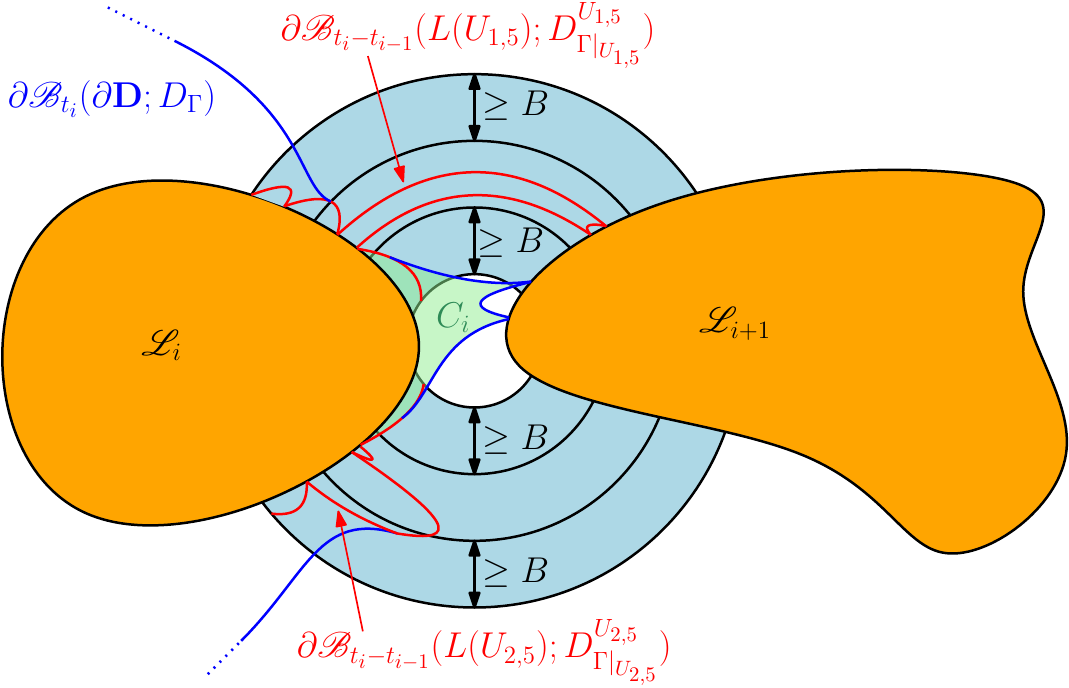}}}$
	\caption{Illustration of the proof of \Cref{prop D tilde does not depend on the metric}. The blue curve is the boundary of the ball from $\partial \BD$ whereas the red curves are the boundaries of the balls starting from the left-hand sides of the rectangles $U_{1,5}$ and $U_{2,5}$ which are drawn in blue. Recall from \Cref{prop:conformal_rectangles} that $U_{1,5}$ and $U_{2,5}$ are divided into three rectangles. Since the distances across the outer rectangles are at least $B$ and since the $D_\Gamma$ distance between $\SCL_i$ and $\SCL_{i+1}$ is smaller than $B$, the red curves and the blue curve agree in the inner rectangles. Finally, by Condition~\eqref{condition b genealogy} of \Cref{prop:conformal_rectangles}, we know that the red curves in each of the inner rectangles have to hit the left boundary. This entails that $C_i$ is included in $\overline{G \cup R_{1,5}\cup R_{2, 5} \cup \mathop{\mathrm{int}}(\SCL_{i+1})}$.}
	\label{fig:genealogy_is_determined}
\end{figure}
Let us end this section by proving \Cref{prop D tilde does not depend on the metric}.
\begin{proof}[Proof of \Cref{prop D tilde does not depend on the metric}]
	\stepn{step:dtilde-indep-setup}{Setup and reduction} Let $x ,y\in \BD_\BQ$ be two distinct points. The converse implication being the content of \Cref{lemma the genealogy of the uniform exploration is contained in the geodesics to the root},  it is enough to show that if $\SCL(x)$ intersects a $D_\Gamma$-geodesic from $\SCL(y)$ to $\partial \BD$, then $\SCL(x) \preceq \SCL(y)$. Let $P$ be a $D_\Gamma$-geodesic from $\SCL(x)$ to $\SCL(y)$.

	By assumption, we can construct a $D_\Gamma$-geodesic $Q$ from $\partial \BD$ to $\SCL(y)$ such that $P\subset Q$ (indeed we know that $\SCL(x)$ intersects a $D_\Gamma$-geodesic $\widetilde{Q}$ from $\partial \BD$ to $\SCL(y)$ and then we replace the part from $\SCL(x)$ to $\SCL(y)$ by $P$). Recall from \Cref{rem:good_event_arbitrary_scale} the definition of the events $E_{z,k}(\Psi)$. Recall from the proof of \Cref{thm:main_result} that we constructed a sequence of loops $\SCL_1 = \SCL(x), \SCL_2, \ldots, \SCL_N=\SCL(y)$ and of Euclidean balls $B_{2^{-k_1}}(z_1), \ldots, B_{2^{-k_{N-1}}}(z_{N-1})$ such that the geodesic $P$ from $\SCL(x)$ to $\SCL(y)$ intersects $\SCL_j$ before intersecting $\SCL_{j+1}$ for all $1\le j \le N-1$. Secondly, for all $1 \le j \le N-1$, letting $A_j\defeq B_{2^{-k_j-1}}(z_j) \setminus \overline{B_{2^{-k_j - R}}(z_j)}$, we know that $\SCL_{j+1}$ and $\SCL_j$ are the only two loops which each make a crossing of $A_j$ and that the event $E_{z_j, k_j}(\Psi)$ occurs for all $1 \le j \le N-1$. Thirdly, we know that $P \cap B_{2^{-k_j-R}}(z_j) \neq \emptyset$ for all $1 \le j \le N-1$ and that $P\subseteq \bigcup_{j=1}^{N-1} B_{2^{-k_j}}(z_j) \cup \bigcup_{j=1}^N \SCL_j$. 
	
	\stepn{step:dtilde-indep-times}{The times $r_i$ and the components $C_i$}

	\substepn{step:dit-claim}{Statement of the claim} Recall from \Cref{section:genealogy-of-loops} the definition of the genealogy of loops. We want to prove that there exist $t_x=r_0<r_1< \cdots< r_{N-1}=t_y$ and connected components $C_i$ for $1 \le i \le N-1$ of respectively $\SCB_{r_i}(\partial \BD; D_\Gamma) \setminus \SCB_{r_{i-1}}(\partial \BD; D_\Gamma)$ for $1 \le i \le N-1$ such that $C_1$ is a descendant of $\SCL(x)$, and for all $2 \le i \le N-1$, the component $C_i$ is a descendant of the component $C_{i-1}$, with $\SCL(y) \subseteq C_{N-1}$.  We use the letters $r_i$ here since $t_1,\ldots,t_N$ already denote times along $P$.

	\substepn{step:dit-definitions}{Definition of $r_i$ and $C_i$} We set for all $i \in [0, N-1]_\BZ$,
	\[
	r_i \defeq D_\Gamma(\partial \BD, \SCL_{i+1}),
	\]
	and let $C_i$ be the connected component of $\SCB_{r_i}(\partial \BD; D_\Gamma) \setminus \SCB_{r_{i-1}}(\partial \BD; D_\Gamma)$ which contains $\SCL_{i+1}$ for $1 \leq i \leq N-1$;  note that $r_0 = D_\Gamma(\partial \BD , \SCL(x)) = t_x$ and $r_{N-1} = D_\Gamma(\partial \BD , \SCL(y)) = t_y$.

	\substepn{step:dit-monotonicity}{Monotonicity of the $r_i$} Note that $r_0 < r_1 < \cdots < r_{N-1}$,  as required.  Indeed,  by our assumption that $\SCL(x)$ is in $Q$, a $D_\Gamma$-geodesic from $\SCL(y)$ to $\partial \BD$, 
		$D_\Gamma(\partial \BD , \ell) = D_\Gamma(\partial \BD , \SCL(x)) + D_\Gamma(\SCL(x) , \ell)$ for every loop $\ell$ which $P$ intersects.  Applying this with $\ell = \SCL_{i+1}$ and using that the $\SCL_i$'s are pairwise distinct,  the strict monotonicity follows.  This is the only place where the standing assumption of the present proof is used.

	\stepn{step:dtilde-indep-locating}{Locating the components $C_i$} Let $i \in [1, N-1] \cap \BZ$. Note that $C_i$ contains every geodesic from $\SCL_i$ to $\SCL_{i+1}$ which is contained in the bounded connected component $G$ of $\BC \setminus (\mathop{\mathrm{int}}(\SCL_i) \cup \mathop{\mathrm{int}}(\SCL_{i+1} )\cup R_{1,5} \cup R_{2,5} )$, where the $R_{j,m}$ for $j=1,2$ and $m=1, \ldots, 8$ are the topological rectangles used in the previous proof (so that $G$ is the set denoted in the same way in the proof of \Cref{thm:main_result},  with $i$ there replaced by $i+1$).  Indeed,  that geodesic is $P|_{[t_i,s_{i+1}]}$,  which is contained in $G$ by the proof of \Cref{thm:main_result},  and every loop $\ell$ which it intersects satisfies $r_{i-1} \leq D_\Gamma(\partial \BD , \ell) \leq r_i$ by the identity displayed above,  so that it is contained in $\SCB_{r_i}(\partial \BD ; D_\Gamma) \setminus \SCB_{r_{i-1}}(\partial \BD ; D_\Gamma)$ and hence,  being connected,  in $C_i$.
	
	Moreover, note that by Condition~\eqref{condition a genealogy} in \Cref{prop:conformal_rectangles}, the $D_\Gamma$-distance between $\SCL_i$ and $\SCL_{i+1}$ is smaller than $B$.

	\stepn{step:dtilde-indep-conclusion}{Conclusion of the proof} Therefore, by Conditions~\eqref{condition b genealogy} and~\eqref{condition c genealogy} in \Cref{prop:conformal_rectangles}, we see that $C_i\subseteq\overline{G \cup R_{1,5}\cup R_{2, 5} \cup \mathop{\mathrm{int}}(\SCL_{i+1})}$. See \Cref{fig:genealogy_is_determined} for more details. In particular, \[\overline{C_i} \cap \SCB_{r_{i-1}}(\partial \BD ; D_\Gamma)\subseteq \SCL_{i}.\]
	This proves that $C_i$ is a descendant of $C_{i-1}$ where by convention $C_0 = \SCL(x)$. Thus, by definition of $\preceq$, we have $\SCL(x) \preceq \SCL(y)$.  Finally,  the ``in particular'' assertion of \Cref{prop D tilde does not depend on the metric} follows,  since $\preceq$ is constructed from the uniform exploration alone.  
More precisely,  the metric $D_{\Gamma}^0$ (and hence $\widetilde{D}_{\Gamma}$) is determined by the uniform exploration.  To see the latter,  note that the paths appearing in the definition of $D_{\Gamma}^0$ refer only to the geodesics to $\partial \BD$.  However,  the argument above allows us to identify the geodesics to $\partial \BD$ via the genealogy induced by $\preceq$ and the latter is determined by the uniform exploration.
\end{proof}

\section{Measurability}
\label{sec:measurability}
By \Cref{thm metric determined by the exploration}, we know that the metric of a geodesic $\CLE_4$ metric coupling in the sense of \Cref{def:weak_axioms} is a measurable function of the associated uniform exploration; in particular, given a uniform exploration of the $\CLE_4$, there is at most one such metric whose metric balls from the boundary agree with this uniform exploration. The following result shows that the metric, and in particular the uniform exploration, is a.s.\ determined by the $\CLE_4$, thus implying that the metric of any geodesic $\CLE_4$ metric coupling coincides a.s.\ with a strong $\CLE_4$ metric in the sense of \Cref{def:axioms}.  This will complete the proofs of Theorems~\ref{thm:convergence_nonsimple_cle} and \ref{thm:uniqueness} and will enable us to prove Corollary~\ref{thm:uniform-exploration-measurability}.
\begin{proposition}\label{lem:measurability}
	Let $D$ be a geodesic $\CLE_4$ metric coupling. Let $\Gamma_U$ be a non-nested $\CLE_4$ in a deterministic simply connected domain $U \subsetneq \BC$. Then $D^U_{\Gamma_U}$ is a.s.\ determined by $\Gamma_U$; that is, $D^U_{\Gamma_U}$ is a.s.\ equal to a measurable function of $\Gamma_U$.
\end{proposition}
The above proposition allows us to upgrade every geodesic $\CLE_4$ metric coupling to a strong $\CLE_4$ metric, as follows.
\begin{proposition}\label{lem:weak-is-strong}
	Let $D$ be a geodesic $\CLE_4$ metric coupling. Then there exists a strong $\CLE_4$ metric $\widetilde{D}$ in the sense of \Cref{def:axioms} such that, for every simply connected domain $U \subsetneq \BC$, we have $D^U_{\Gamma_U} = \widetilde{D}^U_{\Gamma_U}$ a.s.
\end{proposition}
\begin{proof}[Proof of Proposition~\ref{lem:weak-is-strong} assuming Proposition~\ref{lem:measurability}.] By Proposition~\ref{lem:measurability} and the Doob--Dynkin lemma, for each simply connected domain $U \subsetneq \BC$ there is a measurable mapping $\widetilde{D}^U \colon\allowbreak \mathscr{LoopEns}(U) \to \BR^{U_\BQ \times U_\BQ}$ such that $D^U_{\Gamma_U} = \widetilde{D}^U_{\Gamma_U}$ a.s. It remains to check that the family $\widetilde{D} = \{\widetilde{D}^U\}_U$ satisfies the axioms of \Cref{def:axioms}. Axioms~\eqref{it:axiom_geodesic} (geodesic metric) and~\eqref{it:axiom_uniform_exploration} (uniform exploration) are immediate from Axioms~\eqref{it:weak_axiom_geodesic} and~\eqref{it:axiom_uniform_exploration_coupling} of \Cref{def:weak_axioms}, since $\widetilde{D}^U_{\Gamma_U} = D^U_{\Gamma_U}$ a.s.\ and these two axioms only involve the law of the pair $(\Gamma_U, D^U_{\Gamma_U})$. Next, we check Axiom~\eqref{it:axiom_locality} (locality). Let $V \subseteq U$ be a deterministic simply connected subdomain and let $\{V_j\}_j$ be the connected components of $V^\star$. By Axiom~\eqref{it:weak_axiom_locality} of \Cref{def:weak_axioms}, conditionally on~\eqref{eq:weak_axiom_locality}, the pair $(\Gamma_U|_{V_j}, D_{\Gamma_U}^U(\bullet, \bullet; V_j))$ has the law of $(\Gamma_{V_j}, D^{V_j}_{\Gamma_{V_j}})$, and hence also the law of $(\Gamma_{V_j}, \widetilde{D}^{V_j}_{\Gamma_{V_j}})$. Since in the latter pair the second coordinate is a.s.\ equal to a fixed measurable function of the first, the same must be true of the former, i.e.\ a.s.\ $D_{\Gamma_U}^U(\bullet, \bullet; V_j) = \widetilde{D}^{V_j}_{\Gamma_U|_{V_j}}$ for every $j$. Finally, we check Axiom~\eqref{it:axiom_conformal_invariance} (conformal invariance). Let $\phi \colon U \to \phi(U)$ be a deterministic conformal mapping. By Axiom~\eqref{it:weak_axiom_conformal_invariance} of \Cref{def:weak_axioms}, the pairs
	\[
	\left(\phi(\Gamma_U), \left( \widetilde{D}_{\Gamma_U}^U\left(\phi^{-1}(\SCL_1), \phi^{-1}(\SCL_2) \right) \right)_{\SCL_1, \SCL_2 \in \phi(\Gamma_U)} \right) \quad \text{and} \quad \left( \Gamma_{\phi(U)}, \left(\widetilde{D}_{\Gamma_{\phi(U)}}^{\phi(U)}(\SCL_1, \SCL_2) \right)_{\SCL_1, \SCL_2 \in \Gamma_{\phi(U)}} \right)
	\]
	have the same law. In both pairs the second coordinate is a measurable function of the first, and the first coordinates have the same law (by the conformal invariance of the $\CLE_4$), so the two measurable functions agree almost everywhere with respect to the law of a non-nested $\CLE_4$ in $\phi(U)$. This is precisely Axiom~\eqref{it:axiom_conformal_invariance}.
\end{proof}

\begin{proof}[Proof of Theorem~\ref{thm:uniqueness} assuming Proposition~\ref{lem:measurability}.]
	Assertion~\eqref{thm:uniqueness-existence} is exactly Proposition~\ref{lem:weak-is-strong}. As for assertion~\eqref{thm:uniqueness-uniqueness}, let $D$ and $\widetilde{D}$ be two strong $\CLE_4$ metrics and let $U \subsetneq \BC$ be a deterministic simply connected domain. Each of $D$ and $\widetilde{D}$ gives rise to a geodesic $\CLE_4$ metric coupling, so that, by Theorem~\ref{thm metric determined by the exploration}, each of $D^U_{\Gamma_U}$ and $\widetilde{D}^U_{\Gamma_U}$ is given by one and the same measurable function of the associated uniform exploration; since $\Gamma_U$ is itself a measurable function of that exploration, and since the law of the latter does not depend on the metric by Axiom~\eqref{it:axiom_uniform_exploration} (uniform exploration), the couplings $(\Gamma_U, D^U_{\Gamma_U})$ and $(\Gamma_U, \widetilde{D}^U_{\Gamma_U})$ have the same law. Since $D^U_{\Gamma_U}$ and $\widetilde{D}^U_{\Gamma_U}$ are both measurable functions of $\Gamma_U$, it follows that $D^U$ and $\widetilde{D}^U$ agree almost everywhere with respect to the law of a non-nested $\CLE_4$ in $U$, i.e.\ $D^U_{\Gamma_U} = \widetilde{D}^U_{\Gamma_U}$ a.s.
\end{proof}

\begin{proof}[Proof of \Cref{thm:convergence_nonsimple_cle} assuming Proposition~\ref{lem:measurability}.]
	By \cite[Theorem~1.3]{kkmt2026cle4_part1}, the family of the laws of the metrics $\ka_\kappa^{-1} D^U_{\Gamma_U^\kappa}$ for $\kappa \in (4,8)$ is tight with respect to the topology of $\BR^{U_\BQ \times U_\BQ}$, and every subsequential limit as $\kappa \downarrow 4$ is the metric of a weak geodesic $\CLE_4$ metric coupling. By \cite[Theorem~1.4]{kkmt2026cle4_part2}, every weak geodesic $\CLE_4$ metric coupling is a geodesic $\CLE_4$ metric coupling, and by Theorem~\ref{thm metric determined by the exploration} the law of the latter is unique; by Proposition~\ref{lem:weak-is-strong}, it is the law of $(\Gamma_U, D^U_{\Gamma_U})$ for the strong $\CLE_4$ metric $D$ of Theorem~\ref{thm:uniqueness}. Since all subsequential limits agree in law, the convergence in law follows.
\end{proof}
The proof of \Cref{lem:measurability} will follow the same idea as in \cite[Section 7]{ms2019lqg_bm1}, which constructed a metric on $\sqrt{8/3}$-Liouville quantum gravity that was later shown in \cite{ms2021lqg_bm2, ms2021lqg_bm3} to be isometric to the Brownian sphere. However, since in our case the metric has already been constructed, the proof will be much easier.

By Axiom~\eqref{it:weak_axiom_conformal_invariance} (conformal invariance) of \Cref{def:weak_axioms}, in order to prove Proposition~\ref{lem:measurability}, it suffices to consider the case of $U=\BD$. Let $\Gamma$ be a non-nested $\CLE_4$ in $\BD$ and, conditionally on $\Gamma$, let $D$ and $\overline{D}$ be two independent copies of $D^\BD_\Gamma$ (i.e.\ conditionally on $\Gamma$, the metrics $D$ and $\overline{D}$ are independent and each of them has the conditional law of $D^\BD_\Gamma$ given $\Gamma$). We set $\SCB\defeq (\SCB_t)_{t\ge 0} \defeq (\SCB_t(\partial \BD; D))_{t\ge 0}$ and $\overline{\SCB} \defeq (\overline{\SCB}_t)_{t\ge 0} \defeq (\SCB_t(\SCL(0); \overline{D}))_{t\ge 0}$.  We also set for every deterministic or random time $t$, $\SCB_{\le t}\defeq (\SCB_s)_{s\le t}$ and $\overline{\SCB}_{\le t} \defeq (\overline{\SCB}_s)_{s\le t}$. Let $W$ be a uniform random variable on $[0,1]$ which is independent of $(\Gamma, D, \overline{D})$. 
We aim to prove the following result.
\begin{proposition}\label{prop D equals D tilde}
	$D(\partial \BD, \SCL(0)) = \overline{D}(\partial \BD, \SCL(0))$ a.s.
\end{proposition}
One can indeed deduce \Cref{lem:measurability} as follows.
\begin{proof}[Proof of \Cref{lem:measurability} given \Cref{prop D equals D tilde}]
	By conformal invariance, for each fixed $x \in \BD_\BQ$, applying \Cref{prop D equals D tilde} to the image of $(\Gamma, D, \overline{D})$ under a M\"obius transformation of $\BD$ sending $x$ to $0$ (by Axiom~\eqref{it:weak_axiom_conformal_invariance} of \Cref{def:weak_axioms}, the law of $(\Gamma, D)$ is invariant under such a transformation, and hence so is the law of the triple $(\Gamma, D, \overline{D})$, since conditionally on $\Gamma$ the metrics $D$ and $\overline{D}$ are independent copies of $D^\BD_\Gamma$), we deduce that a.s.\ $D(\partial \BD, \SCL(x)) = \overline{D}(\partial \BD, \SCL(x))$. Since $\BD_\BQ$ is countable, a.s.\ this holds simultaneously for all $x \in \BD_\BQ$. Since every loop of $\Gamma$ is of the form $\SCL(x)$ for some $x \in \BD_\BQ$, and since the explored region at time $t$ is determined by $\Gamma$ together with the times at which its loops are discovered (see \Cref{subsec:labelled_cle_4}), we obtain the almost sure equality of metric balls
	\[
	\left(\SCB_t(\partial \BD; D) \right)_{t\ge 0} = \left( \SCB_t(\partial \BD; \overline{D})\right)_{t\ge 0}.
	\]
	But, by \Cref{thm metric determined by the exploration}, $D$ (resp.\ $\overline{D}$) is a function of the uniform exploration $(\SCB_t(\partial \BD; D) )_{t\ge 0}$ (resp.\ $( \SCB_t(\partial \BD; \overline{D}))_{t\ge 0}$), so that $D=\overline{D}$ a.s. Finally, since conditionally on $\Gamma$ the metrics $D$ and $\overline{D}$ are independent with the same conditional law, the almost sure equality $D = \overline{D}$ forces this conditional law to be a.s.\ a Dirac mass; since the metrics take values in the Polish space $\BR^{\BD_\BQ \times \BD_\BQ}$, this means precisely that $D$ is a.s.\ equal to a measurable function of $\Gamma$.
\end{proof}

Let $\tau \defeq W \cdot D (\partial \BD, \SCL(0))$, $\overline{\tau} \defeq \overline{D}(\SCB_\tau, \SCL(0))$, $\overline{\sigma} \defeq W \cdot \overline{D}(\partial \BD, \SCL(0))$ and $\sigma \defeq D(\overline{\SCB}_{\overline{\sigma}}, \partial \BD)$.

Recall from Axiom~\eqref{it:axiom_uniform_exploration_coupling} (uniform exploration) of \Cref{def:weak_axioms} that $D(\partial \BD, \SCL(0))$ is distributed as the time at which the loop surrounding the origin is discovered by the uniform exploration, which is an exponential random variable (see Section~\ref{subsec:labelled_cle_4}, where this time is the first arrival time of a Poisson point process); in particular, $\BE[D(\partial \BD, \SCL(0))] = \BE[\overline{D}(\partial \BD, \SCL(0))] \in (0,\infty)$, so that the following definition makes sense.
We define the probability measures $\nu$ and $\overline{\nu}$ obtained by weighting $\BP$ respectively by $D(\partial \BD, \SCL(0))$ and $\overline{D}(\partial \BD, \SCL(0))$.  In other words, the expectations $\BE_\nu$ and $\BE_{\overline{\nu}}$ under $\nu$ and $\overline{\nu}$ are respectively defined by setting for every non-negative random variable $X$,
\begin{align*}
\BE_\nu[X]&\defeq \frac{1}{\BE[D(\partial \BD, \SCL(0))]} \BE \left[ D(\partial \BD, \SCL(0)) \cdot X \right] \quad \text{and}\\
\BE_{\overline{\nu}}[X]&\defeq \frac{1}{\BE[\overline{D}(\partial \BD, \SCL(0))]} \BE \left[ \overline{D}(\partial \BD, \SCL(0)) \cdot X \right].
\end{align*}
In particular, for every random element $X$, taking values in an arbitrary measurable space, which is a measurable function of $(\Gamma, D, \overline{D})$, and for every non-negative measurable function $f$,
\begin{align*}
	\BE_\nu[f(X,\tau, \overline{\tau})] &=  \frac{1}{\BE[D(\partial \BD, \SCL(0))]} \int_0^\infty \mathrm{d} t \  \BE \left[ \one_{t \le D(\partial \BD, \SCL(0))} f(X,t, \overline{D}(\SCL(0), \SCB_t)) \right]
	\\
	\text{and} \qquad
	\BE_{\overline{\nu}} [f(X,\sigma , \overline{\sigma})] &=  \frac{1}{\BE[\overline{D}(\partial \BD, \SCL(0))]} \int_0^\infty \mathrm{d} t \  \BE \left[ \one_{t \le \overline{D}(\partial \BD, \SCL(0))} f(X,D(\partial \BD, \overline{\SCB}_t), t) \right],
\end{align*}
where both identities are obtained by first integrating over the value $u \in [0,1]$ of $W$, which is independent of $(\Gamma, D, \overline{D})$, and then substituting $t = u D(\partial \BD, \SCL(0))$ (resp.\ $t = u \overline{D}(\partial \BD, \SCL(0))$); we also recall that $\BE[\overline{D}(\partial \BD, \SCL(0))] = \BE[D(\partial \BD, \SCL(0))]$.
The key result of this section is the following identity of laws.
\begin{proposition}\label{prop nu stopped equals nu tilde stopped}
	The $\nu$-law of $(\SCB_{\le \tau}, \overline{\SCB}_{\le \overline{\tau}}, \Gamma)$ is equal to the $\overline{\nu}$-law of $(\SCB_{\le \sigma}, \overline{\SCB}_{\le \overline{\sigma}}, \Gamma)$.
\end{proposition}
\begin{proof}
	\stepn{step:nu-stopped-chain}{The chain of identities} Let $f$ be a non-negative measurable function. We have
	\begin{align*}
		\BE&[D(\partial \BD, \SCL(0))] \cdot\BE_\nu \left[ f\left( \SCB_{\le \tau}, \overline{\SCB}_{\le \overline{\tau}}, \Gamma \right) \right] \\&= 
		\int_0^\infty \mathrm{d}t \ \BE\left[ \one_{t \le D(\partial \BD, \SCL(0))} f(\SCB_{\le t}, \overline{\SCB}_{\le \overline{D}(\SCL(0),\SCB_t)}, \Gamma) \right]\\
		&=\int_0^\infty \mathrm{d}t \ \BE\left[ \one_{t \le D(\partial \BD, \SCL(0))} f(\SCB_{\le t}, (\SCB_s(\SCL(0); D))_{s \le D(\SCL(0),\SCB_t)}, \Gamma) \right]\\
		&= \BE\left[\int_0^{D(\partial \BD, \SCL(0))} \mathrm{d}t \ f(\SCB_{\le t}, (\SCB_s(\SCL(0); D))_{s \le D(\SCL(0),\partial \BD)-t}, \Gamma) \right]
		\\
		&= \BE\left[\int_0^{D(\partial \BD, \SCL(0))} \mathrm{d}t \ f(\SCB_{\le D(\SCL(0),\partial \BD)-t}, (\SCB_s(\SCL(0); D))_{s \le t}, \Gamma) \right]
		\\
		&= \int_0^\infty \mathrm{d}t \ \BE\left[ \one_{t \le \overline{D}(\partial \BD, \SCL(0))}
		f((\overline{\SCB}_s
		(\partial \BD; \overline{D}))_{s \le \overline{D}(\partial \BD, \SCL(0)) -t}, \overline{\SCB}_{\le t}, \Gamma)\right] \\
		&= \int_0^\infty \mathrm{d}t \ \BE \left[  \one_{t \le \overline{D}(\partial \BD, \SCL(0))} 
		f(\SCB_{\le D(\partial \BD, \overline{\SCB}_t)}, \overline{\SCB}_{\le t}, \Gamma)\right]\\
		&=\BE[D(\partial \BD, \SCL(0))] \cdot \BE_{\overline{\nu}}\left[ f(\SCB_{\le \sigma}, \overline{\SCB}_{\le \overline{\sigma}}, \Gamma)\right],
	\end{align*}

	\stepn{step:nu-stopped-justification}{Justification of the equalities} Here, in the first equality, we used the definition of $\nu$  (i.e.\ the first identity displayed above, applied with $X = (\SCB, \overline{\SCB}, \Gamma)$); in the second equality, we used the locality property (\Cref{lem:metric-ball-locality}) together with the fact that, conditionally on $\Gamma$, the metrics $D$ and $\overline{D}$ are independent with the same conditional law. More precisely, on the event $\{t \le D(\partial \BD, \SCL(0))\}$, write $C_t$ for the connected component of $\BD \setminus \SCB_t$ containing the origin. On the one hand, conditionally on $\SCB_{\le t}$, the pair $(\Gamma|_{C_t}, D(\bullet, \bullet; C_t))$ has the law of $(\Gamma_{C_t}, D^{C_t}_{\Gamma_{C_t}})$, so that, conditionally on $(\Gamma, \SCB_{\le t})$, the internal metric $D(\bullet, \bullet; C_t)$ has the conditional law of $D^{C_t}_{\Gamma_{C_t}}$ given $\Gamma_{C_t} = \Gamma|_{C_t}$; moreover, the $D$-balls centered at $\SCL(0)$ agree with the corresponding $D(\bullet, \bullet; C_t)$-balls up until they hit $\partial C_t$. On the other hand, $\SCB_{\le t}$ is a measurable function of $(\Gamma, D)$ and is therefore conditionally independent of $\overline{D}$ given $\Gamma$; consequently the conditional law of $\overline{D}$ given $(\Gamma, \SCB_{\le t})$ is equal to its conditional law given $\Gamma$, and locality applied to $\overline{D}$, with $C_t$ frozen, shows that $\overline{D}(\bullet, \bullet; C_t)$ has the same conditional law given $(\Gamma, \SCB_{\le t})$ as $D(\bullet, \bullet; C_t)$. Consequently, conditionally on $(\Gamma, \SCB_{\le t})$, the growth of the $\overline{D}$-balls from $\SCL(0)$ until they hit $\SCB_t$ has the same law as the growth of the corresponding $D$-balls; in the third equality, we used Fubini's theorem together with the fact that $D(\SCL(0), \SCB_t) = D(\SCL(0), \partial \BD) -t$ for $t \le D(\partial \BD, \SCL(0))$, which follows from Axiom~\eqref{it:weak_axiom_geodesic} (geodesic metric) of \Cref{def:weak_axioms}; in the fourth equality, we performed the change of variables $t \mapsto D(\SCL(0), \partial \BD) - t$; in the fifth equality, we used Fubini's theorem again together with the fact that the law of $(D, \Gamma)$ is the same as the law of $(\overline{D}, \Gamma)$; in the sixth equality, we used the locality property once more, in the same way as in the second equality, together with the identity $\overline{D}(\partial \BD, \overline{\SCB}_t) = \overline{D}(\partial \BD, \SCL(0)) - t$ for $t \le \overline{D}(\partial \BD, \SCL(0))$; and in the last equality, we used the definition of $\overline{\nu}$ (i.e.\ the second identity displayed above).
\end{proof}
\begin{lemma}\label{lemma conditional laws agree}
	The conditional law under $\nu$ of $(\SCB, \overline{\SCB}, \Gamma)$ given $( \SCB_{\le \tau}, \overline{\SCB}_{\le \overline{\tau}}, \Gamma)$ is equal to the conditional law under $\BP$ of $(\SCB, \overline{\SCB}, \Gamma)$ given $(   \SCB_{\le \tau}, \overline{\SCB}_{\le \overline{\tau}}, \Gamma)$. Similarly, the conditional law under $\overline{\nu}$ of $(\SCB, \overline{\SCB}, \Gamma)$ given $( \SCB_{\le \sigma}, \overline{\SCB}_{\le \overline{\sigma}}, \Gamma)$ is equal to the conditional law under $\BP$ of $(\SCB, \overline{\SCB}, \Gamma)$ given $(   \SCB_{\le \sigma}, \overline{\SCB}_{\le \overline{\sigma}}, \Gamma)$.
\end{lemma}
\begin{proof}
	\stepn{step:cond-laws-def-h}{Definition of $h$} Let $f$ be a non-negative bounded measurable function. Let us set, for all $u, v \in [0,1]$,
	\begin{equation}\label{eq def h 0}
		h\left(\SCB_{\le u D(\partial \BD, \SCL(0))}, \overline{\SCB}_{\le v \overline{D}(\SCL(0),\partial \BD)}, \Gamma \right)\defeq \BE\left[\left. f( \SCB,  \overline{\SCB}, \Gamma) \right\vert \SCB_{\le u  D(\partial \BD, \SCL(0))}, \overline{\SCB}_{\le v\overline{D}(\SCL(0), \partial \BD)}, \Gamma \right].
	\end{equation}
We note that there exists a measurable function $h$ such that the equality in ~\eqref{eq def h 0} holds a.s.\ for Lebesgue almost every $(u,v) \in [0,1]^2$ (see \cite[Lemma~8.2]{ms2019lqg_bm1} and the discussion just before the statement of \cite[Lemma~8.2]{ms2019lqg_bm1}).  For the rest of the proof,  we will assume that we are working with a version of the conditional expectation in ~\eqref{eq def h 0} such that ~\eqref{eq def h 0} holds a.s.\ for Lebesgue almost every $(u,v) \in [0,1]^2$.

Note that the right-hand side of~\eqref{eq def h 0} is by definition a measurable function of $(\SCB_{\le u D(\partial \BD, \SCL(0))},\allowbreak \overline{\SCB}_{\le v \overline{D}(\SCL(0),\partial \BD)},\allowbreak \Gamma)$.

Then, since $\overline{D}(\SCL(0), \SCB_{u  D(\partial \BD, \SCL(0))})/\overline{D}(\SCL(0), \partial \BD)$ is a stopping time for the natural filtration of the process $(\SCB_{\le u D(\partial \BD, \SCL(0))},\allowbreak \overline{\SCB}_{\le v \overline{D}(\SCL(0), \partial \BD)},\allowbreak \Gamma)_{v \in [0,1]}$, the optional stopping theorem implies that (after approximating this stopping time from above by dyadic stopping times, and using that the ball processes are c\`adl\`ag by Axiom~\eqref{it:weak_axiom_geodesic} of \Cref{def:weak_axioms}), for all $u \in [0,1]$,
	\begin{equation}\label{eq def h}
	\begin{split}
		&h\left(\SCB_{\le u D(\partial \BD, \SCL(0))}, \overline{\SCB}_{\le \overline{D}(\SCL(0), \SCB_{u  D(\partial \BD, \SCL(0))})}, \Gamma \right)\\
		&\qquad= \BE\left[\left. f( \SCB,  \overline{\SCB}, \Gamma) \right\vert \SCB_{\le u  D(\partial \BD, \SCL(0))}, \overline{\SCB}_{\le \overline{D}(\SCL(0), \SCB_{u D(\partial \BD, \SCL(0))})}, \Gamma \right].
	\end{split}
	\end{equation}

	\stepn{step:cond-laws-stopped}{Proof of~\eqref{eq h equals h stopped at tau}} Let us first show that a.s.,
	\begin{equation}\label{eq h equals h stopped at tau}
		h\left(\SCB_{\le \tau}, \overline{\SCB}_{\le \overline{\tau}}, \Gamma \right)=\BE\left[\left. f( \SCB,  \overline{\SCB}, \Gamma) \right\vert \SCB_{\le \tau}, \overline{\SCB}_{\le \overline{\tau}}, \Gamma \right].
	\end{equation}
	Indeed, for every non-negative measurable function $g$,
	\begin{align*}
		\BE&\left[ h\left(\SCB_{\le \tau}, \overline{\SCB}_{\le \overline{\tau}}, \Gamma \right) g\left(\SCB_{\le \tau}, \overline{\SCB}_{\le \overline{\tau}}, \Gamma \right) \right]\\
		&=
		\int_0^1 \mathrm{d} u \ \BE\left[  h\left(\SCB_{\le u D(\partial \BD, \SCL(0))}, \overline{\SCB}_{\le \overline{D}(\SCL(0), \SCB_{u  D(\partial \BD, \SCL(0))})}, \Gamma \right)\right.\\
		&\qquad\qquad\left.\cdot\, g\left(\SCB_{\le u D(\partial \BD, \SCL(0))}, \overline{\SCB}_{\le \overline{D}(\SCL(0), \SCB_{u  D(\partial \BD, \SCL(0))})}, \Gamma \right)\right]\\
		&=\int_0^1 \mathrm{d}u \ \BE\left[ f(\SCB, \overline{\SCB}, \Gamma) g\left(\SCB_{\le u D(\partial \BD, \SCL(0))}, \overline{\SCB}_{\le \overline{D}(\SCL(0), \SCB_{u  D(\partial \BD, \SCL(0))})}, \Gamma \right)\right] \\
		&=\BE\left[f(\SCB, \overline{\SCB}, \Gamma) g\left(\SCB_{\le \tau}, \overline{\SCB}_{\le \overline{\tau}}, \Gamma \right) \right],
	\end{align*}
	where, in the first and last equalities, we used that $W$ is independent of $(\SCB, \overline{\SCB}, \Gamma)$, and, in the second equality, we used~\eqref{eq def h} together with the fact that $g\left(\SCB_{\le u D(\partial \BD, \SCL(0))}, \overline{\SCB}_{\le \overline{D}(\SCL(0), \SCB_{u  D(\partial \BD, \SCL(0))})}, \Gamma \right)$ is measurable with respect to the corresponding conditioning $\sigma$-algebra. This proves~\eqref{eq h equals h stopped at tau}.

	\stepn{step:cond-laws-nu-case}{The conditional law under $\nu$} Next, let $g$ be a non-negative measurable function. Then, by~\eqref{eq h equals h stopped at tau},
	\begin{align*}
		\BE&[D(\partial \BD, \SCL(0))] \cdot \BE_\nu \left[ \BE\left[\left. f( \SCB,  \overline{\SCB}, \Gamma) \right\vert \SCB_{\le \tau}, \overline{\SCB}_{\le \overline{\tau}}, \Gamma \right] g(\SCB_{\le \tau}, \overline{\SCB}_{\le \overline{\tau}}, \Gamma)\right]\\
		&= \BE\left[ D(\partial \BD, \SCL(0)) h(\SCB_{\le \tau}, \overline{\SCB}_{\le \overline{\tau}}, \Gamma) g(\SCB_{\le \tau}, \overline{\SCB}_{\le \overline{\tau}}, \Gamma) \right]\\
		&=
		\BE\left[ D(\partial \BD, \SCL(0)) \int_0^1 \mathrm{d} u \ h(\SCB_{\le u D(\partial \BD, \SCL(0))}, \overline{\SCB}_{\le \overline{D}(\SCL(0), \SCB_{u D(\partial \BD, \SCL(0))})}, \Gamma)\right.\\
		&\qquad\qquad\left.\cdot\, g(\SCB_{\le u D(\partial \BD, \SCL(0))}, \overline{\SCB}_{\le \overline{D}(\SCL(0), \SCB_{u D(\partial \BD, \SCL(0))})}, \Gamma) \right]
		\\
		&=\int_0^1 \mathrm{d}u \ \BE\left[ D(\partial \BD, \SCL(0))h(\SCB_{\le u D(\partial \BD, \SCL(0))}, \overline{\SCB}_{\le \overline{D}(\SCL(0), \SCB_{u D(\partial \BD, \SCL(0))})}, \Gamma)\right.\\
		&\qquad\qquad\left.\cdot\, g(\SCB_{\le u D(\partial \BD, \SCL(0))}, \overline{\SCB}_{\le \overline{D}(\SCL(0), \SCB_{u D(\partial \BD, \SCL(0))})}, \Gamma)  \right].
	\end{align*}
	Here, we used the definition of $\nu$ and~\eqref{eq h equals h stopped at tau} in the first equality, the independence of $W$ from $(\Gamma, D, \overline{D})$ in the second equality, and Fubini's theorem in the third equality.
	But notice that, for each fixed $u \in (0,1)$, the terminal time $u D(\partial \BD, \SCL(0))$ can be read off from the path $\SCB_{\le u D(\partial \BD, \SCL(0))}$, and we have $D(\partial \BD, \SCL(0)) = uD(\partial \BD, \SCL(0))/u$, so that $D(\partial \BD, \SCL(0))$ is a measurable function of $\SCB_{\le u D(\partial \BD, \SCL(0))}$. Thus, by definition of $h$ and~\eqref{eq def h},
	\begin{align*}
		\int_0^1 &\mathrm{d}u \ \BE\left[ D(\partial \BD, \SCL(0))h(\SCB_{\le u D(\partial \BD, \SCL(0))}, \overline{\SCB}_{\le \overline{D}(\SCL(0), \SCB_{u D(\partial \BD, \SCL(0))})}, \Gamma)\right.\\
		&\qquad\qquad\left.\cdot\, g(\SCB_{\le u D(\partial \BD, \SCL(0))}, \overline{\SCB}_{\le \overline{D}(\SCL(0), \SCB_{u D(\partial \BD, \SCL(0))})}, \Gamma)  \right]\\
		&= \int_0^1 \mathrm{d}u \ \BE\left[ D(\partial \BD, \SCL(0)) f(\SCB, \overline{\SCB}, \Gamma)\right.\\
		&\qquad\qquad\left.\cdot\, g(\SCB_{\le u D(\partial \BD, \SCL(0))}, \overline{\SCB}_{\le \overline{D}(\SCL(0), \SCB_{u D(\partial \BD, \SCL(0))})}, \Gamma) \right] \\
		&= \BE\left[  \int_0^1 \mathrm{d}u \ D(\partial \BD, \SCL(0)) f(\SCB, \overline{\SCB}, \Gamma)\right.\\
		&\qquad\qquad\left.\cdot\, g(\SCB_{\le u D(\partial \BD, \SCL(0))}, \overline{\SCB}_{\le \overline{D}(\SCL(0), \SCB_{u D(\partial \BD, \SCL(0))})}, \Gamma) \right]
		\\
		&= \BE\left[  \int_0^{D(\partial \BD, \SCL(0))} \mathrm{d}t \  f(\SCB, \overline{\SCB}, \Gamma) g(\SCB_{\le t}, \overline{\SCB}_{\le \overline{D}(\SCL(0), \SCB_{t})}, \Gamma) \right].
	\end{align*}
	Undoing the change of variables $t = u D(\partial \BD, \SCL(0))$ and using once more the independence of $W$, the last expression is equal to $\BE[D(\partial \BD, \SCL(0))] \cdot \BE_\nu[f(\SCB, \overline{\SCB}, \Gamma) g(\SCB_{\le \tau}, \overline{\SCB}_{\le \overline{\tau}}, \Gamma)]$.
	Therefore, dividing by $\BE[D(\partial \BD, \SCL(0))]$,
	\[ \BE_\nu \left[ \BE\left[\left. f( \SCB,  \overline{\SCB}, \Gamma) \right\vert \SCB_{\le \tau}, \overline{\SCB}_{\le \overline{\tau}}, \Gamma \right] g(\SCB_{\le \tau}, \overline{\SCB}_{\le \overline{\tau}}, \Gamma)\right]
	=
	\BE_\nu\left[ f(\SCB, \overline{\SCB}, \Gamma) g(\SCB_{\le \tau}, \overline{\SCB}_{\le \overline{\tau}}, \Gamma) \right].
	\]

	\stepn{step:cond-laws-symmetry}{Conclusion by symmetry} By exchanging the roles of $(\SCB, D, \partial \BD)$ and $(\overline{\SCB}, \overline{D}, \SCL(0))$ (which also exchanges $\nu$ with $\overline{\nu}$, $\tau$ with $\overline{\sigma}$, and $\overline{\tau}$ with $\sigma$), we obtain
	\[ \BE_{\overline{\nu}} \left[ \BE\left[\left. f( \SCB,  \overline{\SCB}, \Gamma) \right\vert \SCB_{\le \sigma}, \overline{\SCB}_{\le \overline{\sigma}}, \Gamma \right] g(\SCB_{\le \sigma}, \overline{\SCB}_{\le \overline{\sigma}}, \Gamma)\right]
	=
	\BE_{\overline{\nu}}\left[ f(\SCB, \overline{\SCB}, \Gamma) g(\SCB_{\le \sigma}, \overline{\SCB}_{\le \overline{\sigma}}, \Gamma) \right].
	\]
This completes the proof of the lemma since $f$ and $g$ were arbitrary.
\end{proof}

Combining \Cref{prop nu stopped equals nu tilde stopped} with \Cref{lemma conditional laws agree}, we deduce that $\nu$ and $\overline{\nu}$ induce the same law for the whole explorations.
\begin{corollary}\label{cor nu equals nu tilde}
	The $\nu$-law of $(\SCB, \overline{\SCB}, \Gamma)$ is equal to the $\overline{\nu}$-law of $(\SCB, \overline{\SCB}, \Gamma)$.
\end{corollary}
\begin{proof}
	Let $f$ be a non-negative measurable function. By \Cref{lemma conditional laws agree}, we have
	\begin{align*}
		\BE_\nu\left[ f(\SCB, \overline{\SCB}, \Gamma) \right] &= \BE_\nu\left[ \BE_\nu\left[\left.f(\SCB, \overline{\SCB}, \Gamma) \right\vert \SCB_{\le \tau}, \overline{\SCB}_{\le \overline{\tau}}, \Gamma \right] \right]\\
		&=\BE_\nu\left[ \BE\left[\left.f(\SCB, \overline{\SCB}, \Gamma) \right\vert \SCB_{\le \tau}, \overline{\SCB}_{\le \overline{\tau}}, \Gamma \right] \right]\\
		&=\BE_\nu\left[ h(\SCB_{\le \tau}, \overline{\SCB}_{\le \overline{\tau}}, \Gamma ) \right],
	\end{align*}
	by~\eqref{eq h equals h stopped at tau} where we recall that $h$ was defined in~\eqref{eq def h 0} (note that~\eqref{eq h equals h stopped at tau} holds $\BP$-a.s., and hence also $\nu$- and $\overline{\nu}$-a.s., since $\nu$ and $\overline{\nu}$ are absolutely continuous with respect to $\BP$).
	
	But, by \Cref{prop nu stopped equals nu tilde stopped}, we have
	\[
	\BE_\nu\left[ h(\SCB_{\le \tau}, \overline{\SCB}_{\le \overline{\tau}}, \Gamma ) \right] = \BE_{\overline{\nu}} \left[ h(\SCB_{\le \sigma}, \overline{\SCB}_{\le \overline{\sigma}}, \Gamma ) \right].
	\]
	To conclude, using \Cref{lemma conditional laws agree} once more, we see that it suffices to show that a.s.,
	\begin{equation*}
		h(\SCB_{\le \sigma}, \overline{\SCB}_{\le \overline{\sigma}}, \Gamma ) = \BE \left[ \left. f(\SCB, \overline{\SCB}, \Gamma)\right\vert \SCB_{\le \sigma}, \overline{\SCB}_{\le \overline{\sigma}}, \Gamma  \right].
	\end{equation*}
	This follows from exactly the same reasoning as the one used to obtain~\eqref{eq h equals h stopped at tau}, upon exchanging the roles of $(\SCB, D, \partial \BD)$ and $(\overline{\SCB}, \overline{D}, \SCL(0))$ as at the end of the proof of \Cref{lemma conditional laws agree}.  Indeed, granting this, \Cref{lemma conditional laws agree} gives
	\[
	\BE_{\overline{\nu}} \left[ h(\SCB_{\le \sigma}, \overline{\SCB}_{\le \overline{\sigma}}, \Gamma ) \right] = \BE_{\overline{\nu}}\left[ \BE_{\overline{\nu}}\left[ \left. f(\SCB, \overline{\SCB}, \Gamma) \right\vert \SCB_{\le \sigma}, \overline{\SCB}_{\le \overline{\sigma}}, \Gamma \right] \right] = \BE_{\overline{\nu}}\left[ f(\SCB, \overline{\SCB}, \Gamma) \right],
	\]
	which, combined with the two displays above, completes the proof.
\end{proof}
\begin{proof}[Proof of \Cref{prop D equals D tilde}]
	By \Cref{cor nu equals nu tilde}, we know that the $\nu$-law of $(\SCB, \overline{\SCB})$ is equal to the $\overline{\nu}$-law of $(\SCB, \overline{\SCB})$. But since $D(\partial \BD, \SCL(0))$ is a measurable function of $\SCB$ and since $\overline{D}(\partial \BD, \SCL(0))$ is a measurable function of $\overline{\SCB}$, we deduce that for every measurable non-negative function $f: \BR^2 \to \BR_{\ge 0}$, 
	\[
	\BE_\nu[f(D(\partial \BD, \SCL(0)), \overline{D}(\partial \BD, \SCL(0)))] = \BE_{\overline{\nu}}[f(D(\partial \BD, \SCL(0)), \overline{D}(\partial \BD, \SCL(0)))].
	\]
	Recalling the definitions of $\nu$ and $\overline{\nu}$ and that $\BE[D(\partial \BD, \SCL(0))] = \BE[\overline{D}(\partial \BD, \SCL(0))]$, this means that
	\[
	\BE\left[ D(\partial \BD, \SCL(0)) f(D(\partial \BD, \SCL(0)), \overline{D}(\partial \BD, \SCL(0)))\right] = \BE\left[ \overline{D}(\partial \BD, \SCL(0)) f(D(\partial \BD, \SCL(0)), \overline{D}(\partial \BD, \SCL(0)))\right]
	\]
	for every measurable $f \colon \BR^2 \to \BR_{\ge 0}$, i.e.\ the two measures on $\BR^2$ obtained by weighting the law of $(D(\partial \BD, \SCL(0)), \overline{D}(\partial \BD, \SCL(0)))$ by its first and by its second coordinate coincide. Taking $f(x,y) = \one_{x>y}$ gives $\BE[(D(\partial \BD, \SCL(0)) - \overline{D}(\partial \BD, \SCL(0)))\one_{D(\partial \BD, \SCL(0))>\overline{D}(\partial \BD, \SCL(0))}] = 0$ and hence $D(\partial \BD, \SCL(0)) \le \overline{D}(\partial \BD, \SCL(0))$ a.s.; exchanging the roles of $D$ and $\overline{D}$, we conclude that $D(\partial \BD, \SCL(0)) = \overline{D}(\partial \BD, \SCL(0))$ a.s.
\end{proof}

Finally, we are ready to prove Corollary~\ref{thm:uniform-exploration-measurability}.

\begin{proof}[Proof of Corollary~\ref{thm:uniform-exploration-measurability}]
We start by proving assertion~\eqref{thm:uniform-exploration-measurability-0}.  Let $(\Gamma_U ,  (t_{\ell})_{\ell \in \Gamma_U})$ denote a coupling between a non-nested $\CLE_4$ $\Gamma_U$ on $U$ and a collection of labels $(t_{\ell})_{\ell \in \Gamma_U}$ of the loops in $\Gamma_U$ as in Section~\ref{subsec:gff_labelled_cle_4}. Let also $D_{\Gamma_U}^U$ denote a strong $\CLE_4$ metric on $\Gamma_U$, which exists by Proposition~\ref{lem:weak-is-strong} together with the existence of a geodesic $\CLE_4$ metric coupling established in
\cite[Theorem~1.3]{kkmt2026cle4_part1} and \cite[Theorem~1.4]{kkmt2026cle4_part2}. Then,  Axiom~\eqref{it:axiom_uniform_exploration} of Definition~\ref{def:axioms} implies that there exists a coupling between $\Gamma_U,  D_{\Gamma_U}^U$,  and $(t_{\ell})_{\ell \in \Gamma_U}$ such that $t_{\SCL} =  D_{\Gamma_U}^U(\SCL ,  \partial U)$ for all $\SCL \in \Gamma_U$.  In particular,  Proposition~\ref{lem:measurability} (together with the fact that every strong $\CLE_4$ metric gives rise to a geodesic $\CLE_4$ metric coupling) implies that $(t_{\ell})_{\ell \in \Gamma_U}$ is a measurable function of $\Gamma_U$ in the above coupling.  It follows that there exists a unique collection of labels $(t_{\ell})_{\ell \in \Gamma_U}$ coupled in the same probability space with $\Gamma_U$ such that $(\Gamma_U,  (t_{\ell})_{\ell \in \Gamma_U})$ has the law of a labeled $\CLE_4$ and moreover,  $(t_{\ell})_{\ell \in \Gamma_U}$ is a.s.\ a measurable function of $\Gamma_U$.  Indeed, if $(s_{\ell})_{\ell \in \Gamma_U}$ is another such collection, then $(\Gamma_U, (s_\ell)_{\ell \in \Gamma_U})$ and $(\Gamma_U, (t_\ell)_{\ell \in \Gamma_U})$ have the same law and both collections of labels are measurable functions of $\Gamma_U$, so that the two functions agree almost everywhere with respect to the law of $\Gamma_U$.  Since the uniform exploration $(\SCB_t(\partial U))_{t \ge 0}$ is a measurable function of $(\Gamma_U, (t_\ell)_{\ell \in \Gamma_U})$ (see \Cref{subsec:labelled_cle_4}), it is in turn a.s.\ determined by $\Gamma_U$.  This completes the proof of assertion~\eqref{thm:uniform-exploration-measurability-0}.

Next, we focus on proving assertion~\eqref{thm:uniform-exploration-measurability-1}.  By conformal invariance, we may assume that $U=\BD$. We will use the same notation as in Section~\ref{subsec:gff_labelled_cle_4}. Let $\rho \in (-2, 0)$ and $r= -\rho \lambda$, so that $r \in (0,2\lambda)$. Let $\Gamma_\BD$ and the $\cwBCLE_4(\rho)$ be coupled as in Section~\ref{subsec:gff_labelled_cle_4}. By \Cref{lemma BCLE loops}, the domains encircled by the false (i.e.\@ counterclockwise) loops of the $\cwBCLE_4(\rho)$ are the connected components $O$ of $\BD \setminus \SCB_r(\partial \BD)$ such that $\partial O \cap \partial \BD \neq \emptyset$. Therefore, by assertion~\eqref{thm:uniform-exploration-measurability-0}, we see that the collection of the false loops of the $\cwBCLE_4(\rho)$ is a.s.\ determined by $\Gamma_\BD$. Since, by definition, every arc of the $\cwBCLE_4(\rho)$ separates a true loop from a false loop, the range of the $\cwBCLE_4(\rho)$ is the closure of the union of the boundaries of its false loops, and its true loops are the remaining connected components of the complement of this range.  We therefore conclude that the $\cwBCLE_4(\rho)$ is a.s.\ determined by $\Gamma_\BD$.
\end{proof}

\section{Uniqueness of geodesics}
\label{section:uniqueness-of-geodesics}

In the present section, we prove that there is a.s.\ a unique geodesic connecting any pair of loops. More precisely, we prove \Cref{thm:uniqueness-of-geodesics}. Let $D$ be the (strong) $\CLE_4$ metric (which exists and is unique by \Cref{thm:uniqueness}). By conformal invariance (Axiom~\eqref{it:axiom_conformal_invariance} of \Cref{def:axioms}), together with the fact that a conformal transformation maps admissible paths to admissible paths and hence $D$-geodesics to $D$-geodesics, it is enough to prove the following reformulation of \Cref{thm:uniqueness-of-geodesics}.

\begin{theorem}[Uniqueness of geodesics]\label{thm:uniqueness-geodesics}
Almost surely, for each $\SCL_1, \SCL_2 \in \Gamma_\BH$, there is a unique $D_{\Gamma_\BH}^\BH$-geodesic connecting $\SCL_1$ and $\SCL_2$, in the sense that if $P_1$ and $P_2$ are $D_{\Gamma_\BH}^\BH$-geodesics connecting $\SCL_1$ and $\SCL_2$, then $\{\SCL \in \Gamma_\BH : \SCL \cap P_1 \neq \emptyset\} = \{\SCL \in \Gamma_\BH : \SCL \cap P_2 \neq \emptyset\}$. A similar statement holds with a deterministic connected arc of $\partial\BH$ (resp.\ two deterministic disjoint connected arcs of $\partial\BH$) in place of $\SCL_1$ (resp.\ $\SCL_1$ and $\SCL_2$).
\end{theorem}

The following monotonicity property will be useful.

\begin{lemma}\label{lem:monotonicity}
Let $\ell$ be a random loop in $U$,  independent of $\Gamma_U$,  whose law is absolutely continuous with respect to the Brownian loop measure in $U$. Write $\widetilde\SCL$ for the outer boundary of the cluster obtained by merging $\ell$ together with all of the loops of $\Gamma_U$ which touch $\ell$, and, on the event that $\ell$ is not contained in the interior of a loop of $\Gamma_U$, let
\begin{equation*}
	\widetilde\Gamma_U \defeq \left\{\SCL \in \Gamma_U : \SCL \not\subseteq \overline{\mathop{\mathrm{int}}(\widetilde\SCL)}\right\} \cup \left\{\widetilde\SCL\right\};
\end{equation*}
on the complementary event (on which the merged cluster is not outermost), set $\widetilde\Gamma_U \defeq \Gamma_U$, in which case the statements below are elementary.
Recall that the law of $\widetilde\Gamma_U$ is absolutely continuous with respect to the law of $\Gamma_U$, so that $D_{\widetilde\Gamma_U}^U$ is well defined and almost sure statements about $D_{\Gamma_U}^U$ transfer to $D_{\widetilde\Gamma_U}^U$. Here and below, $D_{\Gamma_U}^U(\SCL, \widetilde\SCL)$ denotes the distance from the loop $\SCL \in \Gamma_U$ to the \emph{subset} $\widetilde\SCL$ of $\overline{U}$, in the sense of \Cref{subsec:setup}. Then, a.s.,
\begin{equation*}
	\begin{cases}
		D_{\widetilde\Gamma_U}^U(\SCL_1, \SCL_2) \le D_{\Gamma_U}^U(\SCL_1, \SCL_2) \wedge \left(D_{\Gamma_U}^U(\SCL_1, \widetilde\SCL) + D_{\Gamma_U}^U(\SCL_2, \widetilde\SCL)\right) & \text{for } \SCL_1, \SCL_2 \in \widetilde\Gamma_U \setminus \{\widetilde\SCL\}; \\
		D_{\widetilde\Gamma_U}^U(\SCL, \widetilde\SCL) \le D_{\Gamma_U}^U(\SCL, \widetilde\SCL) & \text{for } \SCL \in \widetilde\Gamma_U \setminus \{\widetilde\SCL\}.
	\end{cases}
\end{equation*}
In particular, a.s., $D_{\widetilde\Gamma_U}^U(\widetilde\SCL(x), \widetilde\SCL(y)) \le D_{\Gamma_U}^U(\SCL(x), \SCL(y))$ for all $x, y \in U_\BQ$ (where $\widetilde\SCL(x) = \SCL(x)$ if $x$ is not surrounded by $\widetilde\SCL$ and $\widetilde\SCL(x) = \widetilde\SCL$ otherwise).
\end{lemma}

\begin{proof}
Throughout the proof we write $\SCO \defeq \mathop{\mathrm{int}}(\widetilde\SCL)$ and work on the event that $\ell$ is not contained in the interior of a loop of $\Gamma_U$.

\stepn{step:mono-two-ensembles-have}{The two ensembles have the same loops outside $\overline\SCO$}  The loops of $\Gamma_U$ which touch $\ell$ are contained in $\overline\SCO$, and $\widetilde\SCL$ is contained in the closure of the union of $\ell$ and of these loops (whose accumulation points lie on $\ell$). If $\SCL \in \Gamma_U$ is such that $\SCL \not\subseteq \overline{\SCO} $, then, by the loop soup construction of the $\CLE_4$, the cluster of Brownian loops whose outerboundary is $\SCL$ is disjoint from $\overline{\SCO}$. Therefore, $\SCL\cap \overline{\SCO}= \emptyset $. 
Hence
\begin{equation}\label{eq:monotonicity_same_loops}
	\widetilde\Gamma_U \setminus \{\widetilde\SCL\} = \left\{\SCL \in \Gamma_U : \SCL \cap \overline\SCO = \emptyset\right\},
\end{equation}
and the hulls of all of the other loops of $\Gamma_U$, as well as the hull $\overline\SCO$ of $\widetilde\SCL$, are contained in $\overline\SCO$.  Note also that an admissible path for $\widetilde\Gamma_U$ takes its values in the carpet of $\widetilde\Gamma_U$ and therefore never enters $\SCO$.

\stepn{step:mono-two-metrics-agree}{The two metrics agree locally away from $\overline\SCO$}  Let $\SCW$ be the countable collection of open Euclidean balls with rational center and rational radius whose closure is contained in $U$.  We claim that a.s., for every $W \in \SCW$ with $\overline W \cap \overline\SCO = \emptyset$ and every connected component $W_j$ of $W^\star$,
\begin{equation}\label{eq:monotonicity_locality}
	D_{\Gamma_U}^U(\SCL_1, \SCL_2; W_j) = D_{\widetilde\Gamma_U}^U(\SCL_1, \SCL_2; W_j), \qquad \forall \SCL_1, \SCL_2 \in \Gamma_U|_{W_j}.
\end{equation}
Indeed, on the event $\overline W \cap \overline\SCO = \emptyset$, by \Cref{step:mono-two-ensembles-have} the set $W^\star$ and the collection of loops contained in it are the same for the two ensembles, and Axiom~\eqref{it:axiom_locality}, applied to the deterministic simply connected domain $W$ and to each of the two ensembles (the law of $\widetilde\Gamma_U$ being absolutely continuous with respect to that of $\Gamma_U$), identifies both sides of~\eqref{eq:monotonicity_locality} with the same deterministic functional $D_{\bullet|_{W_j}}^{W_j}$ of $\Gamma_U|_{W_j} = \widetilde\Gamma_U|_{W_j}$.  As $\SCW$ is countable, the exceptional set may be chosen not to depend on $W$.

\stepn{step:mono-paths-avoiding}{Paths avoiding $\overline\SCO$}  We claim that a.s.\ the following holds: if $P$ is an admissible path with $P \cap \overline\SCO = \emptyset$ and $\SCL^\prime, \SCL^{\prime\prime}$ are loops hit by $P$, then
\begin{equation}\label{eq:monotonicity_length}
	D_{\widetilde\Gamma_U}^U(\SCL^\prime, \SCL^{\prime\prime}) \le \len\bigl(P; D_{\Gamma_U}^U\bigr).
\end{equation}
Indeed, let $\delta\in (0,  \dist(P, \overline\SCO)) \cap \BQ$. We may cover $P$ by finitely many balls $B_{\delta/8}(q_i)$ with $q_i \in \BQ^2$ and $W_i \defeq B_{\delta/4}(q_i) \in \SCW$, so that $\overline{W_i} \cap \overline\SCO = \emptyset$.  Since $\Gamma_U$ is locally finite, only finitely many loops of diameter at least $\delta/8$ meet $P$, and $P$ meets each of them along a segment of $P$ (cf.\ \Cref{lem:intersection_is_a_segment}).  We may therefore subdivide $P$ into finitely many consecutive pieces $P_1, \ldots, P_n$ whose endpoints lie on loops of $\Gamma_U$ (the points of $P$ on loops being dense in $P$), such that each $P_k$ is contained in some $B_{\delta/8}(q_i)$ and either (i) is contained in a single loop of $\Gamma_U$, or (ii) meets no loop of diameter at least $\delta/8$.  Write $\SCL^{(k)}$ for the loop at the common endpoint of $P_k$ and $P_{k+1}$, with $\SCL^{(0)} = \SCL^\prime$ and $\SCL^{(n)} = \SCL^{\prime\prime}$; by~\eqref{eq:monotonicity_same_loops}, each $\SCL^{(k)}$ is a loop of $\widetilde\Gamma_U$.  In case~(i) we have $\SCL^{(k-1)} = \SCL^{(k)}$.  In case~(ii), any loop whose hull meets $P_k$ meets $P_k$ itself (an admissible path does not enter the interior of a loop), hence has diameter less than $\delta/8$ and is contained in $W_i$; thus $P_k$ lies in a single connected component $W_{i,j}$ of $W_i^\star$, and
\begin{equation*}
	\len\bigl(P_k; D_{\Gamma_U}^U\bigr) \ge D_{\Gamma_U}^U\bigl(\SCL^{(k-1)}, \SCL^{(k)}; W_{i,j}\bigr) = D_{\widetilde\Gamma_U}^U\bigl(\SCL^{(k-1)}, \SCL^{(k)}; W_{i,j}\bigr) \ge D_{\widetilde\Gamma_U}^U\bigl(\SCL^{(k-1)}, \SCL^{(k)}\bigr),
\end{equation*}
where the equality is~\eqref{eq:monotonicity_locality} and the last inequality holds because the chain consisting of the two endpoint loops is admissible in the definition of the length of a path.  Summing over $k$, and using the additivity of lengths together with the triangle inequality for $D_{\widetilde\Gamma_U}^U$, proves~\eqref{eq:monotonicity_length}.

\stepn{step:mono-paths-ending}{Paths ending on $\overline\SCO$}  We claim that a.s.\ the following holds: if $P \colon [0,1] \to \Upsilon_U$ is an admissible path with $P([0,1)) \cap \overline\SCO = \emptyset$ and $P(1) \in \overline\SCO$, and if $\SCL^\prime$ is a loop hit by $P$, then
\begin{equation}\label{eq:monotonicity_endpoint}
	D_{\widetilde\Gamma_U}^U\bigl(\SCL^\prime, \widetilde\SCL\bigr) \le \len\bigl(P; D_{\Gamma_U}^U\bigr).
\end{equation}
Indeed, write $t \defeq \len(P; D_{\Gamma_U}^U)$.  For each $s < 1$, \Cref{step:mono-paths-avoiding} applies to $P|_{[0,s]}$ and shows that every loop hit by $P|_{[0,s]}$ lies at $D_{\widetilde\Gamma_U}^U$-distance at most $t$ from $\SCL^\prime$.  Hence the closed set $\SCB_t(\SCL^\prime; D_{\widetilde\Gamma_U}^U)$ contains the closure of the union of the loops hit by $P$, and therefore $P(1)$ (the points of $P$ on loops being dense in $P$).  Since $\SCB_0(\widetilde\SCL; D_{\widetilde\Gamma_U}^U) = \overline{\SCO} \ni P(1)$, Axiom~\eqref{it:axiom_geodesic} gives the claim.

\stepn{step:mono-conclusion}{Conclusion}  Let $\SCL_1, \SCL_2 \in \widetilde\Gamma_U \setminus \{\widetilde\SCL\}$ and let $P$ be a $D_{\Gamma_U}^U$-geodesic from $\SCL_1$ to $\SCL_2$ (Axiom~\eqref{it:axiom_geodesic}).  If $P \cap \overline\SCO = \emptyset$, then~\eqref{eq:monotonicity_length} gives $D_{\widetilde\Gamma_U}^U(\SCL_1, \SCL_2) \le \len(P; D_{\Gamma_U}^U) = D_{\Gamma_U}^U(\SCL_1, \SCL_2)$; otherwise, applying~\eqref{eq:monotonicity_endpoint} to $P$ and to its time-reversal, each stopped at its first hitting time of $\overline\SCO$, and using the triangle inequality for $D_{\widetilde\Gamma_U}^U$ through $\widetilde\SCL$, gives the same bound.  This proves the first term of the first inequality of the lemma.

Next, let $\SCL \in \widetilde\Gamma_U \setminus \{\widetilde\SCL\}$ and $t > D_{\Gamma_U}^U(\SCL, \widetilde\SCL)$, so that $\SCB_t(\SCL; D_{\Gamma_U}^U)$ meets $\widetilde\SCL$ and there is a loop $\SCL^\ast \in \Gamma_U$ with $D_{\Gamma_U}^U(\SCL, \SCL^\ast) \le t$ whose hull meets $\widetilde\SCL$ (see \cite[Lemma~4.7]{kkmt2026cle4_part2}).  Then $\SCL^\ast \subseteq \overline\SCO$: otherwise, by \Cref{step:mono-two-ensembles-have}, $\SCL^\ast$ would be disjoint from $\overline\SCO \supseteq \widetilde\SCL$, and its hull could meet $\widetilde\SCL$ only if $\widetilde\SCL \subseteq \mathop{\mathrm{int}}(\SCL^\ast)$, which would put $\ell$ in the interior of a loop of $\Gamma_U$.  Hence a $D_{\Gamma_U}^U$-geodesic from $\SCL$ to $\SCL^\ast$ meets $\overline\SCO$, and~\eqref{eq:monotonicity_endpoint}, applied to it stopped at its first hitting time of $\overline\SCO$, gives $D_{\widetilde\Gamma_U}^U(\SCL, \widetilde\SCL) \le D_{\Gamma_U}^U(\SCL, \SCL^\ast) \le t$.  Letting $t \downarrow D_{\Gamma_U}^U(\SCL, \widetilde\SCL)$ proves the second inequality of the lemma, and the second term of the first inequality then follows from the triangle inequality for $D_{\widetilde\Gamma_U}^U$.

Finally, let $x, y \in U_\BQ$.  If $x$ is not surrounded by $\widetilde\SCL$, then $\SCL(x) \not\subseteq \overline\SCO$ (else $x \in \mathop{\mathrm{int}}(\SCL(x)) \subseteq \SCO$), so that $\SCL(x) = \widetilde\SCL(x)$ is a loop of $\widetilde\Gamma_U$; and if $x$ is surrounded by $\widetilde\SCL$, then $\SCL(x) \subseteq \overline\SCO$, since otherwise $\SCL(x) \cap \overline\SCO = \emptyset$ by \Cref{step:mono-two-ensembles-have} and the connected set $\SCO \ni x$ would be contained in $\mathop{\mathrm{int}}(\SCL(x))$, again putting $\ell$ in the interior of a loop.  If neither $x$ nor $y$ is surrounded by $\widetilde\SCL$, the final assertion of the lemma is thus contained in the first inequality; if both are, its left-hand side vanishes; and if exactly one of them, say $x$, is surrounded, then a $D_{\Gamma_U}^U$-geodesic from $\SCL(y)$ to $\SCL(x) \subseteq \overline\SCO$ meets $\overline\SCO$, and~\eqref{eq:monotonicity_endpoint} gives $D_{\widetilde\Gamma_U}^U(\widetilde\SCL(x), \widetilde\SCL(y)) = D_{\widetilde\Gamma_U}^U(\SCL(y), \widetilde\SCL) \le D_{\Gamma_U}^U(\SCL(x), \SCL(y))$.
\end{proof}

Let us recall the FKG inequality for Poisson point processes. Let $(X, \SF)$ be a measurable space. Let $\SCM_\mathrm{pp}(X, \SF)$ be the space of all countable sums of Dirac measures on $X$, equipped with the smallest $\sigma$-algebra which makes the assignments $\mu \mapsto \mu(A)\in \BZ_{\ge 0} \cup \{\infty\}$ measurable for all $A \in \SF$. We shall refer to a measurable mapping $F \colon \SCM_\mathrm{pp}(X, \SF)\to \BR$ as increasing if for every pair $\mu, \nu \in \SCM_\mathrm{pp}(X, \SF)$ such that $\mu(A) \le \nu(A)$ for all $A \in \SF$, we have $F(\mu) \le F(\nu)$. We shall call $F$ decreasing if $-F$ is increasing; note that \Cref{lem:FKG_PPP} below applies verbatim to a pair of decreasing mappings, as one sees by applying it to $-F$ and $-G$.
\begin{lemma}\label{lem:FKG_PPP}
Let $\Xi$ be a Poisson point process with intensity measure given by a $\sigma$-finite measure on $(X, \SF)$. Let $F, G \colon \SCM_\mathrm{pp}(X, \SF) \to \BR$ be bounded, measurable, increasing mappings. Then 
\begin{equation*}
	\BE\lbrack F(\Xi)G(\Xi)\rbrack \ge \BE\lbrack F(\Xi)\rbrack\BE\lbrack G(\Xi)\rbrack. 
\end{equation*}
The same is true with two non-negative, measurable, increasing mappings in place of $F$ and $G$. 
\end{lemma}
\begin{proof}
See, e.g., \cite[Lemma~2.1]{Jan84}. 
\end{proof}
The two lemmas above imply the following FKG inequality for the strong $\CLE_4$ metric $D$.
\begin{lemma}\label{lem:FKG-CLE}
    Let $U \subsetneq \BC$ be a deterministic simply connected domain and let $\Gamma_U$ be a non-nested $\CLE_4$ in $U$. Let $F, G \colon \BR^{U_\BQ \times U_\BQ} \to \BR_{\ge 0}$ be non-negative measurable mappings. Suppose that $F(d) \ge F(d^\prime)$ and $G(d) \ge G(d^\prime)$ for all $d, d^\prime \in \BR^{U_\BQ \times U_\BQ}$ with $d(x, y) \ge d^\prime(x, y)$ for all $x, y \in U_\BQ$. Then $\BE\lbrack F(D_{\Gamma_U}^U) G(D_{\Gamma_U}^U)\rbrack \ge \BE\lbrack F(D_{\Gamma_U}^U)\rbrack \BE\lbrack G(D_{\Gamma_U}^U)\rbrack$. The same is true for any pair of non-negative measurable functionals of $(\Gamma_U, D_{\Gamma_U}^U)$ which are a.s.\ non-increasing under the merging operation of \Cref{lem:monotonicity}, and in particular for the indicator functions of the events $\{D_{\Gamma_U}^U(A, B) > \varepsilon\}$, where $A, B \subseteq \overline U$ are deterministic and $\varepsilon > 0$.
\end{lemma}

\begin{proof}
    This follows immediately from \Cref{lem:monotonicity,lem:FKG_PPP} combined with the fact that $D_{\Gamma_U}^U$ is a.s.\ determined by the critical Brownian loop soup used to construct $\Gamma_U$ since $D_{\Gamma_U}^U$ is a.s.\ determined by $\Gamma_U$ (Proposition~\ref{lem:measurability}).
\end{proof}

Next, let us state the primary ingredient in the proof of \Cref{thm:uniqueness-geodesics}.  Its own proof will be given at the end of the present section.

\begin{proposition}\label{lem:law_cont}
Let $r > 0$. Let $V_r =(0,r) \times (0,1)$. Write $T_r$ (resp.\ $B_r$) for the top (resp.\ bottom) side of $V_r$. Then the law of the random variable $D_{\Gamma_{V_r}}^{V_r}(T_r, B_r)$ is absolutely continuous with respect to the Lebesgue measure. 
\end{proposition}

The proof of \Cref{thm:uniqueness-geodesics} is reminiscent of the proof of geodesic uniqueness for LQG metrics (cf.\ \cite{GeoLQGnSLE,IntroLQGMet}). The idea is as follows. Suppose, for the sake of contradiction, that there are two distinct geodesics connecting a given pair of loops. Then, there exists a dyadic domain such that one of the geodesics passes through it while the other does not. Moreover, the former geodesic enters and exits this dyadic domain through two specific loops. There are only countably many choices for such a dyadic domain and its corresponding entrance and exit loops. However, if we condition on everything not fully contained within the dyadic domain, the length of the latter geodesic becomes deterministic. In contrast, the length of the former geodesic depends on the internal distance between the entrance and exit loops inside the dyadic domain, which, by \Cref{lem:law_cont}, is random. Thus, the conditional probability that the two geodesics have exactly the same length is zero, yielding a contradiction.

\begin{proof}[Proof of \Cref{thm:uniqueness-geodesics} assuming \Cref{lem:law_cont}]
\stepn{step:uniq-geo-dyadic-domains-excursion}{Dyadic domains, excursion arcs, and the $\sigma$-algebra $\SCG_W$}  We shall refer to as a \emph{dyadic square} a square with side length $2^{-n}$ and vertices in $(2^{-n}\BZ)^2$ for some $n \in \BZ$. We shall refer to a \emph{simply connected} domain $W \subseteq \BC$ as \emph{dyadic} if there exists a finite collection of dyadic squares such that $W$ is the interior of the union of their closures. Note that the set of dyadic domains is countable. (Simple connectivity is imposed so that the connected components of $W^\star$ are simply connected, so that $\Gamma_{O}$ and $D^{O}$ are defined for such a component $O$.)

Let us first fix some notation.  Let $W$ be a dyadic domain with $\overline{W} \subseteq \BH$ and $z_1, z_2 \notin \overline W$, and let $O$ be a connected component of $W^\star$, so that $O$ is simply connected.  If $\SCL \in \Gamma_\BH$ is a loop which is not contained in $W$, then each connected component of $\SCL \cap W$ is an open sub-arc of $\SCL$, and we call the closure of such a component an \emph{excursion arc} provided it is contained in $\partial O$.  There are countably many excursion arcs, and, together with $W^\star$ and with the collection of its connected components, they are determined by $\Gamma_\BH \setminus \Gamma_\BH|_{W^\star}$.  Finally, let $\SCG_W$ denote the $\sigma$-algebra generated by~\eqref{eq:weak_axiom_locality} with $V = W$, i.e.\ by the loops of $\Gamma_\BH$ which are not contained in $W$ together with the metric balls started from those loops and from $\partial\BH$ and stopped upon hitting $\partial W^\star$.  By \Cref{def:weak_axioms}, Axiom~\eqref{it:weak_axiom_locality} (locality), conditionally on $\SCG_W$ the pairs $(\Gamma_\BH|_{O}, D_{\Gamma_\BH}^\BH(\bullet, \bullet; O))$, as $O$ ranges over the connected components of $W^\star$, are independent, and the conditional law of the one indexed by $O$ is that of $(\Gamma_{O}, D_{\Gamma_{O}}^{O})$.

We shall use the following elementary observation: any two points of a sub-arc $\alpha$ of a loop $\SCL \in \Gamma_\BH$ are joined by an admissible path of zero $D_{\Gamma_\BH}^\BH$-length (the corresponding sub-arc of $\SCL$), so that splicing into an admissible path a piece which runs along $\alpha$ does not change its $D_{\Gamma_\BH}^\BH$-length.

\stepn{step:uniq-geo-reduction-countably-many}{Reduction to countably many events}  It suffices to show that for each deterministic $z_1, z_2 \in \BH_\BQ$, there is a.s.\ a unique $D_{\Gamma_\BH}^\BH$-geodesic connecting $\SCL(z_1)$ and $\SCL(z_2)$.  Indeed, a.s.\ every loop of $\Gamma_\BH$ is of the form $\SCL(z)$ for some $z \in \BH_\BQ$, and $\BH_\BQ$ is countable.  Given $W$ and $O$ as above and two disjoint excursion arcs $\alpha_1, \alpha_2$ of $\partial O$, we let $E(W, O, \alpha_1, \alpha_2)$ be the event that there exist two $D_{\Gamma_\BH}^\BH$-geodesics $P_1$ and $P_2$ connecting $\SCL(z_1)$ and $\SCL(z_2)$ such that
\begin{itemize}
\item $P_2 \cap \overline{O} = \emptyset$;
\item writing $\sigma_1 \defeq \inf\{t : P_1(t) \in \overline O\}$ and $\sigma_2 \defeq \sup\{t : P_1(t) \in \overline O\}$, we have $\sigma_1 < \sigma_2$, $P_1([\sigma_1, \sigma_2]) \subseteq \overline O$, $P_1(\sigma_1) \in \alpha_1$ and $P_1(\sigma_2) \in \alpha_2$.
\end{itemize}
On the event that there are two $D_{\Gamma_\BH}^\BH$-geodesics connecting $\SCL(z_1)$ and $\SCL(z_2)$ which do not visit the same set of loops, there is a.s.\ such a quadruple $(W, O, \alpha_1, \alpha_2)$ (we use the fact that a metric ball hits a deterministic connected compact set by a loop, as shown in \cite[Lemma 4.7]{kkmt2026cle4_part2}). Since there are countably many dyadic domains $W$ and, given $W$, countably many components $O$ and countably many pairs of excursion arcs, and since all of these are $\SCG_W$-measurable, it suffices to show that for each dyadic domain $W$ with $\overline W \subseteq \BH$ and $z_1, z_2 \notin \overline W$ we have, a.s., simultaneously for every component $O$ of $W^\star$ and every pair of disjoint excursion arcs $\alpha_1, \alpha_2$ of $\partial O$,
\begin{equation*}
\BP\lbrack E(W, O, \alpha_1, \alpha_2) \mid \SCG_W\rbrack = 0 \quad \text{a.s.}
\end{equation*}

\stepn{step:uniq-geo-internal-distance-across}{The internal distance across $O$}  Fix such $W, O, \alpha_1, \alpha_2$ and write
\begin{equation*}
s \defeq D_{\Gamma_\BH}^\BH(\SCL(z_1), \SCL(z_2); \BH \setminus O) - D_{\Gamma_\BH}^\BH(\SCL(z_1), \alpha_1; \BH \setminus O) - D_{\Gamma_\BH}^\BH(\alpha_2, \SCL(z_2); \BH \setminus O).
\end{equation*}
We claim that, on $E(W, O, \alpha_1, \alpha_2)$, we have $D_{\Gamma_\BH}^\BH(\alpha_1, \alpha_2; \overline O) = s$.  Indeed, $P_2$ is an admissible path in $\BH \setminus O$ from $\SCL(z_1)$ to $\SCL(z_2)$, so that $\len(P_2; D_{\Gamma_\BH}^\BH) \ge D_{\Gamma_\BH}^\BH(\SCL(z_1), \SCL(z_2); \BH\setminus O)$, with equality since $P_2$ is a $D_{\Gamma_\BH}^\BH$-geodesic.  Next, $P_1|_{[0,\sigma_1]}$, $P_1|_{[\sigma_1,\sigma_2]}$, and $P_1|_{[\sigma_2,1]}$ are admissible paths in $\BH \setminus O$, in $\overline O$, and in $\BH \setminus O$, respectively, joining the indicated sets, so that $D_{\Gamma_\BH}^\BH(\SCL(z_1),\alpha_1;\BH\setminus O) + D_{\Gamma_\BH}^\BH(\alpha_1,\alpha_2;\overline O) + D_{\Gamma_\BH}^\BH(\alpha_2,\SCL(z_2);\BH\setminus O) \le \len(P_1; D_{\Gamma_\BH}^\BH) = \len(P_2; D_{\Gamma_\BH}^\BH)$.  Conversely, concatenating admissible paths which realize these three internal distances up to $\delta > 0$ (after splicing in sub-arcs of $\alpha_1$ and of $\alpha_2$, which does not change the length, by the observation of \Cref{step:uniq-geo-dyadic-domains-excursion}) yields an admissible path from $\SCL(z_1)$ to $\SCL(z_2)$, so that the same sum is at least $D_{\Gamma_\BH}^\BH(\SCL(z_1), \SCL(z_2)) = \len(P_1; D_{\Gamma_\BH}^\BH)$.  The claim follows.

\stepn{step:uniq-geo-conclusion-via-lem}{Conclusion via \Cref{lem:law_cont}}  Moreover, $s$ is $\SCG_W$-measurable. Since $\alpha_1$ and $\alpha_2$ are sub-arcs of loops which are not contained in $W$, and since $(O; \alpha_1, \alpha_2)$ is thus a conformal rectangle whose two remaining sides are the two arcs of $\partial O$ between $\alpha_1$ and $\alpha_2$, we may let $\rho > 0$ denote its conformal modulus and let $\psi \colon O \to V_\rho$ be a conformal transformation taking $\alpha_1$ to $T_\rho$ and $\alpha_2$ to $B_\rho$; both $\rho$ and $\psi$ are $\SCG_W$-measurable.  Combining the above with Axiom~\eqref{it:axiom_locality} (locality) and Axiom~\eqref{it:axiom_conformal_invariance} (conformal invariance), we obtain
\begin{align*}
\BP\lbrack E(W, O, \alpha_1, \alpha_2) \mid \SCG_W\rbrack &\le \BP\!\left\lbrack D_{\Gamma_\BH}^\BH(\alpha_1, \alpha_2; \overline O) = s \ \middle\vert \ \SCG_W\right\rbrack \\
&= \BP\!\left\lbrack D_{\Gamma_{V_\rho}}^{V_\rho}(T_\rho, B_\rho) = s \ \middle\vert \ \SCG_W\right\rbrack \\
&= 0 \quad \text{a.s.},
\end{align*}
the last equality being \Cref{lem:law_cont} applied with the ($\SCG_W$-measurable, hence conditionally deterministic) parameter $\rho$, together with the fact that $s$ is $\SCG_W$-measurable.
The same reasoning works for geodesics between a loop and a deterministic connected arc of $\partial \BH$ or between two deterministic disjoint connected arcs of $\partial \BH$. This completes the proof of \Cref{thm:uniqueness-geodesics}.
\end{proof}

It remains to prove \Cref{lem:law_cont}. The proof strategy proceeds in stages. First, we state a lemma (\Cref{lem:dist_across_rect_arb_high}) which gives that the distance between the top and bottom sides of the rectangle $V_r$ can be arbitrarily large with positive probability. This lemma will be a consequence of the FKG inequality (cf.\ \Cref{lem:FKG-CLE}). Next, we establish a weaker version of \Cref{lem:law_cont} in \Cref{lem:law_cont_pre}. Using the notation of \Cref{lem:law_cont}, this preliminary lemma states that if $r$ is a random variable which is independent of the $\CLE_4$ and whose law is absolutely continuous with respect to Lebesgue measure, then the distance $D_{\Gamma_{V_r}}^{V_r}(T_r, B_r)$ satisfies the desired property.

The proof of \Cref{lem:law_cont_pre} proceeds as follows. We recall that in a non-nested $\CLE_4$ on the unit disk, the distance from the origin-containing loop to the domain boundary is an exponentially distributed random variable (this is immediate from Axiom~\eqref{it:axiom_uniform_exploration} together with the conformal invariance of the metric, which makes the exploration seen from the origin memoryless). Consequently, removing the origin-containing loop from the unit disk yields a random annulus whose cross-distance is absolutely continuous with respect to Lebesgue measure. We then transition from this random annulus to a random rectangle by conditioning on all loops that intersect a fixed line segment crossing the annulus.

Finally, we deduce \Cref{lem:law_cont} from \Cref{lem:law_cont_pre} via a metric ball exploration. We initially run the metric ball exploration from the bottom boundary for a short duration. Applying \Cref{lem:sle4_cm}, we can deduce that the unexplored region forms a random rectangle whose conformal modulus is absolutely continuous with respect to Lebesgue measure. This observation immediately reduces the proof of \Cref{lem:law_cont} to that of \Cref{lem:law_cont_pre}, the point being that, by the Markov property of the exploration, the conformal modulus of the unexplored region is independent of the $\CLE_4$ inside it, so that \Cref{lem:law_cont_pre} may be applied conditionally on that modulus.

Let us now turn to the first of these steps. The following result states that the distance across a rectangle can be arbitrarily large with positive probability. This is a direct consequence of the FKG inequality (cf.\ \Cref{lem:FKG-CLE}).

\begin{lemma}\label{lem:dist_across_rect_arb_high}
Let $r, N > 0$. Then 
\begin{equation*}
\BP\!\left\lbrack D_{\Gamma_{V_r}}^{V_r}(T_r, B_r) > N\right\rbrack > 0. 
\end{equation*}
\end{lemma}
\begin{proof}

\stepn{step:dist-rect-notation-two-deterministic}{Notation, and two deterministic bounds on the modulus}  To lighten notation, write $V \defeq V_r$, $T \defeq T_r$, and $B \defeq B_r$.  For $0 < s < t $ and $z \in \BC$, we write $A_{s,t}(z) \defeq B_t(z) \setminus \overline{B_s(z)}$ for the corresponding open Euclidean annulus.  We write $Q_{t/2,t}(z)$ for the connected component of $\left(A_{t/2,t}(z) \cap V\right)^\star$ which separates $\partial B_{t/2}(z) \cap V$ from $\partial B_t(z) \cap V$, when such a component exists, and $M_{t/2,t}(z)$ for the extremal distance in $Q_{t/2,t}(z)$ between the part of its boundary coming from $\partial B_{t/2}(z)$ and the part coming from $\partial B_t(z)$.  Note that for all $t \in (0, 1/2)$, the domain $Q_{t/2,t}(z)$ is a conformal annulus if $\overline{B_{t/2}(z)} \subseteq V$ and a conformal rectangle (whose two remaining sides lie on $\partial V$) otherwise, and that in both cases
	\begin{equation}\label{eq:modulus_lower_bound}
		M_{t/2,t}(z) \ge \frac{\log 2}{2\pi},
	\end{equation}
	since $Q_{t/2,t}(z) \subseteq A_{t/2,t}(z)$ and since extremal distance is monotone under inclusion.  In the opposite direction, since the convex rectangle $V$ has interior angles $\pi/2$, the set $A_{2t/3/2,5t/6}(z) \cap V$ contains a sector of $A_{2t/3,5t/6}(z)$ of aperture at least $\pi/2$, so that the extremal distance between $\partial B_{2t/3}(z) \cap V$ and $\partial B_{5t/6}(z) \cap V$ in $A_{2t/3,5t/6}(z) \cap V$ itself is at most a universal constant.

\stepn{step:dist-rect-event-produce}{The event which we shall produce}  We shall show in \Crefrange{step:dist-rect-decoupling-loop-soup}{step:dist-rect-proof-eq-good} that for all $z \in V$ and all $t>0$, with non-zero probability,
	\begin{align}\label{eq:good_annulus_event}
		&\text{no } \SCL \in \Gamma_V \text{ is such that } \overline{\mathop{\mathrm{int}}(\SCL)} \text{ meets both } \partial B_{t/2}(z) \cap V \text{ and } \partial B_{2t/3}(z) \cap V, \notag \\
		&\text{or both} \ \partial B_{5t/6}(z) \cap V \ \text{and}\  \partial B_t(z) \cap V.
	\end{align}
	The first condition is formulated in terms of $\overline{\mathop{\mathrm{int}}(\SCL)}$ so as to rule out simultaneously that a loop of $\Gamma_V$ crosses the annulus and that a loop of $\Gamma_V$ surrounds it; in the latter case $Q_{t/2,t}(z)$ would not exist and the distance between the two arcs would vanish.
	
\stepn{step:dist-rect-decoupling-loop-soup}{Decoupling the loop soup}  Realize $\Gamma_V$ as the outer boundaries of the outermost clusters of a critical Brownian loop soup $\Xi_V$ in $V$.  Write $D \defeq B_{2t}(z) \cap V$ and $\Xi_V = \Xi^{\mathrm{in}} + \Xi^{\mathrm{out}}$, where $\Xi^{\mathrm{in}}$ consists of the loops of $\Xi_V$ contained in $D$; the two parts are independent.  For a cluster $K$ we write $\SCL_K$ for its outer boundary, so that $\overline{\mathop{\mathrm{int}}(\SCL_K)}$ is $K$ together with the domains that it surrounds.  We consider the following five events.
	\begin{itemize}
		\item $E^{\mathrm{out}}$: no loop of $\Xi^{\mathrm{out}}$ intersects $\overline{B_{3t/2}(z)}$.
		\item $E_1$: there is no cluster $K$ of $\Xi^{\mathrm{in}}$ such that $\overline{\mathop{\mathrm{int}}(\SCL_K)}$ meets both $\partial B_{t/2}(z) \cap V$ and $\partial B_{2t/3
			}(z) \cap V$.
		\item $E_2$: there is no cluster $K$ of $\Xi^{\mathrm{in}}$ such that $\overline{\mathop{\mathrm{int}}(\SCL_K)}$ meets both $\partial B_{5t/6}(z) \cap V$ and $\partial B_{t}(z) \cap V$.
		\item $E_3$: there is no cluster $K$ of $\Xi^{\mathrm{in}}$ such that $\overline{\mathop{\mathrm{int}}(\SCL_K)}$ meets both $\partial B_t(z) \cap V$ and $\partial B_{3t/2}(z) \cap V$.
		
		\item $E^{\mathrm{sep}}$: no cluster of $\Xi_V$ separates $B_{2t/3
		}(z) \cap V$ from $\partial V$.
	\end{itemize}

\stepn{step:dist-rect-five-events-step}{The five events of \Cref{step:dist-rect-decoupling-loop-soup} imply~\eqref{eq:good_annulus_event}}  We first claim that, on $E^{\mathrm{out}} \cap E_2$, every cluster $K$ of $\Xi_V$ which intersects $\overline{B_t(z)}$ is a cluster of $\Xi^{\mathrm{in}}$ contained in $B_{3t/2}(z)$.  Indeed, a loop $\ell \in K$ meeting $\overline{B_t(z)}$ lies in $\Xi^{\mathrm{in}}$ by $E^{\mathrm{out}}$; the cluster $K_0$ of $\Xi^{\mathrm{in}}$ containing $\ell$ is contained in $B_{3t/2}(z)$, since otherwise $\overline{\mathop{\mathrm{int}}(\SCL_{K_0})}$, being connected and meeting both $\overline{B_t(z)}$ and the complement of $B_{3t/2}(z)$, would contradict $E_3$; and every loop of $\Xi_V$ meeting $K_0$ meets $\overline{B_{3t/2}(z)}$, hence lies in $\Xi^{\mathrm{in}}$ by $E^{\mathrm{out}}$ and thus in $K_0$.  Hence $K = K_0$.

	Now let $K$ be an outermost cluster of $\Xi_V$ such that $\overline{\mathop{\mathrm{int}}(\SCL_K)}$ meets both $\partial B_{t/2}(z) \cap V$ and $\partial B_{2t/3}(z) \cap V$.  If $K$ meets $\overline{B_t(z)}$, then $K$ is a cluster of $\Xi^{\mathrm{in}}$ by the claim of the previous paragraph, contradicting $E_1$; if not, then the connected set $B_{2t/3}(z) \cap V$, which is disjoint from $K$ and meets $\overline{\mathop{\mathrm{int}}(\SCL_K)}$, lies in a component of $V \setminus K$ surrounded by $K$, contradicting $E^{\mathrm{sep}}$.  Similarly,  suppose that $L$ is an outermost cluster of $\Xi_V$ such that $\overline{\mathop{\mathrm{int}}(\SCL_{L})}$ meets both $\partial B_{5t/6}(z) \cap V$ and $\partial B_t(z) \cap V$.  If $L$ meets $\overline{B_t(z)}$,  then $L$ is a cluster of $\Xi^{\mathrm{in}}$ and this contradicts $E_2$.  If $L$ does not meet $\overline{B_t(z)}$,  then the connected set $B_t(z) \cap V$ lies in a component of $V \setminus {L}$ surrounded by ${L}$,  contradicting $E^{\mathrm{sep}}$.  It follows that if $E^{\mathrm{out}} \cap E_1 \cap E_2 \cap E_3 \cap E^{\mathrm{sep}}$ holds,  we have that ~\eqref{eq:good_annulus_event} holds as well.
	
\stepn{step:dist-rect-each-five-events}{Each of the five events of \Cref{step:dist-rect-decoupling-loop-soup} has positive probability}  A loop of $\Xi_V$ which meets $\overline{B_{3t/2}(z)}$ and is not contained in $B_{2t}(z)$ crosses $A_{3t/2,2t}(z)$, and the number of such loops is a Poisson random variable whose parameter is at most the Brownian loop measure $m_0$ of the set of loops in $\BC$ crossing $A_{3/4,1}(0)$, which is finite and scale invariant.  Hence
	\begin{equation*}
		\BP\lbrack E^{\mathrm{out}}\rbrack \ge e^{-m_0} > 0,
	\end{equation*}
	uniformly in $z$ and $t$.  For $E^{\mathrm{sep}}$: if $\overline{B_t(z)} \not\subseteq V$, then $E^{\mathrm{sep}}$ holds automatically, since the closures of the clusters of $\Xi_V$ are at positive distance from $\partial V$ and hence do not disconnect $B_t(z) \cap V$ from $\partial V$; if $\overline{B_t(z)} \subseteq V$, then a separating cluster surrounds $z$, so that $E^{\mathrm{sep}}$ contains the event $\{\SCL(z) \subseteq B_t(z)\}$, which has positive probability by \cite[Corollary~1.2]{ALS22ExtremalDistance} and conformal invariance.  Finally, $E_1$, $E_2$ and $E_3$ concern only the Brownian loop soup in $D$, so that by scale invariance we may assume that $t = 1/2$ and $z = 0$.  The event $E_1$ (resp.\@ $E_2$) contains the event that the loop $\SCL^0$ of the $\CLE_4$ in $D$ surrounding $0$ satisfies $\overline{B_{1/4}(0)} \subseteq \overline{\mathop{\mathrm{int}}(\SCL^0)} \subseteq B_{1/3}(0)$ (resp.\@  $\overline{B_{5/12}(0)} \subseteq \overline{\mathop{\mathrm{int}}(\SCL^0)} \subseteq B_1(0)$) and this event has positive probability. (Indeed, one can sample an independent Brownian loop $\ell$ that surrounds $B_{1/4}(0)$ and stays in $B_{1/3}(0)$ and add it to the loop soup $\Xi^{\mathrm{in}}$, the resulting loop soup has a law that is absolutely continuous with respect to that of $\Xi^{\mathrm{in}}$. Moreover, conditionally on $\ell$, by the end of the proof of \cite[Lemma 9.6]{CLE}, the probability that there is a crossing by a cluster of $\Xi^{\mathrm{in}}$ of the annulus between $\ell$ and $\partial B_{1/3}(0)$ is strictly less than $1$.) 
	The same argument applied to $A_{1/4,3/4}(0)$ gives $\BP[E_3] > 0$.
	
Suppose now that $D \cap \partial V \neq \emptyset$.  Then,  by the monotonicity property of Brownian loop soup,  we can couple $\Xi^{\mathrm{in}}$ with a Brownian loop soup $\Xi$ in $B_1(0)$ such that $\Xi^{\mathrm{in}} \subseteq \Xi$.  Then,  it is clear that if the event $E_1$ (resp.\@ $E_2$) occurs for $\Xi$,  it has to occur for $\Xi^{\mathrm{in}}$ as well.  Thus,  by arguing as in the previous paragraph,  we obtain that $E_1$ (resp.\@ $E_2$) occurs for $\Xi^{\mathrm{in}}$ with positive probability and the same holds for $E_3$.
	
\stepn{step:dist-rect-proof-eq-good}{Proof of~\eqref{eq:good_annulus_event} by the FKG inequality}  Each of the five events of \Cref{step:dist-rect-decoupling-loop-soup} is a \emph{decreasing} event for the Brownian loop soup $\Xi_V$: adding loops to the configuration can only merge clusters and enlarge the domains that they surround, hence can only destroy $E^{\mathrm{out}}$, $E_1$, $E_2$, $E_3$,  and $E^{\mathrm{sep}}$.  Therefore,  by \Cref{lem:FKG_PPP} applied to the Brownian loop soup and to the indicator functions of these events, which are decreasing, and iterating,
	\begin{equation*}
		\BP\!\left\lbrack E^{\mathrm{out}} \cap E_1 \cap E_2 \cap E_3 \cap E^{\mathrm{sep}}\right\rbrack \ge \BP\lbrack E^{\mathrm{out}}\rbrack\,\BP\lbrack E_1\rbrack\,\BP\lbrack E_2\rbrack\,\BP\lbrack E_3\rbrack\,\BP\lbrack E^{\mathrm{sep}}\rbrack > 0.
	\end{equation*}
	Combining this with \Cref{step:dist-rect-five-events-step,step:dist-rect-each-five-events} proves~\eqref{eq:good_annulus_event}.

\stepn{step:dist-rect-cutting-annulus-along}{Cutting the annulus along deterministic slits}  Next, we claim that there exists $\varepsilon > 0$ such that
\begin{equation}\label{eq:dist_across_rect_arb_high_proof_0}
\BP\!\left\lbrack D_{\Gamma_V}^V(\partial B_{t/2}(z), \partial B_t(z)) > \varepsilon\right\rbrack > 0, \quad \forall z \in V\cap \frac{t}{100} \BZ^2, \ \forall t \in (0, 1/2). 
\end{equation}
(By convention, $D_{\Gamma_V}^V(\partial B_{t/2}(z), \partial B_t(z)) = \infty$ if $\partial B_{t/2}(z) \cap V = \emptyset$ or $\partial B_t(z) \cap V = \emptyset$.) By \Cref{step:dist-rect-proof-eq-good},  we have that the event ~\eqref{eq:good_annulus_event} has positive probability for each $z \in V$ and $t>0$.  On this event,  we have by ~\eqref{eq:modulus_lower_bound} and the discussion as the end of \Cref{step:dist-rect-notation-two-deterministic} that there exists some universal constant $C>1$ such that $C^{-1} \leq M_{t/2,t}(z) \leq C$ (note that if  ~\eqref{eq:good_annulus_event}  occurs,  we have that $A_{2t/3,5t/6}(z) \cap V \subseteq Q_{t/2,t}(z)$).  We work on the event  ~\eqref{eq:good_annulus_event}  for the remainder of \Cref{step:dist-rect-cutting-annulus-along,step:dist-rect-distance-across-single}; it is the two-sided bound on $M_{t/2,t}(z)$,  together with the fact that $\overline{\mathop{\mathrm{int}}(\SCL)}$ meets at most one of the two arcs for every $\SCL \in \Gamma_V$, that will be used.

Suppose first that $\overline{B_t(z)} \subseteq V$ and let $L_1 \defeq z + [t/2,t]$ and $L_2 \defeq z - [t/2,t]$ be two opposite radial slits of $\overline{A_{t/2,t}(z)}$, so that $A_{t/2,t}(z) \setminus (L_1 \cup L_2)$ has exactly two connected components $V^+$ and $V^-$.  Set also $S_j \defeq A_{t/2,t}(z) \setminus L_j$ for $j = 1,2$.  Then $V^\pm$ and $S_j$ are deterministic simply connected domains, so that Axiom~\eqref{it:axiom_locality} (locality) applies to each of them: for $V^\prime$ any one of these four domains and for each connected component $V^\prime_i$ of $(V^\prime)^\star$, we have $D_{\Gamma_V}^V(\bullet, \bullet; V^\prime_i) = D_{\Gamma_V|_{V^\prime_i}}^{V^\prime_i}(\bullet, \bullet)$ and, conditionally on $(V^\prime)^\star$, the loops of $\Gamma_V$ contained in $V^\prime_i$ form a non-nested $\CLE_4$ in $V^\prime_i$.  (Deterministic slits are used because the four sets above are then of the form required by Axiom~\eqref{it:axiom_locality}, with simply connected components of their $\star$-versions, so that both the $\CLE_4$ and the metric $D^\bullet$ are defined on them; neither is the case for $Q_{t/2,t}(z)$ itself, which is doubly connected when $\overline{B_{t/2}(z)} \subseteq V$.)  If instead $\overline{B_{t/2}(z)} \not\subseteq V$, then $A_{t/2,t}(z) \cap V$ is already simply connected and we take a single slit $L_1$ separating it into two components $V^+$ and $V^-$, with $S_1 \defeq (A_{t/2,t}(z) \cap V) \setminus L_1$, and we simply omit $S_2$ and $L_2$ from the discussion below.

\stepn{step:dist-rect-distance-across-single}{The distance across a single annulus}  Let $P$ be an admissible path in $Q_{t/2,t}(z)$ joining the part of $\partial Q_{t/2,t}(z)$ coming from $\partial B_{t/2}(z)$ to the part coming from $\partial B_t(z)$.  If $P$ avoids some $L_j$, then $P$ is an admissible path in $S_j$ joining its two boundary arcs on the two circles; if $P$ meets both $L_1$ and $L_2$, then the portion of $P$ between its last visit to $L_1$ and its first subsequent visit to $L_2$ is an admissible path in $V^+$ or in $V^-$ joining its two sides on $L_1$ and on $L_2$.  Hence
\begin{equation}\label{eq:crossing_min}
	D_{\Gamma_V}^V\left(\partial B_{t/2}(z) \cap V, \partial B_t(z) \cap V\right) \ge X_1 \wedge X_2 \wedge Y^+ \wedge Y^-,
\end{equation}
where $X_j$ (resp.\ $Y^\pm$) denotes the minimum, over the connected components $R$ of $S_j^\star$ (resp.\ of $(V^\pm)^\star$) which separate the two marked boundary arcs of $S_j$ (resp.\ of $V^\pm$), of the internal distance in $R$ between the two boundary arcs of $R$ which they contain.  (Here we use that the left-hand side is the infimum of $D_{\Gamma_V}^V(\SCL^\prime, \SCL^{\prime\prime})$ over loops $\SCL^\prime$ meeting $\partial B_{t/2}(z)$ and $\SCL^{\prime\prime}$ meeting $\partial B_t(z)$, and that a $D_{\Gamma_V}^V$-geodesic joining two such loops, which exists by Axiom~\eqref{it:axiom_geodesic}, contains a path $P$ as above.)

By Axiom~\eqref{it:axiom_locality} (locality) and Axiom~\eqref{it:axiom_conformal_invariance} (conformal invariance), each of the four quantities appearing in~\eqref{eq:crossing_min} is, conditionally on the corresponding region, distributed as $D_{\Gamma_{V_a}}^{V_a}(T_a, B_a)$, where $a > 0$ is the conformal modulus of that region viewed as a conformal rectangle; and each of them is a.s.\ positive, since on the event~\eqref{eq:good_annulus_event} we have, for every $\SCL \in \Gamma_V$, that $\overline{\mathop{\mathrm{int}}(\SCL)}$ meets at most one of $\partial B_{t/2}(z) \cap V$ and $\partial B_t(z) \cap V$.

Therefore, combining~\eqref{eq:crossing_min} with conformal invariance, we obtain that, conditionally on the four regions of \Cref{step:dist-rect-cutting-annulus-along} and on the events~\eqref{eq:good_annulus_event} and $E^{\mathrm{sep}}$, the $D_{\Gamma_V}^V$-distance between $\partial B_{t/2}(z) \cap V$ and $\partial B_t(z) \cap V$ stochastically dominates a random variable which is positive with non-zero probability and whose law does not depend on $t\in (0, 1/2)$ (by scaling) or $z \in V \cap  ({t}/{100}) \BZ^2 $ (since there is a finite number of possibilities that does not depend on $t$, using $E^{\mathrm{sep}}$ when $\overline{B_{2t}(z)} \subset V$).  Choosing $\varepsilon>0$ so that this random variable exceeds $\varepsilon$ with positive probability (which is possible with $\varepsilon$ being universal,  and hence uniformly in $z$ and $t$) and multiplying by the (positive) probability of the event just described,  this proves~\eqref{eq:dist_across_rect_arb_high_proof_0}.

\stepn{step:dist-rect-from-single-annulus}{From a single annulus to a crossing of $V$}  It follows from~\eqref{eq:dist_across_rect_arb_high_proof_0} and \Cref{lem:FKG-CLE} (applied to the decreasing measurable functions $\one_{\{D_{\Gamma_V}^V(\partial B_{t/2}(z), \partial B_t(z)) > \varepsilon\}}$, and iterated over the finitely many $z$'s appearing in the product below) that for each $t  \in (0, 1/2)$,
\begin{multline}\label{eq:dist_across_rect_arb_high_proof_1}
\BP\!\left\lbrack D_{\Gamma_V}^V(\partial B_{t/2}(z), \partial B_t(z)) > \varepsilon \text{ for all } z \in \left(\frac1{100}t\BZ\right)^2 \cap V \text{ with } \dist(z, T \cup B) \ge 2t\right\rbrack \\
\ge \prod_{\substack{z \in \left(\frac1{100}t\BZ\right)^2 \cap V\\\dist(z, T \cup B) \ge 2t}} \BP\!\left\lbrack D_{\Gamma_V}^V(\partial B_{t/2}(z), \partial B_t(z)) > \varepsilon\right\rbrack > 0.
\end{multline}
(Here we have replaced the constraint $\dist(z, T \cup B) \ge t$ by $\dist(z, T \cup B) \ge 2t$, so that every point of $V$ at distance at least $3t$ from $T \cup B$ has a grid point $z$ of the index set within distance $t/100$ of it; the product is over finitely many $z$'s because $V$ is bounded, and each factor is positive by~\eqref{eq:dist_across_rect_arb_high_proof_0}.)
Finally, we observe that, on the event of~\eqref{eq:dist_across_rect_arb_high_proof_1}, we have $D_{\Gamma_V}^V(T, B) > \varepsilon/(100t)$. Indeed, by Axiom~\eqref{it:axiom_geodesic} there is a $D_{\Gamma_V}^V$-geodesic $P$ from $T$ to $B$; since $P$ joins the two horizontal sides of $V$, one may select points $w_1, \ldots, w_m$ visited by $P$ in this order, with $m \ge 1/(4t)$, lying at distance at least $3t$ from $T \cup B$ and at mutual distance at least $3t$, and grid points $z_1, \ldots, z_m$ of the index set of~\eqref{eq:dist_across_rect_arb_high_proof_1} with $\lvert z_j - w_j\rvert \le t/100$.  The annuli $A_{t/2,t}(z_j)$ are then pairwise disjoint and $P$ crosses each of them, so that $\len(P; D_{\Gamma_V}^V) \ge m\varepsilon$. By choosing $t$ to be sufficiently small, namely $t < \varepsilon/(100N)$, we complete the proof of the lemma. 
\end{proof}

Next, let us state the following weaker version of \Cref{lem:law_cont}.

\begin{lemma}\label{lem:law_cont_pre}
For each deterministic subset $S \subseteq \BR$ of Lebesgue measure zero, the subset
\begin{equation*}
\left\{r > 0 : \BP\!\left\lbrack D_{\Gamma_{V_r}}^{V_r}(T_r, B_r) \in S\right\rbrack > 0\right\}
\end{equation*}
has Lebesgue measure zero. 
\end{lemma}

As explained above,  we will obtain \Cref{lem:law_cont} from \Cref{lem:law_cont_pre} by running an appropriate metric ball exploration for a short amount of time and then using the fact that the law of the conformal modulus of the unexplored region is absolutely continuous with respect to the Lebesgue measure.  In order to prove the latter fact,  we will use Lemmas~\ref{lem:sle4_cm} and~\ref{lem:sle4_cm_reverse} below, which prove the analogous property for the regions in the complement of an $\SLE_4$ curve.

\begin{lemma}\label{lem:sle4_cm}
Let $x, y \in (-\infty, 0)$ and $z \in (0, 1)$ be random points. Suppose that, conditionally on $z$, the law of the pair $(x,y)$ is absolutely continuous with respect to the Lebesgue measure on $(-\infty,0)^2$. Given $x$, $y$, and $z$, let $\eta$ be a chordal SLE$_4$ curve in $\BH$ from $x$ to $y$. Then the law of the conformal modulus of the conformal rectangle $(\BH \setminus \eta; 0, z, 1, \infty)$ is absolutely continuous with respect to the Lebesgue measure.
\end{lemma}

Throughout this subsection, whenever $\eta$ is a curve in $\overline\BH$ with both endpoints on $\partial\BH \cap (-\infty,0)$ and $K \subseteq \overline\BH$ is a compact set, we write $(\BH \setminus (\eta \cup K); 0, z, 1, \infty)$ for the conformal rectangle whose underlying domain is the unbounded connected component of $\BH \setminus (\eta \cup K)$.

\begin{lemma}\label{lem:sle4_cm_reverse}
Suppose that we have the exact same setup as in Lemma~\ref{lem:sle4_cm}.  Then,  conditionally on $z$,  if $L_0$ denotes the conformal modulus of $(\BH ; 0,z,1,\infty)$,  we have that the restriction of the Lebesgue measure to $[0,L_0]$ is absolutely continuous with respect to the conditional law of the conformal modulus of $(\BH \setminus \eta ; 0,z,1,\infty)$.
\end{lemma}

Note that $L_0$ is a function of $z$ alone, so that the statement makes sense conditionally on $z$.

Before giving the proofs of Lemmas~\ref{lem:sle4_cm} and~\ref{lem:sle4_cm_reverse},  we state and prove two simple but useful lemmas.

\begin{lemma}\label{lem:conformal-restriction}
    Let $K \subset \overline\BH$ be a compact subset with finitely many connected components, each of which meets $\partial \BH$, such that $\BH \setminus K$ has a unique unbounded connected component $H_K$, and let $u \neq v$ be two points of $\partial\BH \setminus K$ lying on $\partial H_K$. Fix a conformal mapping $\phi$ from $H_K$ onto $\BH$. Let $\eta$ be a chordal SLE$_4$ curve in $\BH$ from $u$ to $v$. Then, the law of $\phi(\eta)$ conditionally on the event that $\eta \cap K = \emptyset$ is mutually absolutely continuous with respect to the law of a chordal SLE$_4$ curve in $\BH$ from $\phi(u)$ to $\phi(v)$.
\end{lemma}
\begin{proof}
    This follows from \cite[Theorem 4.55]{Law09}.  (The case in which $K$ is disconnected follows from the connected case by iterating.)
\end{proof}

\begin{lemma}\label{lem:Newton-Leibniz}
    Let $f \colon [0, 1] \to \BR$ be a continuous and strictly increasing function. Then $f^\prime(x)$ exists for almost every $x \in [0,1]$ and
    \begin{equation*}
        \int_A f^\prime(x) \, \rd x \le \Leb(f(A))
    \end{equation*}
    for all Borel $A \subset [0, 1]$.
\end{lemma}

\begin{proof}
    Write $\mu_f$ for the measure on $[0, 1]$ associated with $f$ (i.e., $\mu_f([a, b]) = f(b) - f(a)$). Since $f$ is monotone, $f^\prime$ exists almost everywhere by Lebesgue's differentiation theorem for monotone functions. By the Lebesgue decomposition theorem, we have $\mu_f = \mu_{ac} + \mu_s$, where $\mu_{ac}$ is absolutely continuous with respect to the Lebesgue measure and $\mu_s$ is singular. Then $\Leb(f(A)) = \mu_f(A) \ge \mu_{ac}(A) = \int_A f^\prime(x) \, \rd x$, where the last equality comes from the fact that $\mu_{ac}= f'(x) \mathrm{d} x$ (see e.g.\ \cite[Theorem 7.14]{Rud87}). For the first equality, note that since $f$ is continuous and strictly increasing it is a homeomorphism onto its image, so that $f(A)$ is a Borel set, and that the Borel measures $A \mapsto \Leb(f(A))$ and $\mu_f$ agree on intervals and hence coincide.
\end{proof}

By composing with an affine map and exhausting by compact subintervals, the same conclusion holds with $[0,1]$ replaced by an arbitrary interval $I \subseteq \BR$ and Borel $A \subseteq I$.  If, moreover, $f$ is absolutely continuous on compact subintervals of $I$, then the inequality above is an equality, since then $\mu_s = 0$.
\begin{proof}[Proof of \Cref{lem:sle4_cm}]
\stepn{step:sle4-cm-setup}{Setup}  Since the assertion of the lemma concerns the law of a quantity which is measurable with respect to $(\eta,z)$, and since by hypothesis the conditional law of $(x,y)$ given $z$ is absolutely continuous with respect to the Lebesgue measure on $(-\infty,0)^2$, we may condition on $z$ and thereby assume throughout the proof that $z \in (0,1)$ is deterministic.  We may then also assume without loss of generality that the law of $(x,y)$ is absolutely continuous with respect to $\one_{\{u < 0,\, v<0\}} \, \rd u \, \rd v$. Fix $0 < a < 1$ and $b < 0$ and a sufficiently small $\varepsilon_0 > 0$. On the event that $a < z-\varepsilon_0$, for $\varepsilon \in (0, \varepsilon_0]$, write $\delta = \delta(\varepsilon) > 0$ for the unique number so that $(\BH \setminus (B_\varepsilon(a) \cup B_\delta(b)); 0, z, 1, \infty)$ and $(\BH; 0, z, 1, \infty)$ have the same conformal modulus. Such a $\delta$ exists for all $\varepsilon \in (0, \varepsilon_0]$ when $\varepsilon_0$ is chosen small enough. Indeed, removing $B_\varepsilon(a) \cap \BH$ (centered on the arc $[0,z]$) strictly increases the modulus, i.e.\ the extremal distance between the boundary arcs $[-\infty,0]$ and $[z,1]$, while removing $B_\delta(b) \cap \BH$ (centered on $[-\infty,0]$) strictly decreases it, in both cases continuously and strictly monotonically and with an effect vanishing as the radius tends to $0$; the intermediate value theorem therefore yields a unique $\delta(\varepsilon)$ once $\varepsilon_0$ is small enough, and $\varepsilon \mapsto \delta(\varepsilon)$ is continuous, strictly increasing, and tends to $0$ at $0$. Write $E(a, b, \varepsilon_0)$ for the event that $a < z-\varepsilon_0$, $b \in (x \wedge y, x \vee y)$, and $\{x, y, 0, z\} \cap (B_\varepsilon(a) \cup B_\delta(b)) = \emptyset$ for all $\varepsilon \in (0, \varepsilon_0]$. Write $\phi_\varepsilon$ for the conformal mapping from the unbounded connected component of $\BH \setminus (B_\varepsilon(a) \cup B_\delta(b))$ onto $\BH$ that fixes $0$, $z$, $1$, and $\infty$; such a map exists by the choice of $\delta(\varepsilon)$, since two conformal rectangles are conformally equivalent if and only if they have the same modulus.

\stepn{step:sle4-cm-randomization-reduction}{Randomization of $\varepsilon$ and reduction}  Given $x$, $y$, $z$, $\eta$, and the event $E(a, b, \varepsilon_0)$, let $\varepsilon$ be sampled uniformly from $(0, \varepsilon_\ast \wedge \varepsilon_0)$, where 
\begin{equation*}
\varepsilon_\ast \defeq \sup\{\varepsilon > 0 : \eta \cap (B_\varepsilon(a) \cup B_{\delta(\varepsilon)}(b)) = \emptyset\}.
\end{equation*}
Since $\varepsilon \mapsto \delta(\varepsilon)$ is increasing, the set on the right-hand side is an interval of the form $(0,\varepsilon_\ast)$.  Moreover $\varepsilon_\ast > 0$ a.s.: a chordal SLE$_4$ curve is simple and meets $\partial \BH$ only at its two endpoints, so that on $E(a,b,\varepsilon_0)$ the curve $\eta$ is at positive distance from both $a$ and $b$.
Write $\eta^\prime \defeq \phi_\varepsilon(\eta)$. Note that, on the event $E(a, b, \varepsilon_0)$, the law of the pair $(\phi_\varepsilon(x), \phi_\varepsilon(y))$ is still absolutely continuous with respect to $\one_{\{u<0,\, v<0\}} \, \rd u \, \rd v$, since $\phi_\varepsilon$ restricted to a neighborhood of $(x \wedge y, x \vee y)$ in $\partial \BH$ is a diffeomorphism onto its image. Thus, by \Cref{lem:conformal-restriction}, the law of $\eta'$ conditionally on the event $E(a, b, \varepsilon_0)$ is absolutely continuous with respect to the law of $\eta$ conditionally on the event $E(a, b, \varepsilon_0)$. Indeed, let $\widetilde{\varepsilon}$ be a uniform random variable in $(0, \varepsilon_0)$. Then, the law of $\phi_\varepsilon(\eta)$ conditionally on $E(a, b, \varepsilon_0)$ is the same as the law of $\phi_{\widetilde{\varepsilon}}$ conditionally on $E(a, b, \varepsilon_0) \cap \{\eta \cap (B_{\widetilde{\varepsilon}}(a) \cup B_{\delta(\widetilde{\varepsilon})}(b)) = \emptyset\}$, so that one can indeed apply \Cref{lem:conformal-restriction}.

Therefore, since there are always $a \in (0, 1) \cap \BQ$, $b \in (-\infty, 0) \cap \BQ$, and $\varepsilon_0 \in \BR_{>0} \cap \BQ$ such that $E(a, b, \varepsilon_0)$ occurs, it suffices to show that, on the event $E(a, b, \varepsilon_0)$, the law of the conformal modulus of the conformal rectangle $(\BH \setminus \eta^\prime; 0, z, 1, \infty)$ is absolutely continuous with respect to the Lebesgue measure. Note that the conformal modulus of $(\BH \setminus \eta^\prime; 0, z, 1, \infty)$ is the same as the conformal modulus of $(\BH \setminus (\eta \cup B_\varepsilon(a)); 0, z, 1, \infty)$. Thus, it suffices to show that, on the event $E(a, b, \varepsilon_0)$, the law of the conformal modulus of $(\BH \setminus (\eta \cup B_\varepsilon(a)); 0, z, 1, \infty)$ is absolutely continuous with respect to the Lebesgue measure.

Indeed, $\phi_\varepsilon$ maps the unbounded connected component of $\BH \setminus (B_\varepsilon(a) \cup B_{\delta(\varepsilon)}(b))$ conformally onto $\BH$ fixing $0$, $z$, $1$, and $\infty$, so the conformal modulus of $(\BH \setminus \eta^\prime; 0,z,1,\infty)$ equals that of $(\BH \setminus (\eta \cup B_\varepsilon(a) \cup B_{\delta(\varepsilon)}(b)); 0,z,1,\infty)$; and on $E(a,b,\varepsilon_0)$, for $\varepsilon < \varepsilon_\ast$, the half-disk $B_{\delta(\varepsilon)}(b) \cap \BH$ is contained in the bounded connected component of $\BH \setminus \eta$ (since $b \in (x \wedge y, x \vee y)$ and $\eta \cap B_{\delta(\varepsilon)}(b) = \emptyset$), so removing it does not affect the unbounded connected component.

\stepn{step:sle4-cm-strict-monotonicity-modulus}{Strict monotonicity of the modulus}

\substepn{step:scm-claim}{The claim and reduction to $\Theta$} We \emph{claim} that, given $\eta$, the assignment
\begin{equation*}
(0, \varepsilon_\ast \wedge \varepsilon_0) \to \BR \colon \varepsilon \mapsto f_\eta(\varepsilon) \defeq \mathop{\mathrm{CM}}(\BH \setminus (\eta \cup B_\varepsilon(a)); 0, z, 1, \infty)
\end{equation*}
is a.s.\ strictly increasing, where $\mathop{\mathrm{CM}}(\bullet)$ denotes the conformal modulus. To lighten notation, write $H_{\eta,\varepsilon} \defeq \BH \setminus (\eta \cup B_\varepsilon(a))$. Recall that $\mathop{\mathrm{CM}}(H_{\eta,\varepsilon}; 0, z, 1, \infty) = d_{H_{\eta,\varepsilon}}([-\infty, 0]_{\partial H_{\eta,\varepsilon}}^\circlearrowleft, [z, 1]_{\partial H_{\eta,\varepsilon}}^\circlearrowleft)$, where $d_{H_{\eta,\varepsilon}}$ denotes the extremal distance in $H_{\eta,\varepsilon}$ and $[-\infty,0]_{\partial H_{\eta,\varepsilon}}^\circlearrowleft$
(resp.\ $[z, 1]_{\partial H_{\eta,\varepsilon}}^\circlearrowleft$) denotes the counterclockwise arc of $\partial H_{\eta,\varepsilon}$ from $-\infty$ to $0$ (resp.\ $z$ to $1$).  Fix $0<\varepsilon^\prime < \varepsilon^{\prime\prime} < \varepsilon_* \wedge \varepsilon_0$ and for all $\varepsilon \in (0 ,  \varepsilon_* \wedge \varepsilon_0)$,  we let $\psi_{\varepsilon}$ denote the conformal transformation mapping the unbounded connected component of $\BH \setminus (\eta \cup B_{\varepsilon}(a))$ onto the rectangle $R_{L_\varepsilon} \defeq (0,L_{\varepsilon}) \times (0,\pi)$ such that $\psi_{\varepsilon}(0)=0,  \psi_{\varepsilon}(\infty) = \ri\pi,  \psi_{\varepsilon}(1) = L_{\varepsilon} +\ri\pi$,  and $\psi_{\varepsilon}(z) = L_{\varepsilon}$,  where $L_{\varepsilon} \defeq \pi f_{\eta}(\varepsilon)$.  (The height $\pi$ matches the normalization of \Cref{lem:extremal_length}; since the extremal distance in $R_L$ between its two vertical sides equals $L/\pi$, we have $L_\varepsilon = \pi f_\eta(\varepsilon)$.  The map $\psi_\varepsilon$ takes $[-\infty,0]^\circlearrowleft_{\partial H_{\eta,\varepsilon}}$ and $[z,1]^\circlearrowleft_{\partial H_{\eta,\varepsilon}}$ onto the left and the right vertical side of $R_{L_\varepsilon}$, respectively.)  Recall the definition of the function
\begin{align*}
\Theta(L) = \sup\left\{2\min\left\{\BP_w\lbrack B_{\tau} \in [0,\ri\pi] \rbrack ,  \BP_w\lbrack B_{\tau} \in [L,L+\ri\pi] \rbrack \right\}:  w \in R_L \right\}
\end{align*}
in Lemma~\ref{lem:extremal_length},  where $B$ is a planar Brownian motion stopped at the first time $\tau$ that it exits $R_L$.  Note also that Lemma~\ref{lem:extremal_length} implies that the supremum in the definition of $\Theta$ is attained when $w = L /2 + \ri \pi / 2$,  and $\Theta(L)$ is a continuous,  strictly decreasing function of $L$.  Thus,  to prove the \emph{claim},  it suffices to show that $\Theta(L_{\varepsilon^\prime}) > \Theta(L_{\varepsilon^{\prime\prime}})$.

\substepn{step:scm-exit-probabilities}{Comparison of exit probabilities} To show the latter,  we note that for all $w \in \BH \cap \partial  B_{\varepsilon^{\prime\prime}}(a)$,  the probability that a planar Brownian motion starting from $w$ exits $H_{\eta,\varepsilon^\prime}$ on $[-\infty,0]_{\partial H_{\eta,\varepsilon^\prime}}^\circlearrowleft$ (resp.\ $[z, 1]_{\partial H_{\eta,\varepsilon^\prime}}^\circlearrowleft$) is positive.  Therefore,  combining with the Markov property of the planar Brownian motion,  we obtain that the following is true.  The probability that a planar Brownian motion starting from $\psi_{\varepsilon^{\prime\prime}}^{-1}(L_{\varepsilon^{\prime\prime}} / 2 + \ri \pi / 2)$ exits $H_{\eta,\varepsilon^{\prime\prime}}$ for the first time on $[-\infty,0]_{\partial H_{\eta,\varepsilon^{\prime\prime}}}^\circlearrowleft$ (resp.\ $[z, 1]_{\partial H_{\eta,\varepsilon^{\prime\prime}}}^\circlearrowleft$) is strictly less than the probability that it exits $H_{\eta,\varepsilon^\prime}$ for the first time on $[-\infty,0]_{\partial H_{\eta,\varepsilon^\prime}}^\circlearrowleft$ (resp.\ $[z, 1]_{\partial H_{\eta,\varepsilon^\prime}}^\circlearrowleft$).  Indeed, since $\varepsilon^\prime < \varepsilon^{\prime\prime}$ we have $H_{\eta,\varepsilon^{\prime\prime}} \subseteq H_{\eta,\varepsilon^\prime}$, and the two marked boundary arcs are the same for both domains, since $B_{\varepsilon^{\prime\prime}}(a)$ is disjoint from $\{0,z\}$ on $E(a,b,\varepsilon_0)$.  Writing $w_0 \defeq \psi_{\varepsilon^{\prime\prime}}^{-1}(L_{\varepsilon^{\prime\prime}}/2 + \ri\pi/2) \in H_{\eta,\varepsilon^{\prime\prime}}$ and $\alpha$ for either of the two marked arcs, the strong Markov property applied at the first hitting time of $\partial B_{\varepsilon^{\prime\prime}}(a)$ gives
\begin{align*}
\BP_{w_0}\lbrack B \text{ exits } H_{\eta,\varepsilon^\prime} \text{ on } \alpha\rbrack &= \BP_{w_0}\lbrack B \text{ exits } H_{\eta,\varepsilon^{\prime\prime}} \text{ on } \alpha\rbrack\\
&\quad + \BE_{w_0}\!\left\lbrack \one_{\{B \text{ exits } H_{\eta,\varepsilon^{\prime\prime}} \text{ on } \partial B_{\varepsilon^{\prime\prime}}(a)\}} \BP_{B_{\tau}}\lbrack B \text{ exits } H_{\eta,\varepsilon^\prime} \text{ on } \alpha\rbrack\right\rbrack,
\end{align*}
and the second term on the right-hand side is strictly positive by the first sentence of this paragraph, together with the fact that $\BP_{w_0}[B \text{ exits } H_{\eta,\varepsilon^{\prime\prime}} \text{ on } \partial B_{\varepsilon^{\prime\prime}}(a)] > 0$.  Thus,  by the conformal invariance of planar Brownian motion, $\Theta(L_{\varepsilon^{\prime\prime}})$ equals twice the minimum of the two exit probabilities from $w_0$ in $H_{\eta,\varepsilon^{\prime\prime}}$ (the supremum in the definition of $\Theta$ being attained at the center of the rectangle), while $\Theta(L_{\varepsilon^\prime})$ is at least twice the corresponding minimum in $H_{\eta,\varepsilon^\prime}$; hence it follows that $\Theta(L_{\varepsilon^\prime}) > \Theta(L_{\varepsilon^{\prime\prime}})$ and so $L_{\varepsilon^\prime} < L_{\varepsilon^{\prime\prime}}$,  i.e.,  $f_{\eta}(\varepsilon^\prime) < f_{\eta}(\varepsilon^{\prime\prime})$.  This proves the claim.  In particular,  since a monotone function is differentiable Lebesgue almost everywhere with a finite derivative, we conclude that $f_{\eta}'(\varepsilon)$ exists and it is finite for Lebesgue almost every $\varepsilon \in (0,\varepsilon_* \wedge \varepsilon_0)$.

\stepn{step:sle4-cm-properties}{Properties of $\Theta$}
Note that it follows from the proof of \cite[Proposition~3.69]{lawler2008conformally} that
\begin{align*}
\Theta(L) = \sum_{\substack{n \in \BN \\ n \ \text{odd}}} (-1)^{(n-1) / 2} \frac{8}{n} \frac{1}{2\pi \cosh(n L / 2)} \quad \text{for all} \quad L>0.
\end{align*}
In particular,  it follows from the above explicit expression and from the local uniform convergence of the differentiated series that $\Theta$ is  Lipschitz on compact subsets of $(0,\infty)$ and hence  we obtain that $\Theta$ is absolutely continuous when restricted to compact subsets of $(0,\infty)$.  Furthermore,  $\Theta'(L)$ exists for all $L>0$ and takes the explicit form
\begin{align*}
\Theta'(L) = -\sum_{\substack{n \in \BN \\ n \ \text{odd}}} (-1)^{(n-1) / 2} \frac{4}{\pi} \frac{e^{n L /2} - e^{-nL/2}}{(e^{nL/2} + e^{-nL /2})^2}.
\end{align*}
From this explicit form,  it is easy to see that $\Theta'$ is continuous.

Moreover, $\Theta'(L) < 0$ for all but countably many $L > 0$, hence for Lebesgue almost every $L>0$.  Indeed, since $\lvert \cosh(nL/2)\rvert \ge \sinh(nc/2)$ for $\Re(L) \ge c > 0$, the series defining $\Theta$ converges locally uniformly on the right half-plane and defines a holomorphic function there, so $\Theta$ and $\Theta'$ are real analytic on $(0,\infty)$; since $\Theta$ is strictly decreasing, $\Theta'$ does not vanish identically, so its zero set is discrete in $(0,\infty)$.

Finally, set $\widehat\Theta(u) \defeq \Theta(\pi u)$ for $u > 0$, so that $\widehat\Theta(f_\eta(\varepsilon)) = \Theta(L_\varepsilon)$.  Then $\widehat\Theta$ inherits all of the above properties of $\Theta$, and $\widehat\Theta^\prime(u) = \pi \Theta^\prime(\pi u) < 0$ for all but countably many $u>0$.

\stepn{step:sle4-cm-two-sided-bounds}{Two-sided bounds on the derivative of $g_\eta \defeq \widehat\Theta \circ f_\eta$}
Consider the function $g_{\eta} = \widehat\Theta \circ f_{\eta}$ and note that $g_{\eta}'(\varepsilon)$ exists and it is equal to $\widehat\Theta'(f_{\eta}(\varepsilon)) f_{\eta}'(\varepsilon)$ for Lebesgue almost every $\varepsilon \in (0,\varepsilon_* \wedge \varepsilon_0)$, by the chain rule, which applies at every $\varepsilon$ at which $f_\eta$ is differentiable since $\widehat\Theta$ is continuously differentiable.  Fix $0 < \widetilde\varepsilon < \varepsilon^\sharp < \varepsilon_\ast \wedge \varepsilon_0$; all the constants below are allowed to depend on $\eta$, on $\widetilde\varepsilon$ and on $\varepsilon^\sharp$, but not on the quantities $\varepsilon$ and $h$ introduced next.  Fix $\varepsilon \in [\widetilde\varepsilon,\varepsilon^\sharp)$ and let $h \in (0,1)$ be such that $\varepsilon + h \le \varepsilon^\sharp$.  Then,  there exist (random) constants $c_1,c_2>0$ which do not depend on $\varepsilon$ and $h$ such that the following is true.  The probability that a planar Brownian motion starting from $\psi_{\varepsilon + h}^{-1}(L_{\varepsilon + h} / 2 + \ri \pi / 2)$ exits $H_{\eta,\varepsilon + h}$ on $\{w \in \BH : \Im(w) \geq \varepsilon/2 \} \cap \partial B_{\varepsilon+h}(a)$ is at least $c_1$,  and the following holds for all $w \in \partial B_{\varepsilon + h}(a)$ such that $\Im(w) \geq \varepsilon/2$.  The probability that a planar Brownian motion starting from $w$ exits $H_{\eta,\varepsilon}$ on $[-\infty,0]_{\partial H_{\eta,\varepsilon}}^\circlearrowleft$
(resp.\ $[z, 1]_{\partial H_{\eta,\varepsilon}}^\circlearrowleft$) is at least $c_2 h$.  (For the first bound, the arc $\{\Im \ge \varepsilon/2\} \cap \partial B_{\varepsilon+h}(a)$ subtends an angle bounded away from $0$ uniformly in $h$, and the configuration converges to that of $H_{\eta,\varepsilon}$ as $h \downarrow 0$, with positive limiting probability.  For the second, such a $w$ lies at distance $h$ from $\partial B_\varepsilon(a)$ and at distance at least $\varepsilon/2$ from $\partial \BH$, so the harmonic measure from $w$ of a fixed boundary arc at positive distance from $w$ is of order $h/\varepsilon$.)  Thus,  combining with the strong Markov property of the planar Brownian motion,  we obtain that the following is true.  The probability that a planar Brownian motion starting from $\psi_{\varepsilon + h}^{-1}(L_{\varepsilon + h} / 2 + \ri \pi / 2)$ exits $H_{\eta,\varepsilon}$ on $[-\infty,0]_{\partial H_{\eta,\varepsilon}}^\circlearrowleft$
(resp.\ $[z, 1]_{\partial H_{\eta,\varepsilon}}^\circlearrowleft$) and before doing that,  it intersects $\BH \cap \partial B_{\varepsilon + h}(a)$,  is at least $c_1 c_2 h$.  Hence,  combining with the definition of $\Theta$ as in \Cref{step:sle4-cm-strict-monotonicity-modulus} ($g_\eta(\varepsilon+h)$ is exactly twice the minimum of the two exit probabilities from that starting point in $H_{\eta,\varepsilon+h}$, while $g_\eta(\varepsilon)$ is at least twice the corresponding minimum in $H_{\eta,\varepsilon}$),  we obtain that $g_{\eta}(\varepsilon) - g_{\eta}(\varepsilon + h) \geq 2c_1 c_2 h$.

Similarly,  there exists a (random) constant $c_3>0$ such that the following holds for all $w \in \BH \cap \partial B_{\varepsilon + h}(a)$.  The probability that a planar Brownian motion starting from $w$ exits $H_{\eta, \varepsilon}$ on $[-\infty,0]_{\partial H_{\eta,\varepsilon}}^\circlearrowleft$
(resp.\ $[z, 1]_{\partial H_{\eta,\varepsilon}}^\circlearrowleft$) is at most $c_3 h$.  (This holds also for $w$ near the two points $a \pm \varepsilon$ at which $\partial B_\varepsilon(a)$ meets $\partial \BH$: there $H_{\eta,\varepsilon}$ has a corner of interior angle $\pi/2$, and the harmonic measure from a point at distance of order $h$ from the corner of a boundary arc at positive distance from it is of order $h^2 \le h$.)  Therefore,  arguing as in the previous paragraph,  we obtain that $g_{\eta}(\varepsilon) - g_{\eta}(\varepsilon + h) \leq 2c_3 h$.  Since $g_\eta$ is decreasing, the last bound shows that $g_\eta$ is Lipschitz, hence absolutely continuous, on $[\widetilde\varepsilon,\varepsilon^\sharp]$; and since $\widehat\Theta$ is a homeomorphism onto its image and $f_\eta = \widehat\Theta^{-1} \circ g_\eta$, the function $f_\eta$ is continuous on $(0,\varepsilon_\ast \wedge \varepsilon_0)$.  It follows that for all $0<\widetilde{\varepsilon} < \varepsilon^\sharp <\varepsilon_* \wedge \varepsilon_0$,  there exists a (random) constant $C>1$ such that $C^{-1} \leq |g_{\eta}'(\varepsilon)| \leq C$ for Lebesgue almost every $\varepsilon \in [\widetilde{\varepsilon},\varepsilon^\sharp]$.  In particular,  since $\widehat\Theta^\prime$ is bounded on the compact set $f_\eta([\widetilde\varepsilon,\varepsilon^\sharp]) \subset (0,\infty)$, we obtain that $f_{\eta}'(\varepsilon) > 0$ for Lebesgue almost every $\varepsilon \in (0,\varepsilon_* \wedge \varepsilon_0)$.

\stepn{step:sle4-cm-conclusion}{Conclusion}
Finally, for each deterministic subset $S \subseteq \BR$ of Lebesgue measure zero, which we may and do assume to be Borel (replacing $S$ by a Borel null set containing it), we have
\begin{multline}\label{eq:sle4_cm_proof_0}
\BP\lbrack\{\mathop{\mathrm{CM}}(\BH \setminus (\eta \cup B_\varepsilon(a)); 0, z, 1, \infty) \in S \} \cap E(a,b,\varepsilon_0)\rbrack \\ = \BE\lbrack\BP\lbrack f_\eta(\varepsilon) \in S \mid \eta \rbrack \one_{E(a,b,\varepsilon_0)} \rbrack = \BE\!\left\lbrack\frac{\Leb(f_\eta^{-1}(S))}{\varepsilon_\ast \wedge \varepsilon_0} \one_{E(a,b,\varepsilon_0)} \right\rbrack,
\end{multline}
where $f_\eta^{-1}(S) \defeq \{\varepsilon \in (0, \varepsilon_\ast \wedge \varepsilon_0) : f_\eta(\varepsilon) \in S\}$, the factor $(\varepsilon_\ast \wedge \varepsilon_0)^{-1}$ coming from the fact that, given $\eta$, the variable $\varepsilon$ is uniform on $(0,\varepsilon_\ast \wedge \varepsilon_0)$.  By \Cref{lem:Newton-Leibniz}, applied on each compact subinterval of $(0,\varepsilon_\ast \wedge \varepsilon_0)$ to the continuous, strictly increasing function $f_\eta$, and letting the subintervals exhaust $(0,\varepsilon_\ast \wedge \varepsilon_0)$,
\begin{equation*}
\int_{f_\eta^{-1}(S)} f_\eta^\prime(\varepsilon) \, \rd\varepsilon \le \Leb(f_\eta(f_\eta^{-1}(S))) = 0 \quad \text{a.s.\ on} \quad E(a,b,\varepsilon_0),
\end{equation*}
where the last equality holds because $f_\eta(f_\eta^{-1}(S)) \subseteq S$ and $\Leb(S) = 0$.
This, together with the positivity of $f_{\eta}'$ Lebesgue almost everywhere, which was established above, implies that $\Leb(f_\eta^{-1}(S)) = 0$ a.s.\ on $E(a,b,\varepsilon_0)$.  Combining this with~\eqref{eq:sle4_cm_proof_0}, we obtain that
\[\BP\lbrack\{\mathop{\mathrm{CM}}(\BH \setminus (\eta \cup B_\varepsilon(a)); 0, z, 1, \infty) \in S\ \cap E(a,b,\varepsilon_0)\rbrack = 0.\]
By \Cref{step:sle4-cm-randomization-reduction}, it follows that on $E(a,b,\varepsilon_0)$ the law of the conformal modulus of $(\BH \setminus \eta; 0,z,1,\infty)$ is absolutely continuous with respect to the Lebesgue measure; taking the union over the countably many rational triples $(a,b,\varepsilon_0)$, whose associated events cover the probability space up to a null set, and averaging over $z$, we conclude that $\BP\lbrack\mathop{\mathrm{CM}}(\BH \setminus \eta; 0,z,1,\infty) \in S\rbrack = 0$.
This completes the proof of the lemma. 
\end{proof}

\begin{proof}[Proof of Lemma~\ref{lem:sle4_cm_reverse}]
\stepn{step:sle4-cm-rev-setup}{Setup}
Suppose that we have the same setup as in the proof of Lemma~\ref{lem:sle4_cm} and let $L_0$ denote the conformal modulus of $(\BH ; 0,z,1,\infty)$.  As there, we condition on $z$ and thus assume that $z$, and hence $L_0$, is deterministic, and we retain the notation $a$, $b$, $\varepsilon_0$, $\delta(\cdot)$, $E(a,b,\varepsilon_0)$, $\varepsilon_\ast$, $H_{\eta,\varepsilon}$, $f_\eta$, $\widehat\Theta$ and $g_\eta = \widehat\Theta \circ f_\eta$ introduced there.  We also set $f_\eta(0) \defeq \lim_{\varepsilon \downarrow 0} f_\eta(\varepsilon) = \mathop{\mathrm{CM}}(\BH \setminus \eta; 0,z,1,\infty)$, the limit existing by monotonicity.

Let $S \subseteq [0,L_0]$ be a Borel subset such that $\Leb(S)>0$.  Then,  by applying the Lebesgue density theorem,  we obtain that there exists $t \in S$ $\cap\,(0,L_0)$ such that $\Leb(S \cap (t,t+\varepsilon)) > 0$ for all $\varepsilon>0$.  (Almost every point of $S$ is a Lebesgue density point of $S$, and any density point $t$ has this property.)  Moreover,  it is easy to see that conditionally on $x,y$,  and $z$,  it is a.s.\ the case that there exists a simple curve $\gamma$ in $\BH$ connecting $x$ to $y$ with $\gamma \cap \partial\BH = \{x,y\}$ such that $\mathop{\mathrm{CM}}(\BH \setminus \gamma; 0, z, 1, \infty) = t$ and $\gamma$ depends only on $x,y,z$,  and $t$; indeed, as $\gamma$ is deformed from the segment $[x,y]$ toward the marked arc $[z,1]$ the modulus of $(\BH \setminus \gamma; 0,z,1,\infty)$ decreases continuously from $L_0$ to $0$, so that the value $t \in (0,L_0)$ is attained.  Since $\gamma$ is a compact subset of $\overline\BH$ meeting $\partial \BH$ only at $x$ and $y$, it lies at positive distance from $a$ and from $b$, so that there exists $\varepsilon_1 \in (0,\varepsilon_0)$ such that $\gamma \cap (\overline{B_{\varepsilon_1}(a)} \cup \overline{B_{\delta(\varepsilon_1)}(b)}) = \emptyset$.  Write $f_\gamma(\varepsilon) \defeq \mathop{\mathrm{CM}}(\BH \setminus (\gamma \cup B_\varepsilon(a)); 0,z,1,\infty)$ for $\varepsilon \in [0,\varepsilon_1]$ and note that the proof of the strict monotonicity of $f_\eta$ given in the proof of Lemma~\ref{lem:sle4_cm} used nothing about the law of $\eta$ and applies verbatim to the deterministic curve $\gamma$; hence $f_\gamma(\varepsilon_1) > f_\gamma(0) = t$.  

Since $\mathrm{Leb}(S \cap (0, f_{\gamma}(\varepsilon_1))) > 0$,  there exists $\delta>0$ such that $\mathrm{Leb}(S \cap (t+\delta ,  f_{\gamma}(\varepsilon_1))) > 0$.  Moreover,   \cite[Lemma~A.1]{CoInCLERiemSph} implies that for each $\widetilde{\delta} \in (0,1)$,  with positive probability,  conditionally on $x,y$,  and $z$,  we have that $\eta$ stays in the $2\widetilde{\delta}$-neighborhood of $\gamma$.  We let $A = A(a,b,\varepsilon_0,\varepsilon_1,\widetilde{\delta})$ denote the event that this occurs together with the event $E(a,b,\varepsilon_0)$,  so that $\BP[A] > 0$.  Note that we can choose $\widetilde{\delta} \in (0,1)$ sufficiently small such that $f_{\eta}(0) < t+\delta < f_{\gamma}(\varepsilon_1) - \delta < f_{\eta}(\varepsilon_1)$.  Indeed,  on $A$ the unbounded connected component of $\BH \setminus B_{2\widetilde{\delta}}(\gamma)$ is contained in the unbounded connected component of $\BH \setminus \eta$ so that the modulus of the conformal rectangle $(\BH \setminus B_{2\widetilde{\delta}}(\gamma) ; 0,z,1,\infty)$ is at least the modulus of $(\BH \setminus \eta ; 0,z,1,\infty)$.  But the former converges to $t = f_{\gamma}(0)$ as $\widetilde{\delta} \downarrow 0$ and so we can choose $\widetilde{\delta}>0$ sufficiently small such that $f_{\eta}(0) < t + \delta$ on $A$.  

Similarly,  we have that $f_{\eta}(\varepsilon_1) \to f_{\gamma}(\varepsilon_1)$ as $\widetilde{\delta} \downarrow 0$,  so that $f_{\eta}(\varepsilon_1) > f_{\gamma}(\varepsilon_1) - \delta$ for $\widetilde{\delta}$ small enough.  In particular,  since $f_{\eta}$ is continuous and strictly increasing,  we can choose $\widetilde{\delta} \in (0,1)$ small enough such that $(t+\delta ,  f_{\gamma}(\varepsilon_1) -\delta) \subseteq f_{\eta}((0,\varepsilon_1))$ and $\mathrm{Leb}(S \cap (t+\delta,f_{\gamma}(\varepsilon_1) - \delta)) > 0$.

\stepn{step:sle4-cm-rev-reduction-eq-sle4}{Reduction to~\eqref{eq:sle4_cm_reverse_proof_0}}
Next,  we sample $\varepsilon$ uniformly from $(0,\varepsilon_1)$ and independently from $\eta$,  and note that
\begin{align*}
&\BP\lbrack\{\mathop{\mathrm{CM}}(\BH \setminus (\eta \cup B_\varepsilon(a)); 0, z, 1, \infty) \in S \cap (t+\delta,f_{\gamma}(\varepsilon_1)-\delta) \} \cap A\rbrack \notag \\
&= \BE\lbrack\BP\lbrack f_\eta(\varepsilon) \in S \cap (t+\delta,f_{\gamma}(\varepsilon_1)-\delta) \mid \eta \rbrack \one_{A} \rbrack = \frac{1}{\varepsilon_1}\BE\lbrack\Leb(f_\eta^{-1}(S \cap (t+\delta,f_{\gamma}(\varepsilon_1)-\delta))) \one_{A} \rbrack,
\end{align*}

It suffices to show that the left-hand side above is positive.  Indeed, exactly as in the proof of Lemma~\ref{lem:sle4_cm} (with $\widetilde\varepsilon$ there taken uniform on $(0,\varepsilon_1)$), \Cref{lem:conformal-restriction} shows that the law of $\eta^\prime = \phi_\varepsilon(\eta)$ conditionally on $y,y,z,$ and on $A$ is mutually absolutely continuous with respect to the law of $\eta$ conditionally on $x,y,z$, and on $A$, and, as was checked there, the conformal modulus of $(\BH \setminus \eta^\prime; 0,z,1,\infty)$ equals that of $(\BH \setminus (\eta \cup B_\varepsilon(a)); 0,z,1,\infty)$.  Hence the positivity of the left-hand side above forces $\BP\lbrack\{\mathop{\mathrm{CM}}(\BH \setminus \eta^\prime; 0,z,1,\infty) \in S\} \cap A\rbrack > 0$ and therefore, by mutual absolute continuity, $\BP\lbrack\mathop{\mathrm{CM}}(\BH \setminus \eta; 0,z,1,\infty) \in S\rbrack > 0$, which is what has to be proved.
Therefore,  it suffices to show that
\begin{align}\label{eq:sle4_cm_reverse_proof_0}
\Leb(f_{\eta}^{-1}(S \cap (t+\delta,f_{\gamma}(\varepsilon_1)-\delta))) >0 \quad \text{a.s.} \quad \text{on} \quad A.
\end{align}

\stepn{step:sle4-cm-rev-proof-eq-sle4}{Proof of~\eqref{eq:sle4_cm_reverse_proof_0}}
To show the latter,  write $S_1 \defeq S \cap (t+\delta,f_{\gamma}(\varepsilon_1)-\delta)$ and recall that the function $g_{\eta}$ defined in the proof of Lemma~\ref{lem:sle4_cm} is strictly decreasing and absolutely continuous when restricted to compact subintervals of $(0,\varepsilon_\ast \wedge \varepsilon_0)$.  Hence,  by the change of variables formula for absolutely continuous strictly monotone functions,
 it holds that
\begin{align*}
\Leb(g_{\eta}(f_{\eta}^{-1}(S_1))) = \Leb(\widehat\Theta(S_1)) = -\int_{f_{\eta}^{-1}(S_1)} g_{\eta}'(x) \, \rd x,
\end{align*}
where the first equality uses that $f_\eta(f_\eta^{-1}(S_1)) = S_1$, which holds because $S_1 \subseteq (t+\delta,f_{\gamma}(\varepsilon_1)-\delta) \subseteq f_\eta((0,\varepsilon_1))$, together with $g_\eta = \widehat\Theta \circ f_\eta$.
Similarly,  we have already shown in the proof of Lemma~\ref{lem:sle4_cm} that $\widehat\Theta$ is strictly decreasing and absolutely continuous when restricted to compact intervals of $(0,\infty)$,  and so we have that
\begin{align*}
\Leb(\widehat\Theta(S_1)) = -\int_{S_1} \widehat\Theta'(x) \, \rd x.
\end{align*}
Note that we have already shown in the proof of Lemma~\ref{lem:sle4_cm} that $-\widehat\Theta'(x)>0$ for Lebesgue almost every $x>0$,  which combined with the fact that $\Leb(S_1) >0$,  gives that $\Leb(\widehat\Theta(S_1)) >0$.  Hence, since the integral of any function over a Lebesgue null set vanishes, the identity $\Leb(g_\eta(f_\eta^{-1}(S_1))) = -\int_{f_\eta^{-1}(S_1)} g_\eta^\prime(x)\, \rd x > 0$ forces  $\Leb(f_{\eta}^{-1}(S_1))>0$ a.s.\ on $A$.  This proves~\eqref{eq:sle4_cm_reverse_proof_0} and hence completes the proof of the lemma since $S \subseteq [0,L_0]$ was arbitrary.
\end{proof}

\begin{proof}[Proof of Lemma~\ref{lem:law_cont_pre}]
\stepn{step:law-cont-pre-setup-events}{Setup and the events $E_r(\varepsilon)$}  Suppose by way of contradiction that there exists a deterministic subset $S \subseteq \BR$ of Lebesgue measure zero such that
\begin{equation}\label{eq:law_cont_pre_proof_0}
\Leb\left(\left\{r > 0 : \BP\!\left\lbrack D_{\Gamma_{V_r}}^{V_r}(T_r, B_r) \in S\right\rbrack > 0\right\}\right) > 0, 
\end{equation}
where $\Leb(\bullet)$ denotes the Lebesgue measure. By replacing $S$ with $S \cap [0, N]$ for some sufficiently large $N > 0$, we may assume without loss of generality that $S \subseteq [0, N]$. (The sets $\{r>0 : \BP\lbrack D_{\Gamma_{V_r}}^{V_r}(T_r,B_r) \in S \cap [0,N]\rbrack>0\}$ increase to the set in~\eqref{eq:law_cont_pre_proof_0} as $N \to \infty$.  We may also assume that $S$ is Borel.  All sets of the form $\{r>0 : \BP\lbrack \bullet \rbrack>0\}$ below are Lebesgue measurable, since the laws of the events in question depend measurably on $r$ by conformal invariance.)

For each $r > 0$, write $L_r$ (resp.\ $R_r$) for the left (resp.\ right) side of $V_r$. For $\varepsilon \in (0,r/2)$ we write $B_\varepsilon(L_r) \defeq \{z \in V_r : \Re(z)<\varepsilon\}$ and $B_\varepsilon(R_r) \defeq \{z \in V_r : \Re(z)>r-\varepsilon\}$ for the two (disjoint) open collars of width $\varepsilon$ along the left and the right side of $V_r$, and $\partial^{\mathrm{in}} B_\varepsilon(L_r) \defeq \{z \in V_r : \Re(z) = \varepsilon\}$ for the inner boundary of the first, and similarly for the second.  For $r, \varepsilon > 0$, write
\begin{equation*}
B_\varepsilon^\star(L_r) \defeq B_\varepsilon(L_r) \setminus \overline{\bigcup_{\SCL \in \Gamma_{V_r} : \SCL \not\subseteq B_\varepsilon(L_r)} \mathop{\mathrm{int}}(\SCL)} \quad \text{and} \quad B_\varepsilon^\star(R_r) \defeq B_\varepsilon(R_r) \setminus \overline{\bigcup_{\SCL \in \Gamma_{V_r} : \SCL \not\subseteq B_\varepsilon(R_r)} \mathop{\mathrm{int}}(\SCL)}; 
\end{equation*}
write $E_r(\varepsilon)$ for the event that the following are true:
\begin{enumerate}
\item\label{it:law_cont_pre_proof_0} $D_{\Gamma_{V_r}}^{V_r}(T_r, B_r) \in S$. 
\item\label{it:law_cont_pre_proof_1} There is a $D_{\Gamma_{V_r}}^{V_r}$-geodesic connecting $T_r$ and $B_r$ that is contained in $V_r \setminus (B_\varepsilon(L_r) \cup B_\varepsilon(R_r))$. 
\item\label{it:law_cont_pre_proof_2} The $D_{\Gamma_{V_r}}^{V_r}$-distance across the left and right sides of $B_\varepsilon^\star(L_r)$ (resp.\ $B_\varepsilon^\star(R_r)$) is greater than $N$.
\end{enumerate}
In condition~\eqref{it:law_cont_pre_proof_2}, the ``$D_{\Gamma_{V_r}}^{V_r}$-distance across the left and right sides of $B_\varepsilon^\star(L_r)$'' means the internal distance $D_{\Gamma_{V_r}}^{V_r}(L_r, \partial^{\mathrm{in}} B_\varepsilon(L_r); \overline{B_\varepsilon(L_r)})$, with the convention $\inf \emptyset = \infty$ of \Cref{subsec:setup}.  Any admissible path in $\overline{V_r}$ joining $\partial^{\mathrm{in}} B_\varepsilon(L_r)$ to $L_r$ contains an admissible path in $\overline{B_\varepsilon(L_r)}$ joining the two sides, and no single loop of $\Gamma_{V_r}$ joins them, since no loop meets $L_r$.

\stepn{step:law-cont-pre-proof-eq-law}{Proof of~\eqref{eq:law_cont_pre_proof_1}}

\substepn{step:lcp-claim}{The claim and the choice of $\varepsilon$} We \emph{claim} that (under the assumption that~\eqref{eq:law_cont_pre_proof_0} holds) there exists $\varepsilon > 0$ such that
\begin{equation}\label{eq:law_cont_pre_proof_1}
\Leb(\{r > 0 : \BP\lbrack E_r(\varepsilon)\rbrack > 0\}) > 0.
\end{equation}
By \Cref{cor:geodesic_avoid_boundary}, a.s., any $D_{\Gamma_{V_r}}^{V_r}$-geodesic connecting $T_r$ and $B_r$ does not touch $L_r$ and $R_r$. More precisely, \Cref{cor:geodesic_avoid_boundary} (applied after conformally mapping $V_r$ onto $\BD$) gives that such a geodesic meets $\partial V_r$ only at its two endpoints, which lie on $T_r$ and on $B_r$ and are a.s.\ not corners of $V_r$; hence the geodesic is at positive distance from $L_r \cup R_r$. This implies that
\begin{equation*}
\bigcup_{\varepsilon > 0} \left\{r > 0 : \BP\lbrack\text{conditions~\eqref{it:law_cont_pre_proof_0} and~\eqref{it:law_cont_pre_proof_1} hold}\rbrack > 0\right\} = \left\{r > 0 : \BP\!\left\lbrack D_{\Gamma_{V_r}}^{V_r}(T_r, B_r) \in S\right\rbrack > 0\right\}, 
\end{equation*}
which implies that there exists $\varepsilon > 0$ such that
\begin{equation}\label{eq:law_cont_pre_proof_2}
\Leb(\{r > 0 : \BP\lbrack\text{conditions~\eqref{it:law_cont_pre_proof_0} and~\eqref{it:law_cont_pre_proof_1} hold}\rbrack > 0\}) > 0.
\end{equation}
(Shrinking $\varepsilon$, we may and do also assume that the set in~\eqref{eq:law_cont_pre_proof_2} intersects $(2\varepsilon,\infty)$ in a set of positive Lebesgue measure, so that the two collars are disjoint for every $r$ under consideration.)
It follows from \Cref{lem:dist_across_rect_arb_high} combined with conformal invariance and the locality property that
\begin{multline}\label{eq:law_cont_pre_proof_3}
\BP\lbrack\text{condition~\eqref{it:law_cont_pre_proof_2} holds} \mid \{\SCL \in \Gamma_{V_r} : \SCL \not\subseteq B_\varepsilon(L_r) \cup B_\varepsilon(R_r)\}\rbrack \\
= \BP\lbrack\text{condition~\eqref{it:law_cont_pre_proof_2} holds} \mid B_\varepsilon^\star(L_r), \ B_\varepsilon^\star(R_r)\rbrack > 0 \quad \text{a.s.}
\end{multline}
Indeed, $B_\varepsilon(L_r)$ and $B_\varepsilon(R_r)$ are deterministic simply connected subdomains of $V_r$, and $B_\varepsilon^\star(L_r)$ and $B_\varepsilon^\star(R_r)$ are the sets attached to them by~\eqref{eq:def-V-star}, measurably with respect to the conditioning $\sigma$-algebra; the first equality is then given by Axiom~\eqref{it:axiom_locality} (locality).  For the positivity: if no admissible path in $\overline{B_\varepsilon(L_r)}$ joins the two sides, then the distance is $\infty$; otherwise every such path crosses a connected component $R$ of $B_\varepsilon^\star(L_r)$ separating them, and, by locality and Axiom~\eqref{it:axiom_conformal_invariance} (conformal invariance), the crossing distance of $R$ is conditionally distributed as $D_{\Gamma_{V_a}}^{V_a}(T_a,B_a)$, where $a>0$ is the conformal modulus of $R$, and this crossing distance exceeds $N$ with positive probability by \Cref{lem:dist_across_rect_arb_high}.  The same applies to $B_\varepsilon^\star(R_r)$, and the two are conditionally independent.

\substepn{step:lcp-fkg}{Conclusion via the FKG inequality} Note that conditionally on $\{\SCL \in \Gamma_{V_r} : \SCL \not \subseteq B_{\varepsilon}(L_r) \cup B_{\varepsilon}(R_r)\}$,  if we have two metrics $d,d'$ on $B_{\varepsilon}^\star(L_r) \cup B_{\varepsilon}^\star(R_r)$ such that $d \geq d'$,  it holds that $\one_{\eqref{it:law_cont_pre_proof_0}, \eqref{it:law_cont_pre_proof_1}}(d) \geq \one_{\eqref{it:law_cont_pre_proof_0}, \eqref{it:law_cont_pre_proof_1}}(d')$ and
$\one_{\eqref{it:law_cont_pre_proof_2}}(d) \geq \one_{\eqref{it:law_cont_pre_proof_2}}(d')$.  Indeed, write $V_r^{\mathrm{mid}} \defeq V_r \setminus (\overline{B_\varepsilon(L_r)} \cup \overline{B_\varepsilon(R_r)})$: the internal distance $D_{\Gamma_{V_r}}^{V_r}(T_r,B_r; V_r^{\mathrm{mid}})$ is measurable with respect to the conditioning $\sigma$-algebra by Axiom~\eqref{it:axiom_locality}, condition~\eqref{it:law_cont_pre_proof_1} states that $D_{\Gamma_{V_r}}^{V_r}(T_r,B_r)$ equals it, and the latter distance can only decrease when loops are added to the collars, so that condition~\eqref{it:law_cont_pre_proof_1} is decreasing and, on it, condition~\eqref{it:law_cont_pre_proof_0} is conditionally deterministic.  The monotonicity of $\one_{\eqref{it:law_cont_pre_proof_2}}$ is immediate.  Therefore,  \Cref{lem:FKG-CLE} implies that
\begin{align*}
\BP\lbrack\eqref{it:law_cont_pre_proof_0}, \ \eqref{it:law_cont_pre_proof_1}, \ \eqref{it:law_cont_pre_proof_2} \mid \cdots\rbrack \geq \BP\lbrack\eqref{it:law_cont_pre_proof_0}, \ \eqref{it:law_cont_pre_proof_1} \mid \cdots\rbrack  \BP\lbrack \eqref{it:law_cont_pre_proof_2} \mid \cdots\rbrack,
\end{align*}
where $\cdots \defeq \{\SCL \in \Gamma_{V_r} : \SCL \not\subseteq B_\varepsilon(L_r) \cup B_\varepsilon(R_r)\}$.  Thus,  it follows that
\begin{align*}
\BE\lbrack\BP\lbrack\eqref{it:law_cont_pre_proof_0}, \ \eqref{it:law_cont_pre_proof_1}, \ \eqref{it:law_cont_pre_proof_2} \mid \cdots\rbrack\rbrack \geq \BE\lbrack\BP\lbrack\eqref{it:law_cont_pre_proof_0}, \ \eqref{it:law_cont_pre_proof_1}  \mid \cdots\rbrack  \BP\lbrack \eqref{it:law_cont_pre_proof_2} \mid \cdots\rbrack\rbrack.
\end{align*}
Recall that by \eqref{eq:law_cont_pre_proof_3},  we have that $\BP\lbrack \eqref{it:law_cont_pre_proof_2} \mid \cdots\rbrack >0$ a.s.   Hence, if $\BP\lbrack\eqref{it:law_cont_pre_proof_0}, \ \eqref{it:law_cont_pre_proof_1}\rbrack >0$,  we obtain that
\begin{align*}
 \BE\lbrack\BP\lbrack\eqref{it:law_cont_pre_proof_0}, \ \eqref{it:law_cont_pre_proof_1}  \mid \cdots\rbrack  \BP\lbrack \eqref{it:law_cont_pre_proof_2} \mid \cdots\rbrack\rbrack>0,
\end{align*}
and so $\BP\lbrack E_r(\varepsilon)\rbrack>0$.  This proves~\eqref{eq:law_cont_pre_proof_1} by combining with~\eqref{eq:law_cont_pre_proof_2}.

\stepn{step:law-cont-pre-random-conformal-rectangle}{A random conformal rectangle in $\BD$}  Now we consider a non-nested $\CLE_4$ in the open unit disk $\BD$. Fix a line segment $I \colon [0, 1] \to \overline\BD$ from $\partial\BD$ to the origin. Write $\tau$ for the first time at which $I$ hits $\SCL(0)$. Write $V$ for the connected component of
\begin{equation*}
\BD \setminus \overline{I([0,\tau]) \cup {}\bigcup_{\SCL \in \Gamma_\BD : \SCL \cap I([0, \tau]) \neq \emptyset} \SCL}
\end{equation*}
whose boundary contains $\partial\BD$. (Removing the slit $I([0,\tau])$ along with the loops meeting it makes $V$ simply connected and produces two prime ends at each of $I(0)$ and $I(\tau)$, denoted $I(0)^\pm$ and $I(\tau)^\pm$ below.) Write $T$ (resp.\ $B$; $L$; $R$) for the top (resp.\ bottom; left; right) sides of $V$ (i.e., $T = \SCL(0) \cap \partial V$ and $B = \partial\BD \cap \partial V$), and $L$ and $R$ for the two sides of the slit, each of which consists of the corresponding side of $I([0,\tau])$ together with the arcs of the loops meeting it which lie on $\partial V$; thus $(V; I(\tau)^+, I(\tau)^-, I(0)^-, I(0)^+)$ is a conformal rectangle with sides $T$, $L$, $B$, $R$ in cyclic order.
Write $\phi_V$ for the conformal map from $V$ onto $V_{r}$, where $r = r(V)>0$ is such that the conformal rectangles $(V; I(\tau)^+, I(\tau)^-, I(0)^-, I(0)^+)$ and $V_r$ 
have the same conformal modulus, taking $T$, $L$, $B$, $R$ to $T_r$, $L_r$, $B_r$, $R_r$, respectively.  Set
\begin{align*}
B_\varepsilon(L; d_V) \defeq \phi_V^{-1}(B_\varepsilon(L_r)) \quad \text{and} \quad B_\varepsilon(R; d_V) \defeq \phi_V^{-1}(B_\varepsilon(R_r)),
\end{align*}
and define $B_\varepsilon^\star(L; d_V)$ and $B_\varepsilon^\star(R; d_V)$ from them as in~\eqref{eq:def-V-star}, so that the event $E(\varepsilon)$ defined below is the image of $E_{r(V)}(\varepsilon)$ under $\phi_V^{-1}$.
	Let the event $E(\varepsilon)$ be defined in the same manner as $E_r(\varepsilon)$ but with $V$ in place of $V_r$ and $D_{\Gamma_\BD}^\BD(\bullet, \bullet; V)$ in place of $D_{\Gamma_{V_r}}^{V_r}$, i.e., $E(\varepsilon)$ is the event that the following are true:
	\begin{enumerate}
\item \label{it:distance_in_the_set} $D_{\Gamma_\BD}^\BD(T, B; V) \in S$. 
\item \label{it:geodesic_away_from_left_right_boundaries} There is a $D_{\Gamma_\BD}^\BD(\bullet, \bullet; V)$-geodesic connecting $T$ and $B$ that is contained in $V \setminus (B_\varepsilon(L; d_V) \cup B_\varepsilon(R; d_V))$. 
\item \label{it:distance_between_left_right_big} The $D_{\Gamma_\BD}^\BD(\bullet, \bullet; V)$-distance across the left and right sides of $B_\varepsilon^\star(L; d_V)$ (resp.\ $B_\varepsilon^\star(R; d_V)$) is greater than $N$.
\end{enumerate}

\stepn{step:law-cont-pre-law}{The law of $r(V)$}
Next,  we note that~\eqref{eq:law_cont_pre_proof_1}  implies that there exists $M>1$ such that if $F$ denotes the set of $r \in (1/M,M)$ for which $\BP\lbrack E_r(\varepsilon)\rbrack >0$, then we have that $\Leb(F) > 0$.  We \emph{claim} that $r(V)$ lies in $F$ with positive probability.  Indeed,  first we note that \cite[Proposition~4.1]{CLE} (see also the discussion just after the statement of \cite[Theorem~1.2]{CLE}) implies that the following is true.  Given the loops in $\Gamma_{\BD}$ intersecting $I((0,\tau))$,  the loop $\SCL(0)$ is an $\SLE_4$ loop in the domain obtained by removing from $\BD$ the loops of $\Gamma_{\BD}$ intersecting $I((0,\tau))$,  rooted at $I(\tau)$ and conditioned to surround $0$.  Let $\eta$ be the counterclockwise parameterization of $\SCL(0)$ and $\widetilde{\eta}$ be its time-reversal.  Let $\sigma$ (resp.\ $\widetilde{\sigma}$) be a stopping time for $\eta$ (resp.\ $\widetilde{\eta}$).  By the above discussion,  on the event that $\eta|_{[0,\sigma]}$ and $\widetilde{\eta}|_{[0,\widetilde{\sigma}]}$ have not hooked up,  the remainder of $\SCL(0)$ is a chordal $\SLE_4$ in the domain obtained by removing from $\BD$ the segment $I([0,\tau])$, the loops of $\Gamma_{\BD}$ intersecting $I((0,\tau))$, and $\eta|_{[0,\sigma]}$ and $\widetilde{\eta}|_{[0,\widetilde{\sigma}]}$.  (More precisely, in the connected component of that set whose boundary contains $\partial\BD$; removing the segment makes this component simply connected, with two prime ends at each of $I(0)$ and $I(\tau)$.)

Consider the conformal transformation $\psi$ which maps this domain onto $\BH$ such that $\psi(I(\tau)^+) = 0,  \psi(I(\tau)^-) = \infty,  \psi(I(0)^-) =1$,  and $\psi(I(0)^+) = z$ for some (random) $z \in (0,1)$.  Recall from \cite[Corollary 1.2]{ALS22ExtremalDistance} that for all $\varepsilon \in (0,1)$,  the probability that $\SCL(0) \subseteq B_{\varepsilon}(0)$ is positive.  Therefore,  we can choose $\sigma$ and $\widetilde{\sigma}$ such that the joint conditional law, given $z$, of the pair of images of $\eta(\sigma)$ and $\widetilde{\eta}(\widetilde{\sigma})$ under $\psi$ is absolutely continuous with respect to the Lebesgue measure on $(-\infty,0)^2$ and the probability that the conformal modulus of $(\BH ; 0,z,1,\infty)$ is larger than $M$ is positive.  (The two images lie in $(-\infty,0)$ because $\eta$ and $\widetilde\eta$ are traced from $I(\tau)^+$ and $I(\tau)^-$, which $\psi$ sends to $0$ and to $\infty$, so that the two partially traced strands lie on the boundary arc of $\BH$ from $\infty$ to $0$.)  Thus,  Lemma~\ref{lem:sle4_cm_reverse} implies that conditionally on the loops of $\Gamma_{\BD}$ that intersect $I((0,\tau))$, and hence on $z$,  we have that the following is true.  The restriction of the Lebesgue measure to $[0,M]$ is absolutely continuous with respect to the conditional law of the conformal modulus of the conformal rectangle $(\BH \setminus \psi(\eta \setminus (\eta([0,\sigma]) \cup \widetilde{\eta}([0,\widetilde{\sigma}]))) ; 0,z,1,\infty)$.  The conformal modulus of this rectangle equals that of $(V; I(\tau)^+, I(\tau)^-, I(0)^-, I(0)^+)$, i.e.\ $1/r(V)$; and $\Leb(\{1/r : r \in F\}) > 0$ with $\{1/r : r \in F\} \subseteq (1/M,M) \subseteq [0,M]$.  Therefore,  we obtain that with positive probability,  $r(V)$ lies in $F$.  This completes the proof of the \emph{claim}.

\stepn{step:law-cont-pre-positive-probability-together}{Positive probability of $E(\varepsilon)$ together with a large crossing of $W_\varepsilon$}  Since $D_{\Gamma_\BD}^\BD(\bullet, \bullet; V) = D_{\Gamma_\BD|_V}^V(\bullet, \bullet)$ a.s.\ (cf.~\Cref{def:axioms}, Axiom~\eqref{it:axiom_locality} (locality)) and $\Gamma_\BD|_V$ is conditionally a non-nested $\CLE_4$ in $V$ given $V$, it follows from~\eqref{eq:law_cont_pre_proof_1} and the fact that $r(V)$ lies in the set of $r>0$ for which $\BP\lbrack E_r(\varepsilon)\rbrack>0$ with positive probability, that $\BP\!\left\lbrack E(\varepsilon)\right\rbrack > 0$. (In detail: $E(\varepsilon)$ is the image of $E_{r(V)}(\varepsilon)$ under $\phi_V^{-1}$, so $\BP\lbrack E(\varepsilon) \mid V\rbrack = \varphi(r(V))$ with $\varphi(r) \defeq \BP\lbrack E_r(\varepsilon)\rbrack$, and $\BP\lbrack E(\varepsilon)\rbrack \ge \BE\lbrack \varphi(r(V)) \one_{\{r(V) \in F\}}\rbrack > 0$ since $\varphi>0$ on $F$ and $\BP\lbrack r(V) \in F\rbrack>0$.) Write
\begin{align*}
W_\varepsilon &\defeq \BD \setminus \overline{\mathop{\mathrm{int}}(\SCL(0)) \cup (V \setminus (B_\varepsilon(L; d_V) \cup B_\varepsilon(R; d_V)))}. 
\end{align*}

Combining Axiom~\eqref{it:axiom_locality} (locality) with Lemma~\ref{lem:dist_across_rect_arb_high},  we obtain that
\begin{align*}
\BP\!\left\lbrack D_{\Gamma_\BD}^\BD(\SCL(0) \cap \partial W_\varepsilon, \partial\BD \cap \partial W_\varepsilon; W_\varepsilon) > N \ \middle\vert \ V, \ B_\varepsilon^\star(L; d_V), \ B_\varepsilon^\star(R; d_V)\right\rbrack > 0 \quad \text{a.s.}
\end{align*}
In particular,  we have that
\begin{align*}
\BP\!\left\lbrack D_{\Gamma_\BD}^\BD(\SCL(0) \cap \partial W_\varepsilon, \partial\BD \cap \partial W_\varepsilon; W_\varepsilon) > N \right\rbrack >0.
\end{align*}
Therefore,  since we have already shown that $\BP\lbrack E(\varepsilon) \rbrack >0$,  we obtain by applying Lemma~\ref{lem:FKG-CLE} as in the proof of~\eqref{eq:law_cont_pre_proof_1} that 
\begin{equation}\label{eq:law_cont_pre_proof_5}
\BP\!\left\lbrack E(\varepsilon) \text{ occurs and } D_{\Gamma_\BD}^\BD(\SCL(0) \cap \partial W_\varepsilon, \partial\BD \cap \partial W_\varepsilon; W_\varepsilon) > N\right\rbrack > 0. 
\end{equation}

\stepn{step:law-cont-pre-conclusion}{Conclusion}
Finally,  suppose that the event $E(\varepsilon)$ occurs and that $D_{\Gamma_\BD}^\BD(\SCL(0) \cap \partial W_\varepsilon, \partial\BD \cap \partial W_\varepsilon; W_\varepsilon) > N$.  We will show that any $D_{\Gamma_{\BD}}^{\BD}$-geodesic $P$ connecting $\SCL(0)$ to $\partial \BD$ must be contained in $\overline V$ and hence is a $D_{\Gamma_\BD}^\BD(\bullet, \bullet; V)$-geodesic connecting $T$ and $B$.  We may and do assume that $P$ does not return to $\SCL(0)$ after time $0$, replacing $P$ if necessary by its restriction to the interval after its last visit to $\SCL(0)$.  Since $P$ takes its values in the carpet $\Upsilon_\BD$ and $\BD \setminus \overline V$ is the union of $\mathop{\mathrm{int}}(\SCL(0))$, of the interiors of the loops meeting $I([0,\tau])$, and of the slit $I([0,\tau])$, the path $P$ can leave $\overline V$ only by meeting $T \cup L \cup R$.  Now $P$ does not meet $L$: the $D_{\Gamma_{\BD}}^{\BD}$-length of $P$ is at most $D_{\Gamma_{\BD}}^{\BD}(T, B; V) \le N$ by condition~\eqref{it:distance_in_the_set} (as $S \subseteq [0,N]$), while if $P$ met $L$, a sub-path of $P$ would be an admissible path in $\overline{B_\varepsilon(L;d_V)}$ joining the inner boundary of that collar to $L$ (no single loop of $\Gamma_\BD$ joins them, the loops meeting $L$ being the loops meeting $I([0,\tau])$, whose intersection with $\overline{B_\varepsilon(L;d_V)}$ is contained in $\partial V$), and hence, by condition~\eqref{it:distance_between_left_right_big}, would have $D_{\Gamma_\BD}^\BD$-length greater than $N$.  The same applies to $R$.  Similarly,  $P$ cannot stay entirely in $\overline{W_{\varepsilon}}$ due to the fact that $D_{\Gamma_\BD}^\BD(\SCL(0) \cap \partial W_\varepsilon, \partial\BD \cap \partial W_\varepsilon; W_\varepsilon) > N$.  Therefore,  it follows that $P \subseteq \overline{V}$ and $P$ does not intersect the closure of $\bigcup_{\SCL \in \Gamma_\BD, \ \SCL \cap I([0, \tau)) \neq \emptyset }\SCL$. Hence, $P$ is a $D_{\Gamma_\BD}^\BD(\bullet, \bullet; V)$-geodesic between $T$ and $B$, and consequently $D_{\Gamma_\BD}^\BD(\SCL(0),\partial\BD) = D_{\Gamma_\BD}^\BD(T,B;V) \in S$ by condition~\eqref{it:distance_in_the_set}.  In particular,  we have that
\begin{equation*}
\BP\!\left\lbrack D_{\Gamma_\BD}^\BD(\SCL(0), \partial\BD) \in S\right\rbrack > 0, 
\end{equation*}
in contradiction to the fact that $D_{\Gamma_\BD}^\BD(\SCL(0), \partial\BD)$ is an exponential random variable (by Axiom~\eqref{it:axiom_uniform_exploration} (uniform exploration)), whose law is therefore absolutely continuous with respect to the Lebesgue measure while $\Leb(S)=0$. This completes the proof of the lemma.
\end{proof}

Before proving \Cref{lem:law_cont}, let us state the following lemma which says that loops are attached at a single point when they are first hit by the metric ball from a segment.

\begin{lemma}[{\cite[Lemma~4.9]{kkmt2026cle4_part2}}]\label{lemma intersection with boundary}
	Let $\alpha \subset\partial \BD$ be a segment. Let $\SCL \in \Gamma_\BD$ and let $t=D^\BD_{\Gamma_\BD}(\alpha, \SCL)$. Then, the intersection $\SCL \cap \SCB^-_t(\alpha; D^\BD_{\Gamma_\BD})$ contains only one point.
\end{lemma}

\begin{figure}[ht!]
    \centering
    \includegraphics[width=.8\linewidth]{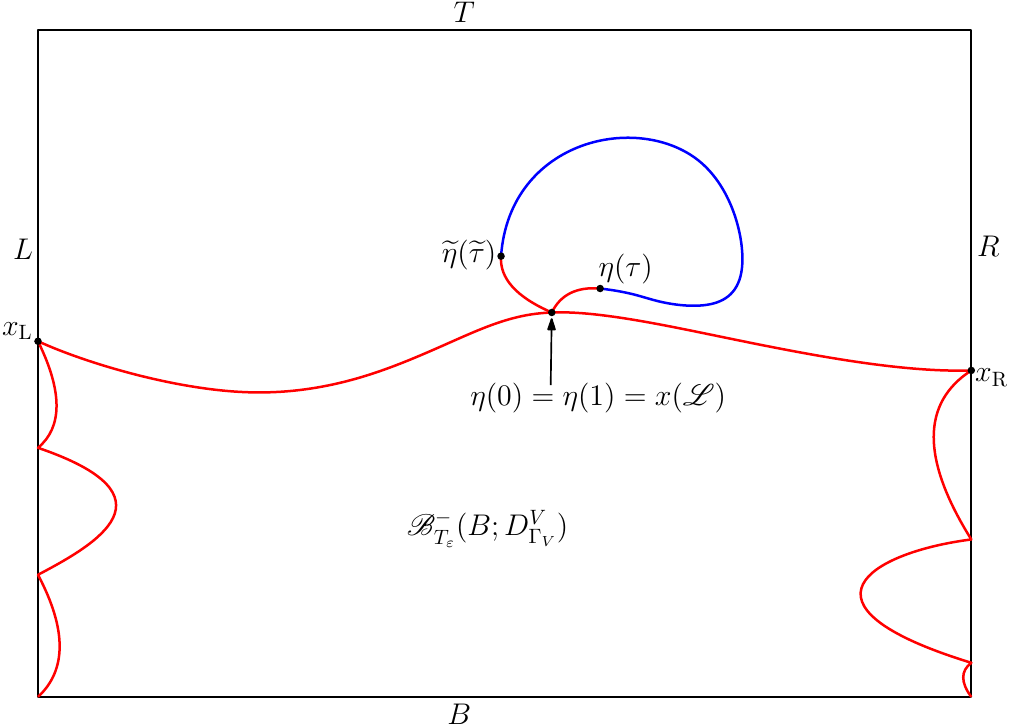}
    \caption{Illustration of the notation in the proof of \Cref{lem:law_cont}.}
    \label{fig:absolute-continuity}
\end{figure}

We conclude this section by proving \Cref{lem:law_cont}.
\begin{proof}[Proof of \Cref{lem:law_cont}]
See \Cref{fig:absolute-continuity} for an illustration.

\stepn{step:law-cont-stopping-time-metric}{A stopping time for the metric exploration}  To lighten notation, write $V \defeq V_r$, $T \defeq T_r$, and $B \defeq B_r$. For each $\SCL \in \Gamma_V$, write $x(\SCL)$ for the unique point of $\SCL \cap \SCB_t^-(B; D_{\Gamma_V}^V)$, which exists thanks to \Cref{lemma intersection with boundary}  (applied to $B$ after conformally mapping $V$ onto $\BD$), where $t \defeq D_{\Gamma_V}^V(\SCL, B)$. Fix $\varepsilon > 0$, which may be chosen to be arbitrarily small. Write 
\begin{equation}\label{eq:law_cont_proof}
T_{\varepsilon} = \inf\left\{t \in (0, D^V_{\Gamma_V}(B, T)) : \text{there exists } \SCL \in \Gamma_V \text{ such that } D_{\Gamma_V}^V(\SCL, B) = t \text{ and } \SCL\not\subseteq \overline{B_\varepsilon(x(\SCL))}
	\right\}, 
\end{equation}
with the convention that $\inf\emptyset = \infty$. On the event that $T_{\varepsilon} < \infty$, write $\SCL \in \Gamma_V$ for the loop described in~\eqref{eq:law_cont_proof}. Let $\eta \colon [0, 1] \to \SCL$ be a counterclockwise parameterization of $\SCL$ with $\eta(0) = \eta(1) = x(\SCL)$ and $\widetilde\eta$ be the time-reversal of $\eta$. Let $\tau \defeq \inf \{t\ge 0: \vert \eta(t)- \eta(0) \vert \ge \varepsilon\}$ and $\widetilde{\tau} \defeq \inf \{t\ge 0: \vert\widetilde{\eta}(t)- \widetilde{\eta}(0) \vert\ge \varepsilon/2 \}$. Note that $\tau$ (resp.\ $\widetilde\tau$) is a stopping time for the filtration generated by 
\begin{equation*}
\left(\SCB_{T_{\varepsilon}}^-(B; D_{\Gamma_V}^V), \eta|_{[0, t]}\right)_{t \ge 0} \quad \left(\text{resp.\ } \left(\SCB_{T_{\varepsilon}}^-(B; D_{\Gamma_V}^V), \widetilde\eta|_{[0, t]}\right)_{t \ge 0}\right).
\end{equation*}
Note that $\diam(\eta([0, \tau])) \ge \varepsilon$ and that $\eta|_{[0, \tau]}$ and $\widetilde\eta|_{[0, \widetilde\tau]}$ have not yet hooked up. Indeed, $\tau<\infty$ since $\SCL \not\subseteq \overline{B_\varepsilon(x(\SCL))}$, while $\eta([0,\tau]) \subseteq \overline{B_\varepsilon(x(\SCL))}$ and $\widetilde\eta([0,\widetilde\tau]) \subseteq \overline{B_{\varepsilon/2}(x(\SCL))}$; if the two strands had hooked up, we would have $\SCL = \eta([0,\tau]) \cup \widetilde\eta([0,\widetilde\tau]) \subseteq \overline{B_\varepsilon(x(\SCL))}$.

\stepn{step:law-cont-conformal-modulus-unexplored}{The conformal modulus of the unexplored region}  Then, given $\SCB_{T_{\varepsilon}}^-(B; D_{\Gamma_V}^V)$, $\eta|_{[0, \tau]}$, $\widetilde\eta|_{[0, \widetilde\tau]}$, 
the conditional law of the remaining strand of $\eta$ is that of a chordal SLE$_4$ curve in $V \setminus (\SCB_{T_{\varepsilon}}^-(B; D_{\Gamma_V}^V) \cup \eta([0, \tau]) \cup \widetilde\eta([0, \widetilde\tau]))$ from $\eta(\tau)$ to $\widetilde\eta(\widetilde\tau)$. Let $U$ be the connected component of $V \setminus (\SCB_{T_{\varepsilon}}^-(B; D_{\Gamma_V}^V) \cup \eta([0, \tau]) \cup \widetilde\eta([0, \widetilde\tau]))$ whose boundary contains $T$. Note that $U$ is simply connected since $\SCB_{T_{\varepsilon}}^-(B; D_{\Gamma_V}^V) \cup \eta([0, \tau]) \cup \widetilde\eta([0, \widetilde\tau]) \cup \partial V$ is connected. Write $L$ (resp.\ $R$) for the left (resp.\ right) boundary of $V$. Write $x_\rL$ (resp.\ $x_\rR$) for the uppermost point of $L \cap \SCB_{T_{\varepsilon}}^-(B; D_{\Gamma_V}^V)$ (resp.\ $R \cap \SCB_{T_{\varepsilon}}^-(B; D_{\Gamma_V}^V)$). Write $\phi$ for the conformal mapping from $U$ onto $\BH$ such that $\phi(x_\rR) = 0$, $\phi(\text{upper-left vertex of } V) = 1$, and $\phi(x_\rL) = \infty$. Write $z \defeq \phi(\text{upper-right vertex of } V) \in (0,1)$ (the points $x_\rR$, the upper-right vertex of $V$, the upper-left vertex of $V$, and $x_\rL$ occur in this counterclockwise order along $\partial U$). We may choose $\tau$ and $\widetilde\tau$ so that the joint conditional law of the pair $(\phi(\eta(\tau)), \phi(\widetilde\eta(\widetilde\tau)))$ given $\SCB_{T_{\varepsilon}}^-(B; D_{\Gamma_V}^V)$, $x(\SCL)$ and $z$ is absolutely continuous with respect to the Lebesgue measure on $(-\infty,0)^2$.
Thus, it follows from \Cref{lem:sle4_cm} that the conditional law of the conformal modulus of the conformal rectangle $(U \setminus \SCB_{T_{\varepsilon}}(B;D^V_{\Gamma_V}); x_\rR, \text{upper-right vertex of } V, \text{upper-left vertex of } V, x_\rL)$ given $\SCB_{T_{\varepsilon}}^-(B; D_{\Gamma_V}^V)$ and $x(\SCL)$ is a.s.\ absolutely continuous with respect to the Lebesgue measure. Combining this with \Cref{lem:law_cont_pre}, we obtain that the conditional law of the $D_{\Gamma_V}^V(\bullet, \bullet; U \setminus \SCB_{T_{\varepsilon}}(B; D_{\Gamma_V}^V))$-distance across the top and bottom sides of the conformal rectangle $U \setminus \SCB_{T_{\varepsilon}}(B; D_{\Gamma_V}^V)$ given $\SCB_{T_{\varepsilon}}^-(B; D_{\Gamma_V}^V)$ and $x(\SCL)$ is a.s.\ absolutely continuous with respect to the Lebesgue measure.

\stepn{step:law-cont-conclusion}{Conclusion}  Thus, for each deterministic subset $S \subseteq \BR$ of Lebesgue measure zero, we conclude that
\begin{align*}
&\BP\!\left\lbrack \{D_{\Gamma_V}^V(T, B) \in S\} \cap \{T_{\varepsilon} < \infty\}\right\rbrack\\
&\qquad= \BP\!\left\lbrack \{T_{\varepsilon} + D_{\Gamma_V}^V(\text{across } V \setminus \SCB_{T_{\varepsilon}}(B; D_{\Gamma_V}^V)) \in S\}\right.\\
&\qquad\qquad\left.\cap \{T_{\varepsilon} < \infty\} \right\rbrack \\
&\qquad\le \BE\!\left\lbrack\one_{\{T_\varepsilon<\infty\}}\BP\lbrack D_{\Gamma_V}^V(\text{across } V \setminus \SCB_{T_{\varepsilon}}(B; D_{\Gamma_V}^V)) \in S - T_{\varepsilon}\right.\\
&\qquad\qquad\qquad\left.\mid \SCB_{T_{\varepsilon}}^-(B; D_{\Gamma_V}^V), \ x(\SCL) \rbrack\right\rbrack \\
&\qquad= 0.
\end{align*}
(Here, $S - T_{\varepsilon} \defeq \{s - T_{\varepsilon} : s \in S\}$.) (The first identity holds since, on $\{T_\varepsilon<\infty\}$, the distance from $B$ to $T$ equals $T_\varepsilon$ plus the distance across the unexplored region $V \setminus \SCB_{T_\varepsilon}(B; D_{\Gamma_V}^V)$, meaning the connected component of that set whose boundary contains $T$, which equals $U \setminus \SCB_{T_\varepsilon}(B; D_{\Gamma_V}^V)$; the final inequality holds since $T_\varepsilon$, and hence the null set $S-T_\varepsilon$, is measurable with respect to the conditioning $\sigma$-algebra.) This completes the proof of \Cref{lem:law_cont} since $\BP[T_{\varepsilon} < \infty] \to 1$ as $\varepsilon \to 0$, so that letting $\varepsilon \to 0$ gives $\BP\lbrack D_{\Gamma_V}^V(T,B) \in S\rbrack = 0$.
\end{proof}

\begin{appendix}

\section{Nested $\CLE_4$ four-arm exponents}
\label{sec:cle_4_four_arm_exponents}

\subsection{Statement of the main result and outline of the proof strategy.}
\label{subsec:outline_of_proof}
The main result that we are going to prove in this appendix is the following estimate on four-arm events for nested $\CLE_4$. Throughout this appendix we set $C_\varepsilon \defeq B_{1/2}(0) \setminus \overline{B_\varepsilon(0)}$. This is the first main input in the proof of Lemma~\ref{prop:bounded_cle_crossings}.

\begin{proposition}\label{prop:nested_cle_4_two_arm_exponent}
Fix $\nu \in (0,1)$.  Then,  there exists $\varepsilon_0 \in (0,1)$ depending only on $\nu$ such that the following holds for all $\varepsilon \in (0,\varepsilon_0)$.  Off an event with probability at most $\varepsilon^{2(1-\nu)}$,  we have that there are at most two crossings of the annulus $C_\varepsilon$ made by loops in $\overline{\Gamma}$;  in particular,  at most one loop of $\overline{\Gamma}$ crosses $C_\varepsilon$.
\end{proposition}

We are also going to prove the following estimate for the probability of the event that there are at least six crossings of $C_\varepsilon$ made by loops in $\overline{\Gamma}$;  recall that the number of such crossings is necessarily even.  This is the second main input in the proof of Lemma~\ref{prop:bounded_cle_crossings}.

\begin{proposition}\label{prop:nested_cle_4_three_arm_exponent}
There exist universal constants $b_0,  \varepsilon_0 \in (0,1)$ such that the following is true for all $\varepsilon \in (0,\varepsilon_0)$.  Off an event with probability at most $\varepsilon^{2+b_0}$,  we have that there are at most four crossings of $C_\varepsilon$ made by loops in $\overline{\Gamma}$;  in particular,  at most two loops of $\overline{\Gamma}$ cross $C_\varepsilon$.
\end{proposition}

\begin{remark}\label{rem:general_outer_radius}
Propositions~\ref{prop:nested_cle_4_two_arm_exponent} and~\ref{prop:nested_cle_4_three_arm_exponent} hold,  with the same proofs,  with $C_\varepsilon$ replaced by $B_{r}(0) \setminus \overline{B_{\varepsilon}(0)}$ for any fixed $r \in (0,1/2]$,  with $\varepsilon_0$ then depending also on $r$.  Indeed,  the only role of the outer radius $1/2$ is to accommodate the annuli $A_k = B_{2^{-3k-1}}(0) \setminus \overline{B_{2^{-3k-2}}(0)}$,  $k \geq 1$,  used below,  and these already lie inside $B_{1/16}(0)$;  for smaller $r$ one simply lets $k$ range over the integers in $[k_r ,  a\log(1/\varepsilon)]$,  where $k_r$ is the least integer with $2^{-3k_r-1} \leq r$ and where $a \in (0,1/100)$ is a small constant appearing in the proofs of Propositions~\ref{prop:nested_cle_4_two_arm_exponent} and~\ref{prop:nested_cle_4_three_arm_exponent}.
\end{remark}

\begin{proof}[Proof of Lemma~\ref{prop:bounded_cle_crossings}.]
\stepn{step:bdd-cross-case}{The case $r_2 \geq d/8$}  Let $\phi$ be the M\"obius transformation of $\BD$ with $\phi(z) = 0$ and $\phi'(z)>0$,  and recall that the law of $\overline\Gamma$ is invariant under $\phi$.  Since $r_2 < d/2 \leq (1-\lvert z\rvert)/2$,  the Koebe distortion theorem gives universal constants $c_1,c_2 \in (0,\infty)$ such that $B_{c_1 r_2 / (1-\lvert z\rvert^2)}(0) \subseteq \phi(B_{r_2}(z))$ and $\phi(\overline{B_{r_1}(z)}) \subseteq \overline{B_{c_2 r_1/(1-\lvert z\rvert^2)}(0)}$. More precisely, by \cite[Corollary 3.25]{lawler2008conformally}, one can take $c_1=1/16$ and $c_2=16$.  As $r_2 \geq d/8$ and $1-\lvert z\rvert^2 \in [d,2]$,  the radius $r \defeq c_1 r_2/(1-\lvert z\rvert^2)$ is bounded below by $c_1 d/16$, which only depends on $d$;  moreover, $r \leq 1/2$ since $c_1 \leq 1$ and $r_2 < d/2 \leq (1-\lvert z\rvert)/2 \leq (1-\lvert z\rvert^2)/2$.  Also $c_2 r_1/(1-\lvert z\rvert^2) = (c_2/c_1) r \cdot (r_1/r_2)$.  Every crossing of $B_{r_2}(z) \setminus \overline{B_{r_1}(z)}$ therefore yields,  under $\phi$,  a crossing of the concentric annulus with outer radius $r$ and inner radius $(c_2/c_1) r (r_1/r_2)$,  and hence also a crossing of the annulus with the \emph{fixed} outer radius $r_* \defeq \min\{ c_1 d / 16 ,  1/2\} \leq r$ and the same inner radius.  Both claims now follow from Remark~\ref{rem:general_outer_radius},  applied to Propositions~\ref{prop:nested_cle_4_two_arm_exponent} and~\ref{prop:nested_cle_4_three_arm_exponent} with the radius $r_*$,  upon taking $\nu$ there smaller than $\epsilon / 2$ and $\delta$ small enough to absorb the constants.

\stepn{step:bdd-cross-case-2}{The case $r_2 < d/8$}  Here conformal invariance alone does not suffice,  since the image $\phi(B_{r_2}(z))$ is then much smaller than $\BD$,  while the bound we are proving is in terms of the ratio $r_1/r_2$ alone.  We therefore insert a comparison of scales.  Write $b_0$ for the universal constant in \Cref{prop:nested_cle_4_three_arm_exponent},  set $\epsilon_1 \defeq \min\{\epsilon , b_0\}/8$,  and let $K_1$ (resp.\ $K_2$) be the smallest (resp.\ largest) integer $k$ with $2^{-k} \leq 2r_2$ (resp.\ $2^{-k} \geq r_2 (r_1/r_2)^{\epsilon_1}$),  so that $K_1 > k_0$ with $k_0 = \lfloor \log(1/(2r_2)) / \log 2 \rfloor$, and $K_2 - K_1 \geq \tfrac{\epsilon_1}{2 \log 2} \log(r_2/r_1)$ provided $r_1/r_2$ is small enough.  We use \Cref{lem:M_good_event_general} in the following form,  which follows from the same proof: given $\alpha \in (0,\infty)$ and $\beta_0 \in (0,1)$,  there are constants $c_0,M>0$ depending only on $F,\alpha$ and $\beta_0$ such that,  for all $L_2 > L_1 \geq k_0$,  off an event of probability at most $c_0 e^{-\alpha(L_2-L_1)}$ there are at least $\beta_0(L_2-L_1)$ values of $k \in [L_1,L_2]_\BZ$ for which $E_{z,k}^M$ occurs.  (Only the length of the window matters,  since $\BP\lbrack E_{z,k}^M \mid \mathcal{F}_{z,k-1}\rbrack$ is bounded from below uniformly in $k$ and $z$.)  Taking $\alpha$ large enough in a way depending only on $\epsilon$ and $b_0$,  we obtain that off an event of probability at most $\tfrac12 (r_1/r_2)^{3}$ there exists $k \in [K_1,K_2]_\BZ$ for which $E_{z,k}^M$ occurs.

Fix such a $k$.  Since $2^{-k-1} \leq r_2$,  every crossing of $B_{r_2}(z) \setminus \overline{B_{r_1}(z)}$ contains a crossing of $B_{2^{-k-1}}(z) \setminus \overline{B_{r_1}(z)}$,  so it suffices to bound the probability of four (resp.\ more than four) crossings of the latter annulus.  As $B_{2^{-k-1}}(z) \subseteq B_{15 \cdot 2^{-k}/16}(z)$,  \Cref{lem:R_N_derivative_bounds_general} together with H\"older's inequality transfers that event to a nested $\CLE_4$ on $B_{2^{-k}}(z)$ coupled with a zero-boundary GFF there,  at the cost of raising its probability to a power $1/q$ with $q>1$ as close to $1$ as we like and multiplying by a constant depending only on $q$ and $M$.  Rescaling $B_{2^{-k}}(z)$ to $\BD$ turns $B_{2^{-k-1}}(z) \setminus \overline{B_{r_1}(z)}$ into precisely the annulus $C_{2^k r_1}$,  and $2^k r_1 \leq (r_1/r_2)^{1-\epsilon_1}$ by the choice of $K_2$.  Hence Propositions~\ref{prop:nested_cle_4_two_arm_exponent} and~\ref{prop:nested_cle_4_three_arm_exponent} bound the probability in question by $(r_1/r_2)^{2(1-\nu)(1-\epsilon_1)/q}$ and by $(r_1/r_2)^{(2+b_0)(1-\epsilon_1)/q}$,  respectively.  Taking $\nu \defeq \epsilon_1$ and then $q$ close enough to $1$,  the first exponent is at least $2-\epsilon$ and the second is at least $2 + b_0/2$;  decreasing $\delta$ so as to absorb the multiplicative constants,  we obtain the two assertions of the lemma,  with the universal constant there given by $b_0 / 2$.
\end{proof}

The first main step of the proof of Proposition~\ref{prop:nested_cle_4_two_arm_exponent} is to reduce the claim of the proposition to the case where we consider crossings of $C_{\varepsilon}$ made by finitely many level lines of a Gaussian free field,  where the number of these level lines is bounded by a large but fixed and deterministic $N \in \BN$.  This will be done in Section~\ref{subsec:finitely_many_crossings}.

In order to achieve this,  we fix $a \in (0,1)$ small and consider the coupling $(\Psi,\overline{\Gamma})$ introduced in Section~\ref{subsec:nested_cle_gff},  where $\Psi$ is a zero-boundary GFF on $\BD$.  More precisely,  for all $k \in \BN_0$,  we consider the annuli $A_k = B_{2^{-3k-1}}(0) \setminus \overline{B_{2^{-3k-2}}(0)}$ and inside every such annulus,  we consider the exploration of crossings of the annulus made by level lines of $\Psi$,  as in Section~\ref{subsubsec:crossings_exploration}.  Then,  since the exploration discovers finitely many crossings a.s.\ and the crossings are a.s.\ determined by the restrictions of $\Psi$ to the corresponding annuli,  we have by the independence across scales property of the GFF (see \cite[Lemma~3.1]{LocMetGFF}) that there exists $N \in \BN$ sufficiently large such that the following holds.  Off an event with probability at most $\varepsilon^3$,  for all sufficiently small $\varepsilon$,  there exists $k \in [1,a\log(1/\varepsilon)]_{\BZ}$ such that at most $N$ of the level line crossings discovered cross $A_k$ between its inner and outer boundaries.  We let $k_f$ be the smallest $k \in [1,a\log(1/\varepsilon)]_{\BZ}$ for which the above occurs.

Suppose that there exists a loop $\ell \in \overline{\Gamma}$ which crosses $C_{\varepsilon}$.  Then,  every time that $\ell$ crosses $C_{\varepsilon}$,  it must first cross $A_{k_f}$.  By using the level line interaction rules and arguing as in the proof of Lemma~\ref{lem:crossings_determined_by_the_field},  we obtain that every time the loop $\ell$ crosses $A_{k_f}$,  it has to merge with one of the level line crossings of $A_{k_f}$.  In particular,  if there are at least four crossings of $C_{\varepsilon}$ made by loops in $\overline{\Gamma}$,  then we can find four distinct level line crossings $\eta_1,\eta_2,\eta_3,\eta_4$ of $A_{k_f}$ such that their continuations after they exit $A_{k_f}$ for the first time cross the annulus $B_{2^{-3k_f - 2}}(0) \setminus \overline{B_{\varepsilon}(0)}$,  and hence,  since $k_f \leq a\log(1/\varepsilon)$,  also the annulus $B_{\varepsilon^{3a\log(2)}/4}(0) \setminus \overline{B_{\varepsilon}(0)}$.  Therefore,  it suffices to prove that the conditional probability of the event that,  for every $j \in \{1,2,3,4\}$,  the path $\eta_j$ crosses $B_{\varepsilon^{3a\log(2)}/4}(0) \setminus \overline{B_{\varepsilon}(0)}$ before intersecting any of the other three paths,  is at most $\varepsilon^{2(1-\nu)}$ for all sufficiently small $\varepsilon$.  This will complete the proof of Proposition~\ref{prop:nested_cle_4_two_arm_exponent}.  We will prove this in Section~\ref{subsec:prob_that_level_lines_survive} and then conclude the proof in Section~\ref{subsec:conlcusion_of_proof}.

The main ingredient of the proof of the latter claim is the following proposition. Its proof will be given in Section~\ref{subsubsec:level_lines_four_arms}.

\begin{proposition}\label{prop:bi_chordal_sle_4_two_arm_exponent}
Let $z_1,z_2,z_3,z_4$ be points on $\partial \BD$ ordered in the counterclockwise way,  fix $\sigma_1 ,  \sigma_3 \in \{-1,+1\}$,  and let $h$ be a GFF on $\BD$ whose boundary conditions are given by $2\sigma_1 \lambda$ on the counterclockwise arc of $\partial\BD$ from $z_1$ to $z_2$,  by $2\sigma_3\lambda$ on the counterclockwise arc of $\partial\BD$ from $z_3$ to $z_4$,  and by $0$ on the two remaining arcs.  For $j \in \{1,3\}$,  let $P_j$ be the level line of $\sigma_j h$ with height $\lambda$ from $z_j$ targeted at $z_{j+1}$,  and let $P_{j+1}$ be the level line of $-\sigma_j h$ with height $-\lambda$ from $z_{j+1}$ targeted at $z_j$.

Then,  there exists a universal constant $C \in (0,\infty)$ such that the following holds for all $\varepsilon \in (0,1)$.  The probability of the event that none of the paths $P_1,P_2,P_3,P_4$ reaches its continuation threshold and that each of them intersects $\partial B_{\varepsilon}(0)$ is at most $C \varepsilon ^2$.
\end{proposition}

Proposition~\ref{prop:bi_chordal_sle_4_two_arm_exponent} will follow 
from Lemma~\ref{lem:bichordal-bulk-4A}, and its proof is given at the end of the appendix.

Let us now briefly explain how to conclude the proof of Proposition~\ref{prop:nested_cle_4_two_arm_exponent} assuming Proposition~\ref{prop:bi_chordal_sle_4_two_arm_exponent}.  Fix $\rho \in (0,a)$ and $\zeta \in (0,\rho/2)$,  and cover the annulus $B_{\varepsilon^{3a\log(2)}/4}(0) \setminus \overline{B_{\varepsilon}(0)}$ by $n_\varepsilon \asymp (\rho/2+\zeta)^{-1}$ successive annuli,  each with a ratio of radii of order $\varepsilon^{\rho/2+\zeta}$.  Proposition~\ref{prop:bi_chordal_sle_4_two_arm_exponent} implies that,  conditionally on $\eta_1,\eta_2,\eta_3,\eta_4$ stopped upon reaching the outer boundary of one of these annuli,  the conditional probability that all four of them reach its inner boundary before intersecting one another is $O(\varepsilon^{\rho})$.  Multiplying these conditional bounds over the $n_\varepsilon$ annuli,  we obtain that the probability that all four paths reach $\partial B_{\varepsilon}(0)$ before intersecting one another is $O(\varepsilon^{\rho n_\varepsilon})$;  since $\rho n_\varepsilon \to 2$ as $\zeta/\rho \to 0$ and then $a \to 0$,  this completes the proof of the claim.

Finally,  we give a brief overview of the proof of Proposition~\ref{prop:nested_cle_4_three_arm_exponent}.  It will be a consequence of Proposition~\ref{prop:nested_cle_4_two_arm_exponent} together with the method introduced in Section~\ref{subsec:conformal_rectangles}.  More precisely,  we fix $\nu \in (0,1)$ and $\varepsilon \in (0,1)$,  and suppose that we are working on the event that there are at least four crossings of $C_{\varepsilon}$ made by loops in $\overline{\Gamma}$;  we let $\ell_1,\ell_2,\ell_3,\ell_4$ be four of them,  labeled in counterclockwise order,  and we set $\ell_5 \defeq \ell_1$.  For all $j \in \{1,2,3,4\}$,  we let $G_j$ denote the connected component of $C_{\varepsilon} \setminus \bigcup_{i=1}^4 \ell_i$ whose boundary contains both the right side of $\ell_j$ and the left side of $\ell_{j+1}$.  Then,  arguing as in Section~\ref{subsec:conformal_rectangles},  we will show that off an event with conditional probability at most $\varepsilon^{b_0}$ as $\varepsilon \to 0$ (for some universal constant $b_0>0$),  we have that the following holds.  For all $j \in \{1,2,3,4\}$,  we can divide $G_j$ into disjoint conformal rectangles $G_{j,1},\ldots,G_{j,n}$ as in Section~\ref{subsec:conformal_rectangles} for $n \asymp \log(1/\varepsilon)$,  such that for some $n_j \in \{1,\ldots,n\}$,  there exists a level line $P_j$ of $\Psi$ which makes a crossing in $\overline{G_{j,n_j}}$ from left to right.  The level lines $P_1,P_2,P_3,P_4$  will prevent any other loop in $\overline{\Gamma}$ from making a crossing of $C_{\varepsilon}$ due to the level line interaction rules. 

Combining with Proposition~\ref{prop:nested_cle_4_two_arm_exponent},  we obtain that off an event with probability at most $\varepsilon^{2(1-\nu) + b_0}$ as $\varepsilon \to 0$,  we have that there exist at most four crossings of $C_{\varepsilon}$ made by loops in $\overline{\Gamma}$.  Therefore,  Proposition~\ref{prop:nested_cle_4_three_arm_exponent} follows upon taking $\nu \defeq b_0 / 4$ and replacing $b_0$ by $b_0/2$,  so that $2(1-\nu)+b_0 = 2 + b_0 / 2$.

\subsection{Reducing the claim to finitely many crossings.}
\label{subsec:finitely_many_crossings}

In this section,  we will show that for fixed $a \in (0,1)$,  with very high probability as $\varepsilon \to 0$,  there is some $k \in [1,a\log(1/\varepsilon)]_{\BZ}$ for which $A_k$ is crossed by at most $N$ level lines of $\Psi$.  As explained in Section~\ref{subsec:outline_of_proof},  every time that a loop in $\overline{\Gamma}$ makes a crossing of $C_{\varepsilon}$,  it has to merge with one of the aforementioned crossings.  This reduces the claim of Proposition~\ref{prop:nested_cle_4_two_arm_exponent} to proving estimates on the probability of four-arm events for finitely many level lines of $\Psi$.

We will use the same exploration as in Section~\ref{subsubsec:crossings_exploration}.  More precisely,  we consider the rectangles $R^{\text{D}},  R^{\text{U}}$,  and $R^{\text{I}}$ as in Section~\ref{subsubsec:crossings_exploration}, and perform the level line explorations of the field $\Psi$.  Note that we have already shown in Lemma~\ref{lem:finite_number_of_crossings} that both of the explorations in $R^{\text{D}}$ and $R^{\text{U}}$ terminate in finitely many steps a.s.  Recall from Section~\ref{subsubsec:crossings_exploration} that each of the crossings discovered in $R^{\text{D}}$ and $R^{\text{U}}$ will continue into $R^{\text{I}}$ and either form a crossing of $R^{\text{I}}$ or make an excursion into $R^{\text{I}}$ back to $R^q$ for $q \in \{D,U\}$.  

We have the following analog of Lemma~\ref{lem:crossings_determined_by_the_field}.

\begin{lemma}\label{lem:all_crossings_determined_by_the_field}
The collection of segments of loops in $\overline{\Gamma}$ which are contained in $A_0$ and whose endpoints lie on $\partial A_0$ is a.s.\ determined by the restriction of the field $\Psi$ to $A_0$.
\end{lemma}

\begin{proof}
The claim in the statement of the lemma follows from the exact same argument used to prove Lemma~\ref{lem:crossings_determined_by_the_field}.
\end{proof}

We can perform the same exploration described above but with the annulus $A_k$ in place of $A_0$ for all $k \in \BN_0$.  Then,  the argument used to prove Lemma~\ref{lem:all_crossings_determined_by_the_field} implies that the segments of loops in $\overline{\Gamma}$ which are contained in $A_k$ and whose endpoints lie on $\partial A_k$ are a.s.\ determined by $\Psi|_{A_k}$.  For all $k \in \BN_0$,  we let $\Gamma_k^{\text{out}}$ (resp.\ $\Gamma_k^{\text{in}}$) denote the collection of segments of loops in $\overline{\Gamma}$ which make an excursion into $A_k$ from the outer (resp.\ inner) boundary of $A_k$.  We also let $\Gamma_k^{\text{I}}$ denote the collection of segments of loops in $\overline{\Gamma}$ which cross $A_k$.

We are now ready to define the main event that we are going to consider.

\begin{definition}

Fix $\delta \in (0,1)$ small and $N \in \BN_0$ large (both deterministic).  For all $k \in \BN_0$,  we let $E_k = E_k(\Psi)$ denote the event that the following hold.

\begin{enumerate}
\item \label{it:number_of_crossings_finite}
The number of crossings in $\Gamma_k^{\text{I}}$ is at most $N$.
\item \label{it:boundary_conditions_bounded}
There exists $c \in 2\lambda\BZ$ such that the boundary conditions of $\Psi$ on the left and right sides of any segment $\ell$ in $\Gamma_k^{\text{I}}$ lie in $[c-N,c+N]$.
\item \label{it:crossings_separated}
For all $\ell,  \ell' \in \Gamma_k^{\text{I}}$ such that $\ell \neq \ell'$,  we have that $\text{dist}(\ell ,  \ell') \geq \delta 2^{-3k}.$
\item \label{it:harmonic_measure_bounded_below}
For all $\ell \in \Gamma_k^{\text{I}}$,  we have that the harmonic measure as seen from $0$ of the part of the left (resp.\ right) side of $\ell$ contained in $\ell \setminus \overline{\left(\bigcup _{\ell' \in \Gamma_k^{\text{out}}} \ell' \right)}$ in the connected component of $\BD \setminus \overline{\left(\bigcup_{\gamma \in \Gamma_k^{\text{out}} \cup \Gamma_k^{\text{I}}} \gamma \right)}$ containing $0$ is at least $\delta$.
\end{enumerate}
\end{definition}

Note that Lemma~\ref{lem:all_crossings_determined_by_the_field} implies that $E_k$ is a.s.\ determined by $\Psi|_{A_k}$.  Hence,  by absolute continuity,  we can define the events $E_k(\widetilde\Psi)$ in the same way for any field $\widetilde\Psi$ on $\BD$ whose boundary conditions are piecewise constant and change only countably many times,  with the segments of loops replaced by the corresponding level lines of $\widetilde\Psi$.  Let us emphasize that condition~\eqref{it:boundary_conditions_bounded} is formulated relative to an unspecified $c \in 2\lambda\BZ$,  so that $E_k(\widetilde\Psi) = E_k(\widetilde\Psi - 2\lambda m)$ for every $m \in \BZ$.  This invariance is essential: the heights of the crossings of $A_k$ perform a random walk in $k$,  so that an absolute bound on them would fail for all but $O(\sqrt{k})$ of the first $k$ scales.

Condition~\eqref{it:number_of_crossings_finite} will allow us to reduce the claim of Proposition~\ref{prop:nested_cle_4_two_arm_exponent} to proving estimates on the probability of four-arm events for at most $N$ level lines of $\Psi$.  Moreover,  condition~\eqref{it:boundary_conditions_bounded} will be useful since it will allow us to compare the law of a field with complicated boundary conditions that can change infinitely many times with the law of a field whose boundary conditions are simpler and whose level lines with appropriate heights have the same joint law as the level lines appearing in the statement of Proposition~\ref{prop:bi_chordal_sle_4_two_arm_exponent}.  As for conditions~\eqref{it:crossings_separated} and~\eqref{it:harmonic_measure_bounded_below},  they will ensure that the level lines discovered by the exploration are sufficiently far away from each other so that Proposition~\ref{prop:bi_chordal_sle_4_two_arm_exponent} can be applied to those level lines.

If we choose the constants $\delta,N$ appearing in the definition of $E_k(\Psi)$ appropriately,  we have that the event $E_k(\Psi)$ occurs for many values of $k$ with very high probability.  This is the content of the following lemma.

\begin{lemma}\label{lem:event_happens_at_many_scales}
Fix $a \in (0,1)$.  Then,  we can choose $\delta,  \varepsilon_0 \in (0,1)$ sufficiently small and $N \in \BN_0$ sufficiently large (depending only on $a$) such that the following holds for all $\varepsilon \in (0,\varepsilon_0)$.  Off an event with probability at most $\varepsilon^3$,  we have that there exists $k \in [1 , a \log(1 / \varepsilon)]_{\BZ}$ such that $E_k(\Psi)$ occurs.  
\end{lemma}

\begin{proof}
First,  we note that by the definition of the exploration,  we have that the crossings in $\Gamma_0^{\text{I}}$ are a.s.\ pairwise disjoint and that a.s.,
\begin{align*}
\overline{\bigcup_{\ell \in \Gamma_{0}^{\text{out}}} \ell } \cap \partial B_{1/4}(0) = \emptyset.
\end{align*}
Moreover,  the boundary conditions of $\Psi$ on the left and right sides of any segment in $\Gamma_0^{\text{I}}$ lie in some compact interval since the exploration of the segments terminates after finitely many steps a.s.  Thus,  since the number of crossings in $\Gamma_0^{\text{I}}$ is a.s.\ finite (by the local finiteness property of nested $\CLE_4$),  we obtain that for all $p \in (0,1)$,  we can choose $\delta \in (0,1)$ sufficiently small and $N \in \BN$ sufficiently large such that $\BP[E_0(\Psi)] \geq p$.  We now transfer this to all scales.  Let us apply \Cref{lem:M_good_event_general} with $F = \overline{B_{1/2}(0)}$ and $z=0$,  with $K \defeq \lfloor 3a\log(1/\varepsilon)\rfloor$,  and with its two parameters (which we denote here by $\alpha$ and $\beta_0$,  so as not to clash with the constants $a$ and $b$ of this appendix) chosen so that $\beta_0 \in (5/6,1)$ and $\alpha$ is large depending on $a$.  It follows that off an event of probability at most $\varepsilon^3$ there are at most $(1-\beta_0)K$ values of $j \in [1,K]_\BZ$ for which $B_{2^{-j}}(0)$ fails to be $M$-good.  Since the number of $j \in [1,K]$ which are divisible by $3$ is at least $K/3 - 2$ and $1-\beta_0 < 1/6$,  we deduce that there are at least $c_1 a \log(1/\varepsilon)$ values of $k \in [1,a\log(1/\varepsilon)]_\BZ$ for which $B_{2^{-3k}}(0)$ is $M$-good,  where $c_1>0$ is universal.  For each such $k$ we have $A_k \subseteq B_{15 \cdot 2^{-3k}/16}(0)$,  and $E_k(\Psi) = E_k(\Psi - \lfloor \Psi_{2^{-3k}}(0)/(2\lambda)\rfloor 2\lambda)$ by the invariance noted after the definition of $E_k$;  hence \Cref{lem:R_N_derivative_bounds_general} together with H\"older's inequality and the scale invariance of the law of the zero-boundary GFF gives $\BP\lbrack E_k(\Psi) \mid \mathcal{F}_{0,3k}\rbrack \geq p^\prime$ a.s.\ on $E_{0,3k}^M$,  where $p^\prime \in (0,1)$ depends only on $p$ and $M$.  Since $E_j(\Psi)$ is $\mathcal{F}_{0,3k}$-measurable for $j<k$,  iterating over the good scales gives that the probability that none of the events $E_k(\Psi)$,  $k \in [1,a\log(1/\varepsilon)]_\BZ$,  occurs is at most $(1-p^\prime)^{c_1 a \log(1/\varepsilon)} = \varepsilon^{c_1 a \log(1/(1-p^\prime))}$,  which is at most $\varepsilon^3$ once $p$ is chosen close enough to $1$ in a way depending only on $a$.
\end{proof}

\begin{figure}[h]
	\centering
	$\vcenter{\hbox{\includegraphics[scale=0.8]{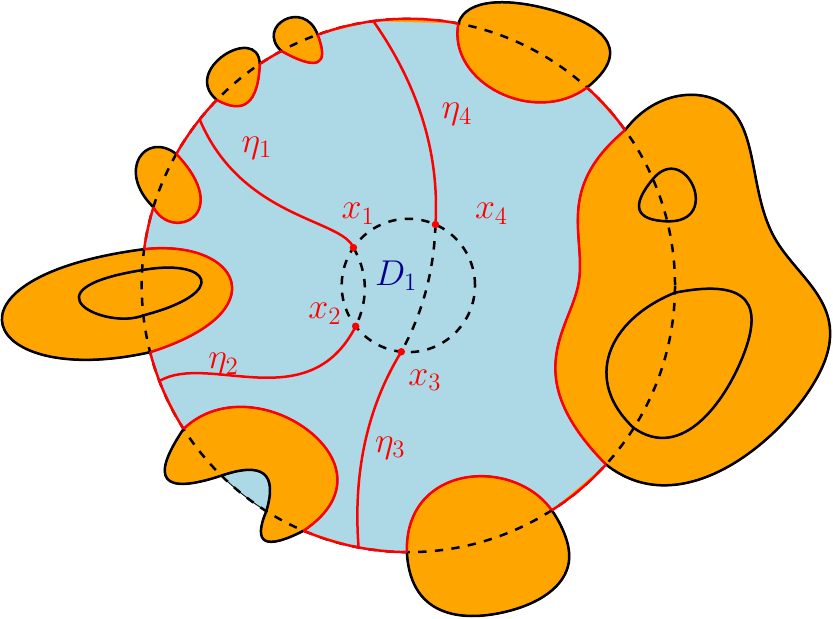}}}$

	\caption{Illustration of the domain $D_1$ in blue with its boundary in red.}
	\label{fig:domain_D1_four_arm_nested_CLE}
\end{figure}

Let $k_f$ denote the smallest integer $k \in [1 ,  a \log(1/\varepsilon)]_{\BZ}$ for which $E_k(\Psi)$ occurs.  Note that Lemma~\ref{lem:event_happens_at_many_scales} implies that $k_f$ exists with very high probability as $\varepsilon \to 0$.  Let also $D_1$ denote the connected component of $\BC \setminus \overline{\bigcup_{\ell \in \Gamma_{k_f}^{\text{out}} \cup \Gamma_{k_f}^{\text{I}}} \ell}$ containing $0$ and note that $D_1$ is a simply connected domain (see Figure \ref{fig:domain_D1_four_arm_nested_CLE}).  Indeed,  a.s.,  Lebesgue-a.e.\ point $z$ of the circle $\partial B_{2^{-3k_f-1}}(0)$ lies on no loop of $\overline\Gamma$,  and the loops of $\overline\Gamma$ surrounding such a point $z$ shrink to $z$.  Any such loop $\gamma$ of sufficiently small diameter must intersect both connected components of $\BC \setminus \partial B_{2^{-3k_f-1}}(0)$ (the interior of $\gamma$ is an open set of diameter at most $\diam(\gamma)$ which contains $z$,  and hence meets both components;  if $\gamma$ avoided the open ball,  then the ball,  being connected,  disjoint from $\gamma$,  and meeting the interior of $\gamma$,  would be contained in that interior,  which has smaller diameter,  while if $\gamma$ avoided the complement of the closed ball,  then the latter would lie in the unbounded complementary component of $\gamma$ and could not meet the interior of $\gamma$);  hence the excursions of $\gamma$ into $B_{2^{-3k_f-1}}(0)$ near $z$ are excursions into $A_{k_f}$ from its outer boundary.  Therefore $\overline{\bigcup_{\ell \in \Gamma_{k_f}^{\mathrm{out}} \cup \Gamma_{k_f}^{\mathrm{I}}} \ell}$ contains a dense subset of the circle $\partial B_{2^{-3k_f-1}}(0)$,  and hence the whole circle;  since every segment in $\Gamma_{k_f}^{\mathrm{out}} \cup \Gamma_{k_f}^{\mathrm{I}}$ has an endpoint on that circle,  this set is connected,  and it separates $0$ from $\infty$.
Moreover,  it is easy to see that $\partial D_1$ is locally connected.  Thus,  there exists a conformal transformation $f_1$ mapping $D_1$ onto $\BD$ such that $f_1(0) = 0$ and $f_1'(0) > 0$,  and $f_1$ extends continuously to a map from $\overline{D_1}$ onto $\overline{\BD}$.  Furthermore,  for the rest of the section,  we will assume that $\Gamma_{k_f}^{\text{I}} \neq \emptyset$,  i.e.,  there are loops in $\overline{\Gamma}$ crossing $A_{k_f}$.

In the following lemma,  we show that the boundary conditions of $\Psi|_{D_1}$ lie on a fixed and deterministic compact interval.  We will use this fact in Sections~\ref{subsec:prob_that_level_lines_survive} and~\ref{subsec:conlcusion_of_proof}.

\begin{lemma}\label{lem:field_has_bounded_boundary_conditions}
Suppose that we have the same setup described above,  and work on the event that $k_f$ exists.  Then,  writing $c \in 2\lambda\BZ$ for the constant furnished by condition~\eqref{it:boundary_conditions_bounded} of $E_{k_f}(\Psi)$,  the boundary conditions of the field $\Psi|_{D_1}$ lie in $[c-N-4\lambda,c+N+4\lambda]$;  in particular,  upon replacing $N$ by $\lceil N+4\lambda \rceil$,  we may and do assume that they lie in $[c-N,c+N]$.
\end{lemma}

\begin{proof}
\stepn{step:bdd-bc-reduction-single-excursion}{Reduction to a single excursion}  Fix $\ell \in \Gamma_{k_f}^{\text{out}}$ such that $\ell \subseteq \partial D_1$ and let $\gamma \in \overline{\Gamma}$ be such that $\ell \subseteq \gamma$.  We will show that the boundary conditions of $\Psi$ on the left and right sides of $\ell$ lie in $[c-N-4\lambda,c+N+4\lambda]$.  This will complete the proof since we already know by~\eqref{it:boundary_conditions_bounded} that the boundary conditions of $\Psi$ on $\partial D_1 \cap \left(\overline{\bigcup_{\ell' \in \Gamma_{k_f}^{\text{I}}} \ell' }\right)$ lie in $[c-N,c+N]$.  All of the intervals appearing below are to be read relative to $c$ in the same way;  the argument is purely relative,  since it only ever compares the boundary conditions on two loops whose nesting generations differ by one.  We will treat the following cases separately,  according to the position of $\gamma$ relative to $\partial B_{2^{-3k_f-2}}(0)$.

\stepn{step:bdd-bc-case-where-does}{The case where $\gamma$ does not meet $\partial B_{2^{-3k_f-2}}(0)$ and does not surround it}
If $\gamma$ is an outermost loop in $\overline{\Gamma}$,  then by the construction of the coupling $(\Psi ,  \overline{\Gamma})$ we have that the boundary conditions of $\Psi$ on $\gamma$ are contained in $\{-2\lambda,0,2\lambda\}$.  In this case $\gamma$ is an outermost loop which meets $\partial B_{2^{-3k_f-1}}(0)$;  the claim therefore holds after increasing $N$ by a universal constant.

Suppose now that $\gamma$ is not an outermost loop in $\overline{\Gamma}$ and let $\gamma'$ be a loop in $\overline{\Gamma}$ that surrounds $\gamma$ such that the nesting generation of $\gamma'$ is given by $n-1$,  where $n$ denotes the nesting generation of $\gamma$.  If $\gamma'$ crosses $A_{k_f}$,  then the boundary conditions of $\Psi$ on $\gamma'$ lie in $[c-N,c+N]$,  and so the construction of the coupling $(\Psi ,  \overline{\Gamma})$ implies that the boundary conditions of $\Psi$ on $\gamma$ lie in $[c-N-2\lambda ,  c+N + 2\lambda]$.  If $\gamma'$ does not cross $A_{k_f}$,  then we must have that $\partial B_{2^{-3k_f - 2}}(0)$ lies in the bounded connected component of $\BC \setminus \gamma'$,  since otherwise we would have that $\gamma'$ disconnects $\gamma$ from $\partial B_{2^{-3k_f-2}}(0)$,  and that contradicts the definition of $\Gamma_{k_f}^{\text{out}}$.  Hence, the loop $\gamma'$ surrounds $B_{2^{-3k_f - 2}}(0)$.  This implies that $\gamma'$ surrounds all the loops in $\overline{\Gamma}$ which cross $A_{k_f}$.  

Let $\tilde{\gamma}$ be a loop in $\overline{\Gamma}$ which crosses $A_{k_f}$ and which is outermost among those loops in the sense that there is no other loop in $\overline{\Gamma}$ which crosses $A_{k_f}$ and surrounds $\tilde{\gamma}$.  We claim that the nesting generation of $\tilde{\gamma}$ is given by $n$,  i.e.,  $\tilde{\gamma}$ and $\gamma$ have the same generation.  Indeed,  suppose that $\tilde{\gamma}$ and $\gamma$ do not have the same generation.  Then,  there exists a loop $\hat{\gamma} \in \overline{\Gamma}$ such that $\hat{\gamma}$ disconnects $\tilde{\gamma}$ from $\gamma$.  Note that $\hat{\gamma} \cap \partial B_{2^{-3k_f - 1}}(0) \neq \emptyset$ since $\hat{\gamma}$ surrounds $\tilde{\gamma}$ and note that $\hat{\gamma}$ does not cross $A_{k_f}$ due to the choice of $\tilde{\gamma}$.  Thus,  we must have that $\hat{\gamma}$ disconnects $\gamma$ from $\partial B_{2^{-3k_f-2}}(0)$ but that contradicts the definition of $\Gamma_{k_f}^{\text{out}}$ together with the fact that $\ell \subseteq \partial D_1$.  It follows that $\tilde{\gamma}$ and $\gamma$ have the same nesting generation.  In particular,  the construction of the coupling $(\Psi ,  \overline{\Gamma})$ implies that the boundary conditions of $\Psi$ on $\gamma'$ lie in $[c-N - 2\lambda ,  c+N + 2\lambda]$,  which also implies that the boundary conditions of $\Psi$ on $\gamma$ lie in $[c-N-4\lambda ,  c+N+4\lambda]$,  as required.

\stepn{step:bdd-bc-case-where-does-2}{The case where $\gamma$ does not meet $\partial B_{2^{-3k_f-2}}(0)$ and surrounds it}
In that case,  the loop $\gamma$ surrounds $B_{2^{-3k_f - 2}}(0)$.  Moreover,  all of the loops in $\overline{\Gamma}$ which cross $A_{k_f}$ are surrounded by $\gamma$;  we let $\tilde{\gamma}$ be a loop in $\overline{\Gamma}$ which crosses $A_{k_f}$ and which is outermost in the same sense as in \Cref{step:bdd-bc-case-where-does}.  Then,  if $n$ is the nesting generation of $\gamma$,  we claim that the nesting generation of $\tilde{\gamma}$ is given by $n+1$.  Indeed,  suppose that this is not the case.  Then,  there exists a loop $\hat{\gamma} \in \overline{\Gamma}$ such that $\hat{\gamma}$ surrounds $\tilde{\gamma}$ and $\hat{\gamma}$ is surrounded by $\gamma$.  Since we have that $\hat{\gamma} \cap \partial B_{2^{-3k_f -2}}(0) = \emptyset$ (by the outermost nature of $\tilde{\gamma}$),  we obtain that $\hat{\gamma}$ surrounds $B_{2^{-3k_f-2}}(0)$,  and hence $\hat{\gamma}$ disconnects $\ell$ from $\partial B_{2^{-3k_f -2}}(0)$.  But this contradicts the fact that $\ell \subseteq \partial D_1$.  It follows that the nesting generation of $\tilde{\gamma}$ is given by $n+1$ and so the boundary conditions of $\Psi$ on $\gamma$ lie in $[c-N-2\lambda ,  c+N+2\lambda]$.  

\stepn{step:bdd-bc-case-where-meets}{The case where $\gamma$ meets $\partial B_{2^{-3k_f-2}}(0)$}
In that case,  we have that $\gamma$ crosses $A_{k_f}$ and thus~\eqref{it:boundary_conditions_bounded} implies that the boundary conditions of $\Psi$ on $\gamma$ lie in $[c-N,c+N]$.  Therefore,  the proof of the lemma is complete by combining with the previous cases.
\end{proof}

\subsection{Bounding the probability that level lines intersect $\partial B_{\varepsilon}(0)$ without intersecting each other.}
\label{subsec:prob_that_level_lines_survive}

\subsubsection{Outline and setup.}
\label{subsubsec:outline_and_setup}

For the rest of the section,  we fix $a \in (0,1/100)$ small and assume that $\delta \in (0,1),  N \in \BN$ are chosen so that the statement of Lemma~\ref{lem:event_happens_at_many_scales} holds.  We also fix points $x_1,x_2,x_3$,  and $x_4$ on $\partial D_1$ which are endpoints of crossings contained in $\Gamma_{k_f}^{\text{I}}$,  labeled in counterclockwise order on $\partial D_1$.   Let $\eta_1,\eta_2,\eta_3,\eta_4$ be the crossings in $\Gamma_{k_f}^{\text{I}}$ which terminate at $x_1,x_2,x_3,x_4$ respectively.  We also assume that there exist $k \in \BZ$ and $\sigma_1 ,  \sigma_3 \in \{-1,+1\}$ such that the following holds for $j \in \{1,3\}$.  The boundary conditions of $\Psi$ on the left (resp.\ right) side of $\eta_j$ are given by $2\lambda k$ (resp.\ $2\lambda(k+\sigma_j)$),  and on the left (resp.\ right) side of $\eta_{j+1}$ are given by $2\lambda(k+\sigma_j)$ (resp.\ $2\lambda k$);  equivalently,  $\eta_j$ (resp.\ $\eta_{j+1}$) corresponds to a level line of $\sigma_j\Psi$ (resp.\ $-\sigma_j\Psi$) with height $\sigma_j(2k+\sigma_j)\lambda$ (resp.\ $-\sigma_j(2k+\sigma_j)\lambda$).

The above assumption corresponds,  after subtracting the constant $2\lambda k$ from $\Psi$,  to the setup of Proposition~\ref{prop:bi_chordal_sle_4_two_arm_exponent} with the same value of $(\sigma_1,\sigma_3)$:  the height $\sigma_j(2k+\sigma_j)\lambda$ of $\eta_j$ becomes $\sigma_j(2k+\sigma_j)\lambda - 2k\sigma_j\lambda = \lambda$ for the field $\sigma_j(\Psi - 2\lambda k)$,  as in the proposition.  The reason for restricting attention to crossings as above is the following.  Suppose that $\ell,\ell' \in \overline{\Gamma}$ are distinct loops so that both of them cross $C_{\varepsilon}$ and so that either $\ell,\ell'$ have the same generation,  or one of the loops surrounds the other and their nesting generations differ by $1$.  Then there are four crossings in $\Gamma_{k_f}^{\text{I}}$ which merge with $\ell$ or with $\ell^\prime$,  and,  since $\ell$ and $\ell^\prime$ are disjoint,  the two which merge with the same loop have adjacent endpoints in the counterclockwise order on $\partial D_1$.  Labeling them $\eta_1,\eta_2,\eta_3,\eta_4$ in counterclockwise order in such a way that $\eta_1,\eta_2$ merge with $\ell^\prime$ and $\eta_3,\eta_4$ merge with $\ell$,  the construction of the coupling $(\Psi,\overline{\Gamma})$ implies that the above assumption holds for some $k \in \BZ$ and $\sigma_1,\sigma_3 \in \{-1,+1\}$.  Let us stress that it is the pairing of the crossings which forces the labeling.  Indeed,  the two arcs of $\partial D_1$ between $\eta_2$ and $\eta_3$ and between $\eta_4$ and $\eta_1$ are the ones lying outside both $\ell$ and $\ell^\prime$ (in the case that $\ell$ surrounds $\ell^\prime$,  outside $\ell^\prime$ but inside $\ell$),  so that the boundary conditions of $\Psi$ on both of them are equal to $2\lambda k$,  exactly as the above assumption requires.  Under a cyclic shift of the labels by one,  the two crossings labeled $\eta_1$ and $\eta_2$ would be arms of different loops,  and the resulting boundary conditions need not satisfy the above assumption.  Similarly,  if the loop $\ell \in \overline{\Gamma}$ makes four crossings of $C_{\varepsilon}$ and $\eta_1,\eta_2,\eta_3,\eta_4$ are the crossings in $\Gamma_{k_f}^{\text{I}}$ which merge with $\ell$,  labeled in counterclockwise order,  we have that the above assumption holds.  We will see in Section~\ref{subsec:conlcusion_of_proof} why it suffices to consider crossings of $C_{\varepsilon}$ made by loops $\ell,\ell'$ in $\overline{\Gamma}$ as above.

Next,  let us fix the hierarchy of the constants,  which we shall use throughout.  We first fix $a$ as above and a ratio $s \in (0,1)$;  next $\delta \in (0,1)$ small and $N \in \BN$ large,  depending only on $a$ and $s$;  next $b \in (0,a/100)$ depending only on $a$;  next $\rho \in (0,a)$ small,  depending only on $a,b,\delta,N$;  and finally $\zeta \defeq s\rho / 2 \in (0,\rho/2)$.  In the proofs of Propositions~\ref{prop:nested_cle_4_two_arm_exponent} and~\ref{prop:nested_cle_4_three_arm_exponent} we shall let $s \to 0$ and then $a \to 0$.  For all $\varepsilon \in (0,1)$,  we let $n_{\varepsilon}$ denote the largest positive integer for which it holds that
\begin{align*}
\left(\frac{\varepsilon^{\rho / 2 + 6a\log(2)}}{64^2}\right) \left(\frac{\varepsilon^{\rho / 2 + \zeta}}{16^2}\right)^{n_{\varepsilon} - 1} \geq \varepsilon,
\end{align*}
and note that 
\begin{align}\label{eqn:number_of_iterations}
n_{\varepsilon} \geq R \defeq \frac{1 - 6a \log(2) - \rho / 2}{\zeta + \rho / 2} \quad \text{for all sufficiently small } \varepsilon .
\end{align}
Moreover,  for all $j \in \{1,2,3,4\}$ and all $n \in \BN$,  we let $\eta_{j,n}$ denote the continuation of the curve $\eta_j$ stopped at the first time that it intersects $\partial B_{(\varepsilon^{\rho / 2 + 6a\log(2)} / 64^2) (\varepsilon^{\rho / 2 + \zeta} / 16^2)^{n-1}}(0)$.  

The main goal of this section is to prove the following lemma.

\begin{lemma}\label{lem:bound_on_number_of_iterations}
Suppose that we have the setup described above.  Then,  for all $\rho$ and $\zeta$ as above,  there exists $\varepsilon_0 \in (0,1)$ depending only on $a,\delta,N,\rho$,  and $\zeta$ such that the following holds for all $\varepsilon \in (0,\varepsilon_0)$.  Let $F_{\varepsilon}$ denote the event that for all $1\leq m \leq n_{\varepsilon}$
and all $j \in \{1,2,3,4\}$,  we have that the path $\eta_{j,m}$ intersects $\partial B_{(\varepsilon^{\rho / 2 + 6a\log(2)} / 64^2) (\varepsilon^{\rho / 2 + \zeta} / 16^2)^{m-1}}(0)$ before intersecting any other path $\eta_{k,m}$ for $k\neq j$.  Then,  we have that
\begin{align*}
\BP[F_{\varepsilon}] \leq \varepsilon^3 + \varepsilon^{\rho n_{\varepsilon} (1-a)}.
\end{align*}
\end{lemma}

Lemma~\ref{lem:bound_on_number_of_iterations} will follow from the following two lemmas.

\begin{lemma}\label{lem:first_stage_cond_prob}
Suppose that we have the same setup described above.  Then,  for all $\rho \in (0,a)$ sufficiently small,  the following holds. Conditional on the event in the statement of Lemma~\ref{lem:event_happens_at_many_scales},  the conditional probability that,  for every $j \in \{1,2,3,4\}$,  the path $\eta_{j,1}$ intersects $\partial B_{\varepsilon^{\frac{\rho}{2} + 6a\log(2)} / 64^2}(0)$ before intersecting any of the other three paths $\eta_{k,1}$,  is at most $\varepsilon^{\rho (1-a)}$ for all $\varepsilon \in (0,\varepsilon_0)$,  where $\varepsilon_0 \in (0,1)$ depends only on $a,\delta,N$,  and $\rho$.
\end{lemma}

\begin{lemma}\label{lem:inductive_stage_cond_prob}
Suppose that we have the same setup described above.  Then,  for all $\rho$ and $\zeta$ as above and all $n \in \BN$,  the following holds.  Conditional on the event that $\eta_{j,m}$ intersects $\partial B_{(\varepsilon^{\rho / 2 + 6a\log(2)} / 64^2) (\varepsilon^{\rho / 2 + \zeta} / 16^2)^{m-1}}(0)$ for the first time before intersecting $\eta_{k,m}$ for all $1 \leq m \leq n$ and all $j,k \in \{1,2,3,4\}$ such that $j \neq k$,  the conditional probability that,  for every $j \in \{1,2,3,4\}$,  the path $\eta_{j,n+1}$ intersects $\partial B_{(\varepsilon^{\rho / 2 + 6a\log(2)} / 64^2) (\varepsilon^{\rho / 2 + \zeta} / 16^2)^n}(0)$ before intersecting any of the other three paths $\eta_{k,n+1}$,  is at most $\varepsilon^{\rho (1-a)}$ for all $\varepsilon \in (0,\varepsilon_0)$,  where $\varepsilon_0 \in (0,1)$ depends only on $a,\delta,N,\rho$,  and $\zeta$.
\end{lemma}

\subsubsection{Proving Lemma~\ref{lem:first_stage_cond_prob}.}

Next,  we prove Lemma~\ref{lem:first_stage_cond_prob}.  First,  we set $\Psi_1 = \Psi \circ f_1^{-1}$,  and for all $j \in \{1,2,3,4\}$,  we set $\tilde{x}_j = f_1(x_j)$ and let $\gamma_j$ be the corresponding level line of either $\Psi_1$ or $-\Psi_1$ starting from $\tilde{x}_j$ and having the same height as $\eta_j$.

We will use Proposition~\ref{prop:bi_chordal_sle_4_two_arm_exponent} in order to prove Lemma~\ref{lem:first_stage_cond_prob}.  We start with the following lemma which proves that with very high probability,  we have that the tips of the level lines $\eta_1,\eta_2,\eta_3,\eta_4$ become sufficiently separated if we draw the curves a bit further.

\begin{lemma}\label{lem:tips_are_separated}
Let $\delta,N$ be the constants in the statement of Lemma~\ref{lem:event_happens_at_many_scales}.  Then,  by possibly taking $\delta \in (0,1)$ to be smaller and $N$ to be larger (in a way that depends only on $a$),  we have that the following holds off an event whose probability is $O(\varepsilon^3)$ as $\varepsilon \to 0$ (at a rate depending only on $a,\delta$,  and $N$).  There exists $n \in [a \log(1/\varepsilon) / 2 , a \log(1/\varepsilon)]_{\BZ}$ such that the event $E_n(\Psi_1)$ occurs.
\end{lemma}

\begin{proof}
First,  we note that Lemma~\ref{lem:field_has_bounded_boundary_conditions} together with \cite[Proposition~3.4]{IG1} implies that for any compact set $K \subseteq \BD$,  the conditional law of $(\Psi_1-c)|_K$ given $D_1$ is absolutely continuous with respect to the corresponding restriction to $K$ of a zero-boundary GFF on $\BD$,  and that the $p$-th moments of the corresponding Radon--Nikodym derivatives are bounded from below and above by constants that depend only on $p,K$,  and $N$.  Since the events $E_n(\Psi_1)$ are invariant under the addition of an element of $2\lambda\BZ$ to the field,  and since $c \in 2\lambda\BZ$,  this comparison may be applied directly to them.  Therefore, taking $K=\overline{B_{1/2}(0)}$ and combining this with the argument used to prove Lemma~\ref{lem:event_happens_at_many_scales}, we complete the proof of the lemma.
\end{proof}

\begin{figure}[h]
	\centering
	$\vcenter{\hbox{\includegraphics[scale=0.8]{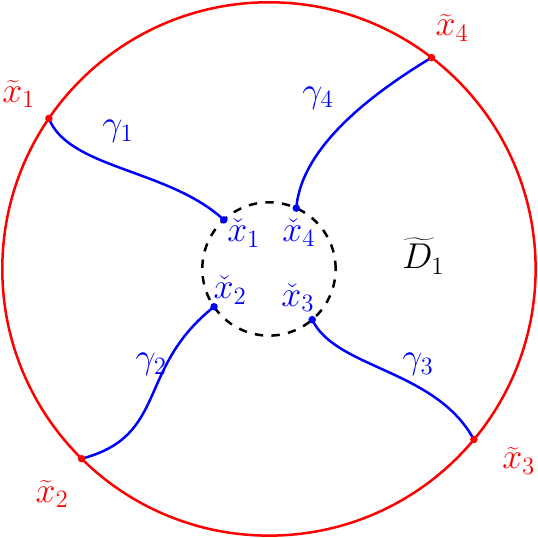}}}$
	\ \
	$\vcenter{\hbox{\includegraphics[scale=0.8]{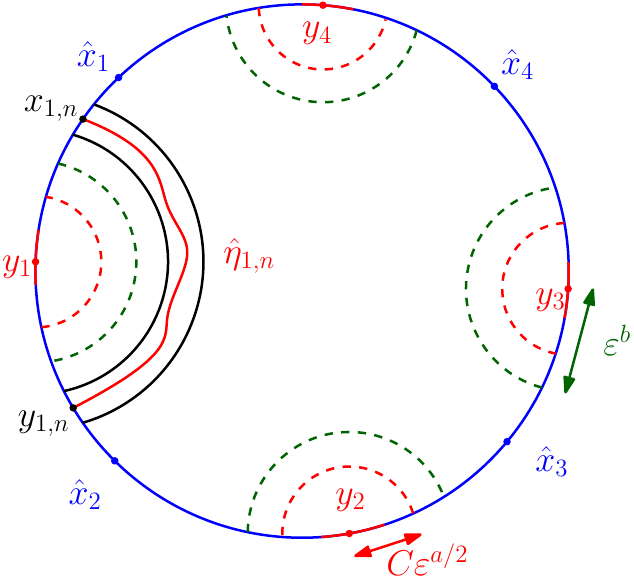}}}$
	\caption{On the left-hand side, illustration of the domain $\widetilde{D}_1$ in white which is the complement of the union of the blue curves $\gamma_1, \ldots, \gamma_4$ in the disk encircled by the red circle. On the right-hand side, we depict the image of $\widetilde{D}_1$ and of its boundary under the map $g_1$. The blue arcs are the images of the curves $\gamma_1, \ldots, \gamma_4$ and the red arcs are the images of the red arcs on the left-hand side. They are included in balls of radius $C\varepsilon^{a/2}$ centered at $y_1, \ldots, y_4$. The domain obtained by removing the green balls of radius $\varepsilon^b$ from  $\BD$ is the domain described in Lemma \ref{lem:comparison_of_fields}. Finally, we draw one instance of a level line $\hat{\eta}_{j,n}$ whose existence is proved and used in the proof of \Cref{lem:level_lines_do_not_hit_small_balls}.}
	\label{fig:domain_D1_tilde_four_arm_nested_CLE}
\end{figure}

For the rest of the section,  we will assume that there exists $n \in [a \log(1/\varepsilon) / 2 , a \log(1/\varepsilon)]_{\BZ}$ such that the event $E_n(\Psi_1)$ occurs,  and we let $\tilde{k}_f$ be the smallest such integer $n$.  Note that the level line interaction rules combined with the definition of $E_{\tilde{k}_f}(\Psi_1)$ imply that the tips of the curves $\gamma_j$ for $j=1,2,3,4$ are distinct when they are stopped at the first time that they intersect $\partial B_{2^{-3\tilde{k}_f - 2}}(0)$.  Let $\tilde{D}_1$ denote the connected component containing $0$ of the complement in $\BD$ of the union of $\gamma_1,\gamma_2,\gamma_3,\gamma_4$ when they are stopped at the first time that they intersect $\partial B_{2^{-3\tilde{k}_f - 2}}(0)$ and let $\check{x}_1, \check{x}_2, \check{x}_3, \check{x}_4$ be their tips. See the left-hand side of Figure \ref{fig:domain_D1_tilde_four_arm_nested_CLE}.
Then,  we let $g_1$ denote the conformal transformation mapping $\tilde{D}_1$ onto $\BD$ such that $g_1(0) = 0$ and $g_1'(0) > 0$.

Set $\hat{x}_j = g_1(\check{x}_j)$ for $j=1,2,3,4$,  and note that condition~\eqref{it:harmonic_measure_bounded_below} in the definition of $E_{\tilde{k}_f}(\Psi_1)$ implies that there exists $\tilde{\delta}_1 \in (0,\delta)$ depending only on $\delta$ such that $|\hat{x}_j - \hat{x}_k| \geq \tilde{\delta}_1$ for all $j,k \in \{1,2,3,4\}$ such that $j \neq k$.  Indeed,  let $\Gamma_{\tilde{k}_f}^{\mathrm{out}}(\Psi_1)$ denote the collection of segments of loops (constructed using the field $\Psi_1$ in the exact same way that the loops in $\overline{\Gamma}$ are constructed using the field $\Psi$) which make an excursion into $A_{\tilde{k}_f}$ from its outer boundary.  Let also $\Gamma_{\tilde{k}_f}^{\mathrm{I}}(\Psi_1)$ denote the collection of segments of loops which cross $A_{\tilde{k}_f}$,  and let $\hat{D}_1$ denote the connected component of $\BD \setminus \overline{\left(\bigcup_{\gamma \in \Gamma_{\tilde{k}_f}^{\mathrm{out}}(\Psi_1)\, \,\cup \, \,\Gamma_{\tilde{k}_f}^{\mathrm{I}}(\Psi_1)} \gamma \right)}$ containing $0$.  Note that $\hat{D}_1 \subseteq \tilde{D}_1$ and that the following holds for all $j,k \in \{1,2,3,4\}$ with $j \neq k$.  The part of the left (resp.\ right) side of $\gamma_j$ that is contained in $\gamma_j \setminus \overline{\left(\bigcup_{\gamma \in \Gamma_{\tilde{k}_f}^{\mathrm{out}}(\Psi_1)\, \,\cup \, \,\Gamma_{\tilde{k}_f}^{\mathrm{I}}(\Psi_1)} \gamma \right)}$ is contained in $\partial \tilde{D}_1$ and its image under $g_1$ is contained in the clockwise (resp.\ counterclockwise) arc of $\partial \BD$ from $\hat{x}_j$ to $\hat{x}_k$.  In particular,  combining the conformal invariance of harmonic measure with its monotonicity under the inclusion $\hat D_1 \subseteq \tilde D_1$,  we obtain that the harmonic measure of the clockwise (resp.\ counterclockwise) arc of $\partial \BD$ from $\hat{x}_j$ to $\hat{x}_k$ in $\BD$ as seen from $0$ is at least $\delta$.  This completes the proof of the claim.

Set also $\tilde{\Psi}_1 = \Psi_1 \circ g_1^{-1}$ and note that the Beurling estimate implies that there exists a universal constant $C \in (0,\infty)$ such that the probability that a complex Brownian motion starting from $0$ exits $\tilde{D}_1$ on $\partial \BD$ is at most $C \varepsilon^{a/2}$.  Indeed,  each of the curves $\gamma_1,\ldots,\gamma_4$ is a connected set joining $\partial \BD$ to $\partial B_{2^{-3\tilde{k}_f-2}}(0)$,  and $\tilde{k}_f \geq a \log(1/\varepsilon)/2$ gives $2^{-3\tilde{k}_f - 2} \leq \tfrac14 \varepsilon^{\frac32 a \log(2)}$,  so that the Beurling estimate bounds the probability in question by a universal constant times $\varepsilon^{\frac34 a \log(2)}$;  the claim follows since $\tfrac34 \log(2) > \tfrac12$.  This is the reason for restricting $\tilde{k}_f$ to lie above $a\log(1/\varepsilon)/2$.  Thus,  combining with \cite[Exercise~2.7]{lawler2008conformally} (which gives that for every connected set $K\subset \overline{\BD}$, the probability that a Brownian motion starting from $0$, stopped when it hits $\partial \BD$, intersects $K$ is at least a universal constant times $\mathrm{diam}(K)$),  we obtain that the following holds by possibly taking $C$ to be larger (in a universal way).  For all $j \in \{1,2,3,4\}$,  there exists $y_j \in \partial \BD$ lying on the counterclockwise arc of $\partial \BD$ from $\hat{x}_j$ to $\hat{x}_{j+1}$ (with the convention that $\hat{x}_5 = \hat{x}_1$) such that the image of the set $\partial \tilde{D}_1 \setminus \bigcup_{j=1}^4 \gamma_j$ under $g_1$ is contained in $\bigcup_{j=1}^4 B_{C \varepsilon^{a/2}}(y_j)$.  Moreover,  possibly decreasing $\tilde\delta_1$ in a way depending only on $\delta$,  we may assume that $\lvert y_j - y_i\rvert \ge \tilde\delta_1$ for $i \neq j$,  since by condition~\eqref{it:harmonic_measure_bounded_below} each of the two boundary arcs of $\partial\BD$ adjacent to $y_j$ has harmonic measure at least $\delta$ as seen from $0$.  Finally,  for all $j \in \{1,2,3,4\}$,  the boundary conditions of $\tilde\Psi_1$ on the counterclockwise arc of $\partial\BD$ from $\hat x_j$ to $\hat x_{j+1}$ are,  outside of $B_{C\varepsilon^{a/2}}(y_j)$,  equal to the boundary conditions of $\Psi$ on the right side of $\eta_j$,  which are also those on the left side of $\eta_{j+1}$;  we denote this common value by $b_j \in \{2\lambda(k-1),2\lambda k,2\lambda(k+1)\}$.

One difficulty in applying Proposition~\ref{prop:bi_chordal_sle_4_two_arm_exponent} is that the boundary conditions of the field $\tilde{\Psi}_1$ are not consistent with the boundary conditions of the field appearing in the statement of Proposition~\ref{prop:bi_chordal_sle_4_two_arm_exponent}.  However,  we will show that Proposition~\ref{prop:bi_chordal_sle_4_two_arm_exponent} can be applied to a random field whose boundary conditions agree with those of $\tilde{\Psi}_1$ outside of small neighborhoods around the points $y_1,y_2,y_3,y_4$,  and that the laws of the two fields,  restricted to the complement of those neighborhoods,  have Radon--Nikodym derivatives with all moments close to $1$.  This is the content of the following lemma.

\begin{lemma}\label{lem:comparison_of_fields}
Fix $b \in (0,a/100)$ deterministic and depending only on $a$.  Suppose that we have the setup described above and let $\hat{\Psi}_1$ denote the GFF on $\BD$ with boundary conditions given by $b_j$ on the counterclockwise arc of $\partial \BD$ from $\hat{x}_j$ to $\hat{x}_{j+1}$ for all $j \in \{1,2,3,4\}$.  Let $\mathcal{Z}_{\varepsilon}$ denote the Radon--Nikodym derivative of the conditional law of $\tilde{\Psi}_1$ restricted to $\BD \setminus \bigcup_{j=1}^4 \overline{B_{\varepsilon^b}(y_j)}$ with respect to the conditional law of $\hat{\Psi}_1$ restricted to $\BD \setminus \bigcup_{j=1}^4 \overline{B_{\varepsilon^b}(y_j)}$,  in both cases given the boundary conditions of $\tilde\Psi_1$ (in particular,  the points $\hat x_j , y_j$ and the values $b_j$).  Then,  for all $p>0$,  we have that
\begin{align*}
\BE[\mathcal{Z}_{\varepsilon}^p ] \to 1 \quad \text{as} \quad \varepsilon \to 0.
\end{align*}
\end{lemma}

\begin{proof}
\stepn{step:comp-fields-overview-proof-strategy}{Overview of the proof strategy}
Write $G_\varepsilon \defeq \BD \setminus \bigcup_{j=1}^4 \overline{B_{\varepsilon^b/2}(y_j)}$.  First,  we note that the Markov property of the GFF implies that we can couple the fields $\tilde{\Psi}_1,\hat{\Psi}_1$ so that they can be expressed as $\tilde{\Psi}_1 = h^0 + \mathcal{H}_{\varepsilon} + \tilde{\mathcal{H}}_1$ and $\hat{\Psi}_1 = h^0 + \mathcal{H}_{\varepsilon} + \hat{\mathcal{H}}_1$ on $G_\varepsilon$,  where $h^0$ is a zero-boundary GFF on $G_\varepsilon$,  $\mathcal{H}_{\varepsilon}$ is a random harmonic function on $G_\varepsilon$ which is independent of $h^0$,  and $\tilde{\mathcal{H}}_1$ (resp.\ $\hat{\mathcal{H}}_1$) is the harmonic function on $\BD$ whose boundary conditions coincide with those of $\tilde{\Psi}_1$ (resp.\ $\hat{\Psi}_1$).  Set $\mathcal{H}_1 = \tilde{\mathcal{H}}_1 - \hat{\mathcal{H}}_1$.

Suppose that we have constructed a smooth function $\Phi_{\varepsilon}$ such that $\Phi_{\varepsilon}$ is equal to $1$ on $\BD \setminus \bigcup_{j=1}^4 \overline{B_{\varepsilon^{b}}(y_j)}$ and vanishes on $\bigcup_{j=1}^4 B_{\varepsilon^b / 2}(y_j)$.  Note that, for $\varepsilon$ small enough that $\varepsilon^b/2 \ge C\varepsilon^{a/2}$, the boundary conditions of $\mathcal{H}_1$ are zero on $\partial \BD \setminus \bigcup_{j=1}^4 \overline{B_{C \varepsilon^{a/2}}(y_j)}$,  and so it holds that $\Phi_{\varepsilon} \mathcal{H}_1 \in H_0^1(G_\varepsilon)$. Therefore,  the proof of \cite[Proposition~3.4]{IG1}  implies that
\begin{align*}
\mathcal{Z}_{\varepsilon} &= \exp\left((\hat{\Psi}_1 ,  \Phi_{\varepsilon} \mathcal{H}_1)_{\nabla} - ||\Phi_{\varepsilon} \mathcal{H}_1||_{\nabla}^2 / 2 \right)=\exp\left((h^0 ,  \Phi_{\varepsilon} \mathcal{H}_1)_{\nabla} - ||\Phi_{\varepsilon} \mathcal{H}_1||_{\nabla}^2 / 2 \right).
\end{align*}
Hence, since $h^0$ is a zero-boundary GFF on $G_\varepsilon$ which is independent of $\mathcal H_1$, for all $p>0$ we have that
\begin{align*}
\BE[\mathcal{Z}_{\varepsilon}^p ] = \BE\left[\exp\left(\frac{p(p-1)}{2} ||\Phi_{\varepsilon} \mathcal{H}_1 ||_{\nabla}^2 \right) \right].
\end{align*}
Therefore,  in order to complete the proof,  we need to find appropriate upper bounds for $||\Phi_{\varepsilon}||_{\infty},  ||\mathcal{H}_1||_{\infty},  ||\Phi_{\varepsilon}||_{\nabla}$,  and $||\mathcal{H}_1||_{\nabla}$,  where all of these norms are taken over $G_\varepsilon$.

\stepn{step:comp-fields-construction}{Construction of $\Phi_{\varepsilon}$}
Fix a smooth function $\chi \colon \BR \to [0,1]$ with $\chi \equiv 0$ on $(-\infty,1/2]$ and $\chi \equiv 1$ on $[1,\infty)$,  and set $\Phi_\varepsilon \defeq \prod_{j=1}^4 \Phi_\varepsilon^j$,  where $\Phi_{\varepsilon}^j(x) \defeq \chi(\varepsilon^{-b}|x - y_j|)$ for $x \in \BC$.  Then $\Phi_\varepsilon$ is smooth,  vanishes on $\bigcup_{j=1}^4 B_{\varepsilon^b/2}(y_j)$,  and equals $1$ on $\BD \setminus \bigcup_{j=1}^4 B_{\varepsilon^b}(y_j)$.  Moreover,  since the four balls $B_{\varepsilon^b}(y_j)$ are pairwise disjoint for small $\varepsilon$ (as $\lvert y_i - y_j \rvert \geq \tilde\delta_1$ for $i \neq j$) and the Dirichlet energy is scale invariant,  we obtain that
\begin{align}\label{eqn:bump_function_dirichlet_energy}
||\Phi_{\varepsilon}||_{\nabla} \lesssim 1,  \quad \text{as} \quad \varepsilon \to 0,
\end{align}
where the implicit constant is universal.  Note also that $||\Phi_{\varepsilon}||_{\infty} \leq 1$ by construction.

\stepn{step:comp-fields-bounding-supremum-dirichlet}{Bounding the supremum and the Dirichlet energy of $\mathcal{H}_1$}
Recall that the boundary conditions of $\mathcal{H}_1$ on $\partial \BD \setminus \bigcup_{j=1}^4 \overline{B_{C \varepsilon^{a/2}}(y_j)}$ are given by the zero function,  and that they lie in $[-2N,2N]$ by Lemma~\ref{lem:field_has_bounded_boundary_conditions};  note that $\mathcal H_1$ is a difference of two harmonic functions whose boundary values lie in a common interval of length $2N$,  so that this bound is unaffected by the recentring by $c$.  Thus,  the reflection principle implies that $\mathcal{H}_1$ can be extended to a harmonic function $\mathcal H_1^*$ on $\BC \setminus \bigcup_{j=1}^4 \overline{B_{C \varepsilon^{a/2}}(y_j)}$ (by possibly taking $C$ to be larger) by setting $\mathcal{H}_1^*(z) = -\mathcal{H}_1(1 / \overline{z})$ for all $z \in \BC \setminus \overline{\BD}$.

For $z \in \BC$ we write $s(z) \defeq \min_{1 \leq i \leq 4} \lvert z - y_i\rvert$.  Since $\mathcal H_1$ vanishes on $\partial\BD \setminus \bigcup_i \overline{B_{C\varepsilon^{a/2}}(y_i)}$ and its boundary values are bounded by $2N$,  the Beurling estimate gives a universal constant $\tilde C \in (0,\infty)$ such that
\begin{align}\label{eqn:bound_on_sup_norm}
\lvert \mathcal{H}_1^*(z)\rvert \leq \tilde{C} N \left(\frac{\varepsilon^{a/2}}{s(z)}\right)^{1/2} \quad \text{whenever} \quad s(z) \geq 2C\varepsilon^{a/2} \quad \text{and} \quad \lvert z \rvert \leq 2.
\end{align}
In particular,  since $s(\cdot) \geq \varepsilon^b/2$ on $G_\varepsilon$ and $b < a/2$,  we have $\lVert\mathcal H_1\rVert_{L^\infty(G_\varepsilon)} \leq 2\tilde CN\varepsilon^{a/4-b/2} \to 0$ as $\varepsilon \to 0$.

Next,  the standard gradient estimate for harmonic functions,  applied on the ball $B_{s(z)/2}(z)$  (on which $\mathcal H_1^*$ is harmonic and on which $s(\cdot) \geq s(z)/2$,  and which is contained in $\overline{B_2(0)}$ when $z \in \BD$,  since $s(z) \leq 2$) gives a universal constant $c^\prime$ such that
\begin{align*}
\lvert \nabla \mathcal{H}_1^*(z)\rvert \leq \frac{4}{s(z)} \sup_{B_{s(z)/2}(z)} \lvert \mathcal{H}_1^*\rvert \leq \frac{c^\prime N \varepsilon^{a/4}}{s(z)^{3/2}} \quad \text{whenever} \quad z \in \BD \quad \text{and} \quad s(z) \geq 4C\varepsilon^{a/2}.
\end{align*}
Partitioning $G_\varepsilon$ according to which of the four points $y_j$ is nearest,  and using that $G_\varepsilon \subseteq \bigcup_j B_{2}(y_j)$ and that $s(\cdot) \geq \varepsilon^b/2$ on $G_\varepsilon$,  we obtain
\begin{align}\label{eqn:dirichlet_energy_bound_harmonic_function}
\int_{G_{\varepsilon}} \lvert \nabla \mathcal{H}_1(z)\rvert^2 \, \rd z \leq 2\pi (c^\prime N)^2 \varepsilon^{a/2} \sum_{j=1}^4 \int_{\varepsilon^b/2}^{2} \frac{\rd r}{r^{2}} \lesssim \varepsilon^{a/2 - b} \quad \text{as} \quad \varepsilon \to 0,
\end{align}
where the implicit constant depends only on $N$,  and $a/2-b>0$ since $b<a/100$.

\stepn{step:comp-fields-conclusion-proof}{Conclusion of the proof}
Finally,  to conclude the proof,  we note that the Cauchy--Schwarz inequality implies that
\begin{align}\label{eqn:cauchy_schwarz_gradient}
||\Phi_{\varepsilon} \mathcal{H}_1||_{\nabla}^2 \leq ||\mathcal{H}_1||_{\infty}^2 ||\Phi_{\varepsilon}||_{\nabla}^2 + ||\mathcal{H}_1||_{\nabla}^2 ||\Phi_{\varepsilon}||_{\infty}^2 + ||\Phi_{\varepsilon}||_{\infty} ||\mathcal{H}_1||_{\infty} (||\Phi_{\varepsilon}||_{\nabla}^2 + ||\mathcal{H}_1||_{\nabla}^2).
\end{align}
Combining~\eqref{eqn:bump_function_dirichlet_energy}--\eqref{eqn:cauchy_schwarz_gradient} with $\lVert\Phi_\varepsilon\rVert_\infty \leq 1$,  we obtain that $\lVert\Phi_\varepsilon\mathcal H_1\rVert_\nabla^2 \lesssim \varepsilon^{a/4-b/2}$,  where the implicit constant depends only on $N$.  Since $a/4-b/2>0$,  the right-hand side is deterministic and tends to $0$,  so the proof is complete in view of the expression for the moments of $\mathcal Z_\varepsilon$ obtained in \Cref{step:comp-fields-overview-proof-strategy}.
\end{proof}

Now,  we let $\tilde{\gamma}_j$ be the level line of $\tilde{\Psi}_1$ (resp.\ $-\tilde{\Psi}_1$,  according to whether $\gamma_j$ is a level line of $\Psi_1$ or of $-\Psi_1$) starting from $\hat{x}_j$ having the same height as $\gamma_j$ for all $j \in \{1,2,3,4\}$;  equivalently,  $\tilde{\gamma}_j$ is the image under $g_1$ of the continuation of $\gamma_j$ beyond the time at which it was stopped.  In the following lemma,  we show that none of the curves $\tilde{\gamma}_1,\tilde{\gamma}_2,\tilde{\gamma}_3,\tilde{\gamma}_4$ gets too close to any of the points in $\{y_1,y_2,y_3,y_4\}$.

\begin{lemma}\label{lem:level_lines_do_not_hit_small_balls}
Suppose that we have the same setup as in Lemma~\ref{lem:comparison_of_fields}. Then,  there exists a deterministic constant $\beta \in (0,b)$ (depending only on $b$ and $N$) such that off an event with probability $O(\varepsilon^{\beta})$ as $\varepsilon \to 0$ (at a rate that depends only on $a,b,\delta$,  and $N$),  we have that
\begin{align*}
\left(\bigcup_{j=1}^4 \tilde{\gamma}_j \right) \cap \left(\bigcup_{i=1}^4 B_{\varepsilon^b}(y_i)\right) = \emptyset.
\end{align*}
\end{lemma}

\begin{proof}
\stepn{step:ll-small-balls-overview}{Overview}
We will prove the claim of the lemma by constructing shields of level lines of $\hat{\Psi}_1$ (which we shall then transfer to $\tilde{\Psi}_1$ by means of Lemma~\ref{lem:comparison_of_fields}) on annuli of the form $\BD \cap (B_{2^{-3n-1}}(y_j) \setminus \overline{B_{2^{-3n-2}}(y_j)})$ for $j=1,2,3,4$,  which disconnect $B_{\varepsilon^b}(y_j)$ from $\{\hat{x}_1,\hat{x}_2,\hat{x}_3,\hat{x}_4\}$.  More precisely,  we will show that with very high probability as $\varepsilon \to 0$,  there are many values of $n$ for which a level line of $\tilde{\Psi}_1$ of a certain height disconnects the inner from the outer boundary of the annulus at scale $n$ before exiting the annulus.  In particular,  the level line will disconnect $B_{\varepsilon^b}(y_j)$ from $\{\hat{x}_1,\hat{x}_2,\hat{x}_3,\hat{x}_4\}$ and thus the level line interaction rules will imply that none of the curves in $\{\tilde{\gamma}_1,\tilde{\gamma}_2,\tilde{\gamma}_3,\tilde{\gamma}_4\}$ enters $B_{\varepsilon^b}(y_j)$.  This will complete the proof of the lemma.

\stepn{step:ll-small-balls-constructing-shields-level}{Constructing shields of level lines of $\hat{\Psi}_1$}

\substepn{step:lls-setup}{Setup at scale $n$} We will follow the same strategy as in the proof of Lemma~\ref{prop:exactly_four_crossings}.  We fix $j \in \{1,2,3,4\}$ and $n \in \BN$ such that
\begin{align*}
n \in \left(\frac{b}{6 \log(2)} \log(1 / \varepsilon) ,  -1 + \frac{b}{3 \log(2)} \log(1 / \varepsilon) \right),
\end{align*}
and note that $\hat{\Psi}_1$ can be decomposed as $\hat{\Psi}_1 = h_{j,n}^0 + \mathcal{H}_{j,n}$,  where $h_{j,n}^0$ (resp.\  $\mathcal{H}_{j,n}$) is a GFF (resp.\ harmonic function) on $\BD \cap B_{2^{-3n}}(y_j)$,  and such that the boundary conditions of $h_{j,n}^0$ are constant and given by $b_j$.  Note that the choice of $n$ implies that the corresponding annulus at scale $n$ is contained in $B_{\varepsilon^{b/2}}(y_j) \setminus \overline{B_{\varepsilon^b}(y_j)}$.  Fix also $\tilde{M}>1$ sufficiently large (to be chosen and depending only on $b$ and $N$) and for $n \in \BN$ as above,  we let $E_{j,n}^{\tilde{M}}$ denote the event that
\begin{align*}
\sup_{z \in \BD \cap B_{2^{-3 n - 1/2}}(y_j)} |\mathcal{H}_{j,n}(z)-\mathcal{H}_{j,n}(y_j)| \leq \tilde{M}.
\end{align*}

\substepn{step:lls-good-scales}{Many good scales $n$} By \cite[Lemma 4.2]{MultiSLE}, we can choose $\tilde{M}>1$ sufficiently large (depending only on $b$ and $N$) such that off an event with probability $O(\varepsilon^3)$ as $\varepsilon \to 0$ (at a rate that depends only on $b,  \tilde{M}$,  and $N$),  we have that $E_{j,n}^{\tilde{M}}$ occurs for at least $b \log(1/\varepsilon) / 100$ values of $n$.

Indeed, by the Markov property of the GFF, we can write $\hat{\Psi}_1= \mathcal{H}+ \Psi^0$ where $\mathcal{H}$ is harmonic on $\BD$ and its boundary conditions are piecewise constant and change values four times. Also, the oscillation of $\mathcal{H}$ is bounded from above by $4\lambda$, since the four boundary values $b_j$ lie in $\{2\lambda(k-1),2\lambda k,2\lambda(k+1)\}$; this is all that is needed, since $E_{j,n}^{\tilde M}$ is defined in terms of a difference of values of $\mathcal H_{j,n}$. Moreover $\Psi^0$ is a zero boundary GFF on $\BD$ which is independent of $\mathcal{H}$. One can also decompose $\Psi^0= (h^0_{j,n}-b_j) + \tilde{\mathcal{H}}_{j,n}$ where $\tilde{\mathcal{H}}_{j,n}$ is harmonic on $\BD \cap B_{2^{-3n}}(y_j)$ with zero boundary conditions on $\partial \BD\cap B_{2^{-3n}}(y_j)$ and independent of $h^0_{j,n}$, so that $\mathcal{H}_{j,n}= \mathcal{H}+ \tilde{\mathcal{H}}_{j,n}-b_j$. By conformally mapping $\BD$ to $\BH$ and applying \cite[Lemma 4.2]{MultiSLE}, we get the analog of $E_{j,n}^{\tilde{M}}$ with $\tilde{\mathcal{H}}_{j,n}$ replacing $\mathcal{H}_{j,n}$, yielding the statement of the previous paragraph.

\substepn{step:lls-crossing}{The crossing event for $h_{j,n}^0$} Fix $u \in (0,\lambda)$ and let $x_{j,n}$ (resp.\  $y_{j,n}$) be the midpoint of the arc of $\partial \BD \cap \left(B_{2^{-3n-1}}(y_j) \setminus \overline{B_{2^{-3n-2}}(y_j)}\right)$ contained in the clockwise (resp.\  counterclockwise) arc of $\partial \BD$ from $y_j$ to $\hat{x}_j$ (resp.\  $\hat{x}_{j+1}$).  Let also $\hat{\eta}_{j,n}$ (resp.\ $\hat{\eta}_{j,n}^{0}$) denote the level line of $\hat{\Psi}_1$ (resp.\  $h_{j,n}^0$) starting from $x_{j,n}$ and targeted at $y_{j,n}$ with height $b_j-u$.  Note that the law of $\hat{\eta}_{j,n}^0$ is given by that of an $\SLE_4(-1 - u / \lambda ;  -1+u/\lambda)$ process in $\BD \cap B_{2^{-3n}}(y_j)$ from $x_{j,n}$ to $y_{j,n}$,  where the force points are located at $x_{j,n}^-$ and $x_{j,n}^+$ respectively.  Note also that both of the weights of the force points lie in $(-2,0)$.  Therefore,  \cite[Lemma~2.5]{miller2017intersections} combined with the translation and scale invariance of the zero boundary GFF implies that there exists a universal constant $q \in (0,1)$ such that with probability at least $q$,  we have that $\hat{\eta}_{j,n}^0$ intersects the counterclockwise arc of $\partial \BD$ from $y_j$ to $\hat{x}_{j+1}$ before exiting $B_{2^{-3n -1}}(y_j) \setminus \overline{B_{2^{-3n-2}}(y_j)}$ for the first time.  Note that the latter event is a.s.\ determined by the restriction of $h_{j,n}^0$ to $\BD \cap (B_{2^{-3n-1}}(y_j) \setminus \overline{B_{2^{-3n-2}}(y_j)})$.

\substepn{step:lls-transfer}{Transferring to $\hat{\Psi}_1$ and iterating} Next, we bound the $p$-th moment of the Radon--Nikodym derivative $\mathcal{Z}_{j,n}$ of the conditional law of $\hat{\Psi}_1$ with respect to that of $h^{0}_{j,n}$ on the event $E^{\tilde{M}}_{j,n}$ for all $p \in \BR$ as in \cite[Lemma~4.1]{GeoLQGnSLE} by a constant $\tilde{c}(p, \tilde{M})$. In particular, we can find $q_0 =q_0(q, \tilde{M})>0$ such that $\BP[1/\mathcal{Z}_{j,n} \le 1/q_0  \vert \hat{\Psi}_1\vert_{\BD \setminus \overline{B_{2^{-3n-1}}(y_j)}}] \ge 1- q/2$ for all $n$ such that $E_{j,n}^{\tilde{M}}$ occurs.

Therefore, for all $n$ such that $E_{j,n}^{\tilde{M}}$ occurs, we have
\begin{align*}
	\BP&\left[\left.\hat{\eta}_{j,n} \text{ reaches the arc from }y_j\text{ to }\hat{x}_{j+1} \right\vert \hat{\Psi}_1\vert_{\BD \setminus \overline{B_{2^{-3n-1}}(y_j)}}\right]  \\
	&= \BE\left[ \left.\mathcal{Z}_{j,n} \one_{\hat{\eta}^0_{j,n} \text{ reaches the arc from }y_j\text{ to }\hat{x}_{j+1}}  \right\vert \hat{\Psi}_1\vert_{\BD \setminus \overline{B_{2^{-3n-1}}(y_j)}}\right] \\
	&\ge q_0 \cdot \frac{q}{2}.
\end{align*}

 Thus, iterating over the at least $b\log(1/\varepsilon)/100$ scales $n$ for which $E_{j,n}^{\tilde M}$ occurs, we obtain that there exists a constant $\beta \in (0,b)$ depending only on $b,q,N$,  and $\tilde{M}$ such that off an event with probability $O(\varepsilon^{\beta})$ as $\varepsilon \to 0$ (at a rate depending only on $b,q,N$,  and $\tilde{M}$),  we have that there exists $n \in \BN$ such that
\begin{align*}
n \in \left(\frac{b}{6 \log(2)} \log(1 / \varepsilon) ,  -1 + \frac{b}{3 \log(2)} \log(1 / \varepsilon) \right),
\end{align*}
and $\hat{\eta}_{j,n}$ intersects the counterclockwise arc of $\partial \BD$ from $y_j$ to $\hat{x}_{j+1}$ before exiting $B_{2^{-3n -1}}(y_j) \setminus \overline{B_{2^{-3n-2}}(y_j)}$ for the first time.

\stepn{step:ll-small-balls-conclusion-proof}{Conclusion of the proof}
Note that the event described in the last paragraph of \Cref{step:lls-transfer} is a.s.\ determined by the restriction of $\hat{\Psi}_1$ to $\BD \setminus \bigcup_{i=1}^4 \overline{B_{\varepsilon^b}(y_i)}$.  Therefore,  we can define the exact same event with the field $\tilde{\Psi}_1$ in place of $\hat{\Psi}_1$.  The proof is then complete by combining the result in the last paragraph of \Cref{step:lls-transfer} with Lemma~\ref{lem:comparison_of_fields} and the level line interaction rules,  and taking a union bound over $j \in \{1,2,3,4\}$.
\end{proof}

Finally,  we are ready to prove Lemma~\ref{lem:first_stage_cond_prob}.

\begin{proof}[Proof of Lemma~\ref{lem:first_stage_cond_prob}.]
\stepn{step:first-stage-two-inclusions-balls}{Two inclusions of balls}  Recall that $\tilde{k}_f$ is the smallest $n \in [(a/2)\log(1/\varepsilon), a \log(1/\varepsilon)]_\BZ$ such that $E_n(\Psi_1)$ occurs. First, we note that since $\tilde{k}_f \leq a \log(1 / \varepsilon)$,  we obtain that
\begin{align*}
\text{dist}(0,\partial \tilde{D}_1) \geq \frac{1}{4} \varepsilon^{3 a \log(2)}.
\end{align*}
Thus,  combining with Koebe's $\frac{1}{4}$-theorem,  we have that
\begin{align*}
|(g_1^{-1})'(0)| \geq \frac{1}{4} \text{dist}(0,\partial \tilde{D}_1) \geq \frac{1}{16} \varepsilon^{3 a \log(2)}. 
\end{align*}
Moreover,  combining with \cite[Theorem~3.23]{lawler2008conformally} gives that
\begin{align*}
|g_1^{-1}(z)| \geq \frac{|z|}{64} \varepsilon^{3 a \log(2)} \quad \text{for all} \quad z \in \BD.
\end{align*}
In particular,  we have that 
\begin{align*}
|g_1^{-1}(z)| \geq \frac{\varepsilon^{\rho / 2}}{64} \varepsilon^{3 a \log(2)} \quad \text{for all} \quad z \in \partial B_{\varepsilon^{\rho / 2}}(0).
\end{align*}
It follows that
\begin{align}\label{eqn:first_inclusion}
B_{\varepsilon^{\frac{\rho}{2} + 3a\log(2)} / 64}(0) \subseteq g_1^{-1}(B_{\varepsilon^{\rho / 2}}(0)).
\end{align}

Similarly,  we have that $k_f \leq a \log(1 / \varepsilon)$ and so 
\begin{align*}
\text{dist}(0,\partial D_1) \geq \frac{1}{4} \varepsilon^{3a\log(2)}.
\end{align*}
Therefore,  combining Koebe's $\frac{1}{4}$-theorem with \cite[Theorem~3.23]{lawler2008conformally} as before and combining with~\eqref{eqn:first_inclusion},  we obtain that 
\begin{align}\label{eqn:second_inclusion}
B_{\varepsilon^{\frac{\rho}{2} + 6a\log(2)} / 64^2}(0) \subseteq f_1^{-1}(g_1^{-1}(B_{\varepsilon^{\rho/2}}(0))).
\end{align}

\stepn{step:first-stage-conclusion}{Conclusion}  Note that~\eqref{eqn:second_inclusion} implies that if all of the paths $\eta_{j,1}$ for $j=1,2,3,4$,  intersect $\partial B_{\varepsilon^{\frac{\rho}{2} + 6a\log(2)} / 64^2}(0)$ before intersecting any of the other paths in $\{\eta_{1,1},\eta_{2,1},\eta_{3,1},\eta_{4,1}\}$,  then we have that all of the paths $\tilde{\gamma}_j$ for $j=1,2,3,4$,  intersect $\partial B_{\varepsilon^{\rho / 2}}(0)$ before intersecting any of the other paths in $\{\tilde{\gamma}_1,\tilde{\gamma}_2,\tilde{\gamma}_3,\tilde{\gamma}_4\}$.  Let now $\hat\gamma_1,\hat\gamma_2,\hat\gamma_3,\hat\gamma_4$ be defined in the same way as $\tilde\gamma_1,\ldots,\tilde\gamma_4$ but with the field $\hat\Psi_1$ of Lemma~\ref{lem:comparison_of_fields} in place of $\tilde\Psi_1$.  Then,  after subtracting the constant $2\lambda k$ from $\hat\Psi_1$,  the boundary conditions of $\hat\Psi_1$ and the heights of $\hat\gamma_1,\ldots,\hat\gamma_4$ are exactly as in \Cref{prop:bi_chordal_sle_4_two_arm_exponent},  with the same value of $(\sigma_1,\sigma_3)$ as in Section~\ref{subsubsec:outline_and_setup};  hence the joint law of $(\hat\gamma_1,\hat\gamma_2,\hat\gamma_3,\hat\gamma_4)$ is that of $(P_1,P_2,P_3,P_4)$,  and \Cref{prop:bi_chordal_sle_4_two_arm_exponent} bounds by $C\varepsilon^{\rho}$ the probability that all four of them reach $\partial B_{\varepsilon^{\rho/2}}(0)$ without intersecting one another.  By Lemma~\ref{lem:level_lines_do_not_hit_small_balls},  off an event of probability $O(\varepsilon^\beta)$ the curves $\tilde\gamma_1,\ldots,\tilde\gamma_4$ do not enter $\bigcup_j B_{\varepsilon^b}(y_j)$,  so that this event is determined by $\tilde\Psi_1$ restricted to $\BD \setminus \bigcup_j\overline{B_{\varepsilon^b}(y_j)}$;  Lemma~\ref{lem:comparison_of_fields} therefore transfers the above bound from $\hat\gamma_1,\ldots,\hat\gamma_4$ to $\tilde\gamma_1,\ldots,\tilde\gamma_4$,  at the cost of a constant factor.  Adding the error terms $O(\varepsilon^3)$ from Lemma~\ref{lem:tips_are_separated} and $O(\varepsilon^\beta)$ from Lemma~\ref{lem:level_lines_do_not_hit_small_balls},  we obtain a bound of the form $C_0(\varepsilon^{\rho} + \varepsilon^3 + \varepsilon^\beta)$,  where $C_0$ depends only on $a,\delta$,  and $N$.  Since $\rho$ was chosen small enough that $\rho < \beta$,  each of the three exponents is strictly larger than $\rho(1-a)$,  and hence this bound is at most $\varepsilon^{\rho(1-a)}$ once $\varepsilon_0$ is taken small enough in a way depending only on $a,\delta,N$,  and $\rho$.  This completes the proof of the lemma.
\end{proof}

\subsubsection{Proving Lemma~\ref{lem:bound_on_number_of_iterations}.}

Next,  we focus on proving Lemma~\ref{lem:bound_on_number_of_iterations}.  First,  we will prove Lemma~\ref{lem:inductive_stage_cond_prob}.  Let $\rho$ and $\zeta$ be as fixed in Section~\ref{subsubsec:outline_and_setup}.  For all $n \in \BN$,  we let $D_{n+1}$ denote the connected component of $D_1 \setminus \bigcup_{j=1}^4 \eta_{j,n+1}$ containing $0$ and let $f_{n+1}$ denote the conformal transformation mapping $D_{n+1}$ onto $\BD$ such that $f_{n+1}(0) = 0$ and $f_{n+1}'(0) > 0$.  Set $\Psi_{n+1}:=\Psi \circ f_{n+1}^{-1}$ and note that the boundary conditions of $\Psi_{n+1}$ on $\partial \BD$ lie in $[c_{n+1}-N,c_{n+1}+N]$ for some $c_{n+1} \in 2\lambda\BZ$,  and they are piecewise constant and change only countably many times.  Let $x_{1,n+1},x_{2,n+1},x_{3,n+1},x_{4,n+1}$ denote the endpoints of $\eta_{1,n+1},\eta_{2,n+1},\eta_{3,n+1},\eta_{4,n+1}$ respectively,  and set $\tilde{x}_{j,n+1}:=f_{n+1}(x_{j,n+1})$ for $j=1,2,3,4$.  As in the $n=0$ case,  we have that the points $\tilde{x}_{1,n+1},\tilde{x}_{2,n+1},\tilde{x}_{3,n+1},\tilde{x}_{4,n+1}$ are ordered counterclockwise on $\partial \BD$.

Next,  we state the analog of Lemma~\ref{lem:tips_are_separated}.

\begin{lemma}\label{lem:tips_are_separated_inductive_step}
Fix $n \in \BN$.  Then,  we have that the following is true off an event with probability at most $\varepsilon^{2\rho}$,  for all $\varepsilon \in (0,\varepsilon_0)$,  where $\varepsilon_0 \in (0,1)$ depends only on $a,\delta,N$,  $\rho$,  and $\zeta$.  There exists $m \in [1 , \zeta \log(1/\varepsilon) / (4 \log(2))]_{\BZ}$ such that $E_m(\Psi_{n+1})$ occurs.
\end{lemma}

\begin{proof}
The statement of the lemma follows from the exact same argument used to prove Lemma~\ref{lem:tips_are_separated}.  The only difference is that there are now only $\asymp \zeta\log(1/\varepsilon)$ scales available rather than $\asymp a \log(1/\varepsilon)$,  so that the failure probability is at most $(1-p)^{\zeta\log(1/\varepsilon)/(4\log 2)}$,  where $p$ is the lower bound on the conditional probability of $E_m(\Psi_{n+1})$ at a given scale.  Since $\zeta = s\rho/2$,  the latter quantity equals $\varepsilon^{s \rho \log(1/(1-p)) / (8\log 2)}$,  and since $p$ may be taken as close to $1$ as we wish upon decreasing $\delta$ and increasing $N$ in a way depending only on $a$ and $s$,  we may arrange that $s\log(1/(1-p)) \geq 16 \log 2$,  so that this is at most $\varepsilon^{2\rho}$;  recall from Section~\ref{subsubsec:outline_and_setup} that $\delta$ and $N$ are chosen after $a$ and $s$ but before $\rho$ and $\zeta$,  so that there is no circularity.
\end{proof}

For the rest of the section,  we will assume that the event in the statement of Lemma~\ref{lem:tips_are_separated_inductive_step} holds and let $\tilde{k}_{f,n+1}$ denote the smallest integer $m \in [1 , \zeta \log(1/\varepsilon) / (4 \log(2))]_{\BZ}$ such that $E_m(\Psi_{n+1})$ occurs.  We also let $\gamma_{1,n+1},\gamma_{2,n+1},\gamma_{3,n+1},\gamma_{4,n+1}$ denote the corresponding level lines starting from $\tilde{x}_{1,n+1},\tilde{x}_{2,n+1},\tilde{x}_{3,n+1},\tilde{x}_{4,n+1}$ respectively,  and note that as in the case $n=0$,  the level line interaction rules imply that the curves $\gamma_{1,n+1},\gamma_{2,n+1},\gamma_{3,n+1},\gamma_{4,n+1}$ are distinct when they are stopped at the first time that they intersect $\partial B_{2^{-3\tilde{k}_{f,n+1} - 2}}(0)$.

Let $\tilde{D}_{n+1}$ denote the connected component containing $0$ of the complement in $\BD$ of the union of $\gamma_{1,n+1},\gamma_{2,n+1},\gamma_{3,n+1},\gamma_{4,n+1}$ when they are stopped at the first time that they intersect $\partial B_{2^{-3\tilde{k}_{f,n+1} - 2}}(0)$.  Then,  we let $g_{n+1}$ denote the conformal transformation mapping $\tilde{D}_{n+1}$ onto $\BD$ such that $g_{n+1}(0) = 0$ and $g_{n+1}'(0) > 0$. Let $\check{x}_{j,n+1}$ be the tip of the curve $\gamma_{j,n+1}$ stopped at the first time that it intersects $\partial B_{2^{-3\tilde{k}_{f,n+1} - 2}}(0)$ for all $j \in \{1, 2, 3, 4\}$.  Set also $\hat{x}_{j,n+1}:=g_{n+1}(\check{x}_{j,n+1})$, $\tilde{\gamma}_{j,n+1}:=g_{n+1}(\gamma_{j,n+1})$ and $\tilde\Psi_{n+1} \defeq \Psi_{n+1} \circ g_{n+1}^{-1}$ for all $j=1,2,3,4$.  Note that combining the Beurling estimate and the fact that $E_{\tilde{k}_{f,n+1}}(\Psi_{n+1})$ occurs as in the case $n=0$,  we obtain that the following hold.  We have that $|\hat{x}_{j,n+1} - \hat{x}_{k,n+1}| \geq \tilde{\delta}_1$ for all $j,k \in \{1,2,3,4\}$ such that $j\neq k$,  and for all $j \in \{1,2,3,4\}$,  there exists a point $\hat{y}_{j,n+1} \in \partial \BD$ lying on the counterclockwise arc of $\partial \BD$ from $\hat{x}_{j,n+1}$ to $\hat{x}_{j+1,n+1}$ (with the convention that $\hat{x}_{5,n+1} = \hat{x}_{1,n+1}$) such that the image of $\partial D_{n+1} \setminus \bigcup_{j=1}^4 \eta_{j,n+1}$ under $g_{n+1} \circ f_{n+1}$ is contained in $\bigcup_{j=1}^4 B_{C \varepsilon^{a/2}}(\hat{y}_{j,n+1})$.  For the latter,  the relevant Beurling bound comes from the fact that the curves $\eta_{1,n+1},\ldots,\eta_{4,n+1}$ reach $\partial B_{r_{n+1}}(0)$ with $r_{n+1} \leq \varepsilon^{\rho/2+6a\log(2)}/64^2$,  while $\dist(0,\partial D_1) \geq \varepsilon^{3a\log(2)}/4$;  hence the harmonic measure from $0$ in $D_{n+1}$ of $\partial D_{n+1} \setminus \bigcup_j \eta_{j,n+1}$ is at most a universal constant times $\varepsilon^{(\rho/2+3a\log(2))/2} \leq \varepsilon^{a/2}$.

\begin{proof}[Proof of Lemma~\ref{lem:inductive_stage_cond_prob}.]
The statement of the lemma follows from the same argument used to prove Lemma~\ref{lem:first_stage_cond_prob}.  In particular,   Lemma~\ref{lem:tips_are_separated} is replaced by Lemma~\ref{lem:tips_are_separated_inductive_step},  and we have statements which are analogous to those of Lemmas~\ref{lem:comparison_of_fields} and~\ref{lem:level_lines_do_not_hit_small_balls}, where the field $\tilde{\Psi}_1$ is replaced by $\tilde{\Psi}_{n+1}$ and the paths $\tilde{\gamma}_1,\tilde{\gamma}_2,\tilde{\gamma}_3,\tilde{\gamma}_4$ are replaced by $\tilde{\gamma}_{1,n+1},\tilde{\gamma}_{2,n+1},\tilde{\gamma}_{3,n+1},\tilde{\gamma}_{4,n+1}$ respectively.  The error term coming from Lemma~\ref{lem:tips_are_separated_inductive_step} is now $\varepsilon^{2\rho}$ in place of $\varepsilon^3$,  and $2\rho > \rho(1-a)$,  so that the conclusion of the previous proof applies verbatim.
\end{proof}

\begin{proof}[Proof of Lemma~\ref{lem:bound_on_number_of_iterations}.]
The statement of the lemma follows from combining Lemmas~\ref{lem:event_happens_at_many_scales}, \ref{lem:first_stage_cond_prob} and~\ref{lem:inductive_stage_cond_prob}:  off an event of probability at most $\varepsilon^3$ the event of Lemma~\ref{lem:event_happens_at_many_scales} occurs,  and on that event the conditional probability of $F_\varepsilon$ is,  by Lemmas~\ref{lem:first_stage_cond_prob} and~\ref{lem:inductive_stage_cond_prob},  at most a product of $n_\varepsilon$ factors each of which is at most $\varepsilon^{\rho(1-a)}$.
\end{proof}

\subsection{Conclusion of the proof.}
\label{subsec:conlcusion_of_proof}

\subsubsection{Proof of Proposition~\ref{prop:nested_cle_4_two_arm_exponent} assuming Proposition~\ref{prop:bi_chordal_sle_4_two_arm_exponent}.}

Now,  we are ready to prove Proposition~\ref{prop:nested_cle_4_two_arm_exponent}.  Suppose that there are at least four crossings of $C_\varepsilon$ made by loops in $\overline\Gamma$;  then either there are two distinct loops $\ell,\ell^\prime \in \overline\Gamma$ which both cross $C_\varepsilon$,  or a single loop $\ell$ makes at least four crossings of $C_\varepsilon$.  As we mentioned in Section~\ref{subsubsec:outline_and_setup},  in the former case it suffices to consider the case that either $\mathop{\mathrm{int}}(\ell) \cap \mathop{\mathrm{int}}(\ell') = \emptyset$ and $\ell,\ell'$ have the same nesting generation,  or one of the loops surrounds the other and their nesting generations differ by one.  Indeed,  suppose that $\mathop{\mathrm{int}}(\ell) \cap \mathop{\mathrm{int}}(\ell') = \emptyset$ and that the two loops have different nesting generations.  Without loss of generality,  we can assume that the nesting generation of $\ell$ is larger than that of $\ell'$.  Then,  there exists a loop $\tilde{\ell} \in \overline{\Gamma}$ such that $\tilde{\ell}$ surrounds $\ell$,  and the loops $\ell'$ and $\tilde{\ell}$ have the same level of nesting.  Thus,  since both $\ell, \ell'$ cross $C_{\varepsilon}$,  the same holds for the loops $\ell' ,  \tilde{\ell}$.  (Indeed,  $\mathop{\mathrm{int}}(\tilde\ell)$ contains the connected set $\ell \cup \mathop{\mathrm{int}}(\ell)$,  which joins the two boundary circles of $C_\varepsilon$,  so that either $\tilde\ell$ crosses $C_\varepsilon$ or $\tilde\ell$ surrounds $\overline{B_{1/2}(0)}$;  the latter is impossible,  since it would force $\ell^\prime \subseteq \mathop{\mathrm{int}}(\tilde\ell)$ and hence $\tilde\ell$ and $\ell^\prime$ would not have the same nesting generation.)  Furthermore,  if $\ell,\ell'$ are such that $\mathop{\mathrm{int}}(\ell) \subseteq \mathop{\mathrm{int}}(\ell')$ and the levels of nesting of $\ell$ and $\ell'$ differ by more than $1$,  then there exists a loop $\tilde{\ell} \in \overline{\Gamma}$ such that $\mathop{\mathrm{int}}(\ell) \subseteq \mathop{\mathrm{int}}(\tilde{\ell}) \subseteq \mathop{\mathrm{int}}(\ell')$ and the levels of nesting of $\tilde{\ell}$ and $\ell'$ differ by $1$.  Hence, if both $\ell$ and $\ell'$ cross $C_{\varepsilon}$,  then the same holds for $\tilde{\ell}$ and $\ell'$.

\begin{proof}[Proof of Proposition~\ref{prop:nested_cle_4_two_arm_exponent}.]
\stepn{step:two-arm-reduction-event}{Reduction to the event $F_\varepsilon$}  Let $\ell,\ell'$ be two loops in $\overline{\Gamma}$ as in the previous paragraph.  We claim that it suffices to show that,  on the event that there are at least four crossings of $C_\varepsilon$ made by loops in $\overline\Gamma$, the event $F_{\varepsilon}$ in the statement of Lemma~\ref{lem:bound_on_number_of_iterations} occurs.  Indeed,  note that by~\eqref{eqn:number_of_iterations},  for all sufficiently small $\varepsilon$ we have that
\begin{align*}
\rho n_\varepsilon(1-a) \geq \rho R (1-a) = \frac{2(1-a) (1-6a \log(2) - \rho / 2)}{1 + 2 \zeta / \rho}.
\end{align*}
Since the right-hand side depends on $\zeta$ and $\rho$ only through the ratio $s = 2\zeta/\rho$ and through $\rho$ itself,  since $\rho < a$,  and since $a \in (0,1/100)$ and $s \in (0,1)$ were arbitrary,  the right-hand side converges to $2$ if we first take $s \to 0$ and then $a \to 0$.  We may therefore apply Lemma~\ref{lem:bound_on_number_of_iterations} together with the fact that the total number of crossings in $\Gamma_{k_f}^{\text{I}}$ is at most $N$ (so that there are at most $4\binom{N}{4}$ counterclockwise labelings of a quadruple $\eta_1,\eta_2,\eta_3,\eta_4$,  and a union bound costs only a constant factor) to complete the proof.  For the rest of the proof,  we fix constants $a,\zeta,$ and $\rho$ as in the statement of Lemma~\ref{lem:bound_on_number_of_iterations} and show that the event $F_{\varepsilon}$ occurs.

\stepn{step:two-arm-four-crossings-their}{The four crossings and their boundary conditions}  Suppose that both $\ell$ and $\ell'$ cross $C_{\varepsilon}$. There exists $k \in \BZ$ such that if the nesting generations of $\ell$ and $\ell'$ are the same, then the boundary conditions of $\Psi$ on the boundaries of the unbounded connected components of $\BC \setminus \overline{\mathop{\mathrm{int}}(\ell)}$ and $\BC \setminus \overline{\mathop{\mathrm{int}}(\ell')}$ are both given by $2\lambda k$ (hence the boundary conditions of $\Psi|_{\mathop{\mathrm{int}}(\ell)}$ and $\Psi|_{\mathop{\mathrm{int}}(\ell')}$ lie in $\{2\lambda (k-1) ,  2\lambda (k+1)\}$),  and if one of the loops surrounds the other (we can assume that $\ell$ surrounds $\ell'$),  it holds that the boundary conditions of $\Psi|_{\mathop{\mathrm{int}}(\ell)}$ and those of $\Psi$ on the boundary of the unbounded connected component of $\BC \setminus \overline{\mathop{\mathrm{int}}(\ell')}$ are both given by $2\lambda k$.  Moreover,  in the latter case,  we have that both the boundary conditions of $\Psi$ on the boundary of the unbounded connected component of $\BC \setminus \mathop{\mathrm{int}}(\ell)$ and the boundary conditions of $\Psi|_{\mathop{\mathrm{int}}(\ell')}$ lie in $\{2\lambda (k-1) ,  2\lambda (k+1)\}$.  

\stepn{step:two-arm-conclusion}{Conclusion}  Note that the construction of the coupling $(\Psi,\overline{\Gamma})$ combined with the level line interaction rules implies that each crossing of $C_{\varepsilon}$ made by $\ell$ or by $\ell'$ has to merge with one of the crossings in $\Gamma_{k_f}^{\text{I}}$.  Therefore,  we obtain that there exist crossings $\eta_1,\eta_2,\eta_3,\eta_4$ of $\Gamma_{k_f}^{\text{I}}$ as above such that $\eta_1$ and $\eta_2$ (resp.\ $\eta_3$ and $\eta_4$) are traced by $\ell'$ (resp.\ $\ell$) while making one crossing from $\partial B_{1/2}(0)$ to $\partial B_{\varepsilon}(0)$ and another crossing from $\partial B_{\varepsilon}(0)$ to $\partial B_{1/2}(0)$.  Also,  in the case that $\ell$ makes four crossings of $C_{\varepsilon}$,  then there exist $k \in \BZ$ and crossings $\eta_1,\eta_2,\eta_3,\eta_4$ as above such that all of $\eta_1,\eta_2,\eta_3,\eta_4$ are traced by $\ell$.  Therefore,  in either case,  the loops must agree with the corresponding level line crossings up until the latter intersect $\partial B_{\varepsilon}(0)$ for the first time.  Since the choice of $n_{\varepsilon}$ implies that
\begin{align*}
B_{\varepsilon}(0) \subseteq B_{(\varepsilon^{\rho / 2 + 6a\log(2)} / 64^2) (\varepsilon^{\rho / 2 + \zeta} / 16^2)^{n_{\varepsilon}-1}}(0),
\end{align*}
we obtain that the event $F_{\varepsilon}$ occurs;  here we use that $\eta_1$ and $\eta_2$ (resp.\ $\eta_3$ and $\eta_4$) are the two strands of the same loop $\ell^\prime$ (resp.\ $\ell$),  so that they do not meet before that loop has completed its crossings of $C_\varepsilon$,  while $\eta_1,\eta_2$ are disjoint from $\eta_3,\eta_4$ since $\ell$ and $\ell^\prime$ are distinct loops.  This completes the proof of the proposition.
\end{proof}

\subsubsection{Proof of Proposition~\ref{prop:nested_cle_4_three_arm_exponent} assuming Proposition~\ref{prop:bi_chordal_sle_4_two_arm_exponent}.}

Finally,  we prove Proposition~\ref{prop:nested_cle_4_three_arm_exponent}.
As explained in Section~\ref{subsec:outline_of_proof},  we will use a strategy similar to that of Section~\ref{subsec:conformal_rectangles}.  We use the same setup and the same constants $a,\rho,\zeta$ as before,  and set $\theta \defeq 1-a/200 \in (1-a/100,1)$;  note that $\theta$ is determined by $a$ alone,  so that the constants $\delta$ and $N$ may be taken to depend on $\theta$ as well.  Then,  we have the following.

\begin{lemma}\label{lem:good_event_many_scales}
By possibly taking $\delta \in (0,1)$ to be smaller and $N \in \BN$ to be larger in the definition of the events $E_n(\Psi)$,  we have that the following is true off an event with probability at most $\varepsilon^3$,  for all sufficiently small $\varepsilon$ (in a way depending only on $a$ and $\theta$).  There exists an integer $n$ such that
\begin{align*}
n \in \left(\frac{\theta \log(1/\varepsilon)}{3\log(2)} ,  -1 + \frac{\log(1/\varepsilon)}{3\log(2)}\right)
\end{align*}
and $E_n(\Psi)$ occurs.
\end{lemma}

\begin{proof}
The statement of the lemma follows from the exact same argument used to prove Lemma~\ref{lem:event_happens_at_many_scales},  the only difference being that the range of scales is now $\asymp (1-\theta)\log(1/\varepsilon) \asymp a \log(1/\varepsilon)$.
\end{proof}

From now on,  we assume that we are working on the event in the statement of Lemma~\ref{lem:good_event_many_scales} and let $\hat{k}_f$ be the smallest integer $n$ as in the lemma statement for which $E_n(\Psi)$ occurs.  We also assume that all of the crossings $\eta_1,\eta_2,\eta_3,\eta_4$ intersect $\partial B_{\varepsilon}(0)$ before intersecting any other curve in $\{\eta_1,\eta_2,\eta_3,\eta_4\}$.  For all $j \in \{1,2,3,4\}$,  we let $U_j$ denote the connected component of 
\begin{align*}
\BD \setminus \overline{\bigcup_{\ell \in \Gamma_{k_f}^{\text{in}}} \ell  \cup \bigcup_{i=1}^4 \eta_i  \cup \bigcup_{\ell' \in \Gamma_{\hat{k}_f}^{\text{out}}} \ell' \cup B_{\varepsilon}(0)}
\end{align*}
whose boundary contains both the right side of $\eta_j$ and the left side of $\eta_{j+1}$,  where we set $\eta_5 = \eta_1$.  Note that arguing as in the proof of Lemma~\ref{lem:field_has_bounded_boundary_conditions},  we have that the boundary conditions of $\Psi|_{U_j}$ lie in $[c-N,c+N]$ for all $j \in \{1,2,3,4\}$.  Moreover,  for all $j \in \{1,2,3,4\}$,  we let $\phi_j$ denote the unique conformal transformation mapping the simply connected domain $U_j$ onto the Euclidean rectangle 
\begin{align*}
V_j:=(0,L_j) \times (0,1)
\end{align*}
such that the right (resp.\ left) side of $\eta_j$ (resp.\ $\eta_{j+1}$) is mapped to $\{0\} \times (0,1)$ (resp.\ $\{L_j\} \times (0,1)$).  Then,  we have the following analog of Lemma~\ref{lem:rectangles_are_not_wide} which states that the rectangles $V_j$ are not too wide.

\begin{lemma}\label{lem:upper_bound_on_width}
Suppose that we have the same setup described above.  Then,  there exists a universal constant $c>0$,  and a constant $\varepsilon_0 \in (0,1)$ depending only on $a,\delta,N$ and $\theta$ such that
\begin{align*}
L_j \leq \frac{c}{(\theta - 3a\log(2)) \log(1/\varepsilon)} \quad \text{for all} \quad \varepsilon \in (0,\varepsilon_0) , \,\,\, j \in \{1,2,3,4\}.
\end{align*}
\end{lemma}

\begin{proof}
It follows from the exact same argument used to prove Lemma~\ref{lem:rectangles_are_not_wide}.
\end{proof}

Note that since we have assumed that $a \in (0,1/100)$ and $\theta > 1-a/100$,  we have $\theta - 3a\log(2) \geq 1/2$,  and hence $L_j \leq 2c / \log(1/\varepsilon)$ for all $\varepsilon \in (0,\varepsilon_0)$.

Next,  for all $j \in \{1,2,3,4\}$,   we set
\begin{align*}
V_{j,n} = (0,L_j) \times \left(\frac{n}{\log(1/\varepsilon)} ,  \frac{n+1}{\log(1/\varepsilon)}\right)
\end{align*}
for all $n = 0,1,\ldots,\lfloor \log(1/\varepsilon) \rfloor - 1$,  and we write $V_{j,n}^m$ for the three rectangles
\begin{align*}
V_{j,n}^m = (0,L_j) \times \left(\frac{n}{\log(1/\varepsilon)} + \frac{m}{3\log(1/\varepsilon)} ,  \frac{n}{\log(1/\varepsilon)} + \frac{m+1}{3\log(1/\varepsilon)} \right) \quad \text{for} \quad m=0,1,2.
\end{align*}
We also fix a universal constant $u \in (0,\lambda)$ and let $a_{j,n}$ denote the midpoint of the left boundary of $V_{j,n}^1$ and let $P_{j,n}$ be the level line of the field $\Psi|_{U_j} \circ \phi_j^{-1}$ starting from $a_{j,n}$,  of height $b_j-u$ and targeted at the midpoint of the right boundary of $V_{j,n}^1$,  where $b_j \in \{2\lambda (k-1) ,  2\lambda k ,  2\lambda (k+1)\}$ is such that the boundary conditions of $\Psi$ on both the right side of $\eta_j$ and the left side of $\eta_{j+1}$ are given by $b_j$.

Next,  as in Lemma~\ref{lem:level_line_hits_right_side_with_positive_prob},  we show in the following lemma that with positive conditional probability,  the curve $P_{j,n}$ intersects the right boundary of $V_{j,n}$ before intersecting either its top or bottom boundaries.

\begin{lemma}\label{lem:cond_prob_crossing}
There exists a universal constant $p \in (0,1)$ such that the following holds for all $\varepsilon \in (0,\varepsilon_0)$,  where $\varepsilon_0$ denotes the constant in the statement of Lemma~\ref{lem:upper_bound_on_width}.  For every fixed $n \in [0,\lfloor \log(1/\varepsilon) \rfloor -1]_{\BZ}$ and $j \in \{1,2,3,4\}$,  and conditional on $U_1,U_2,U_3,U_4$ and on $P_{j,\tilde{n}}$ for $\tilde{n} \neq n$,  we have that the curve $P_{j,n}$ first intersects the right boundary of the rectangle
\begin{align*}
(0,L_j) \times \left(\vert a_{j,n} \vert- \frac{1}{4 \log(1/\varepsilon)} ,  \vert a_{j,n} \vert + \frac{1}{4 \log(1/\varepsilon)} \right)
\end{align*}
before hitting either its top or bottom boundaries,  with conditional probability at least $p$.
\end{lemma}

\begin{proof}
First,  we note that the choice of the height of $P_{j,n}$ combined with \cite[Lemmas~2.7,  2.8]{miller2017intersections} implies that the law of $P_{j,n}$ stopped at the first time that it intersects either the top,  bottom,  or right boundary of 
\begin{align*}
(0,L_j) \times \left(\vert a_{j,n} \vert - \frac{1}{4 \log(1/\varepsilon)} ,  \vert a_{j,n}\vert  + \frac{1}{4 \log(1/\varepsilon)} \right)
\end{align*}
is absolutely continuous with respect to the law of an $\SLE_4(-1-u/\lambda ; -1+u/\lambda)$ in $V_{j,n}$ with the force points located at $a_{j,n}^-$ and $a_{j,n}^+$ respectively,  and stopped at the first time that it exits the aforementioned rectangle.  Also,  the choice of $u$ implies that such an $\SLE_4(-1-u/\lambda ; -1+u/\lambda)$ process intersects both the left and right sides of $V_{j,n}^1$ a.s.  Thus,  recalling that $L_j \leq 2c/\log(1/\varepsilon)$,  the proof of the lemma is complete by arguing as in the proof of Lemma~\ref{lem:level_line_hits_right_side_with_positive_prob}.  
\end{proof}

Now,  we state the final ingredient in the proof of Proposition~\ref{prop:nested_cle_4_three_arm_exponent}.  It follows from Lemma~\ref{lem:cond_prob_crossing} and it states that with very high probability,  we have that,  for each $j$,  there is some $n$ for which $P_{j,n}$ makes a crossing inside the conformal rectangle $\overline{U_j}$ from left to right.

\begin{lemma}\label{lem:many_level_line_shields}
There exists a universal constant $b_0>0$ and a constant $\varepsilon_0 \in (0,1)$ depending only on $a,\delta,N$,  and $\theta$,  such that the following holds for all $\varepsilon \in (0,\varepsilon_0)$.  Conditional on $U_1,U_2,U_3,U_4$,  we have off an event with conditional probability at most $\varepsilon^{b_0}$,  that for all $j \in \{1,2,3,4\}$,  there exists $n_j \in \{0,1,\ldots,\lfloor \log(1/\varepsilon) \rfloor -1\}$ such that the curve $P_{j,n_j}$ intersects the right side of the rectangle 
\begin{align*}
(0,L_j) \times \left(\vert a_{j,n_{j}} \vert - \frac{1}{4 \log(1/\varepsilon)} ,  \vert a_{j,n_{j}}\vert + \frac{1}{4 \log(1/\varepsilon)} \right)
\end{align*}
before intersecting either its top or bottom boundaries.
\end{lemma}

\begin{proof}
The statement of the lemma follows immediately from Lemma~\ref{lem:cond_prob_crossing}:  since that lemma bounds from below the conditional probability of the crossing event for $P_{j,n}$ given all of the other $P_{j,\tilde n}$,  the number of $n$ for which the crossing event occurs stochastically dominates a binomial random variable with parameters $\lfloor \log(1/\varepsilon)\rfloor$ and $p$,  so that the probability that no such $n$ exists is at most $(1-p)^{\lfloor\log(1/\varepsilon)\rfloor} \leq \varepsilon^{b_0}$ with $b_0 \defeq \min\{\tfrac12\log(1/(1-p)) ,  \tfrac12\}$,  the multiplicative constants and the factor $4$ coming from a union bound over $j \in \{1,2,3,4\}$ being absorbed by decreasing $\varepsilon_0$.
\end{proof}

Now,  we are ready to prove Proposition~\ref{prop:nested_cle_4_three_arm_exponent}.  

\begin{proof}[Proof of Proposition~\ref{prop:nested_cle_4_three_arm_exponent}.]
\stepn{step:three-arm-construction-shields}{Construction of the shields}  First,  we note that combining Lemmas~\ref{lem:bound_on_number_of_iterations} and~\ref{lem:many_level_line_shields},  we obtain that off an event with probability at most $\varepsilon^3 + \varepsilon^{\rho n_{\varepsilon} (1-a) + b_0}$ as $\varepsilon \to 0$ (at a rate that depends only on $a,\delta,N,\rho$,  and $\theta$),  we have that the following holds.  Suppose that the event $F_{\varepsilon}$ in the statement of Lemma~\ref{lem:bound_on_number_of_iterations} occurs.  Then,  for all $j \in \{1,2,3,4\}$,  there exists $n_j \in \{0,1,\ldots,\lfloor \log(1/\varepsilon) \rfloor - 1\}$ such that the curve $\phi_j^{-1}(P_{j,n_{j}})$ intersects $\eta_{j+1}$ before intersecting $\partial U_j \setminus (\eta_j \cup \eta_{j+1})$ for the first time.  Suppose that we are working on the above event.

\stepn{step:three-arm-conclusion}{Conclusion}  Suppose that there are at least four crossings of $C_\varepsilon$ made by loops in $\overline{\Gamma}$.  Then,  as explained in the proof of Proposition~\ref{prop:nested_cle_4_two_arm_exponent},  the event $F_{\varepsilon}$ occurs.  In particular,  there exist level lines $\eta_1,\eta_2,\eta_3,\eta_4$ as in the definition of $F_{\varepsilon}$,  such that at least four crossings of $C_\varepsilon$ made by loops in $\overline{\Gamma}$ merge with $\eta_1,\eta_2,\eta_3$, and $\eta_4$ respectively.  We claim that there are no other crossings of $C_\varepsilon$ made by loops in $\overline{\Gamma}$.  Indeed,  suppose that there exists $\ell \in \overline{\Gamma}$ such that $\ell \cap (\cup_{i=1}^4 \eta_i) = \emptyset$ and $\ell$ crosses $C_\varepsilon$.  Then,  there exists $j \in \{1,2,3,4\}$ such that $\ell$ makes a crossing in $\overline{U_j}$ between the two connected components of $\partial U_j \setminus (\eta_j \cup \eta_{j+1})$.  But then,  $\ell$ has to cross $P_{j,n_j}$ and this contradicts the level line interaction rules since level lines do not cross each other.  It follows that there are at most four crossings of $C_\varepsilon$ made by loops in $\overline{\Gamma}$.  Therefore,  we conclude the proof of the proposition as in the proof of Proposition~\ref{prop:nested_cle_4_two_arm_exponent} by first taking $s \to 0$ and then $a \to 0$,  and noting that $b_0$ is independent of $a,\zeta,\rho$,  and $\theta$;  more precisely,  we choose $a$ and $s$ small enough that $\rho n_\varepsilon(1-a) \geq 2 - b_0/2$,  and then the probability bound $\varepsilon^3 + \varepsilon^{\rho n_\varepsilon(1-a)+b_0}$ is at most $2\varepsilon^{2+b_0/2}$ since $b_0<1<2$,  which gives the statement of the proposition upon replacing $b_0$ by $b_0/4$.
\end{proof}

\subsubsection{Proof of Proposition~\ref{prop:bi_chordal_sle_4_two_arm_exponent}.}
\label{subsubsec:level_lines_four_arms}

Before we give the proof of Proposition~\ref{prop:bi_chordal_sle_4_two_arm_exponent},  we give some background on bi-chordal $\SLE_4$ processes.

Let $D$ be a simply connected domain and fix prime ends $x_1,x_2,x_3,x_4 \in \partial D$ ordered in the counterclockwise way.  A bi-chordal $\SLE_4$ on $D$ with link pattern $\{\{x_1,x_2\},\{x_3,x_4\}\}$ is a probability measure on pairs of simple curves $(\eta_{12},\eta_{34})$ such that $\eta_{12}$ (resp.\ $\eta_{34}$) connects $x_1$ to $x_2$ (resp.\ $x_3$ to $x_4$) in $D$ and the curves $\eta_{12},\eta_{34}$ do not intersect each other,  and such that the following is true: The law of $\eta_{12}$ (resp.\ $\eta_{34}$) given $\eta_{34}$ (resp.\ $\eta_{12}$) is that of a chordal $\SLE_4$ from $x_1$ to $x_2$ (resp.\ from $x_3$ to $x_4$) on the connected component of $D \setminus \eta_{34}$ (resp.\ $D \setminus \eta_{12}$) whose boundary contains $x_1$ and $x_2$ (resp.\ $x_3$ and $x_4$).

The existence of bi-chordal $\SLE_4$ was proved in \cite{mjconfigurational} using the Brownian loop soup measure.  The uniqueness of bi-chordal $\SLE_4$ follows from the exact same argument used to prove \cite[Theorem~4.1]{IG2}.  In particular,  the authors in \cite{IG2} proved the uniqueness of bi-chordal $\SLE_{\kappa}$ for all $\kappa \in (0,4)$ (for a fixed domain and a link pattern),  using the theory of flow lines introduced in \cite{IG1}.  However,  the exact same argument can be applied in the case of bi-chordal $\SLE_4$ where the results in \cite{IG1} are replaced by the corresponding results in the context of level lines proved in \cite{LevelLineGFFI}.  See also \cite[Theorem~1.2]{beffara2021uniqueness} for another proof for all $\kappa \in (0,4]$.

\begin{proposition}\label{prop:conditional_law_of_strands}
Suppose that we have the same setup as in Proposition~\ref{prop:bi_chordal_sle_4_two_arm_exponent}.  Let $E$ denote the event that $P_1$ and $P_2$ merge and that $P_3$ and $P_4$ merge,  in each case before either of the two curves in question reaches its continuation threshold.  Then,  the conditional law given $E$ of the pair $(P_1,P_3)$ (equivalently,  of the pair $(P_2,P_4)$,  since on $E$ these are the time-reversals of $P_1$ and $P_3$,  respectively) is that of a bi-chordal $\SLE_4$ with link pattern $\{\{z_1, z_2\}, \{z_3, z_4\}\}$.
\end{proposition}

\begin{proof}
Let $E_1$ (resp.\ $E_2$) denote the event that the curves $P_1,P_2$ (resp.\ $P_3,P_4$) merge before either of them reaches its continuation threshold.  Let also $E = E_1 \cap E_2$.  Note that \cite[Theorem~1.1.6]{LevelLineGFFI} implies that on the event $E_1$ (resp.\ $E_2$),  the curve $P_2$ (resp.\ $P_4$) is equal to the time-reversal of $P_1$ (resp.\ $P_3$) a.s.  Moreover,  if $E_1$ occurs,  the level line interaction rules imply that both $P_3$ and $P_4$ stay to the left of $P_1$.  In particular,  since both of the force points of $P_3$ and $P_4$ with weight at most $-2$ lie on the boundary of the connected component of $\BD \setminus P_1$ lying to the right of $P_1$,  we obtain that none of the curves $P_3,P_4$ hit their continuation threshold.  Hence $E_1 \subseteq E_2$,  and by symmetry $E_2 \subseteq E_1$;  it follows that $E = E_1 = E_2$.  Also,  \cite[Lemmas~2.3,~2.5]{miller2017intersections} implies that $\BP[E] > 0$.

Next,  we note that the boundary conditions of $h$ combined with the choice of the paths $P_1$, $P_2$, $P_3$ and $P_4$, together imply that the conditional law of $P_3$ given $(P_1,P_2)$ and on the event $E_1$,  is that of a chordal $\SLE_4$ from $z_3$ to $z_4$ on the connected component of $\BD \setminus P_1$ lying to the left of $P_1$.  Similarly,  we have that the conditional law of $P_1$ given $(P_3,P_4)$ and on the event $E_2$,  is that of a chordal $\SLE_4$ from $z_1$ to $z_2$ on the connected component of $\BD \setminus P_3$ lying on the left of $P_3$.  Note that since $E = E_1 = E_2$,  we have that $E$ is a.s.\ determined by $P_1$ (resp.\ $P_3$).  It thus follows from combining with the uniqueness of the law of bi-chordal $\SLE_4$ that the conditional law of $(P_1,P_3)$ (resp.\ $(P_2,P_4)$) given $E$ is that of a bi-chordal $\SLE_4$ on $\BD$ with link pattern $\{\{z_1,z_2\} ,  \{z_3,z_4\}\}$.
\end{proof}

\begin{proof}[Proof of Proposition~\ref{prop:bi_chordal_sle_4_two_arm_exponent}.]
Let $E$ be the event in the statement of Proposition~\ref{prop:conditional_law_of_strands}.  The event in the statement of the proposition is contained in $E$,  and on it the two curves $P_1 \cup P_2$ and $P_3 \cup P_4$ both intersect $\overline{B_{\varepsilon}(0)}$ without intersecting each other.  Therefore,  by Proposition~\ref{prop:conditional_law_of_strands} and Lemma~\ref{lem:bichordal-bulk-4A},  the probability of that event is at most $\BP[E]$ times the probability that both curves of a bi-chordal $\SLE_4$ with link pattern $\{\{z_1,z_2\},\{z_3,z_4\}\}$ intersect $B_{2\varepsilon}(0)$,  which is at most $C\varepsilon^2$ for a universal constant $C$.
\end{proof}
\end{appendix}

\bibliographystyle{alpha}
\bibliography{references}

\end{document}